\documentclass[reqno, 12pt]{amsart}
\usepackage{marginnote,amsmath,geometry,fancyhdr,amssymb,etoolbox,wrapfig,setspace,amsthm,tikz,babel,inputenc,mathtools,bbm,tikz-cd,mathrsfs,graphicx,bm,listings}
\usepackage{longtable}
\usepackage{comment}
\usetikzlibrary{calc,intersections,decorations.markings}
\usepackage[shortlabels]{enumitem}
\usepackage[colorlinks=true, pdfstartview=FitV, linkcolor=purple, citecolor=blue,
filecolor=violet, urlcolor=violet]{hyperref}
\usepackage{cleveref} 

\definecolor{darkGreen}{rgb}{0.01,0.50,0.32}

\newtheorem{theorem}{Theorem}[section]

\newtheorem{lemma}[theorem]{Lemma}
\newtheorem{proposition}[theorem]{Proposition}
\newtheorem{corollary}[theorem]{Corollary}

\newtheorem{construction}[theorem]{Construction}

\theoremstyle{definition}
\newtheorem{definition}[theorem]{Definition}
\newtheorem{remark}[theorem]{Remark}
\newtheorem{examplex}[theorem]{Example}
\newenvironment{example}
  {\pushQED{\qed}\examplex}
  {\popQED\endexamplex}
\Crefname{examplex}{Example}{Examples}

\newcommand{\CS}{\mathcal{S}}

\newcommand{\CC}{\mathcal{C}}

\newcommand{\weezer}{\leftexp{=}{\kern-0.23em\operatorname{W}}^{\kern-0.21em =}}

\newcommand{\CI}{\mathcal{I}}

\newcommand{\CL}{\mathcal{L}}
\newcommand{\CE}{\mathcal{E}}
\newcommand{\CF}{\mathcal{F}}

\newcommand{\CX}{\mathcal{X}}

\newcommand{\CO}{\mathcal{O}}
\newcommand{\CH}{\mathcal{H}}

\newcommand{\R}{\mathbb{R}}

\newcommand{\Q}{\mathbb{Q}}
\newcommand{\G}{\mathbb{G}}
\newcommand{\Z}{\mathbb{Z}}

\newcommand{\A}{\mathbb{A}}
\renewcommand{\P}{\mathbb{P}}
\newcommand{\ve}{\varepsilon}

\newcommand{\del}{\partial}
\newcommand{\ds}{\displaystyle}

\newcommand{\spn}{\op{span}}

\DeclareMathOperator{\Hom}{Hom}

\DeclareMathOperator{\coker}{coker}
\DeclareMathOperator{\supp}{supp}

\newcommand{\lp}{\left(}
\newcommand{\rp}{\right)}

\newcommand{\la}{\langle}
\newcommand{\ra}{\rangle}

\newcommand{\hide}[1]{}

\newcommand{\wh}[1]{\widehat{#1}}

\newcommand{\inc}{\xhookrightarrow{}}
\newcommand{\wt}{\widetilde}

\newcommand{\op}[1]{\operatorname{#1}}
\newcommand{\Gm}{\mathbb{G}_{\op{m}}}

\newcommand{\Db}{\op{D}^{\op{b}}}
\newcommand{\D}{\op{D}}

\newcommand{\Cox}{{\op{Cox}}}
\newcommand{\Cone}{\op{Cone}}
\newcommand{\Conv}{\op{Conv}}

\newcommand{\Spec}{\op{Spec}}

\newcommand{\floor}[1]{\lfloor #1 \rfloor}

\renewcommand{\bar}{\overline}
\newcommand{\Kbar}{\bar{K}}
\newcommand{\Relint}{\op{Relint}}
\newcommand{\Int}{\op{Int}}
\newcommand{\face}{\preceq}

\DeclareMathOperator{\codim}{codim}
\newcommand{\us}{\op{us}}
\renewcommand{\ss}{\op{ss}}

\DeclareMathOperator{\Tor}{Tor}
\DeclareMathOperator{\Cl}{Cl}
\newcommand{\bs}[1]{\boldsymbol{#1}}
\newcommand{\wrel}{\prec}

\DeclareMathOperator{\HHL}{H}
\DeclareMathOperator{\cells}{cells}
\newcommand{\ceil}[1]{\left\lceil #1 \right\rceil}
\DeclareMathOperator{\proj}{proj}
\DeclareMathOperator{\omax}{\op{max}}
\DeclareMathOperator{\omin}{\op{min}}

\tikzset{arr /.tip = {Latex[round,width=1.5mm,length=2mm]},labelnode/.style={draw=none,fill=none,circle, inner sep=0mm}}

\tikzset{double line with arrow/.style args={#1,#2}{decorate,decoration={markings,%
mark=at position 0 with {\coordinate (ta-base-1) at (0,1pt);
\coordinate (ta-base-2) at (0,-1pt);},
mark=at position 1 with {\draw[#1] (ta-base-1) -- (0,1pt);
\draw[#2] (ta-base-2) -- (0,-1pt);
}}}}
\newcommand{\lst}[1]{\lstinline{#1}}
\definecolor{codegreen}{rgb}{0,0.6,0}
\definecolor{codegray}{rgb}{0.5,0.5,0.5}
\definecolor{codepurple}{rgb}{0.58,0,0.82}
\definecolor{backcolor}{rgb}{0.95,0.95,0.92}
\definecolor{m2blue}{RGB}{0,0,245}
\definecolor{m2green}{RGB}{68,136,51}
\definecolor{m2purple}{RGB}{135,42,141}
\definecolor{m2teal}{RGB}{60,136,137}
\lstdefinestyle{mystyle}{
    backgroundcolor=\color{backcolor},   
    commentstyle=\color{codegreen},
    keywordstyle=[1]{\color{m2blue}},
    keywordstyle=[2]{\color{m2purple}},
    keywordstyle=[3]{\color{m2green}},
    keywordstyle=[4]{\color{m2teal}},
    numberstyle=\tiny\color{codegray},
    stringstyle=\color{codepurple},
    basicstyle=\ttfamily\footnotesize,
    breakatwhitespace=false,         
    breaklines=true,
    breakautoindent=false,
    breakindent=0pt,
    keepspaces=true,                 
    numbers=left,                    
    numbersep=5pt,                  
    showspaces=false,                
    showstringspaces=false,
    showtabs=false,                  
    tabsize=2,
    morekeywords=[1]{method,interiorPoint,matrix,submatrix,map,apply,floor,intersection,ceiling,min,max,mutableMatrix,numrows,numcols,any,dim,entries,join,append,scan,flatten,transpose,isMember,unique,loadPackage,positions,promote,product,select,length,rays,contains,toList,table},
    morekeywords=[2]{for,from,list,then,return,and,not,do,to,continue,else,in,new,if},
    morekeywords=[3]{Matrix,ZZ,Polyhedron,List,HashTable,SimplicialComplex,Cone,Sequence,Boolean},
    morekeywords=[4]{SimplicialComplexes,HH,true,false}
}
\title{Comparing grade restriction windows to Cox windows}
\author{David Favero, Mykola Sapronov, Dash Stevanovich}
\begin{document}

\maketitle

\begin{abstract}
We introduce a new class of toric varieties which admit an HKKN stratification compatible with the combinatorics of the GKZ fan.  This class includes all smooth projective toric varieties of Picard ranks 1 and 2.  We show that grade restriction windows \cite{HL, BFK} arising from our HKKN stratifications agree with the intrinsic window categories appearing in \cite{CoxCategory}.  On the other hand, we provide an example of a smooth projective toric DM stack where this intrinsic window is never a grade restriction window despite the many possible choices.  As an application, we demonstrate that reordering and/or taking a subset of an HKKN stratification sometimes changes the  grade restriction window and sometimes does not.
\end{abstract}

\tableofcontents

\section{Introduction}
Grade restriction windows originated in a physics manuscript of Herbst--Hori--Page \cite{herbst2008phases} as part of their analysis of $B$-branes for abelian gauged linear sigma models.  These were quickly used by physicists Herbst--Walcher to create powerful comparisons between derived categories of toric varieties \cite{herbst2012unipotence}. Concurrently, Segal \cite{segal2011equivalences} provided new insight into  work of Orlov \cite{orlov2010derived} by demonstrating that both derived categories and matrix factorization categories could be compared using geometric invariant theory (GIT).
These works were  followed up by Halpern-Leistner \cite{HL} and Ballard--Favero--Katzarkov \cite{BFK} who generalized Segal's GIT result to a much broader geometric setting.  GIT is now established as a powerful way to compare and connect derived categories and birational geometry, unifying some of the core literature on equivalences, semi-orthogonal decompositions, and exceptional collections. Still, at the technical heart of these results, remains the notion of grade restriction windows from high-energy theoretical physics.

Grade restriction windows were formalized to mathematics roughly as follows.  Let $G$ act on $X$ and choose a one-parameter subgroup $\lambda: \mathbb G_m \to G$.  Consider the contracting locus
\[
S_\lambda \coloneq \{ x \in X \ | \ \lim_{t \to 0} \lambda(t) \cdot x \text{ exists} \}.
\]
Given an interval $[a,b] \subseteq \mathbb R$, a \emph{grade restriction window} is the full subcategory of the derived category of the quotient stack
\[
W_{\lambda, [a,b]} \subseteq D([X/G])
\]
such that the derived restriction of objects to a $\lambda$-fixed point has $\lambda$-weights in $[a,b]$.  The key result in \cite{BFK, HL}  is that restriction induces an equivalence of categories
\begin{equation} \label{eq: window equiv}
    D([X \backslash (G \cdot S_\lambda)/G]) \cong W_{\lambda, [a,a+T^{>0}_{\lambda})}
\end{equation}
where $T_\lambda^{>0}$ is the sum of the positive $\lambda$-weights of the conormal bundle to $S_\lambda$. 

Going back to work of Kempf \cite{kempf1978instability} and Hesselink \cite{hesselink1979desingularizations},  results of Kirwan \cite{kirwan2020cohomology} and Ness \cite{ness1984stratification} demonstrate that   the unstable locus from GIT is stratified by $G \cdot S_\lambda$ as above for a certain ordered set of one-parameter subgroups $\bs{\lambda} = \lambda_1, \dots, \lambda_t$. This is called a \emph{Hesselink--Kampen--Kempf--Ness (HKKN)-stratification}. 
Hence, iterating \eqref{eq: window equiv} realizes the derived category of a GIT quotient as a full subcategory of the derived category of the quotient stack
\[
D([U/G]) \cong W_{\bs{\lambda}} \subset D([X/G]).
\]

The grade restriction windows coming from HKKN stratifications  are specific instances in what has now been imagined more broadly as window subcategories.  That is, a \emph{window}  $W \subseteq D([X/G])$ is a full subcategory of the derived category of the quotient stack whose restriction to an open substack $D([U/G])$ is an equivalence.  Indeed, this is the perspective which has been taken by \v{S}penko--Van den Bergh \cite{vspenko2017non} and Halpern-Leistner--Sam \cite{halpern2020combinatorial} who described ``magic windows'' which exist in the particularly nice, quasi-symmetric case.  
In this paper, we focus our attention on comparing grade restriction windows to a new type of window  discovered by Ballard--Berkesch--Brown--Cranton Heller--Erman--Favero--Hanlon--Huang--Ganatra \cite{CoxCategory} as part of their birational solution to King's conjecture. Therein they describe what we call the \emph{Cox window} 
\[
D(X_\Sigma) \cong W_{\op{Cox}} \subset D([X/G])
\]
for a smooth projective toric variety $X_\Sigma$. The Cox window is intrinsic to a toric variety. It comes simply from the toric data and is not decorated by all the additional choices of grade restriction windows.

In summary, we now have two methods of embedding $D(X_\Sigma)$ into $D([X/G])$; one via grade restriction windows, and the other via Cox windows.  The former involves a number of choices and the latter seems more intrinsic.  On the other hand, both of these methods are intimately related to birational geometry and geometric invariant theory.  It is therefore natural to wonder what the relationship between the two constructions is (see \cite[Remark 7.9]{CoxCategory}).  

Our work here provides a comprehensive analysis of the comparison between these two.  For one, we demonstrate that for a particular class of toric varieties which we call arrangeable  (see \Cref{def:arrangeable}), the two windows agree. 
\begin{theorem}[see \Cref{thm:window_equality}]
\label{intro thm: window_equality}
 For arrangable GIT problems, there exists an HKKN stratification such that the grade restriction window and Cox window agree.    
\end{theorem}
While limiting as a whole, arrangeable GIT problems still include the important special cases of all projective toric varieties of picard rank 1 and 2 as well as higher picard rank examples like blow-ups of projective space.  On the other hand, without the arrangability assumption the theorem does not always hold.  That is, we construct a smooth projective toric Deligne-Mumford stack where no matter what choice of HKKN stratification and shift parameters one chooses, the two constructions cannot agree.
\begin{theorem}
[=\Cref{thm:no_stratification_making_windows_equal_exist}]
\label{intro thm: counterexample}
There exists a smooth projective toric Deligne-Mumford stack whose Cox window is not a grade restriction window.
\end{theorem}

To prove our results, we consider objects in the Cox window that restrict to nef line bundles and check whether or not they obey the grade restriction rules.  Here, we benefit from the work of Cranton Heller \cite{CH25} (upgraded to stacks in \Cref{prop:FMonAmplesGivesF}), who gives a topological interpretation of these objects.    As a consequence, we obtain a precise topological criterion for when a given grade restriction window agrees with the Cox window.  That is, it reduces to the contractibility of certain locally-closed subsets of $\mathbb R^n$, see \Cref{thm:GRW_is_Cox_iff}.  
To analyze the homotopy type of these regions, we cover them with closed convex regions and use the nerve theorem due to Leray \cite{Leray1945} and Borsuk \cite{Borsuk1948} (see also \cite{BjornerTopologicalMethods}).  In other words,  the contractibility of these regions is reduced to combinatorics and a careful analysis of these combinatorics leads to our theorems.

Once all these tools are locked in, we consider examples in \S\ref{section:examples}.  One such example is \Cref{intro thm: counterexample}.  Using our topological criterion, we also analyze several different HKKN stratifications for the same GIT problem.  We observe that some of these satisfy our topological criterion while others do not.  In particular, we observe that closely related HKKN stratifications may or may not give the same grade restriction window.  
\begin{theorem}[see \Cref{cor:may_or_may_not}\ref{cor:may_or_may_not_reorder_MAY},\ref{cor:may_or_may_not_reorder_NOT}]
\label{intro thm: reordering}
Grade restriction windows from reordered HKKN stratifications sometimes agree and sometimes do not.    
\end{theorem}
\begin{theorem}[see \Cref{cor:may_or_may_not}\ref{cor:may_or_may_not_sub_MAY},\ref{cor:may_or_may_not_sub_NOT}]
\label{intro thm: substratification}
A HKKN stratification which is an ordered subset of another HKKN stratification sometimes yields the same grade restriction window and sometimes does not.   
\end{theorem}
The appearance of a condition like arrangability in our results is not ad hoc.  Rather, it seems to arise due to a fundamental inconsistency between the Gelfand--Kapranov--Zelevinky combinatorics of the secondary fan \cite{GKZ} (see also \cite[Chapters 14-15]{CLS}) and those of HKKN stratifications.  We show that arrangability rectifies this problem, allowing for an HKKN stratification which matches the toric combinatorics.  
\begin{theorem}
[=\Cref{thm:existence_of_wall_stratifications}]
\label{intro thm: arrangable poset} Let $\Gamma$ be a chamber of an arrangeable secondary fan $\Sigma_{GKZ}$.  There is a natural poset ordering on the walls of $\Gamma$ such that any extension to a total ordering  induces an HKKN stratification of $V(B_\Gamma)$.
\end{theorem}

The paper is organized as follows.  In \S\ref{sec: toric geometry}, we lay out the GIT and GKZ combinatorics underlying our discussion.  HKKN stratifications and their relationship with arrangeability are discussed in \S\ref{section:stratifications}, proving \Cref{intro thm: arrangable poset}.  In \S\ref{section: windows}, we define grade restriction windows and Cox windows and introduce Cranton Heller's generating set for the Cox window (which we upgrade to stacks).  This sets us up for \S\ref{section:window_comparison} where we dig into the topological interpretation of these generating objects.  In \S\ref{section:main_thm} we prove \Cref{intro thm: window_equality} and in \S\ref{section:examples} we give several examples which demonstrate Theorems~\ref{intro thm: counterexample},~\ref{intro thm: reordering}, and~\ref{intro thm: substratification}.
\section*{Acknowledgments}
All three authors were supported by the National Science Foundation under DMS award numbers 2302262 and 2412039.  Dash Stevanovich was additionally supported by a GAANN Fellowship under award P200A240046.
This project greatly benefited from the creative atmosphere and collegial environment at the Homological mirror symmetry and multigraded commutative algebra workshop at the American Institute for Mathematics in 2025. The computer algebra system Macaulay2 \cite{M2} provided valuable assistance. We are also grateful for fruitful conversations with Matthew Ballard, Daniel Erman, Kimoi Kemboi, Mahrud Sayrafi, and Jay Yang.

\section*{AI Declaration}
AI tools were not used for this work. 

\section{Toric Geometry Background}\label{sec: toric geometry}
In this section we collect the basic facts about toric geometry and toric geometric invariant theory that we will need in later sections. Much of this is restating standard results from references like \cite{CLS} in the language of toric stacks. In \S\ref{subsection:toric_GIT} we recall the toric GIT setup and introduce the secondary fan. We describe a standard combinatorial method for constructing the secondary fan in \S\ref{subsection:combo_description_of_secondary_fan}. In \S\ref{subsection:stacky_fans} we describe the stacky fan data associated to toric stacks and discuss the effect of ray-dropping. In \S\ref{subsection:nef_divisors} we describe nefness in the stacky setting and recall some of the polytope combinatorics one uses to anaylze divisors in toric geometry. We introduce a distinguished class of 1-PSs and describe their relationship to the walls of the secondary fan in \S\ref{subsection:1PS_and_walls}.

Throughout we fix an algebraically closed field of characteristic zero $k$ with multiplicative group $\Gm \coloneqq \Spec k[t,t^{-1}]$. We let $D(\CX)$ denote the bounded derived category of coherent sheaves on $\CX$. For a $\wh{G}$-graded ring $S$ we let $D(S)$ denote the bounded derived category of finitely-generated $\wh{G}$-graded $S$-modules. All functors between derived categories are taken to be derived unless otherwise stated. In section \hyperref[section:notational_index]{NI}
there is a notational index summarizing the main notation used frequently throughout the paper.

\subsection{Overview of Toric Geometric Invariant Theory} \label{subsection:toric_GIT}
Roughly speaking, given a group $G$ acting on a space $X$, Geometric Invariant Theory (GIT) provides a way of parameterizing the ``nice" quotients that can be formed by taking quotients of open subsets of $X$ by $G$. Much of the study of simplicial toric varieties, and more generally smooth toric Deligne-Mumford (DM) stacks, is described well in terms of toric GIT. In our study of these stacks, we will primarily take this toric GIT perspective, and ignore the underlying fans/stacky fans when we can. For a more thorough treatment of toric GIT, see chapters 14 and 15 of \cite{CLS}.

Throughout we will assume that we have an embedding of a linear algebraic group $G\subset \Gm^n$ acting diagonally on $\A^n$. Applying $\wh{(-)}\coloneqq \Hom(-,\Gm)$ to this inclusion of linear algebraic groups we get a surjective map of Abelian groups $\pi:\Z^n\to \wh{G}$. We define $M\coloneqq \ker \pi$ to obtain the exact sequence
\begin{equation}\label{eqn:ses_for_lattices}
    \begin{tikzcd}[column sep = 0.5cm] 0 \arrow[r] & M\arrow[r,"i"] & \Z^n \arrow[r,"\pi"] & \wh{G}\arrow[r] & 0.
\end{tikzcd}
\end{equation}
By applying $-\otimes \R$ to the above short exact sequence of abelian groups, we obtain a short exact sequence of vector spaces:
\begin{equation}\label{eqn:ses}
\begin{tikzcd}[column sep = 0.5cm] 0 \arrow[r] & M_\R\arrow[r,"i_\R"] & \R^n \arrow[r,"\pi_\R"] & \wh{G}_\R\arrow[r] & 0.
\end{tikzcd}
\end{equation}
We define the \emph{weights} of the embedding $G\subset \Gm^n$ to be $\beta_\rho\coloneqq \pi_\R(e_\rho)$ for $\rho\in [n]$. Throughout we will let $S\coloneqq k[x_1,\hdots, x_n]$ denote the \emph{Cox ring} and interpret the $G$ action on $\A^n = \Spec S$ equivalently as a $\wh{G}$-grading on $S$ via $\deg(x_\rho) = \pi(e_\rho)$. 

Our construction of toric GIT will rely on a description via \emph{one-parameter subgroups}, i.e. injective algebraic group homomorphisms $\lambda:\Gm \to G$. To any such 1-PS $\lambda$ we dualize via $\wh{(-)}$ to get
\[\wh{\lambda}:\wh{G}\to \Z,\]
as $\Hom(\Gm,\Gm) \cong \Z$. Similarly, if we have an abelian group homomorphism $f:\wh{G}\to \Z$, we can dualize to a morphism
\[\wh{f}:\Gm \to \wh{\wh{G}}\]
of linear algebraic groups. Because $\Gm$ is divisible, $\wh{\lp -\rp}$ is exact. This, with the fact that $G$ is a subtorus of $\Gm^n$ implies that $\wh{\wh{G}}\cong G$. So, the above dualization $\wh{f}:\Gm \to G$ can be interpreted as a 1-PS. In fact, this dualization induces a bijection
\[\Hom(\Gm, G) \cong \Hom_{\Z}(\wh{G}, \Z).\]
As such, throughout we will conflate 1-PSs and their duals. When clear, we will also conflate $\wh{\lambda}:\wh{G}\to\Z$ with its associated $\R$-tensored vector space map $\wh{\lambda}_\R : \wh{G}_\R\to \R$. With this in mind, we think of 1-PSs primarily in terms of their associated rational hyperplanes in $\wh{G}_\R$.

\begin{definition} Let $\lambda:\Gm\to G$ be a 1-PS. We define the following subsets of $I_\lambda^{\star}\subset [n]$ and $H_\lambda^{\star}\subset \wh{G}_\R$ respectively:
\[\def\arraycolsep{1cm}\begin{array}{cc}
    I_\lambda^{=0} \coloneqq \{i\in [n]\ | \ \wh{\lambda}(\beta_i) = 0\}, & H_{\lambda}^{=0} = \{\alpha\in \wh{G}_\R \ | \ \wh{\lambda}_\R(\alpha) = 0\}\\
    I_\lambda^{\le 0} \coloneqq \{i\in [n]\ | \ \wh{\lambda}(\beta_i)  \le 0\}, & H_{\lambda}^{\le 0} = \{\alpha\in \wh{G}_\R \ | \ \wh{\lambda}_\R(\alpha) \le 0\}\\
    I_\lambda^{\ge 0} \coloneqq \{i\in [n]\ | \ \wh{\lambda}(\beta_i)  \ge 0\}, & H_{\lambda}^{\ge 0} = \{\alpha\in \wh{G}_\R \ | \ \wh{\lambda}_\R(\alpha) \ge 0\}\\
    I_\lambda^{<0} \coloneqq \{i\in [n]\ | \ \wh{\lambda}(\beta_i)  < 0\}, & H_{\lambda}^{<0} = \{\alpha\in \wh{G}_\R \ | \ \wh{\lambda}_\R(\alpha) < 0\}\\
    I_\lambda^{>0} \coloneqq \{i\in [n]\ | \ \wh{\lambda}(\beta_i)  > 0\}, & H_{\lambda}^{>0} = \{\alpha\in \wh{G}_\R \ | \ \wh{\lambda}_\R(\alpha) > 0\}.
\end{array}\]
\end{definition}
It will be useful to us to both consider these collections $I_\lambda^\star$ as subsets of $[n]$ and as ideals of $S$. For $I\subset [n]$, we will let $\la I\ra \coloneqq \la x_i \ | \ i\in I\ra\subset S $ and $V(I) \coloneqq V(\la I \ra )$. In particular we let $\la \emptyset \ra = \la 0\ra$. With this notation in hand, we construct the following ideal of $S$.
\begin{definition} Given $\chi\in \wh{G}_\R$, we define the following ideal in $S$:
\[B_\chi \coloneqq \bigcap_{\wh{\lambda}_\R(\chi) > 0}\la I_\lambda^{>0}\ra \]
where this intersection is taken over all 1-PSs $\lambda$ so that $\wh{\lambda}_\R(\chi) > 0$. In the language of GIT, $V(B_\chi)$ is the \emph{unstable locus} of the linearization with respect to $\chi$. Sometimes we will denote $V(B_\chi)$ as $(\A^n)^{\us}(\chi)$ for this reason. For any $\chi$ so that $B_\chi\ne 0$, define
\[\sigma_\chi = \overline{\{\alpha\in \wh{G}_\R \ | \ B_\chi = B_\alpha\}}.\]
\end{definition}
The following is a well-known consequence of these definitions (see e.g. \cite{GKZ,CLS}).
\begin{proposition}\label{prop:def_of_secondary_fan} For any $\chi\in \wh{G}_\R$ such that $B_\chi \ne 0$, $\sigma_\chi$ is a strongly convex rational polyhedral cone. Moreover, $\{\sigma_\chi \ | \ \chi\in\wh{G}_\R, \ B_\chi \ne 0\}$ form a fan.
\end{proposition}
We denote the fan guaranteed by the above proposition by $\Sigma_{GKZ}$ and refer to it as the \textit{secondary fan} or the \textit{GKZ fan} after Gelfand, Kapranov, and Zelevinsky \cite{GKZ}. In light of \Cref{prop:def_of_secondary_fan}, for any cone $\sigma\in \Sigma_{GKZ}$ we may define $B_\sigma$ to be $B_\chi$ for any $\chi\in \Relint \sigma$. It is a fundamental result of toric GIT that associated to the maximal cones $\Gamma$ of $\Sigma_{GKZ}$, which we call \emph{chambers}, one obtains a smooth toric DM stack. To a chamber $\Gamma \in \Sigma_{GKZ}$, define the open set
\[U_\Gamma \coloneqq \A^n \setminus V(B_\Gamma).\]
The quotient stack $\CX_\Gamma\coloneqq [U_\Gamma/ G]$ is a smooth toric DM stack. 

In this paper we will focus exclusively on \emph{projective} GIT problems. In all but \S\ref{section:stratifications}, we will focus on GIT problems which are additionally \emph{primitive geometric} (see \cite[Proposition 15.1.6]{CLS}). By these we mean the following.

\begin{definition}\label{def:projective_GIT_problem} We say that the GIT problem associated to an inclusion $G\subset \Gm^n$ is \emph{projective} if for each $\rho =1 ,\hdots, n$, $\beta_\rho\ne 0$ and $|\Sigma_{GKZ}|$ is a strongly convex cone in $\wh{G}_\R$. 
\end{definition}

\begin{definition}\label{def:primitive_geometric} We say that the GIT problem associated to an inclusion $G\subset \Gm^n$ is \emph{primitive geometric} if the lattice points $u_{\rho}\coloneqq i^{\vee}(e_{\rho})\in M^\vee, \rho \in [n]$ are non-zero, primitive and generate distinct rays.
\end{definition}

\begin{remark}\label{rmk:primitive_geometric} 
Any projective, primitive geometric GIT problem arises from the torus-invariant divisor short exact sequence of a projective toric variety (see \cite[Example 14.2.14]{CLS}). 
Although we make the assumption that the GIT problems appearing in this paper are primitive geometric (with the exception of \S\ref{section:stratifications}), we believe this assumption is unnecessary. In particular, none of our proofs rely on this assumption.
However, one of the main categories appearing in this paper, the Cox windows of \cite{CoxCategory}, were only proven to exist in the primitive geometric setting.  
\end{remark}

\subsection{A combinatorial description of the secondary fan}\label{subsection:combo_description_of_secondary_fan} In the previous subsection we introduced the secondary fan from the GIT perspective. In this subsection we describe the methods by which to view $\Sigma_{GKZ}$ as a purely combinatorial gadget arising from the weights $\beta_\rho\in \wh{G}_\R$. We closely follow chapter 15 of \cite{CLS} in this subsection.

Strictly speaking, to make sense of the cones in the secondary fan as \textit{rational}, one should first define the lattice in which ray generators for cones in this fan live. In our case this lattice is $\wh{G}/\Tor(\wh{G})$ where $\Tor(\wh{G})$ is the torsion subgroup of $\wh{G}$. Observe that $(\wh{G}/\Tor(\wh{G}))_\R \cong \wh{G}_\R$. Later we will construct cones and fans in the dual vector space $\wh{G}_\R^*$. As $\Hom_\Z(\Tor(\wh{G}),\Z) = 0$, we have
\[\Hom_\Z(\wh{G}/\Tor(\wh{G}),\Z) \cong \Hom_\Z(\wh{G},\Z)\]
so the dual lattice to $\wh{G}/\Tor(\wh{G})$ is $\Hom(\wh{G},\Z)$. Similarly one computes that
\[\Hom_\Z(\wh{G},\Z)_\R \cong \Hom_\R(\wh{G}_\R,\R) = \wh{G}_\R^\ast.\]
Thus we have dual lattices $\wh{G}/\Tor(\wh{G})$, $\Hom_\Z(\wh{G},\Z)$ with respective dual vector spaces $\wh{G}_\R, \wh{G}_\R^\ast$.

In light of the above we interpret the cones over $\beta_\rho$ as rational cones.
\begin{definition} For $J\subset [n]$, we define the rational cone $\beta_J\coloneqq \Cone(\beta_\rho \ | \ \rho\in J)$ in $\wh{G}_\R$.
\end{definition}
This leads to the following well-known characterization of $B_\chi$. We provide a proof for clarity and because the convex-geometric ideas are similar to those of proofs found in \S\ref{section:stratifications}.
\begin{lemma}\label{lemma:B_chi_via_betas} For any $\sigma\in \Sigma_{GKZ}$,
\[B_\sigma = \la x^J \ | \ \sigma\subset \beta_J\ra.\]
\end{lemma}
\begin{proof} Let $\chi\in \Relint \sigma$ so that by \Cref{prop:def_of_secondary_fan}, $B_\sigma = B_\chi$. We will first show that $B_\chi = \la x^J \ | \ \chi\in \beta_J\ra$. Suppose by way of contradiction $\la x^J \ | \ \chi\in \beta_J\ra$ was not contained in $B_\chi$. Then, there is a $J\subset [n]$ so that $\chi\in \beta_J$ and a 1-PS $\lambda$ so that $\wh{\lambda}_\R(\chi) > 0$ but $J\cap I_\lambda^{>0} = \emptyset$. This would imply $J\subset I_{\lambda}^{\le 0}$. Thus $\chi \in \beta_J\subset H_{\lambda}^{\le 0}$. This contradicts $\wh{\lambda}_\R(\chi) > 0$. So we have that $\la x^J \ | \ \chi\in \beta_J \ra\subset B_\chi$.

To show the opposite inclusion, see that minimal generators of $B_\chi$ are of the type $x^J$ where $J\subset [n]$ has the property that for every $\lambda$ such that $\wh{\lambda}_\R(\chi) > 0$, $I_\lambda^{>0}\cap J \ne\emptyset$. Suppose by way of contradiction that $x^J$ is a minimal generator of $B_\chi$ but $\chi\notin \beta_J$. Then there is a 1-PS $\lambda$ so that the hyperplane $H_{\lambda}^{=0}$ separates $\chi$ and $\beta_J$. That is, $\wh{\lambda}_\R(\chi) > 0$ but $\beta_J\subset H_{\lambda}^{\le 0}$. Thus $J\cap I_{\lambda}^{>0} = \emptyset$. This contradicts our choice of $J$, so we conclude that $\chi\in \beta_J$. This implies $B_\chi\subset \la x^J \ | \ \chi\in \beta_J\ra$ thus
\[B_\chi = \la x^J \ | \ \chi\in \beta_J\ra.\]

Now by \Cref{prop:def_of_secondary_fan}, it follows that $\chi\in \beta_J$ if and only if $\Relint \sigma\subset \beta_J$. As each $\beta_J$ is closed, we conclude that $\chi\in \beta_J$ if and only if $\sigma\subset \beta_J$. The statement follows.
\end{proof}
As an immediate consequence of \Cref{lemma:B_chi_via_betas}, one obtains the following well-known characterization of $\sigma_\chi$ and chambers of $\Sigma_{GKZ}$.
\begin{corollary}\label{cor:gamma_hits_beta_J_implies_contained} For any $\chi\in |\Sigma_{GKZ}|$, recall $\sigma_\chi$, the minimal cone of $\Sigma_{GKZ}$ containing $\chi$. We have
\[\sigma_\chi = \bigcap_{J\subset [n], \ \chi\in \beta_J}\beta_J = \bigcap_{J\subset [n], \sigma_\chi\subset \beta_J}\beta_J.\]
As a consequence, if $\Gamma\in \Sigma_{GKZ}$ is a chamber, then $\beta_J \cap \Int \Gamma \ne \emptyset$ implies $\Gamma\subset \beta_J$.
\end{corollary}

\subsection{Stacky fans and ray dropping}\label{subsection:stacky_fans}
Suppose as before that we have a projective, geometrically primitive GIT problem with associated secondary fan $\Sigma_{GKZ}$. When $\Gamma$ is a chamber of $\Sigma_{GKZ}$, we associate to it the global quotient stack $\CX_\Gamma = [(\A^n\setminus V(B_\Gamma))/G]$ as described in \S\ref{subsection:toric_GIT}. There is an alternative description of $\CX_\Gamma$ using stacky fans that we make use of here. Following \cite{BCS04} we make the following definition.
\begin{definition}\label{def:stacky_fan} A \emph{stacky fan} is the data of a fan $\Sigma$ in $N_\R$ and a homomorphism $\alpha:\Z^{\Sigma(1)}\to N$ with $\alpha(e_\rho) = b_\rho u_\rho$ where $b_\rho > 0$ and $u_\rho$ is the primitive ray generator of $\rho\in \Sigma(1)$. Then $\CX_{\Sigma, \alpha}$ is defined to be the quotient stack
\[\CX_{\Sigma,\alpha} \coloneqq \left[ \lp \A^{\Sigma(1)}\setminus V(B_\Sigma)\rp / G_{\Sigma,\alpha}\right]\]
where $B_\Sigma$ is the irrelevant ideal of $\Sigma$ (see \cite[Section 5.1]{CLS}) and
\[G_{\Sigma,\alpha} \coloneqq \Hom(\coker \beta^\vee, \Gm)\subset \Gm^{\Sigma(1)}.\]
\end{definition}
By \cite[Proposition 3.2]{BCS04}, when $\Sigma$ is simplicial, $\CX_{\Sigma, \beta}$ is a smooth toric DM stack. If $b_\rho = 1$ for each $\rho$, $\Hom(G_{\Sigma,\alpha},\Gm) = \Cl(\CX_{\Sigma,\alpha}) = \Cl(X_\Sigma)$, so $\CX_{\Sigma,\alpha} = \CX_\Sigma$, the usual construction of a quotient stack associated with  toric variety $X_{\Sigma}$, appearing in \cite[Chapter 14]{CLS}. In fact, since we assume that our GIT problems are geometrically primitive, this is precisely our setting. Nevertheless, since we believe that geometric primitivity assumption could be lifted, we give a general construction of a stacky fan from $\Gamma \in \Sigma_{GKZ}$. 

Given a chamber $\Gamma \in \Sigma_{GKZ}$, we first need to associate a fan $\Sigma_\Gamma$ in $N_\R$ to $\Gamma$. Here, the construction is precisely as in  \cite[Chapters 14, 15]{CLS}, but in brief, we define the fan
\begin{equation}\label{eq:fan_from_chamber}
\Sigma_\Gamma \coloneqq \{\Cone(u_\rho) \ | \ \rho\notin J) \ | \ J\subset [n], \ \Gamma\subset \beta_J\}
\end{equation}
in $N_\R$ where $i^\vee:\Z^{n}\to N$ is the dual of $i:M\to \Z^n$ in \eqref{eqn:ses_for_lattices}. We then define
\[I_\Gamma^{\emptyset} \coloneqq [n]\setminus \Sigma_\Gamma(1).\]
The set $I_\Gamma^{\emptyset}$ is often referred to as the set of \emph{dropped rays} of $\Gamma$, because it correspond to rays which occur in the fans associated to other chambers, but not in $\Sigma_\Gamma$. 

We then extend the exact sequence of \eqref{eqn:ses_for_lattices} to the following exact square which appears in the proof of Corollary 4.23 of \cite{FaveroKellyFCY}:
\begin{equation}\label{eqn:exact_square}\begin{tikzcd} 0 \arrow[r] & 0\arrow[r] \arrow[d] & \Z^{I_\Gamma^\emptyset} \ar[-,double line with arrow={-,-}]{r} \arrow[d] & \Z^{I_\Gamma^\emptyset}\arrow[r] \arrow[d] & 0 \\
0 \arrow[r]& M \arrow[r,"i"] \ar[-,double line with arrow={-,-}]{d} & \Z^n \arrow[r,"\pi"] \arrow[d,"\proj_\Gamma"]& \wh{G}\arrow[r] \arrow[d,"\proj_{\wh{G}_\Gamma}"]& 0 \\
0\arrow[r] & M \arrow[r,"i_\Gamma"] & \Z^{\Sigma_\Gamma(1)}\arrow[r,"\pi_\Gamma"] & \wh{G}_\Gamma\arrow[r] & 0
\end{tikzcd}
\end{equation}
where $\wh{G}_\Gamma$ is defined to be the cokernel of $\Z^{I_\Gamma^\emptyset}\to \wh{G}$.

\begin{definition}\label{def:stacky_fan_from_chamber}
Consider the stacky fan $\Sigma_{\Gamma}$ as in \eqref{eq:fan_from_chamber} and $\alpha:= i^{\vee}_{\Gamma}$ as in  \eqref{eqn:exact_square}. We let $\CX_{\Sigma_{\Gamma},i^{\vee}_{\Gamma}}$ denote the associated quotient stack.
\end{definition}

We now explain the stacky isomorphism between $\CX_{\Gamma}$ and $\CX_{\Sigma_{\Gamma},i^{\vee}_{\Gamma}}$. Under the natural isomorphism $\A^n \cong \A^{\Sigma_\Gamma(1)}\times \A^{I_\Gamma^\emptyset}$ one gets
\begin{equation*}
    \A^n\setminus V(B_\Gamma) \cong (\A^{\Sigma_\Gamma(1)}\setminus V(B_{\Sigma_\Gamma}))\times \Gm^{I_\Gamma^\emptyset},
\end{equation*}
which descends to an isomorphism of the associated global quotient stacks:
\begin{equation}\label{eqn:stack_iso}
[(\A^n\setminus V(B_\Gamma))/G] \simeq
[((\A^{\Sigma_\Gamma(1)}\setminus V(B_{\Sigma_\Gamma}))\times \Gm^{I_\Gamma^\emptyset})/G.] 
\end{equation}
Applying $\Hom(-,\Gm)$ to the right vertical exact sequence in \eqref{eqn:exact_square}, we set
\[G_\Gamma \coloneqq \ker (G \inc \Gm^n \twoheadrightarrow \Gm^{I_\Gamma^\emptyset}).\]
In particular, as $\wh{G}_\Gamma$ is also equal to the cokernel of $i_\Gamma$ by \eqref{eqn:exact_square}, it follows that 
\[G_\Gamma = G_{\Sigma_\Gamma, i_\Gamma^\vee}.\]
Combining the stack isomorphism \eqref{eqn:stack_iso} with Corollary 4.23 of \cite{FaveroKellyFCY} we get the following.
\begin{proposition}[Corollary 4.23 of \cite{FaveroKellyFCY}] \label{prop:stack_isomorphism} Let $\Gamma$ be as above. Then there is a natural isomorphism of stacks
\[\CX_\Gamma \simeq \CX_{\Sigma_\Gamma, i_\Gamma^\vee}.\]
\end{proposition}

All stacks considered in this paper are quotient stacks. The category of (quasi-)coherent sheaves on a quotient stack $[X/G]$ are naturally equivalent to $G$-equivariant (quasi-)coherent sheaves on $X$. We implicitly use this equivalence throughout this paper and write $\CO_{[X/G]}(\chi)$ for a $G$-equivariant line bundle on $X$.

\begin{corollary}\label{cor:image_of_line_bundles_under_stack_iso} Under the natural isomorphism of stacks of \Cref{prop:stack_isomorphism}, $G$-equivariant line bundles $\CO_{\CX_\Gamma}(\chi)$ for $\chi\in \wh{G}$ are sent to line bundles $\CO_{\CX_{\Sigma_\Gamma, i_\Gamma^\vee}}(\proj_{\wh{G}_\Gamma}(\chi))$ where $\proj_{\wh{G}_\Gamma}$ is as in \eqref{eqn:exact_square}.
\end{corollary}

\subsection{Nef divisors and polytopes}\label{subsection:nef_divisors} In this subsection we recall some facts about nef divisors and their polytopes. Letting $X_{\Sigma_\Gamma}$ denote the coarse moduli space of $\CX_{\Sigma_\Gamma,i_\Gamma^\vee}$, observe that $\Cl(X_{\Sigma_{\Gamma}})$ is a finite index subgroup of $\wh{G}_{\Gamma}$. In view of this, we make the following definition.
\begin{definition}\label{def:nefdivisorsonastack}
We call a character $d \in \wh{G}_{\Gamma}$ and the associated equivariant line bundle $\CO_{\CX_{\Sigma_{\Gamma},i^{\vee}_\Gamma}}(d)$ \emph{nef} if $d$ is a $\Q$-Cartier nef class of $X_{\Sigma_{\Gamma}}$. Equivalently, $d$ is nef if $d \in \text{Nef}(X_{\Sigma_{\Gamma}})$.
\end{definition}

\begin{proposition}\label{prop:every nefisisomorphictooneinthechamber}
For any nef line bundle $\CO_{\CX_{\Sigma_{\Gamma},i^{\vee}_\Gamma}}(d)$ there exists $d' \in \wh{G} \cap \Gamma$ such that $\proj_{\wh{G}_{\Gamma}}(d')=d$.
\end{proposition}

\begin{proof}

By \cite[Proposition 15.1.3]{CLS}, \[
\Gamma
\cong
\operatorname{Nef}(X_{\Sigma_{\Gamma}})\times \mathbb{R}_{\ge 0}^{I_{\emptyset}}.
\]
In particular, $(\proj_{\wh{G}_{\Gamma}})_{\R}$ restricts to a surjective map from $\Gamma$ to $\operatorname{Nef}(X_{\Sigma_{\Gamma}})$. Hence, $\proj_{\wh{G}_{\Gamma}}$ supplies a surjection from $\Gamma \cap \wh{G}$ to $\operatorname{Nef}(X_{\Sigma_{\Gamma}}) \cap \wh{G}_{\Gamma}$, which is precisely the set of all nef characters.
\end{proof}

\begin{definition}\label{def:divisor_polytope} Let $I\subset [n]$ and let $D = \sum_{\rho\in I}a_\rho e_\rho\in \Z^I$ be any formal sum. We define the polytope
\[P_D\coloneqq \{x\in M_\R \ | \ \la x, u_\rho\ra \ge -a_\rho\text{ for each }\rho \in I\}.\]
\end{definition}
This definition is meant to be a slight generalization of the polytope of sections that allows us to talk about $P_D$ without specifying an underlying stack. That is, when we take $D= \sum_{\rho\in [n]}a_\rho e_\rho\in \Z^n$, for any chamber $\Gamma$, $P_{\proj_\Gamma(D)}$ is the usual polytope of sections associated to the Weil divisor $\proj_\Gamma(D)$ on the toric variety $X_{\Sigma_\Gamma}$. We recall a well-known fact. 
\begin{proposition}\label{prop:nef_polytopes} Let $\Gamma$ be a chamber of a secondary fan $\Sigma_{GKZ}$ defined by a toric GIT problem. Then for any $d\in \Gamma\cap \wh{G}$ and any lift $D\in \Z^{n}$ of $d$,
\[P_D = P_{\proj_\Gamma(D)}.\]
\end{proposition}
\begin{proof} It follows from \cite[Lemma 14.4.6]{CLS} and \cite[Lemma 14.4.4(a)]{CLS} that for any lift $D$ of $d\in \Gamma$, the normal fan to $P_D$ is refined by $\Sigma_\Gamma$. In particular, one only needs to use the facets defined by $u_\rho$ for $\rho\in \Sigma_\Gamma(1)$ to recover $P_D$.
\end{proof}

The following elementary convex geometry lemma is proven in the course of \cite[Proposition 4.15]{CoxCategory}, but as we require this statement in a slightly more general setting later (see \Cref{prop:nefcohomologyvanishes}), we reprove it for completeness. This proof exactly mirrors the one appearing in \cite{CoxCategory}.
\begin{lemma}\label{lemma:lattice_points_of_shifted_polytope} Fix some $D=\sum_{\rho\in [n]}a_\rho e_\rho\in \Z^n$ and any chamber $\Gamma$ of $\Sigma_{GKZ}$. Then for any $\theta\in M_\R$,
\[(P_{\proj_\Gamma(D)}+ \theta)\cap M = P_{\proj_\Gamma(D-\ceil{i_\R(\theta)})}\cap M\]
where $i_\R$ is the inclusion in \eqref{eqn:ses}. It follows, for instance by taking $\Gamma$ to be any chamber in the moving cone, that
\[(P_{D}+ \theta)\cap M = P_{D-\ceil{i_\R(\theta)}}\cap M.\]
\end{lemma}
\begin{proof} In one direction, suppose $x\in (P_{\proj_\Gamma(D)} +\theta)\cap M$. Then, $x-\theta\in P_{\proj_\Gamma(D)}$, so for any $\rho\in \Sigma_\Gamma(1)$, $\la x-\theta, u_\rho\ra \ge -a_\rho$. Rearranging,
\[\la x,u_\rho\ra \ge -a_\rho + \la \theta,u_\rho\ra.\]
As $x\in M$ and $u_\rho\in N$, we have that $\la x,u_\rho\ra\in \Z$, so the above inequality implies
\[\la x,u_\rho\ra \ge \ceil{-a_\rho+\la \theta,u_\rho\ra} = -a_\rho+\ceil{\la \theta,u_\rho\ra}\]
Thus $x\in P_{\proj_\Gamma(D-\ceil{i_\R(\theta)})} \cap M$.

In the other direction, suppose $x\in P_{\proj_\Gamma(D-\ceil{i_\R(\theta)})}\cap M$, then for any $\rho\in \Sigma_{\Gamma}(1)$,
\begin{align*}
    \la x-\theta, u_\rho\ra &= \la x,u_\rho\ra- \la \theta,u_\rho\ra\\
    &\ge -a_\rho+\ceil{\la \theta,u_\rho\ra}-\la \theta,u_\rho\ra &\text{as }x\in P_{D+\ceil{i_\R(\theta)}}^\Gamma\\
    &\ge -a_\rho. & \text{as }\ceil{\la \theta, u_\rho\ra} \ge \la \theta,u_\rho\ra
\end{align*}
Thus $x-\theta\in P_{\proj_\Gamma(D)}$, or equivalently $x\in P_{\proj_\Gamma(D)}+\theta$.
\end{proof}

\subsection{1-parameter subgroups and walls of the secondary fan}\label{subsection:1PS_and_walls}
Suppose as above that we have an inclusion $G\subset \Gm^n$ which induces a projective GIT problem.

We have already seen that the cones $\beta_J$ dictate the combinatorics of the secondary fan $\Sigma_{GKZ}$. It turns out that the walls (codimension 1 cones) of the secondary fan are in particularly close correspondence with the cones $\beta_J$ of codimension 1 due to the following result of \cite{Schmitt}. We include a proof of the following statement for completeness.
\begin{lemma}[Theorem 3.2 of \cite{Schmitt}]\label{lemma:schmitt} The union of walls of the secondary fan equals the union of $\beta_J$ of codimension 1: 
\[\bigcup_{\codim \tau = 1}\tau = \bigcup_{\codim \beta_J = 1}\beta_J.\]
\end{lemma}
\begin{proof} We proceed by showing the complements in $|\Sigma_{GKZ}|$ of the two sides are equal. For one direction, suppose $\chi\in |\Sigma_{GKZ}|$ is such that if $\chi \in \beta_J$, then $\beta_J$ is full-dimensional. In particular, if $\chi\in \beta_J$, no face of $\beta_J$ contains $\chi$. So, $\chi$ is in the interior of each of the $\beta_J$ cones which contain it. It follows that the intersection of the $\beta_J$ cones which contain $\chi$ itself contains an open, and thus full-dimensional, set. By \Cref{cor:gamma_hits_beta_J_implies_contained}, this implies that the minimal cone of $\Sigma_{GKZ}$ containing $\chi$ is full-dimensional, showing the first containment.

To show the opposite containment, suppose $\chi$ is in the interior of a chamber $\Gamma$ of $\Sigma_{GKZ}$. If $\chi\in \beta_J$, then $\Gamma\subset \beta_J$ by \Cref{cor:gamma_hits_beta_J_implies_contained}. As $\Gamma$ is full-dimensional, this implies that every $\beta_J$ containing $\chi$ is full-dimensional, as desired.
\end{proof}
We now explain the relationship between walls and one-parameter subgroups.  To begin, there is an obvious relationship between one-parameter subgroups and $\beta_J$ for which $\codim \beta_J = 1$.  Namely, as $\beta_J$ is a codimension 1 rational polyhedral cone of $\wh{G}_\R$, its perpendicular $\beta_J^\perp$ is a dimension 1 subspace of $\wh{G}_\R^\ast$.  Thus there is a primitive ray generator $\lambda_J\in \Hom_\Z(\wh{G},\Z)$ so that $\beta_J^\perp = \Cone(\lambda_J,-\lambda_J)$. The relationship is given by $\spn \beta_J = H_{\pm \lambda_J}^{=0}$. That is, the span of $\beta_J$ is precisely the hyperplane obtained from pairing with the one-parameter subgroup $\lambda_J$.  We record this set of 1-PSs in the following.
\begin{definition} We define the following set of 1-PSs:
\[\Lambda \coloneqq \bigcup_{\codim \beta_J = 1}\{\lambda_J,-\lambda_J\}.\]
Later we will consider the subset of $\lambda\in \Lambda$ containing internal walls of $\Sigma_{GKZ}$, so we additionally define
\[\Lambda^{\circ} \coloneqq \Lambda\setminus \lp |\Sigma_{GKZ}|^\vee \cup -|\Sigma_{GKZ}|^\vee\rp\]
\end{definition}
The precise relationship between the walls $\tau \in  \Sigma_{GKZ}(\dim \wh{G}_\R-1)$ of the secondary fan and the set of 1-PSs $\Lambda$ is given by the following.
\begin{corollary}\label{cor:schmitt} \begin{enumerate}[(i),itemsep=10pt]
    \item\label{cor:schmitt_i} For any $\lambda\in \Lambda$, $\beta_{I_{\lambda}^{=0}}$ is a union of walls of $\Sigma_{GKZ}$. 
    \item\label{cor:schmitt_ii} For any wall $\tau\in \Sigma_{GKZ}(\dim \wh{G}_\R-1)$, there is a 1-PS $\lambda\in \Lambda$, unique up to sign, so that $\spn \tau = H_\lambda^{=0}$.
    \item\label{cor:schmitt_iii} If $\xi\in \wh{G}_\R^\ast$ and $\codim \beta_{I_\xi^{=0}} = 1$, then there is a unique $\lambda\in \Lambda$ so that $\Cone(\xi) = \Cone(\lambda)$.
    \item\label{cor:schmitt_iv} If $J\subset [n]$ is such that $\codim \beta_{J} = 1$ and $0 \ne \xi\in\wh{G}_\R^\ast$ is such that $J\subset I_\xi^{=0}$, then there is a unique $\lambda\in \Lambda$ so that $\xi$ is a multiple of $\lambda$.
    \item\label{cor:schmitt_v} If $\lambda\in \Lambda$ and $0\ne \xi\in \wh{G}_\R^\ast$ are such that $I_{\lambda}^{=0}\subset I_\xi^{=0}$, then $\xi$ is a multiple of $\lambda$. 
\end{enumerate}
\end{corollary}
\begin{proof} Statements \ref{cor:schmitt_i}, \ref{cor:schmitt_ii} are immediate from \Cref{lemma:schmitt} and the construction of $\Lambda$. Statements \ref{cor:schmitt_iii}, \ref{cor:schmitt_iv}, \ref{cor:schmitt_v} follow from \Cref{lemma:schmitt}.
\end{proof}
Throughout, we will be more concerned with the 1-PSs in $\Lambda$ which define walls of the secondary fan than the walls themselves. By \Cref{cor:schmitt}, we view these 1-PSs which define walls interchangeably with linear functionals $\wh{\lambda}_\R:\wh{G}_\R\to \R$ whose associated hyperplanes contain enough weights $\beta_i$ to span the hyperplane. \Cref{cor:schmitt} details the variety of ways in which we will describe $\lambda\in \Lambda$.

It is a well-known fact that for any rational polyhedral cone $\sigma$, $(\sigma^\vee)^\vee = \sigma$ (see e.g. \cite[Proposition 1.2.4]{CLS}). Unpacking the meaning of the dual cone, this is $\sigma = \bigcap_{f\in \sigma^\vee}H_f^{\ge 0}$. It turns out that for cones of the type $\beta_J$ for $J\subset [n]$, we can restrict this description from an intersection over all functionals in $\beta_J^\vee$ to one over the finite subset of functionals which are also in $\Lambda$. This is made precise in the following lemma.
\begin{lemma}\label{lemma:beta_cones_determined_by_lambdas} If $\codim |\Sigma_{GKZ}| =0$, then for any $J\subset [n]$,
\[\beta_J = \bigcap_{\lambda\in \Lambda\cap \beta_J^\vee}H_{\lambda}^{\ge 0}.\]
\end{lemma}

\begin{proof} That $\beta_J$ is contained in this intersection follows from the fact that this is an intersection over $H_{\lambda}^{\ge 0}$ for $\lambda \in \beta_J^\vee$. To show the opposite inclusion, first suppose that $\beta_J$ is full-dimensional. Let $\beta_{J'}$ be an arbitrary facet of $\beta_J$ with associated linear functional $f: \wh{G}_\R\to \R$. By construction $J'\subset I_f^{=0}$, so $\beta_{J'}\subset \beta_{I_f^{=0}}$. Then, as $\beta_{I_f^{=0}}$ has codimension at least 1 and $\beta_{J'}$ has codimension 1, we get that $\beta_{I_f^{=0}}$ has codimension 1. By \Cref{cor:schmitt}\ref{cor:schmitt_iii} we get that there is a $\lambda\in \Lambda$ so that $\Cone(f) = \Cone(\lambda)$. In particular $H_f^{\ge 0} = H_{\lambda}^{\ge 0}$, so we conclude that the right hand side is contained in $\beta_J$. 

Now suppose $\codim \beta_J = \ell \ge 1$. Then, as $|\Sigma_{GKZ}| = \beta_{[n]}$ is full-dimensional, we can choose $\rho_1$ so that $\beta_{\rho_1}\notin \spn \beta_J$. In particular, we find that $\beta_J$ is a face of $\beta_{J\cup \{\rho_1\}}$ and $\dim \beta_{J\cup \{\rho_1\}} = \dim \beta_J+1$. Inductively then, we can choose $\rho_1,\hdots, \rho_\ell$ so that for $J' \coloneqq J \cup \{\rho_1,\hdots, \rho_\ell\}$, $\beta_{J'}$ is full-dimensional and $\beta_J$ is a face of $\beta_{J'}$. 

As $\beta_J$ is a face of $\beta_{J'}$, $\beta_J$ is equal to the intersection of a finite number of facets of $\beta_{J'}$. Because $\beta_{J'}$ is full-dimensional, we have already proven that the facets of $\beta_{J'}$ are given by functionals in $\Lambda$. That is to say, there is a subset $C\subset \Lambda\cap \beta_{J}^\vee$ such that
\begin{align*}
    \beta_J &= \beta_{J'}\cap \bigcap_{\lambda\in C}H_{\lambda}^{=0} \supset \beta_{J'}\cap \bigcap_{\lambda\in \Lambda \cap \beta_{J}^\vee} H_{\lambda}^{\ge 0} = \bigcap_{\xi \in \Lambda\cap \beta_{J'}^\vee}H_\xi^{\ge 0}\cap \bigcap_{\lambda\in \Lambda\cap \beta_J^\vee} H_\lambda^{\ge 0} \supset \bigcap_{\lambda\in \Lambda\cap \beta_J^\vee}H_\lambda^{\ge 0}.
\end{align*}
The first containment follows from the fact that for any $\lambda\in C$, $- \lambda\in \Lambda\cap \beta_J^\vee$ and $H_\lambda^{=0} = H_{\lambda}^{\ge 0}\cap H_{-\lambda}^{\ge 0}$. The last inclusion follows from $\beta_{J'}^\vee \subset \beta_J^\vee$. Having shown the opposite containment, we are done.
\end{proof}
\section{HKKN Stratifications}\label{section:stratifications}
In this section we first review the notion of a Hesselink-Kampen-Kempf-Ness (HKKN) stratification and describe the well-known Hilbert-Mumford algorithm for constructing such stratifications. Next, we introduce the notion of arrangeable GIT problems and prove that such problems admit a new form of particularly nice HKKN stratifications which are not captured by the Hilbert-Mumford procedure. As a note, in this section, we do not assume that our GIT problems are primitive geometric as is required in the following sections.

\subsection{A review of HKKN stratifications} In this subsection, we provide a review of HKKN stratifications. Throughout this subsection, we will closely follow the descriptions from \S2.1 of \cite{HL} and \S2.1 of \cite{BFK}.

Let $G\subset \Gm^n$ be an algebraic torus acting on $\A^n$ induced by the natural action of $\Gm^n$ on $\A^n$. Let $X\subseteq \A^n$ be a $G$-equivariant open subvariety so that $[X/G]$ is a smooth toric DM stack. 
\begin{definition} Let $\lambda:\Gm \to G$ be a one-parameter subgroup. We define the \emph{fixed locus} of $X$ with respect to $\lambda$ to be
\[X_{\lambda}^{=0} \coloneqq \{x\in X \ | \ \lambda(t)\cdot x = x,\ \forall t\in \Gm\}.\]
We define the \emph{contracting locus} of $X$ with respect to $\lambda$ to be
\[X_{\lambda}^{\le 0} \coloneqq \{x\in X \ | \ \lim_{t\to \infty}\lambda(t)\cdot x\in X_{\lambda}^{=0}\}.\]
We also define the following coordinate projection map:
\begin{align*}
    \pi_\lambda: (\A^n)_\lambda^{\le 0} &\to (\A^n)_\lambda^{=0}\\
    x&\mapsto \lim_{t\to\infty}\lambda(t)\cdot x.
\end{align*}
\end{definition}

\begin{definition} Let $X\subset \A^n$ be a $G$-invariant open subvariety. We call a closed, nonempty subvariety $S\subset X$ a \emph{Hesselink-Kampen-Kempf-Ness} or \emph{HKKN stratum} if there is a one-parameter subgroup $\lambda:\Gm \to G$ so that $S = X_\lambda^{\le 0}$.

Let $Y\subset \A^n$ be a closed $G$-equivariant subvariety. We call a collection of locally-closed pairwise-disjoint subvarieties $S_1,\hdots, S_k\subset \A^n$ an \emph{HKKN stratification} of $Y$ if $Y = \bigsqcup_{i=1}^kS_i$
and $S_i\subset \A^n\setminus \sqcup_{j < i}S_j$ is a closed HKKN stratum for each $i$.
\end{definition}

\begin{remark}\label{rmk:1PS_sign_difference} Before continuing we should remark that we take a different sign convention for our stratifications than those which appear in our primary references \cite{HL, BFK}. Namely, we define the contracting locus as the $x$ so that $\lambda(t)\cdot x\in X_\lambda^{=0}$ as $t\to\infty$, whereas the contracting locus in the literature is generally defined instead with the limit taken as $t\to 0$. The effect of this is that our stratifications differ by a sign from the existing literature. Namely, if $\lambda_1,\hdots, \lambda_k$ induces an HKKN stratification of $Y$ in our convention, then $-\lambda_1,\hdots, -\lambda_k$ would induce an HKKN stratification of $Y$ in the convention of \cite{HL,BFK}.
\end{remark}

If $S\subset X$ is an HKKN stratum witnessed by a one-parameter subgroup $\lambda$, then by assumption, $S = X_\lambda^{\le 0}$ is closed in $X$. However, as $X\subset \A^n$ is an open subvariety, $\A^n\setminus X$ is of codimension at least 1. This implies
\[X_{\lambda}^{\le 0} = (\A^n)_{\lambda}^{\le 0}\cap X.\]
The assumption that $S$ be nonempty implies that $X_{\lambda}^{=0}$ is nonempty. In fact, $\lambda$ induces an HKKN stratum $X_\lambda^{\le 0}$ of $X$ if and only if $X_{\lambda}^{=0}$ is nonempty and $X_\lambda^{\le 0} = (\A^n)_\lambda^{\le 0}\cap X$. This observation directly implies the following lemma.
\begin{lemma}\label{lemma:HKKN_strata_determined_by_lambda} 
The data of an HKKN stratification of a closed $G$-equivariant subvariety $Y\subset \A^n$  is equivalent to an ordered choice of 1-PSs $\lambda_1,\hdots, \lambda_k$ such that the following three properties hold:
\begin{enumerate}[(a),itemsep=10pt]
    \item\label{item:HKKN_condition_1} $\ds Y = \bigcup_{i=1}^k (\A^n)_{\lambda_i}^{\le 0}$.
    \item\label{item:HKKN_condition_2} $\ds (\A^n)_{\lambda_i}^{=0}\setminus \bigcup_{j < i}(\A^n)_{\lambda_j}^{\le 0} \ne \emptyset$ for all $2 \le i \le k$.
    \item\label{item:HKKN_condition_3} $\ds \lp \A^n\setminus \bigcup_{j < i}(\A^n)_{\lambda_j}^{\le 0}\rp_{\lambda_i}^{\le 0} = (\A^n)_{\lambda_i}^{\le 0}\setminus \bigcup_{j < i}(\A^n)_{\lambda_j}^{\le 0}$ for all $2 \le i \le k$.
\end{enumerate}
\end{lemma}
In light of the above, when $\lambda_1,\hdots, \lambda_k$ satisfy the conditions of \Cref{lemma:HKKN_strata_determined_by_lambda}, we will say that they \emph{induce} an HKKN stratification, or simply refer to $\lambda_1,\hdots, \lambda_k$ as an HKKN stratification of $Y$. When we have a fixed HKKN stratification $\lambda_1,\hdots, \lambda_k$ we define
\[X_i\coloneqq \A^n\setminus \bigcup_{j<i}(\A^n)_{\lambda_j}^{\le 0}.\]
The conditions in \Cref{lemma:HKKN_strata_determined_by_lambda} are somewhat cumbersome but can be combinatorialized nicely. 
Observe that
\begin{align*}
    (\A^n)_\lambda^{=0} &= V(I_\lambda^{\ne 0}),\\
    (\A^n)_\lambda^{\le 0} &= V(I_{\lambda}^{>0}).
\end{align*}
Thus we can recharacterize when $\lambda_1,\hdots,\lambda_k$ introduce a stratification of the unstable locus $Y = (\A^n)^{\rm us}(\chi)$ as follows.
\begin{lemma}\label{lemma:HKKN_strat_combo} Let $G$ be as above and let $\chi\in \wh{G}$. Then $\lambda_1,\hdots, \lambda_k$ induce an HKKN stratification of $(\A^n)^{\us}(\chi)$ if and only if the following three properties hold.\vspace{0.5cm}
\begin{enumerate}[(i), itemsep = 0.5cm]
    \item\label{lemma:HKKN_strat_combo_i} $B_\chi = \bigcap_{i=1}^k \la I_{\lambda_i}^{>0}\ra$.
    \item\label{lemma:HKKN_strat_combo_ii} For all $j < i$, $I_{\lambda_j}^{>0}\not\subset I_{\lambda_i}^{\ne 0}$ or equivalently $I_{\lambda_i}^{=0}\not\subset I_{\lambda_j}^{\le 0} $.
    \item\label{lemma:HKKN_strat_combo_iii} For any subset $K\subset I_{\lambda_i}^{\le 0}$, if for every $j < i$, $I_{\lambda_j}^{>0} \not\subset I_{\lambda_i}^{>0}\sqcup K$, then $I_{\lambda_j}^{>0}\not\subset I_{\lambda_i}^{\ne 0}\cup K$ for every $j < i$.
\end{enumerate}
\end{lemma}
\begin{proof} It suffices to show that properties \ref{lemma:HKKN_strat_combo_i}, \ref{lemma:HKKN_strat_combo_ii}, \ref{lemma:HKKN_strat_combo_iii} are equivalent respectively to properties \ref{item:HKKN_condition_1}, \ref{item:HKKN_condition_2}, \ref{item:HKKN_condition_3} of \Cref{lemma:HKKN_strata_determined_by_lambda}. So see that \ref{lemma:HKKN_strat_combo_i} and \ref{item:HKKN_condition_1} are equivalent, first note that \ref{item:HKKN_condition_1} is equivalent to
\[V(B_{\chi}) = \bigcup_{i=1}^kV(I_{\lambda_i}^{>0}) = V\lp \bigcap_{i=1}^k\la I_{\lambda_i}^{>0}\ra\rp.\]
As $B_{\chi}$ and $\bigcap_{i=1}^k\la I_{\lambda_i}^{>0}\ra$ are radical, we get that the above is equivalent to condition \ref{lemma:HKKN_strat_combo_i}.

To see that \ref{lemma:HKKN_strat_combo_ii} and \ref{item:HKKN_condition_2} are equivalent, first see that \ref{item:HKKN_condition_2} is equivalent to
\[V(I_{\lambda_i}^{\ne 0})\setminus \bigcup_{j < i}V(I_{\lambda_j}^{>0})\ne \emptyset.\]
which is equivalent to saying that for each $j < i$, $V(I_{\lambda_i}^{\ne 0})\not\subset V(I_{\lambda_j}^{>0})$. Again, by the ideal-variety correspondence, this is equivalent to $I_{\lambda_j}^{>0}\not\subset I_{\lambda_i}^{\ne 0}$ for each $j<i$ which is exactly condition \ref{lemma:HKKN_strat_combo_ii}.

To see that \ref{lemma:HKKN_strat_combo_iii} and \ref{item:HKKN_condition_3} are equivalent, first see that \ref{item:HKKN_condition_3} is equivalent to
\[\pi_{\lambda_i}^{-1} \lp V(I_{\lambda_i}^{\ne 0})\setminus \bigcup_{j<i}V(I_{\lambda_j}^{>0})\rp = V(I_{\lambda_i}^{>0})\setminus \bigcup_{j<i}V(I_{\lambda_j}^{>0})\]
The left hand side is always a subset of the right, so the above is equivalent to the right side being contained in the left. Thus \ref{item:HKKN_condition_3} holds if and only if for any point $x\in V(I_{\lambda_i}^{>0})\setminus \bigcup_{j<i}V(I_{\lambda_j}^{>0})$, $\pi_{\lambda_i}(x)\in V(I_{\lambda_i}^{\ne 0})\setminus \bigcup_{j < i}V(I_{\lambda_j}^{>0})$. For a point $x\in \A^n$ let \[K(x) \coloneqq \{i\in [n] \ : \ x_i=0\}.\] 
Then $x\in V(I_{\lambda_i}^{>0})\setminus \bigcup_{j<i}V(I_{\lambda_j}^{>0})$ if and only if $I_{\lambda_i}^{>0}\subset K(x)$ and for every $j < i$, $I_{\lambda_j}^{>0}\not\subset K(x)$. Similarly, $K(\pi_{\lambda_i}(x)) = K(x)\cup I_{\lambda_i}^{\ne 0}$ so $\pi_{\lambda_i}(x)\in V(I_{\lambda_i}^{\ne 0})\setminus \bigcup_{j < i}V(I_{\lambda_j}^{>0})$ if and only if for every $j <i$, $I_{\lambda_j}^{>0}\not\subset I_{\lambda_i}^{\ne 0}\cup K(x)$. Thus \ref{lemma:HKKN_strat_combo_iii} implies \ref{item:HKKN_condition_3} by considering $K = K(x)\setminus I_{\lambda_i}^{>0}$. 

Similarly, to show \ref{item:HKKN_condition_3} implies \ref{lemma:HKKN_strat_combo_iii}, consider a subset $K\subset I_{\lambda_i}^{\le 0}$ such that for all $j < i$, $I_{\lambda_j}^{>0}\not\subset I_{\lambda_j}^{>0}\sqcup K$. Define $x\in \A^n$ via
\[x_i = \begin{cases} 0 & i\in K\sqcup I_{\lambda_i}^{>0} \\
1 & \text{else}
\end{cases}\]
Then, $K(x) = K\sqcup I_{\lambda_i}^{>0}$ so $I_{\lambda_j}^{>0}\not\subset I_{\lambda_i}^{>0}\sqcup K$ for $j < i$. Thus by our equivalent characterization of \ref{item:HKKN_condition_3}, we have that for all $j < i$, $I_{\lambda_j}^{>0}\not\subset I_{\lambda_i}^{\ne 0}\cup K(x) = I_{\lambda_i}^{\ne 0}\cup K$. Thus \ref{item:HKKN_condition_3} implies \ref{lemma:HKKN_strat_combo_iii}.
\end{proof}

\begin{example}\label{example:H2_stratifications} Let $G=\Gm^2$ be embedded in $\Gm^4$ with the dual map $\pi:\Z^4 \to \Z^2$ given by
\[\pi = \begin{pmatrix} 1 & -2 & 1 & 0 \\ 0 & 1 & 0 & 1
\end{pmatrix}\]
The associated GIT problem has the following secondary fan:
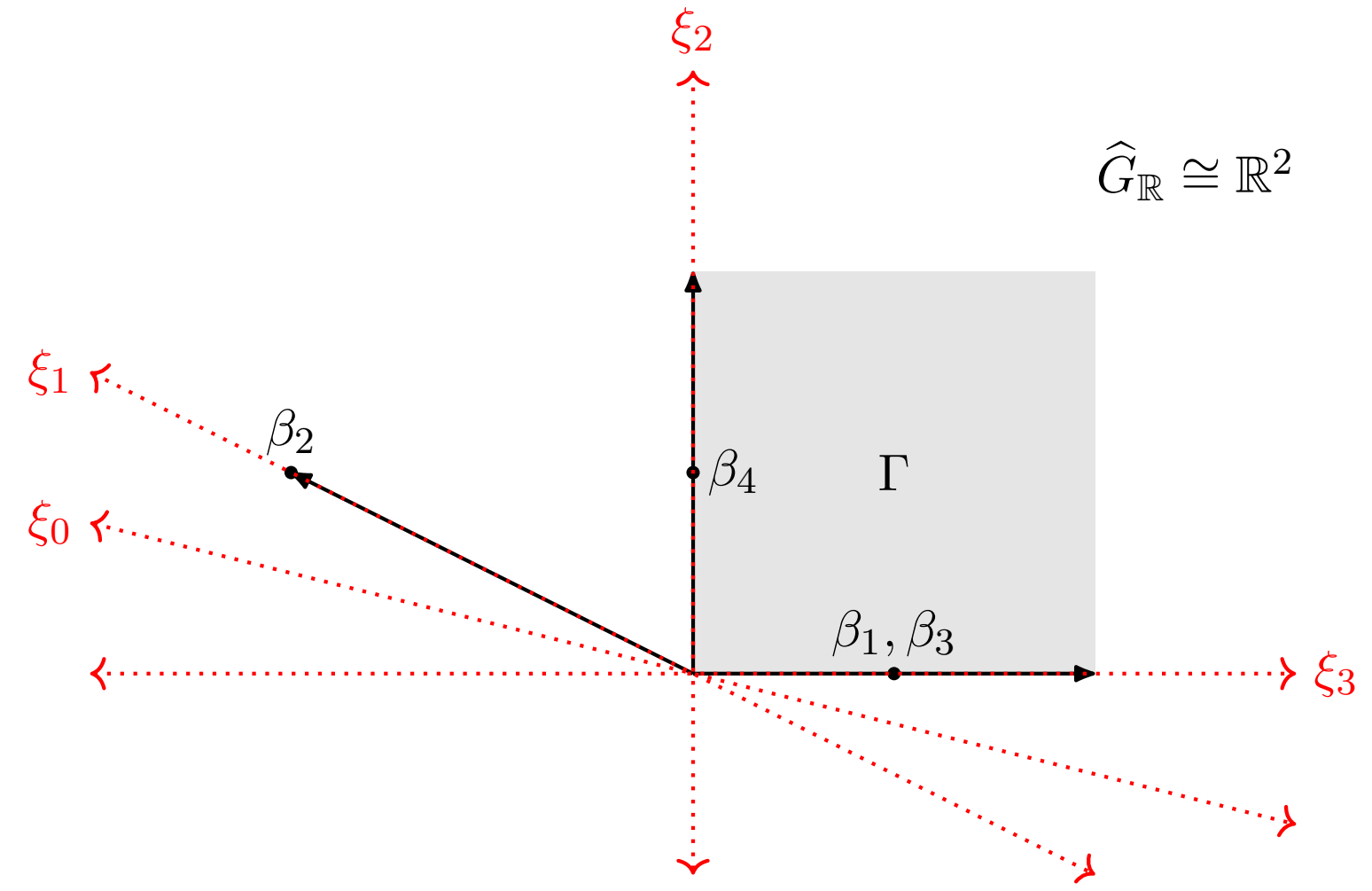
\[\begin{tikzpicture}
\def\fanscale{1.5}
\def\boxscale{1}
\def\xshift{0}
\def\yshift{0}
\def\xmin{-3};
\def\xmax{3};
\def\ymin{-1};
\def\ymax{3};

\coordinate (origin) at (\xshift,\yshift);
\coordinate (beta13) at (\xshift+\fanscale*1,\yshift);
\coordinate (beta2) at (\xshift+\fanscale*-2,\yshift+\fanscale*1);
\coordinate (beta4) at (\xshift, \yshift+1*\fanscale*1);
\coordinate (ray13) at (\xshift+2*\fanscale,\yshift);
\coordinate (ray2) at (\xshift-2*\fanscale,\yshift+1*\fanscale);
\coordinate (ray4) at (\xshift+0*\fanscale,\yshift+2*\fanscale);
\coordinate (chamberur) at (\xshift+2*\fanscale, \yshift+2*\fanscale);
\coordinate (boxll) at (\xshift+\xmin*\fanscale, \yshift+\ymin*\fanscale);
\coordinate (boxul) at (\xshift+\xmin*\fanscale, \yshift+\ymax*\fanscale);
\coordinate (boxlr) at (\xshift+\xmax*\fanscale, \yshift+\ymin*\fanscale);
\coordinate (boxur) at (\xshift+\xmax*\fanscale, \yshift+\ymax*\fanscale);
\def\xiOneSlope{-0.5}
\coordinate (xi11) at (\xshift+\xmin*\fanscale, \yshift+\xiOneSlope*\xmin*\fanscale);
\coordinate (xi12) at (\xshift+\fanscale*\ymin/\xiOneSlope,\yshift+\ymin*\fanscale);
\coordinate (xi21) at (\xshift, \yshift+\ymax*\fanscale);
\coordinate (xi22) at (\xshift, \yshift+\ymin*\fanscale);
\coordinate (xi31) at (\xshift+\xmin*\fanscale, \yshift);
\coordinate (xi32) at (\xshift+\xmax*\fanscale, \yshift);
\def\xiZeroSlope{-0.25}
\coordinate (xi01) at (\xshift+\xmin*\fanscale, \yshift+\xiZeroSlope*\xmin*\fanscale);
\coordinate (xi02) at (\xshift+\xmax*\fanscale, \yshift+\xiZeroSlope*\xmax*\fanscale);

% \draw[fill=none, draw=black] (boxll) -- (boxul) -- (boxur) -- (boxlr) -- cycle;
% chamber shading
\draw[draw = none, fill = gray, opacity= 0.2] (origin) -- (ray13) -- (chamberur) -- (ray4) -- cycle;

% betas
\node[circle,fill,inner sep=1pt] at (beta13) {};
\node[circle,fill,inner sep=1pt] at (beta2) {};
\node[circle,fill,inner sep=1pt] at (beta4) {};
\node[draw = none, fill = none] (betas13label) at ($(beta13)+(0,0.3)$) {\small $\beta_1,\beta_3$};
\node[draw = none, fill = none] (beta3label) at ($(beta4)+(0.3,0)$) {\small $\beta_4$};
\node[draw = none, fill = none] (beta1label) at ($(beta2)+(0,0.3)$) {\small $\beta_2$};

% chamber label
\node[draw = none, fill = none] at (\xshift+1*\fanscale,\yshift+1*\fanscale) {\small $\Gamma$}; 

% rays
\draw[-arr,black,thick] (origin) -- (ray13);
\draw[-arr,black,thick] (origin) -- (ray2);
\draw[-arr,black,thick] (origin) -- (ray4);

% 1-PSs
\draw[red,thick,dotted,<->] (xi11) -- (xi12);
\draw[red,thick,dotted,<->] (xi21) -- (xi22);
\draw[red,thick,dotted,<->] (xi31) -- (xi32);
\draw[red,thick,dotted,<->] (xi01) -- (xi02);

% 1-PS labels
\node[draw=none,fill=none,red] at ($(xi01)+(-0.3,0)$) {\small $\xi_0$};
\node[draw=none,fill=none,red] at ($(xi11)+(-0.3,0)$) {\small $\xi_1$};
\node[draw=none,fill=none,red] at ($(xi21)+(0,0.3)$) {\small $\xi_2$};
\node[draw=none,fill=none,red] at ($(xi32)+(0.3,0)$) {\small $\xi_3$};

\node[draw=none,fill=none,black] at (\xshift+\fanscale*2.5,\yshift+\fanscale*2.5) {\small $\wh{G}_\R\cong \R^2$};
\end{tikzpicture}\]
We have also labeled four hyperplanes $H_{\xi_i}^{=0}$ by their associated 1-PSs $\xi_0,\hdots, \xi_3$.  These 1-PSs are given a choice of sign so that $\Gamma \subset H_{\xi_i}^{\ge 0}$. By \Cref{lemma:B_chi_via_betas} one can compute that $B_{\Gamma} = \la x_1,x_3\ra\cap \la x_2,x_4\ra$. It follows that $\CX_\Gamma = \CH_2$, the Hirzebruch surface of type 2. We will use this GIT problem and chamber as a running example throughout the remainder of the paper.

We will now illustrate how to use \Cref{lemma:HKKN_strat_combo} to show that the ordering $\xi_1,\xi_2,\xi_3$ induces an HKKN stratification of $V(B_\Gamma)$. Explicitly, condition \ref{lemma:HKKN_strat_combo_i} is satisfied because
\[\la I_{\xi_1}^{>0}\ra \cap \la I_{\xi_2}^{>0}\ra \cap\la I_{\xi_3}^{>0}\ra = \la x_1,x_3,x_4\ra \cap \la x_1,x_3\ra \cap \la x_2,x_4\ra = \la x_1,x_3\ra\cap \la x_2,x_4\ra = B_\Gamma.\]
To check condition \ref{lemma:HKKN_strat_combo_ii}, see that
\begin{align*}
    I_{\xi_2}^{=0} &= \{4\}\not\subset \{2\} = I_{\xi_1}^{\le 0}\\
    I_{\xi_3}^{=0} &= \{1,3\}\not\subset \{2\} = I_{\xi_1}^{\le 0}\\
    I_{\xi_3}^{=0} &= \{1,3\}\not\subset \{2,4\} = I_{\xi_2}^{\le 0}
\end{align*}
To check condition \ref{lemma:HKKN_strat_combo_iii} for $\xi_2$, first note that the only subsets $K$ of $I_{\xi_2}^{\le 0}$ for which $I_{\xi_1}^{>0}\not\subset I_{\xi_2}^{>0}\sqcup K$ are $K= \emptyset, \{2\}$. Then, $I_{\xi_1}^{>0}\not\subset I_{\xi_2}^{\ne 0}\cup K$ for both of these choices as shown above in proving that condition \ref{lemma:HKKN_strat_combo_ii} held. To check condition \ref{lemma:HKKN_strat_combo_iii} for $\xi_3$, first note that $I_{\xi_2}^{>0}\subset I_{\xi_1}^{>0}$. Then, the only subsets $K$ of $I_{\xi_3}^{\le 0} = I_{\xi_3}^{=0}$ such that $I_{\xi_2}^{>0}\not\subset I_{\xi_3}^{>0}\sqcup K$ are $K= \emptyset, \{1\},\{3\}$. However,
\begin{align*}
    \{1,3\} = I_{\xi_2}^{>0}&\not\subset I_{\xi_3}^{\ne 0}\cup\emptyset =\{2,4\}\\
    \{1,3\} = I_{\xi_2}^{>0}&\not\subset I_{\xi_3}^{\ne 0}\cup\{1\} =\{1,2,4\}\\
    \{1,3\} = I_{\xi_2}^{>0}&\not\subset I_{\xi_3}^{\ne 0}\cup\{3\} =\{2,3,4\}
\end{align*}
The above non-containments for $I_{\xi_2}^{>0}$ imply the corresponding non-containments for $I_{\xi_1}^{>0}$ as $I_{\xi_1}^{>0}\supset I_{\xi_2}^{>0}$. So, condition \ref{lemma:HKKN_strat_combo_iii} is satisfied, so we conclude that $\xi_1,\xi_2,\xi_3$ is an HKKN stratification of $V(B_\Gamma)$.

By analyzing the conditions of \Cref{lemma:HKKN_strat_combo} in the same way as above, one can show that the following orderings of the 1-PSs $\xi_i$ induce HKKN stratifications of $V(B_\Gamma)$:
\begin{itemize}
    \item $\xi_0,\xi_1, \xi_2,\xi_3$
    \item $\xi_1,\xi_3,\xi_2$
    \item $\xi_3,\xi_1,\xi_2$
    \item $\xi_2,\xi_3$
\end{itemize}
The above list is certainly non-exhaustive. For instance, replacing $\xi_0$ with any other 1-PS $\xi$ such that $|\Sigma_{GKZ}|\subset H_\xi^{>0}$ yields another HKKN stratification.

It is also useful to note which orderings cannot occur by analyzing the violation of condition \ref{lemma:HKKN_strat_combo_ii}. Roughly speaking, condition \ref{lemma:HKKN_strat_combo_ii} enforces that in an HKKN stratification, the hyperplanes associated to the ordered 1-PSs should not get ``further" from the chamber going along in the stratification. For example, see that
\[I_{\xi_1}^{=0}= \{2\}\subset \{2,4\} = I_{\xi_2}^{\le 0}.\]
So, in any HKKN stratification involving both $\xi_1$ and $\xi_2$, condition \ref{lemma:HKKN_strat_combo_ii} implies that $\xi_1$ must precede $\xi_2$. Similarly, because $I_{\xi_0}^{=0} = \emptyset$, any HKKN stratification including $\xi_0$ must begin with $\xi_0$.
\end{example}

\subsection{Hilbert-Mumford HKKN stratifications}\label{subsection:hilbert_mumford}
In this subsection we outline, without proof, an HKKN stratification which exists for any projective toric variety which we call the \textit{Hilbert-Mumford HKKN stratification}. We will closely follow the exposition of this procedure as detailed in \cite[\S2.1]{HL}. For a complete treatment with proofs, see \cite{dolgachev1998variation}.

As in the previous subsection, let $G\subset \Gm^n$ act on $\A^n$ via the natural action of $\Gm^n$ on $\A^n$. Let $X\subset \A^n$ be an open $G$-equivariant subvariety so that $[X/G]$ is a projective toric DM stack. Fix a linearization $\CL = \CO_X(1)$, viewing $\CL\in \wh{G}$. In particular, $X = (\A^n)^{\ss}(\CL)$.

Fix a conjugation-invariant norm $|\cdot|$ on the one-parameter subgroups, interpreting $|\cdot|$ as a distance on the vector space $\wh{G}_\R^\ast$. With this norm in hand we define the numerical invariant
\[\mu_\CL(\lambda) \coloneqq \frac{\wh{\lambda}(\CL)}{|\lambda|}.\]
Viewing $\CL\in \wh{G}_\R$, $\mu_\CL(\lambda)$ measures the distance of $\CL$ from the hyperplane $H_\lambda^{=0}$. 

We inductively define 1-PSs $\lambda_i$ and strata $S_i$ as follows. First choose $\lambda_1$ to be the 1-PS which maximizes $\mu_\CL$ and define $S_1 \coloneqq V(I_{\lambda_1}^{>0})$. Given a list of 1-PSs $\lambda_1,\hdots, \lambda_{i}$, define
\[S_i \coloneqq V(I_{\lambda_i}^{>0})\setminus \bigcup_{j < i}S_j.\]
We choose $\lambda_{i+1}$ to be a 1-PS $\lambda$ which maximizes $\mu_L$ among the $\lambda$ such that $V(I_{\lambda}^{\ne 0})$ is not contained in the union of previous strata: $\bigcup_{j \le i}S_j$. Terminate this process when $\bigcup_{j \le k}S_j = \A^n\setminus X$. By work of Dolgachev and Hu, this process terminates and $\lambda_1, \hdots, \lambda_k$ induce an HKKN stratification of $X$ (see \cite{dolgachev1998variation}). We call this stratification $\lambda_1,\hdots, \lambda_k$ the \textit{Hilbert-Mumford HKKN stratification} of $\A^n\setminus X$ associated to $\CL, |\cdot|$.

This construction offers a high degree of flexibility: varying the norm and linearization generally yields distinct Hilbert-Mumford HKKN stratifications. See \cite{VariationOfStratifications} for a complete combinatorial treatment of this flexibility when fixing the norm and varying the linearization. 

\begin{example}\label{example:H2_HM_stratifications}
We continue \Cref{example:H2_stratifications} by describing a Hilbert-Mumford HKKN stratification of $V(B_\Gamma)$. Let $\xi_i$ be as in \Cref{example:H2_stratifications}. Choose $\CL = (1,4)^T\in \wh{G}$ to be our ample linearization and fix the standard Euclidean norm on $\wh{G}_\R^\ast$ dual to our chosen identification $\wh{G}_\R\cong \R^2$. Then, $\lambda_1$ is the 1-PS which maximizes $\mu_\CL$, which in this case is $\lambda_1 = (4,1) = \xi_0$. Note that $I_{\lambda_1}^{>0} = [4]$, so in order for $\lambda_i$ for $i>1$ to be such that $V(I_{\lambda_i}^{\ne 0})$ is not contained in $S_1 = V(I_{\lambda_1}^{>0}) = \{0\}$, we must have that $I_{\lambda_i}^{=0}\ne\emptyset$. This, with the fact that $\wh{\lambda}_i(\CL) > 0$ for each $i > 1$, implies that $\lambda_i\in \{\xi_1,\xi_2,\xi_3\}$ for each $i >1$. 

Noting that $\xi_1 = (1,2)$, $\xi_2 = (1,0)$, and $\xi_3 = (0,1)$, one can directly compute that $\mu_\CL(\xi_3) > \mu_\CL(\xi_2) > \mu_\CL(\xi_3)$. After verifying that the fixed loci are not contained in previously defined strata, one can check that setting $\lambda_2 = \xi_1$, $\lambda_3 = \xi_3$, and $\lambda_4 = \xi_2$ yields a Hilbert-Mumford HKKN stratification $\lambda_1,\hdots, \lambda_4$ of $V(B_\Gamma)$. It is also worth verifying that after selecting $\lambda_4$, all 1-PSs $\xi$ such that $\wh{\xi}(\CL) > 0$ and $I_\xi^{=0}\ne \emptyset$ have been exhausted, thus the Hilbert-Mumford procedure indeed terminates after selecting $\lambda_4$.
\end{example}
\subsection{Arrangeable GIT problems and wall stratifications}\label{subsection:wall_stratifications}
In this subsection we will introduce a class of GIT problems which admit nice HKKN stratifications which are the focus of our main theorem. Throughout this subsection we let $G\subset \Gm^n$ define a GIT problem with secondary fan $\Sigma_{GKZ}$ with set of associated wall 1-PSs $\Lambda$. To make sense of this class of problems, we first need notation to collect the sets of indices $\rho$ whose associated weights are scalar multiples of each other.
\begin{definition}\label{def:beta_class_definition} Let $G\subset \Gm^n$ define a toric GIT problem. Then, for any $\rho\in [n]$, let
\[[\rho] \coloneqq \{\rho'\in [n] \ | \ \beta_{\{\rho'\}} = \beta_{\{\rho\}}\}.\]
\end{definition}

We are now in a position to describe the main class of GIT problems of interest to us.
\begin{definition}\label{def:arrangeable}
Let $G\subset \Gm^n$ define a projective GIT problem with secondary fan $\Sigma_{GKZ}$ with weights $\beta_1,\hdots, \beta_n\in \wh{G}_\R$.
We call $\Sigma_{GKZ}$ and its associated GIT problem an \emph{arrangeable fan} and \emph{arrangeable GIT problem} respectively, if the following two properties hold for any $\lambda\in \Lambda^\circ$: \vspace{5pt}
\begin{enumerate}[(i),itemsep = 0.5cm]
    \item\label{def:arrangeable_property_i} $H_{\lambda}^{=0}\cap |\Sigma_{GKZ}| = \beta_{I_{\lambda}^{=0}}$.
    \item\label{def:arrangeable_property_ii} $|\{[\rho] \ | \ \rho\in I_{\lambda}^{=0}\}|= \dim \wh{G}_\R-1$.
    %For any $\rho\in I_{\lambda}^{=0}$, $\beta_{\{\rho\}}$ is a ray of $H_{\lambda}^{=0}\cap |\Sigma_{GKZ}|$
\end{enumerate}
\end{definition}
Property \ref{def:arrangeable_property_i} of arrangeability requires that the vector space span of any internal wall in the secondary fan only hits other walls, while missing the interior of chambers. This will be made more precise in \Cref{lemma:arrangeable_implies_nonsplitting}. 

Recall that for any $\lambda\in \Lambda$, $\dim \beta_{I_{\lambda}^{=0}} =\dim \wh{G}_\R -1$ by \Cref{cor:schmitt}\ref{cor:schmitt_i}. Property \ref{def:arrangeable_property_ii} of arrangeability states that for $\lambda\in \Lambda^\circ$, $I_{\lambda}^{=0}$ contains minimally many weights up to repeats and positive rescaling. In particular, property \ref{def:arrangeable_property_ii} implies that $\beta_{I_{\lambda}^{=0}}$ is a simplicial cone where each $\rho\in I_\lambda^{=0}$ defines a ray. With this in mind, it follows that for any $\rho_1,\hdots, \rho_\ell\in I_\lambda^{=0}$ such that $[\rho_i]\ne [\rho_j]$ for any $i\ne j$, that
\[\dim \beta_{\{\rho_1,\hdots, \rho_\ell\}} = \dim \beta_{[\rho_1]\cup \hdots \cup [\rho_\ell]} = \ell.\]
\begin{example}\label{example:H2_is_arrangeable} We continue with the $\CH_2$ secondary fan $\Sigma_{GKZ}$ from \Cref{example:H2_stratifications}. From the depiction of $\Sigma_{GKZ}$, we see that the only $\xi$ in $\Lambda^\circ$ are $\pm \xi_2$. As the conditions of \Cref{def:arrangeable} are equivalent for $\lambda$ and $-\lambda$, to verify that $\Sigma_{GKZ}$ is arrangeable, we need only verify that conditions \ref{def:arrangeable_property_i}, \ref{def:arrangeable_property_ii} hold for $\xi_2$.

To check condition \ref{def:arrangeable_property_i}, observe that
\[H_{\xi_0}^{=0}\cap |\Sigma_{GKZ}| = \Cone((0,1)^T) = \beta_{\{4\}} = \beta_{I_{\xi_0}^{=0}}.\]
To check condition \ref{def:arrangeable_property_ii}, note
\[|\{[\rho]\ : \ \rho\in I_{\xi_0}^{=0}\}| = |\{[4]\}| = 1 = \dim \wh{G}_\R - 1.\]
So we conclude that $\Sigma_{GKZ}$ is arrangeable.
\end{example}
\begin{example}\label{example:non_arrangeable_examples} In \Cref{fig:arrangeable_non_examples} we depict three different secondary fans $\Sigma_1,\Sigma_2,\Sigma_3$ in $\wh{G}_\R = \R^3$. To see that $\Sigma_1$ is arrangeable, we consider the three distinct 1-PSs (up to a sign) associated to the internal walls of $\Sigma_1$. For each such 1-PS $\lambda$, the associated hyperplane $H_\lambda^{=0}$ spans the entirety of the secondary fan, so they each satisfy property \ref{def:arrangeable_property_i}. Each such $\lambda$ also has only 2 distinct weights in $I_{\lambda}^{=0}$, so property \ref{def:arrangeable_property_ii} is satisfied. 

To see that $\Sigma_2$ is not arrangeable, consider the 1-PS $\lambda$ associated to the red internal wall. Here we see that $H_{\lambda}^{=0}$ extends beyond $\beta_{I_{\lambda}^{=0}}$ into the northeast chamber, thus violating property \ref{def:arrangeable_property_i} of arrangeability.

To see that $\Sigma_3$ is not arrangeable, consider the 1-PS $\lambda$ associated to the red internal walls. We see that there are 3 weights in $I_{\lambda}^{=0}$, which exceeds $\dim \wh{G}_\R-1 = 2$, so $\lambda$ violates property \ref{def:arrangeable_property_ii} of arrangeability.
\begin{figure}
    \centering
\[\begin{tikzpicture}[beta/.style={circle,fill=black,inner sep=1.5pt}]
%%% Fan 1
\def\rootthree{1.73205}
\def\fanscale{1.5}
\def\xshift{4}

\coordinate (beta1) at (-1*\fanscale,0);
\coordinate (beta2) at (0,\rootthree*\fanscale);
\coordinate (beta3) at (1*\fanscale,0);
\coordinate (beta4) at ($0.5*(beta1)+0.5*(beta2)$);
\coordinate (beta5) at ($0.5*(beta1)+0.5*(beta3)$);
\coordinate (beta6) at ($0.6667*(beta5)+0.3333*(beta2)$);
\coordinate (beta7) at ($0.5*(beta2)+0.5*(beta3)$);

\draw[-,black,thick] (beta1) -- (beta2);
\draw[-,black,thick] (beta2) -- (beta3);
\draw[-,black,thick] (beta1) -- (beta3);
\draw[-,black,thick] (beta4) -- (beta5);
\draw[-,black,thick] (beta4) -- (beta3);
\draw[-,black,thick] (beta5) -- (beta2);

\draw[draw=none,fill=gray,opacity=.2] (beta1) -- (beta2) -- (beta3) -- cycle;

\node[beta] at (beta1) {};
\node[beta] at (beta2) {};
\node[beta] at (beta3) {};
\node[beta] at (beta4) {};
\node[beta] at (beta5) {};

\node[draw=none, fill=none] at (1*\fanscale,\rootthree*\fanscale) {\small arrangeable};
\node[draw=none, fill=none] at (-0.5*\fanscale,\rootthree*\fanscale) {\small $\Sigma_1$};

%% Fan 2
\coordinate (beta12) at ($(beta1)+(\xshift*\fanscale,0)$);
\coordinate (beta22) at ($(beta2)+(\xshift*\fanscale,0)$);
\coordinate (beta32) at ($(beta3)+(\xshift*\fanscale,0)$);
\coordinate (beta42) at ($(beta4)+(\xshift*\fanscale,0)$);
\coordinate (beta52) at ($(beta5)+(\xshift*\fanscale,0)$);
\coordinate (beta62) at ($(beta6)+(\xshift*\fanscale,0)$);
\coordinate (beta72) at ($(beta7)+(\xshift*\fanscale,0)$);

\draw[draw=none,fill=gray,opacity=.2] (beta12) -- (beta22) -- (beta32) -- cycle;

\draw[-,black,thick] (beta12) -- (beta22);
\draw[-,black,thick] (beta22) -- (beta32);
\draw[-,black,thick] (beta12) -- (beta32);
\draw[-,red,thick] (beta12) -- (beta62);
\draw[-,black,thick] (beta22) -- (beta62);
\draw[-,black,thick] (beta32) -- (beta62);

\draw[-,dashed,red,thick] (beta62) -- (beta72);
% \draw[-,dashed,red,thick] (beta62) -- (beta52);
% \draw[-,dashed,red,thick] (beta62) -- (beta72);

\node[beta] at (beta12) {};
\node[beta] at (beta22) {};
\node[beta] at (beta32) {};
\node[beta] at (beta62) {};

\node[draw=none, fill=none] at (1*\fanscale+\xshift*\fanscale,\rootthree*\fanscale) {\small violates \ref{def:arrangeable_property_i}};
\node[draw=none, fill=none] at (-0.5*\fanscale+\xshift*\fanscale,\rootthree*\fanscale) {\small $\Sigma_2$};

%% Fan 3
\coordinate (beta13) at ($(beta1)+(2*\xshift*\fanscale,0)$);
\coordinate (beta23) at ($(beta2)+(2*\xshift*\fanscale,0)$);
\coordinate (beta33) at ($(beta3)+(2*\xshift*\fanscale,0)$);
\coordinate (beta43) at ($(beta4)+(2*\xshift*\fanscale,0)$);
\coordinate (beta53) at ($(beta5)+(2*\xshift*\fanscale,0)$);
\coordinate (beta63) at ($(beta6)+(2*\xshift*\fanscale,0)$);
\coordinate (beta73) at ($(beta7)+(2*\xshift*\fanscale,0)$);

\draw[-,black,thick] (beta13) -- (beta23);
\draw[-,black,thick] (beta23) -- (beta33);
\draw[-,black,thick] (beta13) -- (beta33);
\draw[-,black,thick] (beta23) -- (beta53);
\draw[-,black,thick] (beta33) -- (beta43);
\draw[-,black,thick] (beta43) -- (beta53);
\draw[-,black,thick] (beta53) -- (beta73);
\draw[-,black,thick] (beta43) -- (beta73);

\draw[draw=none,fill=gray,opacity=.2] (beta13) -- (beta23) -- (beta33) -- cycle;

\draw[-,red,thick] (beta13) -- (beta73);

\node[beta,red] at (beta13) {};
\node[beta] at (beta23) {};
\node[beta] at (beta33) {};
\node[beta] at (beta43) {};
\node[beta] at (beta53) {};
\node[beta,red] at (beta63) {};
\node[beta,red] at (beta73) {};

\node[draw=none, fill=none] at (1*\fanscale+2*\xshift*\fanscale,\rootthree*\fanscale) {\small violates \ref{def:arrangeable_property_ii}};
\node[draw=none, fill=none] at (-0.5*\fanscale+2*\xshift*\fanscale,\rootthree*\fanscale) {\small $\Sigma_3$};

\end{tikzpicture}
\]
\caption{In this figure we draw a cross section of three secondary fans $\Sigma_1,\Sigma_2,\Sigma_3$ in $\wh{G}_\R = \R^3$ with the weights $\beta_\rho$ for each problem being depicted as dots. The first fan is arrangeable, but the other two are not. For an explanation, see \Cref{example:non_arrangeable_examples}.}
\label{fig:arrangeable_non_examples}
\end{figure}
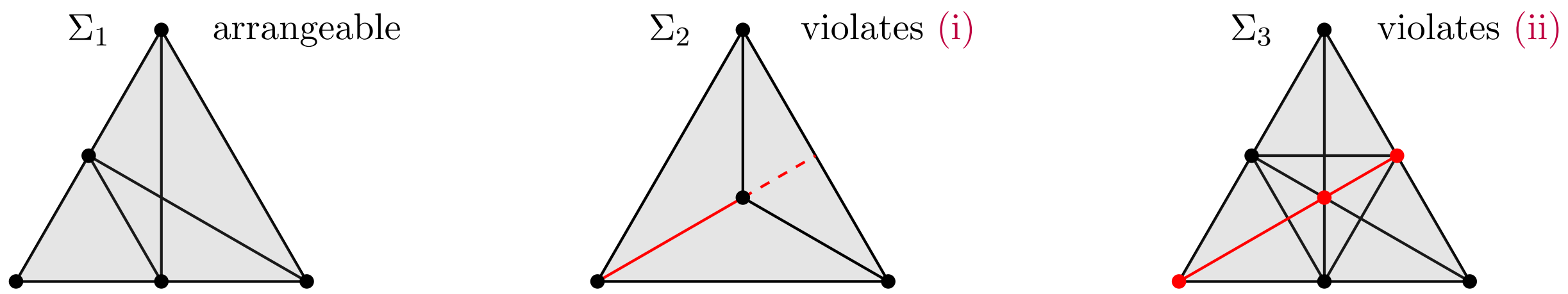
\end{example}
\begin{remark}\label{rmk:all_rank_2_problems_are_arrangeable} Observe that if $\dim \wh{G}_\R = 1$ or $2$, then any projective GIT problem is automatically arrangeable. If $\dim \wh{G}_\R = 3$, it is easy to show that a projective GIT problem is arrangeable if and only if none of the weights $\beta_i$ are in the interior of $|\Sigma_{GKZ}|$. Arrangeability becomes a more restrictive condition for higher rank GIT problems. However, we have examples of arrangeable GIT problems in arbitrarily high rank as shown in the following example.
\end{remark}

\begin{example}\label{example:blowups_of_Pn}
Choose $m \ge 2$ and $0 \le k \le m+1$. Let $G\cong \Gm^{k+1}$ so that $\wh{G}\cong \Z^{k+1}$. Let $e_1,\hdots, e_{k+1}$ be the standard basis for $\Z^{k+1}$. We define a GIT problem $G\subset \Gm^{m+k+1}$ via the weights $\beta_\rho\in \wh{G}$:
\[\begin{array}{cccc} \beta_1 = e_1,& \hdots,& \beta_k=e_k, & \beta_{k+1}=\hdots = \beta_{m+1} = e_{k+1},\\ \beta_{m+2} = e_1+e_{k+1},& \hdots, & \beta_{m+k+1} = e_k+e_{k+1}.
\end{array}\]
\begin{figure}
\begin{tikzpicture}
\def\fanscale{6}
\def\trih{0.866} % height of equilateral triangle
\def\op{0.3}
\def\lineop{0.5}
\def\thickness{thick}
\def\dashsty{dashed}
\def\pslength{1.5cm}

\coordinate (beta1) at (-0.5*\fanscale, 0*\fanscale);
\coordinate (beta3) at (0.5*\fanscale, 0*\fanscale);
\coordinate (beta2) at (0*\fanscale, \trih*\fanscale);
\coordinate (betam2) at ($(beta1)!0.5!(beta3)$);
\coordinate (betam3) at ($(beta2)!0.5!(beta3)$);
\coordinate (center) at ($(betam2)!0.333!(beta2)$);

\draw[draw=none, fill=gray,opacity=\op] (beta1) -- (beta2) -- (betam3) -- (center) -- (betam2) -- cycle;
\draw[draw=none, fill=gray,opacity=\op] (betam2) -- (betam3) -- (beta3) -- cycle;
\draw[draw=none,fill = red, opacity=\op] (betam2) -- (center) -- (betam3) -- cycle;

\tikzset{beta/.style={circle,inner sep =1.5pt, fill = black, draw =black}}
\node[beta] at (beta1) {};
\node[beta] at (beta2) {};
\node[beta] at (beta3) {};
\node[beta] at (betam2) {};
\node[beta] at (betam3) {};

\draw[draw=black,\thickness,opacity=\lineop] (beta1) -- (beta2);
\draw[draw=black,\thickness,opacity=\lineop] (beta1) -- (beta3);
\draw[draw=black,\thickness,opacity=\lineop] (beta2) -- (beta3);
\draw[draw=black,\thickness,opacity=\lineop] (betam2) -- (beta2);
\draw[draw=black,\thickness,opacity=\lineop] (beta1) -- (betam3);
\draw[draw=black,\thickness,opacity=\lineop] (betam2) -- (betam3);

\node[draw=none,fill=none] at ($(beta1)+(-0.4,0)$) {\small $\beta_1$};
\node[draw=none,fill=none] at ($(beta2)+(0,+0.3)$) {\small $\beta_2$};
\node[draw=none,fill=none] at ($(beta3)+(1.6,0)$) {\small $\beta_3 = \hdots = \beta_{m+1}$};
\node[draw=none,fill=none] at ($(betam2)+(0,-0.3)$) {\small $\beta_{m+2}$};
\node[draw=none,fill=none] at ($(betam3)+(0.4,0.3)$) {\small $\beta_{m+3}$};

\node[draw=none,fill=none,red] at ($(betam2)!0.5!(betam3)!0.45!(center)$) {$\Gamma$};

\end{tikzpicture}
    \caption{The secondary fan $\Sigma_{GKZ}$ in $\wh{G}_\R\cong \R^3$ of the GIT problem described in \Cref{example:blowups_of_Pn} in the case of $k=2$ and $m \ge 2$. For the highlighted chamber $\Gamma$, $\CX_\Gamma = X$, the blowup of $\P^m$ at $k$ torus fixed points.}
    \label{fig:GKZ_for_blowups_of_Pn}
\end{figure}
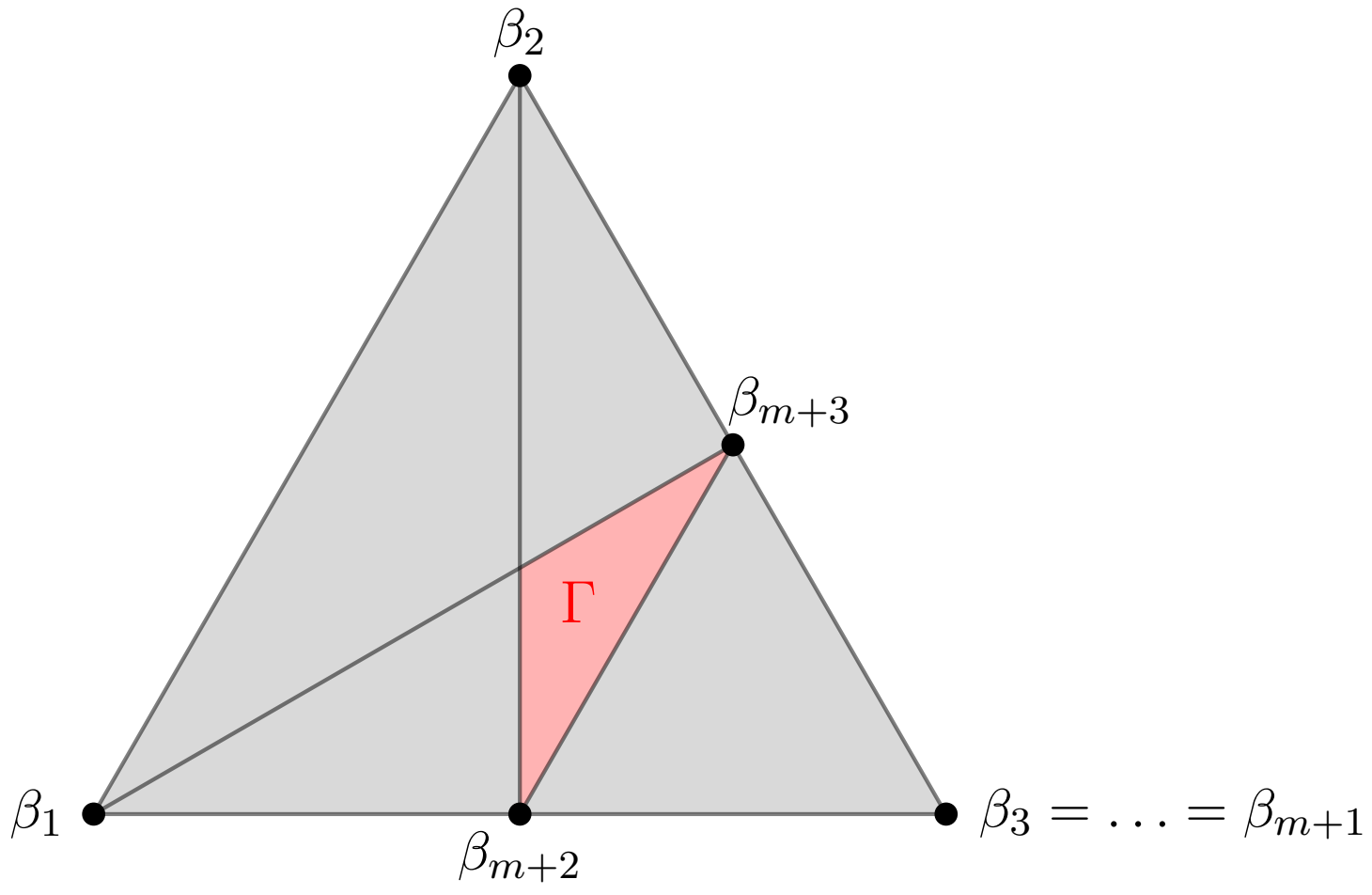
Let $X$ be the space obtained by iteratively blowing up $\P^m$ at $k$ of the $m+1$ torus fixed points for some $0 \le k \le m+1$. Then, one can check that $X$ is a GIT quotient of $\A^{m+k+1}$ by $G$ and thus is associated to a chamber $\Gamma$ of the secondary fan $\Sigma_{GKZ}$ associated to $G\subset \Gm^{m+k+1}$. See \Cref{fig:GKZ_for_blowups_of_Pn} for the associated secondary fan and chamber in the case of $k=2$ and $m \ge 2$.

We will now describe how to check that $\Sigma_{GKZ}$ is arrangeable. If $k = 0$, it is easy to verify that $\Sigma_{GKZ}$ is arrangeable, so instead assume $k\ge 1$. In the first case, suppose $k \le m$. For $J\subset [k]$, define the 1-PS $\lambda_J$ via
\[\wh{\lambda}_J(e_j) = \begin{cases} 0 & j \in J \\
-1 & j\in [k]\setminus J\\
1 & j = k+1
\end{cases}\]
One can check that $\lambda_J$ is a well-defined 1-PS and in fact

\[\Lambda^\circ = \bigcup_{J\subsetneq [k]}\{\lambda_J, -\lambda_J\}.\]
Thus to show that $\Sigma_{GKZ}$ is arrangeable, it suffices to show that $\lambda_J$ satisfies properties \ref{def:arrangeable_property_i} and \ref{def:arrangeable_property_ii} of \Cref{def:arrangeable} for each $J\subsetneq [k]$. First observe that
\[I_{\lambda_J}^{=0} = J\cup (m+([k]\setminus J))\]
so $|I_{\lambda_J}^{=0}| = k$. In particular, this is strictly stronger than property \ref{def:arrangeable_property_ii}. Because $\beta_{I_{\lambda}^{=0}}\subset H_{\lambda}^{=0}\cap |\Sigma_{GKZ}|$ for any $\lambda$, to show \ref{def:arrangeable_property_i} it suffices to show the opposite inclusion. So, let $x\in H_{\lambda_J}^{=0}\cap |\Sigma_{GKZ}|$ with coordinate expression $x = \sum_{j\in [k+1]}x_je_j$. Then we have that
\begin{equation}\label{eqn:blown_up_Pn}
    0 = \wh{\lambda}_J(x) = x_{k+1}-\sum_{i\in [k]\setminus J}x_i.
\end{equation}
Thus
\begin{align*}
    x &= \sum_{i\in [k+1]}x_ie_i\\
    &= \sum_{j\in J}x_j e_j+\sum_{i\in [k]\setminus J}x_ie_i+x_{k+1}e_{k+1} \\
    &= \sum_{j\in J}x_j e_j+\sum_{i\in [k]\setminus J}x_ie_i+\lp \sum_{i\in [k]\setminus J}x_i\rp e_k & \text{by eqn. \eqref{eqn:blown_up_Pn}} \\
    &= \sum_{j\in J}x_j e_j+\sum_{i\in [k]\setminus J}x_i(e_i+e_k).
\end{align*}
Because $|\Sigma_{GKZ}| = \Cone(e_1,\hdots, e_{k+1})$, we also have that $x_i \ge 0$ for each $i$, so the above implies that $x\in \beta_{I_{\lambda_J}^{=0}}$. So, property \ref{def:arrangeable_property_ii} holds and we conclude that $\Sigma_{GKZ}$ is arrangeable. 

If instead $k = m+1$, $\Lambda^\circ$ has the same description, except $\lambda_{[k]}\in |\Sigma_{GKZ}|^\vee$. So, the same analysis yields that for each $\lambda\in \Lambda^\circ$, $\lambda$ satisfies properties \ref{def:arrangeable_property_i} and \ref{def:arrangeable_property_ii} of \Cref{def:arrangeable}.
\end{example}

We now seek to construct HKKN stratifications in arrangeable GIT problems. As a first step, we show the following lemma.
\begin{lemma}\label{lemma:arrangeable_implies_nonsplitting} Let $\sigma$ be any cone of an arrangeable secondary fan $\Sigma_{GKZ}$. Then, for every $\lambda\in \Lambda$, $\sigma\subset H_{\lambda}^{\ge 0}$ or $\sigma\subset H_{-\lambda}^{\ge 0}$ (or both).
\end{lemma}
\begin{proof} Suppose, by way of contradiction, there exists $\lambda\in \Lambda$ such that $\sigma\not\subset H_{\lambda}^{\ge 0}$ and $\sigma\not\subset H_{-\lambda}^{\ge 0}$. Equivalently, $\sigma \cap H_{\lambda}^{<0}\ne \emptyset$ and $\sigma \cap H_{\lambda}^{>0} \ne \emptyset$. 
Because $H_\lambda^{<0},H_\lambda^{>0}$ are open, these imply that $(\Relint \sigma)\cap H_\lambda^{<0}\ne\emptyset$ and $(\Relint \sigma)\cap H_\lambda^{>0}\ne \emptyset$. 
So, we can choose $x_-\in (\Relint \sigma)\cap H_\lambda^{<0}$ and $x_+\in (\Relint \sigma)\cap H_\lambda^{>0}$. Then, on the straight line path from $x_-$ to $x_+$ there must be a unique $x_0\in (\Relint \sigma) \cap H_\lambda^{=0}$. If this $x_0$ were in $\beta_{I_{\lambda}^{=0}}$, then by \Cref{cor:schmitt}\ref{cor:schmitt_i} there is a wall $\tau\subset H_{\lambda}^{=0}$ containing $x_0$. However, $x_0\in \Relint \sigma$ thus $\sigma$ is the unique minimal cone (via containment) containing $x_0$, which implies $\sigma \subseteq \tau$.
However, $\tau\subset H_{\lambda}^{=0}$ and $\sigma\not\subset H_{\lambda}^{=0}$, so this cannot be the case. So instead it must be that $x_0\notin \beta_{I_{\lambda}}^{=0}$, thus violating condition \ref{def:arrangeable_property_i} of arrangeability.

\end{proof}
\begin{remark}\label{rmk:non_splitting_eq_to_arr_property_1} In fact, one can show that condition \ref{def:arrangeable_property_i} of arrangeability is equivalent to the conclusion of \Cref{lemma:arrangeable_implies_nonsplitting} holding for every $\sigma\in \Sigma_{GKZ}$. 
\end{remark}
The main use of \Cref{lemma:arrangeable_implies_nonsplitting} is the following lemma.

\begin{lemma}\label{lemma:arrangeable_implies_plus_ideal_equals_irr} Let $\Sigma_{GKZ}$ be an arrangeable GIT problem and $\Gamma$ be a chamber of $\Sigma_{GKZ}$. Then
\[B_{\Gamma} = \bigcap_{\lambda\in \Lambda \cap \Gamma^\vee} \la I_\lambda^{>0}\ra.\]
\end{lemma}
\begin{proof} Recall by \Cref{lemma:B_chi_via_betas},
\[B_\Gamma = \la x^J \ | \ \Gamma \subset \beta_J\ra.\]
So, take any $J\subset [n]$ so that $\Gamma \subset \beta_J$. Take any $\lambda\in \Gamma^\vee$ and suppose by way of contradiction that $J\subset I_\lambda^{\le 0}$. Then this would imply 
\[\Gamma \subset \beta_J\subset H_\lambda^{\le 0}\]
which contradicts $\lambda \in \Gamma^\vee$ as $\Gamma$ is full-dimensional. Thus we have that $J\cap I_\lambda^{>0}\ne \emptyset$, so $x^J\in \la I_\lambda^{>0}\ra$ and 
\[B_\Gamma \subset \bigcap_{\lambda\in \Lambda\cap \Gamma^\vee}\la  I_\lambda^{>0}\ra.\]

To show the opposite inclusion, see that the minimal generators of the right hand side are exactly the square-free monomials $x^J$ such that for each $\lambda\in \Lambda\cap \Gamma^\vee$, there exists some $\rho\in I_\lambda^{>0}\cap J$.   Thus by \Cref{cor:gamma_hits_beta_J_implies_contained}, it suffices to show for any such $J$, one has $\Gamma \subset \beta_J$. Suppose this is not the case. By \Cref{lemma:beta_cones_determined_by_lambdas}, there is some $\lambda \in \Lambda\cap \beta_J^\vee$ so that $H_{\lambda}^{<0} \cap \Gamma \ne \emptyset$. However, as all cones of arrangeable secondary fans are contained on one side of the half space for any $\lambda\in \Lambda$ by \Cref{lemma:arrangeable_implies_nonsplitting}, this implies that $\Gamma \subset H_{\lambda}^{\le 0}$. Thus $-\lambda \in \Lambda \cap \Gamma^\vee$, so by assumption, $J\cap I_{-\lambda}^{>0} = J \cap I_{\lambda}^{<0}$ is nonempty. However, this contradicts $\lambda\in \beta_J^\vee$, so we conclude instead that $\Gamma \subset \beta_J$ as desired. 
\end{proof}

Geometrically, \Cref{lemma:arrangeable_implies_plus_ideal_equals_irr} states that the contracting loci of the 1-PSs appearing in $\Gamma^\vee \cap \Lambda$ cut out the unstable locus $V(B_\Gamma)$. In this way, the unordered set $\Gamma^\vee \cap \Lambda$ satisfies the first of the three conditions required to induce an HKKN stratification from \Cref{lemma:HKKN_strat_combo}. So, our goal in the remainder of the subsection will be to place an ordering on $\Gamma^\vee \cap \Lambda$ that satisfies the second and third conditions.

\begin{definition}\label{def:wall_partial_order} Let $\Gamma$ be a chamber of a secondary fan $\Sigma_{GKZ}$. We define a relation $\wrel$ on pairs $\lambda,\xi\in \Gamma^\vee \cap \Lambda$ as follows.
\begin{enumerate}[(i)]
    \item\label{def:wall_partial_order_i} If $I_\xi^{>0}\subset I_\lambda^{\ne 0}$ or equivalently $I_{\lambda}^{=0}\subset I_{\xi}^{\le 0}$, then $\lambda \wrel \xi$.
    \item\label{def:wall_partial_order_ii} If $\lambda\in |\Sigma_{GKZ}|^\vee$ and $\xi\notin |\Sigma_{GKZ}|^{\vee}$, then $\lambda \wrel \xi$.
\end{enumerate}
\end{definition}
When we refer to the relation $\wrel$ it will be clear from context which chamber $\Gamma$ defines the relation. This relation is exactly constructed so that if $(\Gamma^\vee\cap \Lambda, \wrel)$ is indeed a poset, any extension to a total order automatically satisfies condition \ref{lemma:HKKN_strat_combo_ii} of \Cref{lemma:HKKN_strat_combo}. As such, we wish to show $(\Gamma^\vee\cap \Lambda,\wrel)$ is a poset (see \Cref{lemma:wall_relation_is_poset}). To do so, we will first prove the following simple convex geometry lemma.

\begin{lemma} \label{lemma:hyperplane_and_halfspace_containment} If $C$ is a cone in $\wh{G}_\R$ and $\lambda,\xi$ are two one-parameter subgroups such that $H_{\lambda}^{=0}\cap C \subset H_{\xi}^{\le 0} \cap C$, then either
\[H_{\lambda}^{\le 0}\cap C \subset H_{\xi}^{\le 0}\cap C\quad \text{or}\quad H_{\lambda}^{\ge 0}\cap C \subset H_{\xi}^{\le 0}\cap C.\]
\end{lemma}
\begin{proof} Suppose this was not the case. Then there are points
\begin{align*}
    x_- &\in H_{\lambda}^{\le 0}\cap C \cap H_{\xi}^{>0},\\
    x_+ &\in H_{\lambda}^{ \ge 0}\cap C\cap H_{\xi}^{>0}.
\end{align*}
Let $\gamma:[0,1]\to \wh{G}_\R$ be a straight line path from $x_-$ to $x_+$. As $\wh{\lambda}\circ \gamma$ is monotonic and $\wh{\lambda}(\gamma(0)) \le 0$ and $\wh{\lambda}(\gamma(1)) \ge 0$, there must be some $t_0\in [0,1]$ so that $\wh{\lambda}(\gamma(t_0)) = 0$. As $C$ is convex, we have that $\gamma(t_0)\in C\cap H_{\lambda}^{=0}$, thus by assumption, $\wh{\xi}(\gamma(t_0)) \le 0$. However, by our choice of $x_-,x_+$, we have that $\wh{\xi}(\gamma(0)), \wh{\xi}(\gamma(1)) > 0$. The linearity of $\wh{\xi}\circ \gamma$ would imply $\wh{\xi}(\gamma(t_0)) > 0$, a contradiction.
\end{proof}
In the arrangeable setting, we use the above to get the following lemma.
\begin{lemma}\label{lemma:arrangeable_implies_negative_halfspace_containment} Let $\Sigma_{GKZ}$ be an arrangeable secondary fan, $\Gamma$ a chamber of $\Sigma_{GKZ}$, and $\lambda,\xi\in \Gamma^\vee \cap \Lambda$. If $I_{\lambda}^{=0}\subset I_{\xi}^{\le 0}$, then $H_{\lambda}^{\le 0}\cap |\Sigma_{GKZ}|\subset H_{\xi}^{\le 0}\cap |\Sigma_{GKZ}|$. As an immediate consequence, $I_{\lambda}^{\le 0}\subset I_\xi^{\le 0}$. 
\end{lemma}

\begin{proof} We immediately get the following inclusions and equalities:
\[H_{\lambda}^{=0}\cap |\Sigma_{GKZ}| = \beta_{I_{\lambda}^{=0}}\subset H_{\xi}^{\le 0}\cap |\Sigma_{GKZ}|.\]
The equality is by condition \ref{def:arrangeable_property_i} of arrangeability and the inclusion is by the assumption that $I_{\lambda}^{=0}\subset I_\xi^{\le 0}$. As $\Gamma$ is full-dimensional and $\Gamma\subset H_\xi^{\ge 0}$, it follows that $\Gamma\not\subset H_\xi^{\le 0}\cap |\Sigma_{GKZ}|$. Then, because $\Gamma\subset H_{\lambda}^{\ge 0}\cap |\Sigma_{GKZ}|$ and $\Gamma\not\subset H_{\xi}^{\le 0}\cap |\Sigma_{GKZ}|$, we conclude the desired result directly from \Cref{lemma:hyperplane_and_halfspace_containment}.
\end{proof}
We are now in a position to show that $(\Gamma^\vee \cap \Lambda,\wrel)$ is a poset.
\begin{lemma}\label{lemma:wall_relation_is_poset} Let $\Gamma$ be a chamber of an arrangeable secondary fan $\Sigma_{GKZ}$. Then, $(\Gamma^\vee\cap \Lambda, \wrel)$ of \Cref{def:wall_partial_order} is a poset.
\end{lemma}
\begin{proof} First observe that as $I_{\lambda}^{>0}\subset I_\lambda^{\ne 0}$, we have $\lambda \wrel \lambda$ for every $\lambda\in \Gamma^\vee \cap \Lambda$. Thus $\wrel$ is reflexive. 

Before showing that $\wrel$ is transitive and antisymmetric, we will show that for any relation $\xi \wrel \lambda$, if $\lambda\in |\Sigma_{GKZ}|^\vee$, then $\xi = \lambda$. Suppose we had such a relation and $\lambda\in |\Sigma_{GKZ}|^\vee$. Then, by the definition of $\wrel$, it cannot be that $\xi \wrel \lambda$ is of type \ref{def:wall_partial_order_ii}, so instead we must have $I_\lambda^{>0}\subset I_{\xi}^{\ne 0}$. Taking complements, this implies $I_\xi^{=0}\subset I_\lambda^{\le 0}$. However, $\lambda\in |\Sigma_{GKZ}|^\vee$, so $I_\lambda^{\le 0} = I_{\lambda}^{=0}$. Thus $I_\xi^{=0} = I_\lambda^{=0}$, which implies $\xi = \pm \lambda$ by \Cref{cor:schmitt}\ref{cor:schmitt_v}. That $\xi,\lambda\in \Gamma^\vee$ and $\Gamma$ is full dimensional together imply that in fact $\xi = \lambda$.

We will now show that $\wrel$ is transitive. Suppose $\lambda_1 \wrel \lambda_2\wrel \lambda_3$. If $\lambda_3\in |\Sigma_{GKZ}|^\vee$, then by the above observation we would get $\lambda_1 = \lambda_2 = \lambda_3$. If $\lambda_2\in |\Sigma_{GKZ}|^\vee$, then again by the above $\lambda_1 = \lambda_2$ so $\lambda_1 \wrel \lambda_3$. Assume then that $\lambda_2,\lambda_3\notin |\Sigma_{GKZ}|^\vee$. If $\lambda_1\in |\Sigma_{GKZ}|^\vee$, then $\lambda_1\wrel \lambda_3$. The only remaining case is that each $\lambda_1,\lambda_2,\lambda_3\notin |\Sigma_{GKZ}|^\vee$. Then both relations are of type \ref{def:wall_partial_order_i}, so 
\[I_{\lambda_1}^{= 0}\subset I_{\lambda_2}^{\le 0},\quad I_{\lambda_2}^{=0}\subset I_{\lambda_3}^{\le 0}.\]
By \Cref{lemma:arrangeable_implies_negative_halfspace_containment}, the second inclusion implies $I_{\lambda_2}^{\le 0}\subset I_{\lambda_3}^{\le 0}$. Thus $I_{\lambda_1}^{= 0}\subset I_{\lambda_3}^{\le 0}$, so $\lambda_1\wrel \lambda_3$ as desired.

We will now show that $\wrel$ is antisymmetric, completing the proof. Suppose $\lambda\wrel \xi$ and $\xi \wrel \lambda$. If either of $\lambda$ or $\xi$ are in $|\Sigma_{GKZ}|^\vee$, then again $\lambda = \xi$ as desired. So, assume that both $\lambda,\xi\notin |\Sigma_{GKZ}|^\vee$. Then both relations are of type \ref{def:wall_partial_order_i}, so
\[I_{\lambda}^{=0}\subset I_{\xi}^{\le 0},\quad I_{\xi}^{=0}\subset I_{\lambda}^{\le 0}.\]
By \Cref{lemma:arrangeable_implies_negative_halfspace_containment} this implies
\[H_{\lambda}^{\le 0}\cap |\Sigma_{GKZ}|\subset H_{\xi}^{\le 0}\cap |\Sigma_{GKZ}| \subset H_{\lambda}^{\le 0}\cap |\Sigma_{GKZ}|\]
so we have that $H_{\lambda}^{\le 0}\cap |\Sigma_{GKZ}| = H_{\xi}^{\le 0}\cap |\Sigma_{GKZ}|$. Taking dual cones we find
\[\Cone(|\Sigma_{GKZ}|^\vee \cup \{\lambda\}) = \Cone(|\Sigma_{GKZ}|^\vee \cup \{\xi\}).\]
Because $\lambda\in \Gamma^\vee$, it follows that $-\lambda\notin |\Sigma_{GKZ}|^{\vee}$, so $\pm \lambda\notin |\Sigma_{GKZ}|^\vee$. It follows that $\Cone(\lambda)$ is a ray of the above cone which was not a ray of $|\Sigma_{GKZ}|^\vee$. Similarly, $\Cone(\xi)$ is a ray of the above cone which was not a ray of $|\Sigma_{GKZ}|^\vee$, so $\Cone(\lambda) = \Cone(\xi)$. As $\lambda,\xi\in \Lambda$, they are primitive, thus $\lambda=\xi$.
\end{proof}

\Cref{lemma:wall_relation_is_poset} and our work to this point together suggest that to obtain an HKKN stratification of $V(B_\Gamma)$, one should extend the poset $(\Gamma^\vee\cap \Lambda, \wrel)$ to a total order. The remaining work in showing that such an extension is indeed a wall stratification is showing that the third condition of \Cref{lemma:HKKN_strat_combo} is satisfied. The main technical ingredient in that demonstration is the following lemma.
\begin{lemma}\label{lemma:internally_simplicial_implies_close_previous_wall} 
Let $\Gamma$ be a chamber of an arrangeable secondary fan and let $\lambda\in \Gamma^\vee \cap \Lambda^{\circ}$. Then, for any $\rho_0\in I_{\lambda}^{=0}$, there is a 1-PS $\xi\in \Gamma^\vee \cap\Lambda$ such that $\xi \wrel \lambda$ and
\[I_{\xi}^{>0}\cap I_\lambda^{\le 0} = [\rho_0].\]
\end{lemma}
 \begin{proof} Because $\lambda\in \Lambda^{\circ}$,
 \[\sigma \coloneqq |\Sigma_{GKZ}|\cap H_{\lambda}^{\le 0}\]
is a dimension $k$ cone. Assume first that $\dim \wh{G}_\R \ge 3$. By condition \ref{def:arrangeable_property_ii} of arrangeability, it follows that  $\beta_{I_{\lambda}^{=0}\setminus [\rho_0]}$ is a dimension $\dim \wh{G}_\R-2$ face of $\sigma$. By generalized facet-ray duality (see e.g. \cite[Proposition 1.2.10]{CLS}), we have that
\[\tau \coloneqq \sigma^{\vee} \cap \bigcap_{\rho \in I_{\lambda}^{=0}\setminus [\rho_0]}H_{\beta_\rho}^{=0}\]
is thus a dimension 2 face of $\sigma^\vee$. If $\dim \wh{G}_\R=2$, simply define
\[\tau \coloneqq \sigma^\vee.\]
The proof of both cases proceeds the same way hereafter. By construction, $\lambda\in \tau$, and as $\lambda\in \sigma^\vee(1)$ and $\tau \face \sigma^\vee$, we have that $\lambda\in \tau(1)$. As $\dim \tau = 2$, $\tau$ is simplicial and thus has exactly two rays, so let $\mu \ne \lambda$ be the generator of the other ray of $\tau$. As $\mu \in \tau(1)$ and $\tau\face \sigma^\vee$, we have that $\mu \in \sigma^\vee(1)$. That is, $\mu$ defines a facet of $\sigma = |\Sigma_{GKZ}|\cap H_{\lambda}^{\le 0}$ which is not $\lambda$. Thus $\mu$ defines a facet of $|\Sigma_{GKZ}|$.

We now want to parametrize the 1-PSs in $\tau$ (up to positive scaling) by using our description of the rays of $\tau$ so that we can find a wall of $\Sigma_{GKZ}$ in $\tau$ corresponding to a 1-PS which is in some sense ``close" to $\lambda$.
For $t\in [0,1]$ we define a functional $\wh{\xi}_t\in \wh{G}_\R^\ast$ via
\[\wh{\xi}_t \coloneqq (1-t)\wh{\mu}+t\wh{\lambda}.\]
As $\mu$ corresponds to an external wall of $\Sigma_{GKZ}$, for any $\rho\in I_{\lambda}^{<0} $, $\wh{\xi}_0(\beta_\rho) = \wh{\mu}(\beta_\rho) \ge 0$. However for $\rho\in I_{\lambda}^{<0}$, $\wh{\xi}_1(\beta_\rho) = \wh{\lambda}(\beta_\rho) < 0$. 

Thus, there is a unique $t_\rho\in [0,1)$ so that $\wh{\xi}_{t_\rho}(\beta_\rho) = 0$. See also that for fixed $\rho\in I_{\lambda}^{<0}$, $\wh{\xi}_{\bullet}(\beta_\rho)$ is decreasing. Define $\rho_-$ to be a $\rho\in I_{\lambda}^{<0}$ which maximizes $t_\rho$ among $\rho\in I_\lambda^{<0}$.
For any $\rho\in I_{\lambda}^{<0}$, as $\wh{\xi}_{\bullet}(\beta_\rho)$ is decreasing and $t_\rho \le t_{\rho_-}$, we have that
\[\wh{\xi}_{t_{\rho_-}}(\beta_\rho) \le \wh{\xi}_{t_\rho}(\beta_\rho) = 0,\]
thus $I_{\lambda}^{<0}\subset I_{\xi_{t_{\rho_-}}}^{\le 0}$. So we find that
\[I_{\xi_{t_{\rho_-}}}^{>0} \cap I_{\lambda}^{\le 0} = I_{\xi_{t_{\rho_-}}}^{>0} \cap I_{\lambda}^{= 0}\subset I_{\xi_{t_{\rho_-}}}^{\ne 0} \cap I_{\lambda}^{=0} = [\rho_0].\]
As $\xi_0 = \mu$ corresponds to an external wall of $\Sigma_{GKZ}$ distinct from $\lambda$, we have that $[\rho_0]\subset I_{\xi_0}^{>0}$. As $\wh{\xi}_{1}(\beta_{\rho_0}) = 0$, we have that for any $t < 1$, $\wh{\xi}_t(\beta_{\rho_0}) > 0$. In particular $[\rho_0]\subset I_{\xi_{t_{\rho_-}}}^{>0}$ so we have that in fact
\[I_{\xi_{t_{\rho_-}}}^{>0} \cap I_{\lambda}^{\le 0} = [\rho_0].\]
From our definition of $\rho_-$, one can deduce that $t_{\rho_-}\in \Q$ so that $\Cone(\xi_{t_{\rho_-}})$ is rational. We define $\xi$ to be the primitive ray generator of $\Cone(\xi_{t_{\rho_-}})$. Recall
\begin{equation}\label{eqn:proof_of_hard_HKKN}
    I_{\xi}^{=0} = I_{\xi_{t_{\rho_-}}}^{=0} \supset [\rho_-]\sqcup (I_{\lambda}^{=0}\setminus [\rho_0]).
\end{equation}
As $\lambda\in \Lambda^\circ$, property \ref{def:arrangeable_property_ii} of arrangeability implies $\codim \beta_{I_{\lambda}^{=0}\setminus [\rho_0]} = 2$. Because $\rho_-\notin I_\lambda^{=0}$, 
\[\codim \beta_{[\rho_-]\sqcup I_{\lambda}^{=0}\setminus [\rho_0]} = 1,\]
so $\codim \beta_{I_{\xi}^{=0}} \le 1$. Then $\xi \ne 0$ implies that in fact $\codim \beta_{I_{\xi}^{=0}}=1$ so $\xi\in \Lambda$ by \Cref{cor:schmitt}\ref{cor:schmitt_iv}. 

To complete the proof we need only show that $\xi\wrel \lambda$. First note that $\xi\in \tau\face \sigma^\vee$, so $\Gamma\subset \sigma\subset H_{\xi}^{\ge 0}$. That is, $\xi \in \Gamma^\vee \cap \Lambda$. If $\xi\in |\Sigma_{GKZ}|^\vee$, we are done because $\lambda\notin |\Sigma_{GKZ}|^\vee$ implies $\xi \wrel \lambda$. So, assume $\xi\notin |\Sigma_{GKZ}|^\vee$. Then, by condition \ref{def:arrangeable_property_ii} of arrangeability, $I_\xi^{=0}$ is a union of $\dim \wh{G}_\R-1$ distinct sets $[\rho]$ for $\rho\in I_\xi^{=0}$. Because $\lambda\in \Lambda^\circ$, it also follows from condition \ref{def:arrangeable_property_ii} of arrangeability that $[\rho_-]\sqcup I_\lambda^{=0}\setminus [\rho_0]$ is a union of $\dim \wh{G}_\R -1$ distinct sets $[\rho]$. Because $I_\xi^{=0}$ is partitioned by the sets $[\rho]$, this implies that the inclusion from eqn. \eqref{eqn:proof_of_hard_HKKN} is in fact an equality. It follows that
\[I_{\xi}^{=0} = [\rho_-]\sqcup I_{\lambda}^{=0}\setminus [\rho_0] \subset I_{\lambda}^{\le 0},\]
so $\xi \wrel \lambda$ as desired.
\end{proof}
From \Cref{lemma:internally_simplicial_implies_close_previous_wall} we obtain the following lemma which will be useful for determining when a collection of 1-PSs induces an HKKN stratification by analyzing condition \ref{lemma:HKKN_strat_combo_iii} of \Cref{lemma:HKKN_strat_combo}.
\begin{lemma}\label{lemma:HKKN_condition_iii_technical_step} Let $\Gamma$ be a chamber of an arrangeable secondary fan $\Sigma_{GKZ}$. Let $\lambda\in \Gamma^\vee \cap \Lambda^{\circ}$. Let $K\subset I_{\lambda}^{\le 0}$ be such that for every $\xi \wrel \lambda$, $I_\xi^{>0}\not\subset I_\lambda^{>0}\sqcup K$. If $\xi\in \Gamma^\vee \cap \Lambda$ is such that $I_{\xi}^{>0}\subset I_{\lambda}^{\ne 0}\cup K$, then $\lambda \wrel \xi$.
\end{lemma}
\begin{proof} Let $\lambda$ and $K$ be as in the statement. For each $\rho\in I_{\lambda}^{=0}$, by \Cref{lemma:internally_simplicial_implies_close_previous_wall}, there exists $\xi \wrel \lambda$ so that $I_{\xi}^{>0}\cap I_{\lambda}^{\le 0} = [\rho]$. So $I_\xi^{>0}\not\subset I_\lambda^{>0}\sqcup K$ implies that $[\rho] \not\subset K$. Now suppose there were some $\xi\in \Gamma^\vee \cap \Lambda$ so that $I_\xi^{>0}\subset I_\lambda^{\ne 0}\cup K$. If $I_\xi^{>0}\cap I_\lambda^{=0} \ne \emptyset$, then there is some $\rho\in I_\xi^{>0}\cap I_\lambda^{=0}$ and $[\rho]\subset I_\xi^{>0}\cap I_\lambda^{=0}$. Because $I_\xi^{>0}\subset I_\lambda^{\ne 0}\cup K$, this would imply $[\rho]\subset K$, contradicting $[\rho]\not\subset K$. Thus $I_\xi^{>0}\subset I_\lambda^{\ne 0}$, so $\lambda \wrel \xi$. 
\end{proof}

With \Cref{lemma:HKKN_condition_iii_technical_step}, we are now in a position to show that any extension of $(\Gamma^\vee \cap \Lambda, \wrel)$ to a total ordering induces an HKKN stratification of $V(B_\Gamma)$.
\begin{theorem}\label{thm:existence_of_wall_stratifications} Let $\Gamma$ be a chamber of an arrangeable secondary fan $\Sigma_{GKZ}$. Then any extension of the poset $(\Gamma^\vee\cap \Lambda, \wrel)$ to a total ordering $\Gamma^\vee \cap \Lambda = \{\lambda_1,\hdots, \lambda_k\}$ induces an HKKN stratification of $V(B_\Gamma)$.
\end{theorem}
\begin{proof} Let $\lambda_1,\hdots, \lambda_k$ be any totalization of $(\Gamma^\vee \cap \Lambda, \wrel)$ guaranteed by the Szpilrajn extension theorem \cite{Szpilrajn}, i.e. $\Gamma^\vee\cap \Lambda = \{\lambda_1,\hdots, \lambda_k\}$ and if $\lambda_i \wrel \lambda_j$ for any distinct $\lambda_i,\lambda_j\in \Gamma^\vee \cap \Lambda$, then $i < j$. To show that $\lambda_1,\hdots, \lambda_k$ induces an HKKN stratification of $V(B_\Gamma)$, we need only show $\lambda_1,\hdots, \lambda_k$ satisfy properties \ref{lemma:HKKN_strat_combo_i}, \ref{lemma:HKKN_strat_combo_ii}, and \ref{lemma:HKKN_strat_combo_iii} of \Cref{lemma:HKKN_strat_combo}. Property \ref{lemma:HKKN_strat_combo_i} is satisfied by \Cref{lemma:arrangeable_implies_plus_ideal_equals_irr}. Property \ref{lemma:HKKN_strat_combo_ii} is satisfied by condition \ref{def:wall_partial_order_i} of the partial order $\wrel$ and the fact that the ordering $\lambda_1,\hdots,\lambda_k$ extends $(\Gamma^\vee \cap \Lambda, \wrel)$. To see that property \ref{lemma:HKKN_strat_combo_iii} is satisfied, let $K\subset I_{\lambda_i}^{\le 0}$ and for all $j < i$ suppose $I_{\lambda_j}^{>0}\not\subset I_{\lambda_i}^{>0}\sqcup K$. If $\lambda_i\in |\Sigma_{GKZ}|^\vee$, then $I_{\lambda_i}^{\ne 0} = I_{\lambda_i}^{>0}$, so for any $j < i$ we get
\[I_{\lambda_j}^{>0}\not\subset I_{\lambda_i}^{>0}\sqcup K = I_{\lambda_i}^{\ne 0}\sqcup K\]
as desired. In the other case, suppose $\lambda_i\notin |\Sigma_{GKZ}|^\vee$. Then as $\lambda_j\wrel \lambda_i$ implies $j < i$, we have that for every $\lambda_j \wrel \lambda_i$, $I_{\lambda_j}^{>0}\not\subset I_{\lambda_i}^{>0}\sqcup K$. Thus by \Cref{lemma:HKKN_condition_iii_technical_step}, if $I_{\lambda_j}^{>0}\subset I_{\lambda_i}^{\ne 0}\cup K$, then $\lambda_i \wrel \lambda_j$, so $i < j$. So for any $j < i$, $I_{\lambda_j}^{>0}\not\subset I_{\lambda_i}^{\ne 0}\cup K$ as desired. Thus $\lambda_1,\hdots,\lambda_k$ satisfies property \ref{lemma:HKKN_strat_combo_iii} of \Cref{lemma:HKKN_strat_combo} so we conclude that $\lambda_1,\hdots, \lambda_k$ is an HKKN stratification of $V(B_\Gamma)$.

\end{proof}
For chambers $\Gamma$ of non-arrangeable problems, the relation $\wrel$ on $\Gamma^\vee \cap \Lambda$ is generally not a poset as illustrated in the following example.
\begin{example}\label{example:relation_not_a_poset} Consider the GIT problem defined by $G\cong \Gm^3$ including into $\Gm^7$ via weights
\[\begin{array}{cccc} \beta_0 =\beta_1 = (0,0,1), & \beta_2 =\beta_3 = (0,2,1), & \beta_4 = (1,1,1),\\ \beta_5 = (2,1,1), & \beta_6 = (3,0,1).
\end{array}\]
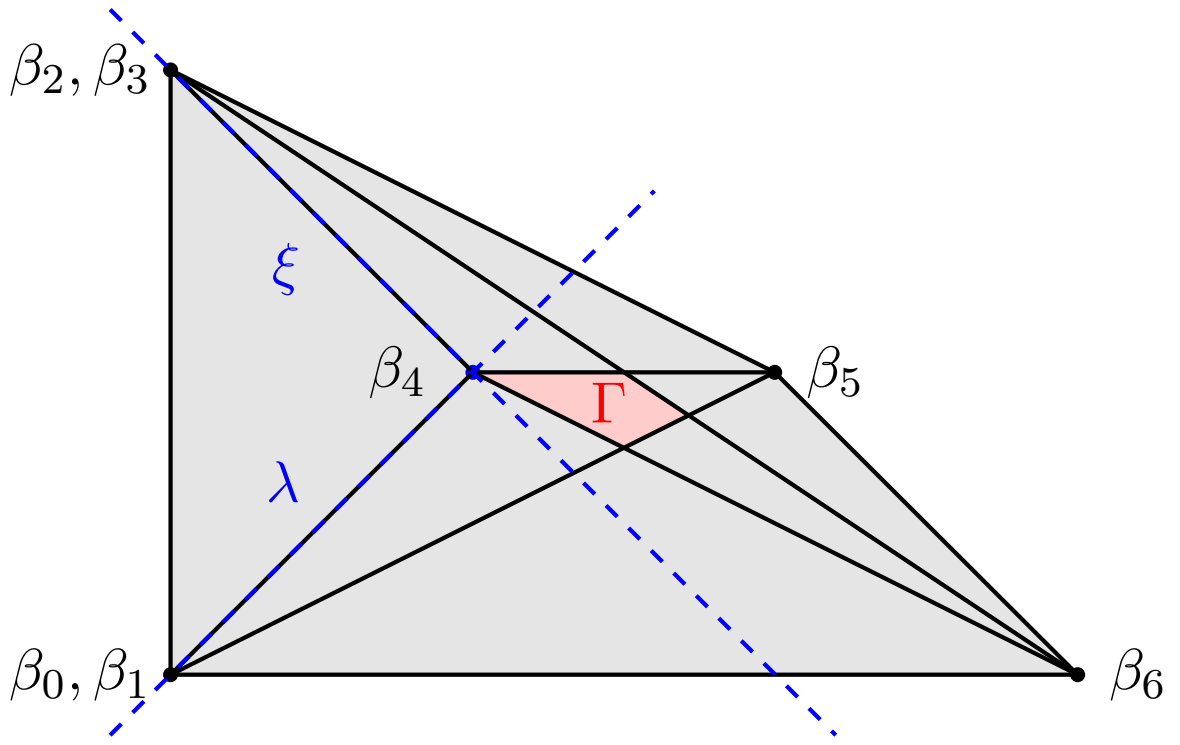
\begin{figure}
\begin{tikzpicture}
%%% Fan 1

\def\fanscale{2}
\def\xshift{4}

\coordinate (beta2) at (0,2*\fanscale);
\coordinate (beta6) at (3*\fanscale,0);
\coordinate (beta0) at (0,0);
\coordinate (beta3) at ($(beta2)$);
\coordinate (beta4) at (1*\fanscale,1*\fanscale);
\coordinate (beta1) at ($(beta0)$);
\coordinate (beta5) at ($0.666*(beta6)+0.5*(beta2)$);

\coordinate (alpha0) at (1.5*\fanscale,1*\fanscale);
\coordinate (alpha1) at (1.5*\fanscale,0.75*\fanscale);
\coordinate (alpha2) at (1.714*\fanscale,0.857*\fanscale);

\draw[draw=none,fill=gray,opacity=.2] (beta1) -- (beta2) -- (beta5) -- (beta4)--(beta6)--cycle;

\draw[draw=none,fill=gray,opacity=.2] (beta5) -- (alpha0) -- (alpha2) -- (alpha1)--(beta6)--cycle;

\draw[draw=none,fill=red,opacity=.2] (beta4) -- (alpha0) -- (alpha2) -- (alpha1)--cycle;

\node[circle,fill,inner sep=1pt,black] at (beta0) {};
%\node[circle,fill,inner sep=1pt] at (beta1) {};
\node[circle,fill,inner sep=1pt,black] at (beta2) {};
%\node[circle,fill,inner sep=1pt] at (beta3) {};
\node[circle,fill,inner sep=1pt,blue] at (beta4) {};
\node[circle,fill,inner sep=1pt] at (beta5) {};
\node[circle,fill,inner sep=1pt] at (beta6) {};

\draw[-,black,thick] (beta1) -- (beta4);
\draw[-,black,thick] (beta2) -- (beta4);
\draw[-,black,thick] (beta0) -- (beta2);
\draw[-,black,thick] (beta1) -- (beta5);
\draw[-,black,thick] (beta4) -- (beta6);
\draw[-,black,thick] (beta0) -- (beta6);
\draw[-,black,thick] (beta5) -- (beta6);
\draw[-,black,thick] (beta2) -- (beta6);
\draw[-,black,thick] (beta4) -- (beta5);
\draw[-,black,thick] (beta2) -- (beta5);

\draw[-,blue,thick,dashed]($(beta1)!-0.2!(beta4)$) -- ($(beta1)!1.6!(beta4)$);
\draw[-,blue,thick,dashed]($(beta2)!-0.2!(beta4)$) -- ($(beta2)!2.2!(beta4)$);

\node[draw = none, fill = none, black] at ($(beta0)+(-0.6,0)$) {\small $\beta_0,\beta_1$};

\node[draw = none, fill = none, black] at ($(beta2)+(-0.6,0)$) {\small $\beta_2,\beta_3$};

\node[draw = none, fill = none, black] at ($(beta4)+(-0.5,0)$) {\small $\beta_4$};

\node[draw = none, fill = none, black] at ($(beta5)+(0.4,0)$) {\small $\beta_5$};

\node[draw = none, fill = none, black] at ($(beta6)+(0.4,0)$) {\small $\beta_6$};
\node[draw = none, fill = none, red] at ($(beta4)+(0.9,-0.2)$) {\small $\Gamma$};

\node[draw=none,fill=none,blue, above left] at ($(beta1)!0.5!(beta4)$) {\small $\lambda$};
\node[draw=none,fill=none,blue, below left] at ($(beta2)!0.5!(beta4)$) {\small $\xi$};
\end{tikzpicture}
    \caption{Secondary fan $\Sigma_{GKZ}$ and chamber $\Gamma$ associated to the weights described in \Cref{example:relation_not_a_poset} demonstrating that $(\Gamma^\vee \cap \Lambda, \wrel)$ need not define a poset when $\Sigma_{GKZ}$ is not arrangeable. The hyperplanes associated to the 1-PSs $\lambda,\xi$ from \Cref{example:relation_not_a_poset} are also depicted.}
    \label{fig:relation_not_a_poset}
\end{figure}
The associated secondary fan is depicted in \Cref{fig:relation_not_a_poset}. Let $\Gamma$ be the highlighted chamber in the figure. Let $\lambda$ be the primitive 1-PS such that $I_{\lambda}^{=0} = \{0,1,4\}$ and $\lambda\in \Gamma^\vee$ and let $\xi$ be the primitive 1-PS such that $I_\xi^{=0} = \{2,3,4\}$ and $\xi\in \Gamma^\vee$. Then $\lambda,\xi\in \Gamma^\vee \cap \Lambda$, but $I_{\lambda}^{>0} = I_\xi^{>0} = \{5,6\}$. This implies that there are two type \ref{def:wall_partial_order_i} relations $\lambda \wrel \xi$ and $\xi \wrel \lambda$. Thus $(\Gamma^\vee\cap \Lambda, \wrel)$ is not antisymmetric and thus not a poset.
\end{example}

\section{Windows}\label{section: windows}
In this section we introduce the two window categories appearing in our main statement and provide a condition for when the two constructions are equal. In subsection \S\ref{subsection:grade_restriction_windows} we define grade restriction windows. In subsection \S\S\ref{subsection:cox_windows} we define Cox windows. In \S\ref{subsection:nef_generation} we construct generators for $\D(\CX_\Gamma)$ and in \S\ref{subsection:cox_generator_images} we describe the images of these generators under the Cox window functor.

Before continuing we recall the definition of a window. Let $X\subset \A^n= \Spec S$ be a $G$-equivariant open subset. 
\begin{definition}\label{def:window} A full subcategory $W\subset \Db_{\wh{G}}(S)$ is a \emph{window} for $X\subset \Spec S$ if the restriction $r_X^\ast:W\to \D([X/G])$ is an equivalence.
\end{definition}
\subsection{Grade restriction windows}\label{subsection:grade_restriction_windows} In this subsection we introduce the grade restriction windows, largely following the treatment of grade restriction rules from \cite{BFK} but our definition of a grade restriction window more clearly mirrors the definition in \cite{HL}.

As before, let $G\subset \Gm^n$ be a reductive linear algebraic group. Let $\CF_\bullet\in \Db([\A^n/G])\cong \Db_{\wh{G}}(S)$. Let $\lambda$ be a 1-PS and $V$ be an open subvariety of $(\A^n)^{=0}_\lambda$ and let $i:V\inc \A^n$ be the associated inclusion. Choose a free replacement $\CE_\bullet$ of $\CF_\bullet$. Because open immersions are flat, the derived pullback need not be derived. We define the subcomplex $(i^*\CE_\bullet)_{\lambda = k}$ of $i^*\CE_\bullet$ to be the subcomplex consisting of the terms $\CO_V(-\chi)$ such that $\wh{\lambda}(\chi) = k$. This indeed defines a complex because all morphisms between terms $\CO_V(-\chi)$, $\CO_V(-\alpha)$ where $\wh{\lambda}(\chi) \ne \wh{\lambda}(\alpha)$ are 0 on $V$, as $V\subset (\A^n)^{=0}_\lambda$. Further, $i^*\CE_\bullet$ splits $G$-equivariantly as a direct sum of these subcomplexes:
\begin{equation}\label{eqn:splitting}
    i^*\CE_\bullet \cong \bigoplus_{k\in \Z}(i^*\CE_\bullet)_{\lambda =k}.
\end{equation}
We define the \emph{$\lambda$-weights of $\CF_\bullet$ along $V$} to be the set of $k\in \Z$ so that $(i^*\CE_\bullet)_{\lambda = k}$ is not quasi-isomorphic to 0. We will denote the $\lambda$-weights of $\CF_\bullet$ along $V$ by $\mu(\CF_\bullet, \lambda, V)$. There is a more general definition of $\lambda$-weights for complexes on arbitrary smooth stacks $[X/G]$ (see e.g. \cite[Definition 3.1.1]{BFK}); however, in this work we are only concerned with computing the $\lambda$-weights of complexes of free $\wh{G}$-graded $S$-modules (line bundles on a stack $[\A^{n}/G]$), so we use the above simplified description.

We are now in a position to introduce one of the two main categories of interest in this paper, the grade restriction windows. Before doing so, we define the \emph{positive width} and \emph{negative width} respectively of a 1-PS $\lambda$ to be
\[T_\lambda^{>0}\coloneqq \sum_{\rho\in I_{\lambda}^{>0}}\wh{\lambda}(\beta_\rho),\quad \quad T_\lambda^{<0}\coloneqq \sum_{\rho\in I_{\lambda}^{<0}}\wh{\lambda}(\beta_\rho).\]
\begin{definition}\label{def:grade_restriction_windows} Let $G$ be as above and $X\subset \A^n$ an open $G$-equivariant subvariety. Let $\bs{\lambda} = \lambda_1,\hdots, \lambda_\ell$ induce an HKKN stratification of $\A^n\setminus X$. For each $i =1,\hdots, \ell$ choose an integer $w_i\in \Z$. We define the \emph{grade restriction window} associated to $\bs{\lambda}$, $\bs{w}$ to be the full subcategory of $\Db([\A^n/G])$
\begin{equation*}
W_{\bs{\lambda}}^{\bs{w}} \coloneqq \big\{\CF_\bullet \in \Db([\A^n/G]) \ | \ \mu(\CF_\bullet, \lambda_{i}, (X_{i})^{=0}_{\lambda_i}) \subset [w_i, w_i+T_{\lambda_i}^{>0})\text{ for }i=1,\hdots, \ell\big\} .
\end{equation*}
Often we take $\bs{w}=\bs{0}$, so we let $W_{\bs{\lambda}}$ denote $W_{\bs{\lambda}}^{\bs{0}}$. We say $\CF_\bullet$ \emph{satisfies the $\lambda_i$ grade restriction rule} if $\mu(\CF^{\bullet},\lambda_i,(X_{i})^{=0}_{\lambda_i})\subset [w_i,w_i+T^{>0}_{\lambda_i})$.
\end{definition}
The following theorem of Halpern-Leistner, explains why $W_{\bs{\lambda}}^{\bs{w}}$ are referred to as grade restriction \textit{windows}.
\begin{theorem}[Theorem 1.1, \cite{HL}]\label{thm:DHL} Let $G$ be as above and $X\subset \A^n$ an open $G$-equivariant subvariety and let $\bs{\lambda} = \lambda_1,\hdots, \lambda_\ell$ induce an HKKN stratification of $\A^n\setminus X$. Then for any $\bs{w}\in \Z^\ell$, $W_{\bs{\lambda}}^{\bs{w}}$ is a window for $\Db([X/G])$, i.e. the restriction functor $r_X^*:W_{\bs{\lambda}}^{\bs{w}}\to \Db([X/G])$ induces an equivalence of categories.
\end{theorem}
\begin{remark}\label{rmk:weights_sign_discrepancy} Our description of the weights follows the description from \cite{BFK} and is negative of the convention taken by Halpern-Leistner in \cite{HL}. However, because our stratifications are negative of those used by Halpern-Leistner (see \Cref{rmk:1PS_sign_difference}), we have the exact statement of the main theorem of \cite{HL} in our setup.
\end{remark}

Generally, computing the $\lambda$-weights of an arbitrary complex $\CF_\bullet$ can be challenging. However, if $\CF_\bullet$ is in addition $T\coloneqq \Gm^n$-equivariant, there is a concrete, combinatorial recipe for this computation. Suppose $\Gamma$ is a chamber of a $\Sigma_{GKZ}$ defined by $G\subset \Gm^n$. Let $\lambda_1,\hdots, \lambda_k$ be a collection of 1-PSs which induce an HKKN stratification of $V(B_\Gamma)$. Recall
\[
(X_{i})^{=0}_{\lambda_i} =  \A^n\setminus \bigcup_{j < i}(\A^n)_{\lambda_j}^{\le 0}.
\]
Since $\bigcup_{j < i}(\A^n)_{\lambda_j}^{\le 0}$ is a union of $T$-orbits, $X_{i}$ is $T$-equivariant and so is the fixed locus of $\lambda_{i}$ in $X_{i}$. For any $L\subset [n]$, we define the following distinguished point $p_L$ in $\A^n$ as
\[(p_{L})_t \coloneqq \begin{cases} 1 & t\in L \\
0 & t \in [n]\setminus L
\end{cases}\]
Each $T$-orbit contains exactly one distinguished point $p_L$. Because $(X_{i})^{=0}_{\lambda_i}$ is $T$-equivariant, we can write $(X_{i})^{=0}_{\lambda_i}$ as a finite disjoint union of $T$-orbits:
\begin{equation}\label{eq:splittingintoTorbits}
  (X_{i})^{=0}_{\lambda_i} = \bigsqcup\limits_{j=1}^{r} T \cdot p_{L_{j}}  
\end{equation}
for $L_1,\hdots, L_r\subset [n]$. Finally, abusing notation, let $p_{L_j}$ denote the embedding of $p_{L_j}$ in $(X_{i})^{=0}_{\lambda_i}$.

\begin{proposition}\label{prop:derivednakayama}
Suppose $\CF \in \Db((X_{i})^{=0}_{\lambda_i} )$ admits a $T$-equivariant structure. Then  $\CF = 0$ if and only if $ L{p^{*}_{L_{j}}} \CF = 0$ for each $j = 1,\dots,r$.
\end{proposition}

\begin{proof} 
\begin{align*}
 \CF = 0 & \iff \supp (\CF) = \varnothing \\
 & \iff \forall j, p_{L_{j}} \notin \supp(\CF)  & \text{since $\CF$ is $T$-equivariant and the $T \cdot p_{L_j}$ cover $(X_{i})^{=0}_{\lambda_i}$}  \\
 & \iff \forall j, Lp^*_{L_j}\CF=0.
\end{align*}
The last equivalance is the application of the Nakayama's Lemma and an induction argument on the length of the complex as in \cite[Lemma 3.9]{HuybrechtsFM}.
\end{proof}

\subsection{Cox windows}\label{subsection:cox_windows}
In this section we review the construction of the Cox category and the Cox windows of \cite[Section 7.2]{CoxCategory}. Throughout this section we consider a projective GIT problem $G\subset \Gm^n$ with associated secondary fan $\Sigma_{GKZ}$. In order to define Cox windows, following the authors of \cite{CoxCategory}, we additionally assume the action is primitive geometric (see \Cref{rmk:primitive_geometric})

\begin{definition}\label{def:zonotope}
We define the half-open zonotope $Z_{G}\subseteq \wh{G}_{\R}$ associated to the inclusion $G\subset \Gm^n$ 
as the image of $[0,1)^{n}$ under the map $\pi_{\R}$ in (\ref{eqn:ses}).
\end{definition}

\begin{definition}\label{defn:BT}
The \emph{Bondal--Thomsen collection} for the
inclusion $G\subset \Gm^n$ is the set of characters $\Theta \subseteq \wh{G}$ that lie in the preimage of $-Z_{G}$ under the natural map:
\begin{align*}
    \rho:\wh{G}&\to\wh{G}_{\R} \\
    \chi &\mapsto \chi\otimes 1
\end{align*}
\end{definition}

\begin{example}\label{example:H2_BT_degrees} Following weights from \Cref{example:H2_stratifications}, we draw the associated zonotope $Z_{G}$. The lattice points in $Z_{G}$ correspond to elements in $-\Theta$. In the chosen basis for $\wh{G}_\R$, that is 
\[
\Theta = \{(0,0),(-1,0),(0,-1),(-1,-1),(1,-1)\}.
\]

\[\begin{tikzpicture}
\def\fanscale{1.5}
\def\xshift{0}
\def\yshift{0}

\coordinate (origin) at (\xshift,\yshift);
\coordinate (beta13) at (\xshift+\fanscale*1,\yshift);
\coordinate (beta2) at (\xshift+\fanscale*-2,\yshift+\fanscale*1);
\coordinate (beta4) at (\xshift, \yshift+1*\fanscale*1);
\coordinate (ray13) at (\xshift+3*\fanscale,\yshift);
\coordinate (ray2) at (\xshift-3*\fanscale,\yshift+1.5*\fanscale);
\coordinate (ray4) at (\xshift+0*\fanscale,\yshift+3*\fanscale);
\coordinate (zonul) at (\xshift-2*\fanscale, \yshift+2*\fanscale);
\coordinate (zonum) at (\xshift,\yshift+2*\fanscale);
\coordinate (zonur) at (\xshift+2*\fanscale,\yshift+1*\fanscale);
\coordinate (zonlr) at (\xshift+2*\fanscale, \yshift);
\coordinate (d1) at (origin);
\coordinate (d2) at (\xshift+1*\fanscale, \yshift+0*\fanscale);
\coordinate (d3) at (\xshift-1*\fanscale, \yshift+1*\fanscale);
\coordinate (d4) at (\xshift+0*\fanscale, \yshift+1*\fanscale);
\coordinate (d5) at (\xshift+1*\fanscale, \yshift+1*\fanscale);

% betas
\def\betaop{0.5}
\node[circle,fill,inner sep=1pt,opacity=\betaop] at (beta13) {};
\node[circle,fill,inner sep=1pt,opacity=\betaop] at (beta2) {};
\node[circle,fill,inner sep=1pt,opacity=\betaop] at (beta4) {};
\node[draw = none, fill = none,opacity=\betaop] (betas13label) at ($(beta13)+(0,-0.3)$) {\small $\beta_1,\beta_3$};
\node[draw = none, fill = none,opacity=\betaop] (beta3label) at ($(beta4)+(0.3,0)$) {\small $\beta_4$};
\node[draw = none, fill = none,opacity=\betaop] (beta1label) at ($(beta2)+(-0.3,-0.2)$) {\small $\beta_2$};

% zonotope
\draw[fill=red,opacity=0.1,draw = none] (origin) -- (beta2) -- (zonul) -- (zonum) -- (zonur) -- (zonlr) -- cycle;

% rays
\draw[-arr,black] (origin) -- (ray13);
\draw[-arr,black] (origin) -- (ray2);
\draw[-arr,black] (origin) -- (ray4);

% border
\draw[red,thick] (beta2) -- (origin) -- (zonlr);
\draw[red,thick,dotted] (beta2) -- (zonul) -- (zonum) -- (zonur) -- (zonlr);

% BT degrees
\node[circle,fill,inner sep=1.5pt,red] at (d1) {};
\node[circle,fill,inner sep=1.5pt,red] at (d2) {};
\node[circle,fill,inner sep=1.5pt,red] at (d3) {};
\node[circle,fill,inner sep=1.5pt,red] at (d4) {};
\node[circle,fill,inner sep=1.5pt,red] at (d5) {};

\node[draw=none,fill=none,black] at (\xshift+\fanscale*2,\yshift+\fanscale*2) {\small $\wh{G}_\R\cong \R^2$};
\node[draw=none,fill=none,red] at (\xshift+\fanscale*-1,\yshift+\fanscale*1.5) {\small $Z_G$};
\end{tikzpicture}\]
\end{example}

\begin{definition}
The \emph{Cox category}, denoted $K_{\Theta}$, (following \cite{CoxCategory}), is the full subcategory of $\D([\A^{n}/G])$ generated by line bundles $\CO(\chi)$ for $\chi\in \Theta$. Equivalently, it is a full subcategory of $D_{\wh{G}}(S)$ generated by $S(\chi)$ with $\chi\in \Theta$.
\end{definition}

For each chamber $\Gamma$ in $\Sigma_{GKZ}$, the authors of \cite{CoxCategory} construct a window which is itself a subcategory of $K_{\Theta}\subset \Db([\A^n/G])$. We will provide an outline of this construction closely following that of \cite{CH25}. For a more detailed exposition with proofs see \cite{HHL}. First consider the cell structure on $\R^n$ whose codimension 1 skeleton is the union of all integral translates of coordinate hyperplanes in $\Z^n$. We then restrict this cell structure to the vector space $i_\R(M_\R)\subset \R^n$ (see eqn. \ref{eqn:ses}). This cell structure descends further to the real torus $i_\R(M_\R)/i(M)$. Let $E$ denote the set of cells in the quotient $i_\R(M_\R)/i(M)$, $E_j$ the set of $j$-dimensional cells in $E$, and $\ve$ the associated incidence function with $\ve(\delta,\tau) = \pm 1$ depending on the chosen orientation when $\tau$ is a boundary cell of $\delta$.

From the data of this cell structure $E$, we will construct a complex of sheaves in $\Db(\CX_\Gamma \times [\A^n/G])$ which we call the \emph{HHL complex} after work of Hanlon, Hicks, Lazarev. For any $\gamma\in i_\R(M_\R)$, $\pi(\ceil{\gamma})$ is invariant up to shifts by elements of $i(M)$ by the exact sequence in eqn. \eqref{eqn:ses}, thus $\pi(\ceil{-})$ is a well-defined stratification $i_\R(M_\R)/i(M)\to \wh{G}$. Further, $\pi(\ceil{-})$ is constant on cells $\delta\in E$ and thus descends to a labeling $\pi(\ceil{-}):E\to \wh{G}$. Similarly, $\pi(\floor{-}):E\to \wh{G}$ is a well-defined assignment. 
Let $\HHL_\Gamma\in \Db(\CX_\Gamma \times [\A^n/G])$ be the complex whose component in homological degree $j$ is given by
\[\bigoplus_{\delta\in E_{j}} \CO_{\CX_\Gamma}(\pi(\floor{\delta}))\boxtimes S(-\pi(\ceil{\delta})).\]
The differential from a term labeled by $\delta$ to the term labeled by $\tau$ is nonzero only if $\tau$ is a boundary cell of $\delta$. For such $\delta, \tau$, choose lifts $\wt{\delta},\wt{\tau}$ in the cell structure on $i_\R(M_\R)$ so that $\wt{\tau}$ is a boundary of $\wt{\delta}$. Then the morphism from the term labeled by $\delta$ to the term labeled by $\tau$ is given by
\[\ve(\delta,\tau)x^{\floor{\wt{\tau}}-\floor{\wt{\delta}}}\otimes x^{\ceil{\wt{\delta}}-\ceil{\wt{\tau}}} .\]
Note that the twists appearing in the HHL complex belong to $\Theta$. With the description of the HHL resolution $\HHL_\Gamma$ in hand, we are now ready to define the Cox window functor as a Fourier-Mukai transform. Let $p,q$ be the projection morphisms from $\CX_\Gamma\times [\A^n/G]$ onto its two factors:
\[\begin{tikzcd}
    & \CX_\Gamma\times [\A^n/G] \arrow[ld,"p"'] \arrow[rd,"q"]& \\ \CX_\Gamma & & {[\A^n/G]}
\end{tikzcd}\]
\begin{definition}\label{def:cox_window_functor} Let $\Gamma$ be a chamber of $\Sigma_{GKZ}$ associated to the inclusion $G\subset \Gm^n$. Then the \emph{Cox window functor} $\Phi_\Gamma$ is defined to be the (stacky/equivariant) Fourier-Mukai transform (see e.g. \cite[Section~1.2]{Ploog2007}) with kernel $\HHL_\Gamma$:
\begin{equation}\label{def:FMfunctor}
\Phi_\Gamma \coloneqq Rq^{G}_*\lp Lp^*(-)\otimes_L \HHL_\Gamma\rp: \Db(\CX_\Gamma)\to \Db([\A^n/G]).
\end{equation}
We define the \emph{Cox window} associated to $\Gamma$ to be the essential image of $\Phi_\Gamma$. We denote the Cox window by $W_\Gamma^{\Cox}$.
\end{definition}
\begin{theorem}[Theorem 7.6, \cite{CoxCategory}]\label{thm:cox_window_exists} For any chamber $\Gamma$ of $\Sigma_{GKZ}$, $W_\Gamma^{\Cox}$ is a subcategory of  $K_{\Theta}$ and a window for $\CX_\Gamma$.
\end{theorem}

\subsection{Generation of Cox Windows}\label{subsection:nef_generation} 
The goal of this section is to recall 
a generation statement for $\CX_{\Sigma_{\Gamma},i^{\vee}_\Gamma} \simeq \CX_{\Gamma}$.

\begin{proposition}\label{prop:amplesgenerate}
The collection of all nef line bundles  on $\CX_{\Sigma_{\Gamma},i^{\vee}_\Gamma}$ (see \Cref{def:nefdivisorsonastack})
generates $D(\CX_{\Sigma_{\Gamma},i^{\vee}_\Gamma})$.
\end{proposition}

\begin{proof}
By \cite[Lemma 4.12]{CoxCategory}, there exist characters
$d_1,\ldots,d_\ell$ in $\wh{G}_{\Gamma}$ such that
\[
\pi^*\bigl(\operatorname{Perf}(X_{\Sigma_{\Gamma}})\bigr)
\otimes
\bigoplus_{i=1}^{\ell}\mathcal{O}_{\CX_{\Sigma_{\Gamma},i^{\vee}_\Gamma}}(d_i)
\]
generates $D(\CX_{\Sigma_{\Gamma},i^{\vee}_\Gamma})$. Let $d$ be any ample line bundle on the coarse moduli space $X_{\Sigma_{\Gamma}}$. Pick a sufficiently large $k$, so that $k\cdot d - d_{i} \in \text{Nef}(X_{\Sigma_{\Gamma}})$ for $i = 1, \dots,l$. The collection of line bundles $\{\CO_{X_{\Sigma_{\Gamma}}}(\lambda(k\cdot d)),\lambda\geq1\}$ generates $\operatorname{Perf}(X_{\Sigma_{\Gamma}})$ by \cite[Theorem~4]{Orlov2009}. Since $b_{\rho}=1$ for all $\rho \in \Sigma_{\Gamma}(1)$, we have that $\pi^{*}(\CO_{X_{\Sigma_{\Gamma}}}(d))\cong \CO_{\CX_{\Sigma_{\Gamma},i^{\vee}_\Gamma}}(d)$. Then, the collection of line bundles:
\[
\{\CO_{\CX_{\Sigma_{\Gamma},i^{\vee}_\Gamma}}(\lambda(k\cdot d)-d_i)\mid \lambda\ge1,i=1,\dots,l\}
\]
generates  $D(\CX_{\Sigma_{\Gamma},i^{\vee}_\Gamma})$ by \cite[Lemma 4.12]{CoxCategory} and lies in $\text{Nef}(X_{\Sigma_{\Gamma}})$ by construction, as desired.
\end{proof}

\begin{corollary}\label{prop:chamber_characters_generate}
The collection of line bundles: 
\[
\{ \CO_{\CX_{\Gamma}}(d), d \in \Gamma\cap \wh{G}\} 
\] 
generates $D(\CX_{\Gamma})$.
\end{corollary}
\begin{proof}
This follows directly by applying \Cref{prop:amplesgenerate} and \Cref{prop:every nefisisomorphictooneinthechamber} to \Cref{cor:image_of_line_bundles_under_stack_iso}.
\end{proof}

Here we state a useful vanishing result about this generating set of line bundles.

\begin{proposition}\label{prop:nefcohomologyvanishes}
Let $d \in \Gamma$ and $\chi \in \Theta$. Then, one has
\[ H^p(\CX_\Gamma,\CO_{\CX_\Gamma}(d+\chi)) = \begin{cases} S_{d+\chi} & \text{if } p = 0, \\
0 & \text{if }p \neq 0, \end{cases}\]
where $S_{d+\chi}$ is the vector space of degree $d+\chi$ elements of $S$.
\end{proposition}

\begin{proof}
Pick any lift $D=\sum_{\rho\in [n]}a_\rho e_\rho$ of $d$. Since $\chi \in \Theta$, by one of the equivalent characterizations of the Bondal-Thomsen collection (see \cite[Definition 2.12]{CoxCategory}), there exists $\theta \in M_{\R}$ such that the formal linear combination $D_{\chi}:=\sum_{\rho \in [n]}
\left\lfloor
\left\langle -\theta, u_{\rho} \right\rangle
\right\rfloor e_\rho$ is a lift of $\chi$. 

We first prove vanishing of higher sheaf cohomology. By the identification of line bundles in  \Cref{cor:image_of_line_bundles_under_stack_iso} arising from the stack isomorphism of \Cref{prop:stack_isomorphism}, we have for all $p$:
\[H^p(\CX_\Gamma,\CO_{\CX_\Gamma}(d+\chi)) = H^p(\CX_{\Sigma_\Gamma, i_\Gamma^\vee},\CO_{\CX_{\Sigma_\Gamma, i_\Gamma^\vee}}(\proj_{\wh{G}_{\Gamma}} (d+\chi))). \]
Since $D + D_{\chi}$ is a lift for $d+\chi$, by commutativity of \eqref{eqn:exact_square}, \[ \proj_{\Gamma}(D + D_{\chi})=\sum_{\rho\in \Sigma_{\Gamma}(1)}a_\rho e_\rho + \sum_{\rho \in \Sigma_{\Gamma}(1)}
\left\lfloor
\left\langle -\theta, u_{\rho} \right\rangle
\right\rfloor e_\rho\]
is a lift for $\proj_{\wh{G}_{\Gamma}} (d+\chi)$. In particular, observe that $\proj_{\wh{G}_{\Gamma}} (\chi)$ is a Bondal-Thomsen element for $\wh{G}_{\Gamma}$. Since $d\in \Gamma$, $\proj_{\wh{G}_{\Gamma}}(d)$ is nef (see \Cref{prop:every nefisisomorphictooneinthechamber}).  Hence, the higher vanishing follows directly by \cite[Proposition 4.15]{CoxCategory}.  
To show the equality of zero sections, we use the following set of equalities.
\begin{align*}
H^0(\CX_{\Sigma_\Gamma, i_\Gamma^\vee},\CO_{\CX_{\Sigma_\Gamma, i_\Gamma^\vee}}(\proj_{\wh{G}_{\Gamma}} (d+\chi))) & = k\left\langle (P_{\proj_{\Gamma}(D)}+\theta)\cap M\right\rangle  & \text{ by \cite[Proposition 4.15]{CoxCategory}} \\
&= k\left\langle (P_{D} + \theta)\cap M\right\rangle.  & \text{ by \Cref{prop:nef_polytopes} since $d \in \Gamma$} \\
&= k\left\langle (P_{D+D_{\chi}})\cap M\right\rangle  & \text{ \Cref{lemma:lattice_points_of_shifted_polytope} and $\floor{-x}= -\ceil{x}$}  
\end{align*}
Finally, the bijection between $S_{d+\chi}$ and lattice points of $P_{D+D_{\chi}}$ is standard (see \cite[Proposition 5.4.1]{CLS}).
\end{proof}

\subsection{Images of nef line bundles in Cox Window}\label{subsection:cox_generator_images} In this subsection, we describe the images of line bundles $\CO_{\CX_{\Gamma}}(d)$, $d \in \Gamma\cap \wh{G}$ under the Cox window functor $\Phi_\Gamma$.

We now define one of the central constructions of our paper, closely following the work of \cite{CH25}, where it was introduced.  As before, let $d \in \Gamma\cap \wh{G}$ and $D = \sum_{\rho\in [n]}a_\rho e_\rho$ a lift of $d$ with respect to $\pi$ in \eqref{eqn:ses}.

As before, given a formal linear combination $D = \sum_{\rho\in [n]}a_\rho e_\rho$, we can consider a polytope $P_{D} \subseteq M_{\R}$. It is homeomorphic to its image $i_\R(P_D) = i_\R(M_\R) \cap \R^n_{\ge -D}\subset \R^n$, so we will write $P_D$ for both when the ambient space is clear. Following \cite{CH25}, we can associate a cellular complex of $\Z^{n}$-graded $S$-modules supported on $P_D\subset \R^{n}$.

Consider the CW structure  on $\R^n$ whose codimension 1 skeleton is the union of all integral translates of the coordinate planes in $\Z^n$. Concretely, every cell is of the form:
\[
\delta = \prod_{i \in J}(b_{i}-1,b_{i})\times\prod_{j \in [n]\setminus J}\{b_{j}\}
\]
for some $\bs{b}=(b_1,\dots,b_n)\in\Z^n$ and $J\subseteq[n]$.
We can restrict the cells of this CW complex to the linear space $i_{\R}(M_{\R})$ via the embedding from \Cref{eqn:ses}. Since the boundary of $\R^n_{\ge -D}$ is contained in the coordinate planes, $P_D$ is a subcomplex of $i_{\R}(M_{\R})$. We let $\cells^i(P_D)$ denote the $i$-dimensional cells of $P_D$ and $\cells(P_D)$ denote all cells of $P_D$. Now for any cell $\delta$ of $P_D$, observe that by the definition, $\ceil{-}$ is constant on $\delta$, so we may define $\ceil{\delta}\in \Z^n$. We then label every cell $\delta\subset P_D$ with the monomial $x^{\ceil{\delta}+D}$. The monomials in the labeling are well-defined because by definition, $\delta + D \subseteq \R_{\geq0}^{n}$, and $\ceil{a_{i}+\delta_{i}}\le\ceil{a_{i}}+\ceil{\delta_{i}}= a_{i} +\ceil{\delta_{i}}$. Moreover, by
\cite[Lemma 2.3]{CH25} this labeling is a least common multiple labeling in the sense of \cite{BS98}.
Consequently, we have the following classical construction.

\begin{construction}
Choose an incidence function $\ve$ on $P_{D}$. Define the cellular complex of free $S$-modules $F^{D}_{\bullet}$ associated to the labeled cell complex $P_{D}$ as in \cite[\S1]{BS98}.
\end{construction}

Let us describe the resulting complex concretely. Given the regular CW complex $P_{D}$ and a choice of an incidence function $\ve$, we start with an oriented chain complex  $C_{\bullet}(P_D,k)$:
\begin{equation}\label{eq:cellcomplexonPd}
 C_i(P_D,k) = \bigoplus_{\delta\in \cells^{i}(P_D)}k\cdot e_\delta   
\end{equation}
with the differential
\[\del_i(e_\delta) = \bigoplus_{\tau\in  \cells^{i-1}(P_D)}\ve(\delta,\tau)e_\tau.
\]
This complex is the ordinary cellular homology complex for $P_{D}$. Using the ceiling labeling on the cells, we can lift the cellular homology complex to a complex of $\Z^{n}$-graded $S$-modules as follows. Set
\[
F^{D}_i \coloneqq \bigoplus_{\delta\in \cells^{i}(P_D)}S(-(\ceil{\delta}+D)).
\]

For $\delta\in \cells^{i}(P_D)$ and $\tau\in \cells^{i-1}(P_D)$, the differential between the associated summands is given by
\begin{equation*}
\begin{tikzcd}[column sep = 2cm] \del_i\coloneqq S(-(\ceil{\delta}+D))\arrow[r,"\varepsilon(\delta{,}\tau)x^{\ceil{\delta}-\ceil{\tau}}"] & S(-(\ceil{\tau}+D))
\end{tikzcd}
\end{equation*}
Note that $\tau$ being a boundary cell of $\delta$ ensures that $\ceil{\delta}-\ceil{\tau}\in \Z_{\ge 0}^n$ so that the above differential is well-defined.

\begin{remark}\label{rmk:FisinKtheta}
Observe that since we add $D$ to every monomial label of cells of $P_{D}$, every summand of $F^{D}_{\bullet}$ is twisted by $-D$. If we untwist, the resulting complex, denoted $F^{D}_{\bullet}(D)$, has summands of the form $S(-\ceil{\delta})$. In this paper, we will largely view $F^{D}_{\bullet}(D)$ as a complex of $\wh{G}$-graded $S$-modules by regrading the $\Z^n$-grading via $\pi:\Z^{n} \to \wh{G}$. We want to make a useful observation about the terms of $F^{D}_{\bullet}(D)$ as a complex of $\wh{G}$-graded $S$-modules. Since $\delta \subset i_{\R}(M_{\R})$, we can write $\ceil{\delta} = m + \epsilon$, where $m \in i_{\R}(M_{\R})$ and $\epsilon \in [0,1)^{n}$. Then, $\pi_{\R}(-\ceil{\delta}) =\pi_{\R}(-m+(-\epsilon))= -\pi_{\R}(\epsilon)$ using eqn.~\eqref{eqn:ses}. Hence, by \Cref{defn:BT}, $\pi(-\ceil{\delta}) \in \Theta$ and $F^{D}_{\bullet}(D) \in K_{\Theta}$.
\end{remark}

\begin{example}\label{example:H2_cellular_resolution} We continue with our running example of $\CH_{2}$. We pick an ample divisor $d = \CO(1,1)$ and take a lift $D = (0,0,1,1)$. We choose the same GIT problem, chamber, and lift as in \S4 of \cite{CH25} so that the interested reader may explore Cranton Heller's more detailed treatment of this example. The associated labeled regular CW complex $P_{D}$ is depicted in \Cref{fig:labeled_cellular_complex}. The resulting (twisted) complex $F^{D}_{\bullet}(D)$ is depicted in \Cref{fig:full_cellular_resolution}. The entire (untwisted) complex $F^{D}_{\bullet}$ is found in \S4 of \cite{CH25}.
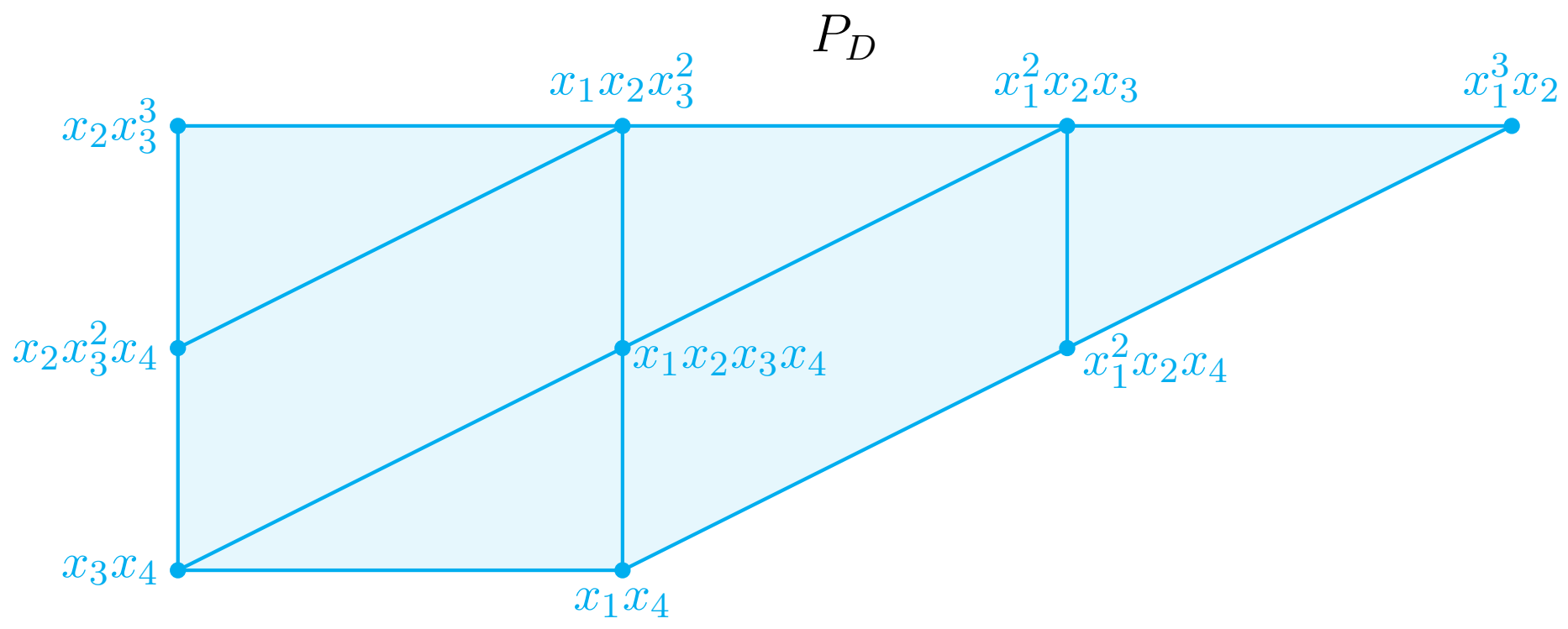
\begin{figure}
    \centering
\[\begin{tikzpicture}
\def\polyscale{3.5}
\def\xshift{0}
\def\yshift{0}
\def\guideop{1} % guideline opacity
\def\op{0.1} % zoomed in 2 cell opacity
\def\cellscale{1.75} % half scale of polyscale
\def\cellop{0.4} % zoomed cell opacity
\def\labelop{1} % opacity of cell labels
\def\guidestyle{dotted}
\def\cellcolor{cyan}
\def\twocellcolor{cyan}
\def\labelcolor{black}
\def\guidecolor{black}
\def\polylabelcolor{cyan}

\coordinate (v00) at (\xshift+0*\polyscale,\yshift+0*\polyscale);
\coordinate (v01) at (\xshift+0*\polyscale, \yshift+0.5*\polyscale);
\coordinate (v02) at (\xshift+0*\polyscale,\yshift+1*\polyscale);
\coordinate (v10) at (\xshift+1*\polyscale,\yshift+0*\polyscale);
\coordinate (v11) at (\xshift+1*\polyscale,\yshift+0.5*\polyscale);
\coordinate (v12) at (\xshift+1*\polyscale,\yshift+1*\polyscale);
\coordinate (v21) at (\xshift+2*\polyscale,\yshift+0.5*\polyscale);
\coordinate (v22) at (\xshift+2*\polyscale,\yshift+1*\polyscale);
\coordinate (v32) at (\xshift+3*\polyscale,\yshift+1*\polyscale);

\draw[-,\cellcolor,thick] (v00) -- (v02);
\draw[-,\cellcolor,thick] (v10) -- (v12);
\draw[-,\cellcolor,thick] (v21) -- (v22);
\draw[-,\cellcolor,thick] (v00) -- (v10);
\draw[-,\cellcolor,thick] (v02) -- (v32);
\draw[-,\cellcolor,thick] (v01) -- (v12);
\draw[-,\cellcolor,thick] (v00) -- (v22);
\draw[-,\cellcolor,thick] (v10) -- (v32);

% fill
\draw[draw=none,fill=\twocellcolor,opacity=\op] (v00) -- (v02) -- (v32) -- (v10) -- cycle;

% labels
\node[draw=none,fill=none,left,\polylabelcolor] at (v00) {\small $x_3x_4$};
\node[draw=none,fill=none,left,\polylabelcolor] at (v01) {\small $x_2x_3^2x_4$};
\node[draw=none,fill=none,left,\polylabelcolor] at (v02) {\small $x_2x_3^3$};
\node[draw=none,fill=none,below,\polylabelcolor] at (v10) {\small $x_1x_4$};
\node[draw=none,fill=none,\polylabelcolor] at ($(v11)+(0.85,-0.1)$) {\small $x_1x_2x_3x_4$};
\node[draw=none,fill=none,\polylabelcolor] at ($(v21)+(0.7,-0.1)$) {\small $x_1^2x_2x_4$};
\node[draw=none,fill=none,above,\polylabelcolor] at (v12) {\small $x_1x_2x_3^2$};
\node[draw=none,fill=none,above,\polylabelcolor] at (v22) {\small $x_1^2x_2x_3$};
\node[draw=none,fill=none,above,\polylabelcolor] at (v32) {\small $x_1^3x_2$};

\def\nodewidth{1pt}
\node[inner sep = \nodewidth, draw=\cellcolor, circle, thick, fill=\cellcolor,opacity=1] at (v00) {};
\node[inner sep = \nodewidth, draw=\cellcolor, circle, thick, fill=\cellcolor,opacity=1] at (v01) {};
\node[inner sep = \nodewidth, draw=\cellcolor, circle, thick, fill=\cellcolor,opacity=1] at (v02) {};
\node[inner sep = \nodewidth, draw=\cellcolor, circle, thick, fill=\cellcolor,opacity=1] at (v10) {};
\node[inner sep = \nodewidth, draw=\cellcolor, circle, thick, fill=\cellcolor,opacity=1] at (v11) {};
\node[inner sep = \nodewidth, draw=\cellcolor, circle, thick, fill=\cellcolor,opacity=1] at (v12) {};
\node[inner sep = \nodewidth, draw=\cellcolor, circle, thick, fill=\cellcolor,opacity=1] at (v21) {};
\node[inner sep = \nodewidth, draw=\cellcolor, circle, thick, fill=\cellcolor,opacity=1] at (v22) {};
\node[inner sep = \nodewidth, draw=\cellcolor, circle, thick, fill=\cellcolor,opacity=1] at (v32) {};

\node[draw=none,fill=none] at (\xshift+1.5*\polyscale,\yshift+1.2*\polyscale) {$P_D$};
\end{tikzpicture}\]
    \caption{Labeled cell complex $P_D$ on $\CH_2$ for $D = (0,0,1,1)$. Every cell $\delta$ is labeled by the monomial $x^{\ceil{\delta}+D}$.}
    \label{fig:labeled_cellular_complex}
\end{figure}

\begin{figure}
\[\hspace{-1.8cm}\begin{tikzpicture}
\def\polyscale{6.5}
\def\cellcolor{cyan}
\def\fillcolor{cyan}
\def\labelcolor{black}
\def\cellop{0.1}
\def\labelop{1}
\def\twocellop{0.03}
\def\arrowop{0.2}
\def\diffop{1}
\def\nodewidth{2pt}
\def\xshift{0}
\def\yshift{0}
\def\thickness{ultra thick}
\def\nodesize{\tiny}
\def\diffthick{thin}
\coordinate (nodeshift) at (-0.02*\polyscale,0.02*\polyscale);

\tikzset{every node/.style={inner sep=0pt}}

\coordinate (v00) at (\xshift+0*\polyscale,\yshift+0*\polyscale);
\coordinate (v01) at (\xshift+0*\polyscale,\yshift+0.5*\polyscale);
\coordinate (v02) at (\xshift+0*\polyscale,\yshift+1*\polyscale);
\coordinate (v10) at (\xshift+1*\polyscale,\yshift+0*\polyscale);
\coordinate (v11) at (\xshift+1*\polyscale,\yshift+0.5*\polyscale);
\coordinate (v12) at (\xshift+1*\polyscale,\yshift+1*\polyscale);
\coordinate (v21) at (\xshift+2*\polyscale,\yshift+0.5*\polyscale);
\coordinate (v22) at (\xshift+2*\polyscale,\yshift+1*\polyscale);
\coordinate (v32) at (\xshift+3*\polyscale,\yshift+1*\polyscale);
\coordinate (v31) at (\xshift+3*\polyscale,\yshift+0.5*\polyscale);

% 1 cell middle coords
\coordinate (w0001) at ($(v00)!0.5!(v01)$); % wijkl is between vij and vkl
\coordinate (w0102) at ($(v01)!0.5!(v02)$);
\coordinate (w0011) at ($(v00)!0.5!(v11)$);
\coordinate (w0112) at ($(v01)!0.5!(v12)$);
\coordinate (w0212) at ($(v02)!0.5!(v12)$);
\coordinate (w0010) at ($(v00)!0.5!(v10)$);
\coordinate (w1011) at ($(v10)!0.5!(v11)$);
\coordinate (w1112) at ($(v11)!0.5!(v12)$);
\coordinate (w1021) at ($(v10)!0.5!(v21)$);
\coordinate (w1122) at ($(v11)!0.5!(v22)$);
\coordinate (w1222) at ($(v12)!0.5!(v22)$);
\coordinate (w2122) at ($(v21)!0.5!(v22)$);
\coordinate (w2132) at ($(v21)!0.5!(v32)$);
\coordinate (w2232) at ($(v22)!0.5!(v32)$);

% 2 cell middle coords
\def\twocellx{0.25} % xshift, yshift of the centers of the two cell labels
\def\twocelly{0.38}
\coordinate (f1) at ($(v01)!0.5!(v11)$);
\coordinate (f2) at ($(v11)!0.5!(v21)$);
\coordinate (f3) at ($(v21)!0.5!(v31)$);
\coordinate (e1) at ($(f1)+(-1*\polyscale*\twocellx,\polyscale*\twocelly)$);
\coordinate (e2) at ($(f2)+(-1*\polyscale*\twocellx,\polyscale*\twocelly)$);
\coordinate (e3) at ($(f3)+(-1*\polyscale*\twocellx,\polyscale*\twocelly)$);
\coordinate (g1) at ($(f1)+(1*\polyscale*\twocellx,-1*\polyscale*\twocelly)$);

\draw[-,\cellcolor,\thickness,opacity=\cellop] (v00) -- (v02);
\draw[-,\cellcolor,\thickness,opacity=\cellop] (v10) -- (v12);
\draw[-,\cellcolor,\thickness,opacity=\cellop] (v21) -- (v22);
\draw[-,\cellcolor,\thickness,opacity=\cellop] (v00) -- (v10);
\draw[-,\cellcolor,\thickness,opacity=\cellop] (v02) -- (v32);
\draw[-,\cellcolor,\thickness,opacity=\cellop] (v01) -- (v12);
\draw[-,\cellcolor,\thickness,opacity=\cellop] (v00) -- (v22);
\draw[-,\cellcolor,\thickness,opacity=\cellop] (v10) -- (v32);

\draw[draw=none,fill=\fillcolor,opacity=\twocellop] (v00) -- (v02) -- (v32) -- (v10) -- cycle;

\node[inner sep = \nodewidth, draw=\cellcolor, circle, thick, fill=\cellcolor,opacity=\cellop] at (v00) {};
\node[inner sep = \nodewidth, draw=\cellcolor, circle, thick, fill=\cellcolor,opacity=\cellop] at (v01) {};
\node[inner sep = \nodewidth, draw=\cellcolor, circle, thick, fill=\cellcolor,opacity=\cellop] at (v02) {};
\node[inner sep = \nodewidth, draw=\cellcolor, circle, thick, fill=\cellcolor,opacity=\cellop] at (v10) {};
\node[inner sep = \nodewidth, draw=\cellcolor, circle, thick, fill=\cellcolor,opacity=\cellop] at (v11) {};
\node[inner sep = \nodewidth, draw=\cellcolor, circle, thick, fill=\cellcolor,opacity=\cellop] at (v12) {};
\node[inner sep = \nodewidth, draw=\cellcolor, circle, thick, fill=\cellcolor,opacity=\cellop] at (v21) {};
\node[inner sep = \nodewidth, draw=\cellcolor, circle, thick, fill=\cellcolor,opacity=\cellop] at (v22) {};
\node[inner sep = \nodewidth, draw=\cellcolor, circle, thick, fill=\cellcolor,opacity=\cellop] at (v32) {};

% 0 cell labels
\node[draw=none,fill=none,color=\labelcolor,opacity=\labelop] (v00label) at ($(v00)+(nodeshift)$) {\nodesize $(0,0)$};
\node[draw=none,fill=none,color=\labelcolor,opacity=\labelop] (v01label) at ($(v01)+(nodeshift)$) {\nodesize $(1,-1)$};
\node[draw=none,fill=none,color=\labelcolor,opacity=\labelop] (v02label) at ($(v02)+(nodeshift)$) {\nodesize $(0,0)$};
\node[draw=none,fill=none,color=\labelcolor,opacity=\labelop] (v10label) at ($(v10)+(nodeshift)$) {\nodesize $(0,0)$};
\node[draw=none,fill=none,color=\labelcolor,opacity=\labelop] (v11label) at ($(v11)+(nodeshift)$) {\nodesize $(1,-1)$};
\node[draw=none,fill=none,color=\labelcolor,opacity=\labelop] (v12label) at ($(v12)+(nodeshift)$) {\nodesize $(0,0)$};
\node[draw=none,fill=none,color=\labelcolor,opacity=\labelop] (v21label) at ($(v21)+(nodeshift)$) {\nodesize $(1,-1)$};
\node[draw=none,fill=none,color=\labelcolor,opacity=\labelop] (v22label) at ($(v22)+(nodeshift)$) {\nodesize $(0,0)$};
\node[draw=none,fill=none,color=\labelcolor,opacity=\labelop] (v32label) at ($(v32)+(nodeshift)$) {\nodesize $(0,0)$};

% 1 cell labels
\node[draw=none,fill=none,color=\labelcolor,opacity=\labelop] (w0001label) at ($(w0001)+(nodeshift)$) {\nodesize $(1,-1)$};
\node[draw=none,fill=none,color=\labelcolor,opacity=\labelop] (w0102label) at ($(w0102)+(nodeshift)$) {\nodesize $(0,-1)$};
\node[draw=none,fill=none,color=\labelcolor,opacity=\labelop] (w0010label) at ($(w0010)+(nodeshift)$) {\nodesize $(-1,0)$};
\node[draw=none,fill=none,color=\labelcolor,opacity=\labelop] (w0011label) at ($(w0011)+(nodeshift)$) {\nodesize $(1,-1)$};
\node[draw=none,fill=none,color=\labelcolor,opacity=\labelop] (w0112label) at ($(w0112)+(nodeshift)$) {\nodesize $(0,-1)$};
\node[draw=none,fill=none,color=\labelcolor,opacity=\labelop] (w0212label) at ($(w0212)+(nodeshift)$) {\nodesize $(-1,0)$};
\node[draw=none,fill=none,color=\labelcolor,opacity=\labelop] (w1011label) at ($(w1011)+(nodeshift)$) {\nodesize $(1,-1)$};
\node[draw=none,fill=none,color=\labelcolor,opacity=\labelop] (w1112label) at ($(w1112)+(nodeshift)$) {\nodesize $(0,-1)$};
\node[draw=none,fill=none,color=\labelcolor,opacity=\labelop] (w1021label) at ($(w1021)+(nodeshift)$) {\nodesize $(1,-1)$};
\node[draw=none,fill=none,color=\labelcolor,opacity=\labelop] (w1122label) at ($(w1122)+(nodeshift)$) {\nodesize $(0,-1)$};
\node[draw=none,fill=none,color=\labelcolor,opacity=\labelop] (w1222label) at ($(w1222)+(nodeshift)$) {\nodesize $(-1,0)$};
\node[draw=none,fill=none,color=\labelcolor,opacity=\labelop] (w2122label) at ($(w2122)+(nodeshift)$) {\nodesize $(0,-1)$};
\node[draw=none,fill=none,color=\labelcolor,opacity=\labelop] (w2132label) at ($(w2132)+(nodeshift)$) {\nodesize $(0,-1)$};
\node[draw=none,fill=none,color=\labelcolor,opacity=\labelop] (w2232label) at ($(w2232)+(nodeshift)$) {\nodesize $(-1,0)$};

% 2-cell labels
\node[draw=none,fill=none,color=\labelcolor,opacity=\labelop] (e1label) at ($(e1)+(nodeshift)$) {\nodesize $(-1,-1)$};
\node[draw=none,fill=none,color=\labelcolor,opacity=\labelop] (e2label) at ($(e2)+(nodeshift)$) {\nodesize $(-1,-1)$};
\node[draw=none,fill=none,color=\labelcolor,opacity=\labelop] (e3label) at ($(e3)+(nodeshift)$) {\nodesize $(-1,-1)$};
\node[draw=none,fill=none,color=\labelcolor,opacity=\labelop] (f1label) at ($(f1)+(nodeshift)$) {\nodesize $(0,-1)$};
\node[draw=none,fill=none,color=\labelcolor,opacity=\labelop] (f2label) at ($(f2)+(nodeshift)$) {\nodesize $(0,-1)$};
\node[draw=none,fill=none,color=\labelcolor,opacity=\labelop] (g1label) at ($(g1)+(nodeshift)$) {\nodesize $(1,-1)$};

% differentials
% 2-cell to 1-cell
\draw[->,\labelcolor,\diffthick,opacity=\arrowop] (g1label) -- (w0010label) node[midway, opacity=\diffop] {\tiny $-x_2$};
\draw[->,\labelcolor,\diffthick,opacity=\arrowop] (g1label) -- (w0011label) node[midway, opacity=\diffop] {\tiny $1$};
\draw[->,\labelcolor,\diffthick,opacity=\arrowop] (g1label) -- (w1011label) node[midway, opacity=\diffop] {\tiny $-1$};

\draw[->,\labelcolor,\diffthick,opacity=\arrowop] (e1label) -- (w0102label) node[midway, opacity=\diffop] {\tiny $x_1$};
\draw[->,\labelcolor,\diffthick,opacity=\arrowop] (e1label) -- (w0212label) node[midway, opacity=\diffop] {\tiny $-x_4$};
\draw[->,\labelcolor,\diffthick,opacity=\arrowop] (e1label) -- (w0112label) node[midway, opacity=\diffop] {\tiny $-x_3$};

\draw[->,\labelcolor,\diffthick,opacity=\arrowop] (e2label) -- (w1112label) node[midway, opacity=\diffop] {\tiny $x_1$};
\draw[->,\labelcolor,\diffthick,opacity=\arrowop] (e2label) -- (w1222label) node[midway, opacity=\diffop] {\tiny $-x_4$};
\draw[->,\labelcolor,\diffthick,opacity=\arrowop] (e2label) -- (w1122label) node[midway, opacity=\diffop] {\tiny $-x_3$};

\draw[->,\labelcolor,\diffthick,opacity=\arrowop] (e3label) -- (w2122label) node[midway, opacity=\diffop] {\tiny $x_1$};
\draw[->,\labelcolor,\diffthick,opacity=\arrowop] (e3label) -- (w2232label) node[midway, opacity=\diffop] {\tiny $-x_4$};
\draw[->,\labelcolor,\diffthick,opacity=\arrowop] (e3label) -- (w2132label) node[midway, opacity=\diffop] {\tiny $-x_3$};

\draw[->,\labelcolor,\diffthick,opacity=\arrowop] (f1label) -- (w0001label) node[midway, opacity=\diffop] {\tiny $x_1$};
\draw[->,\labelcolor,\diffthick,opacity=\arrowop] (f1label) -- (w0011label) node[midway, opacity=\diffop] {\tiny $-x_3$};
\draw[->,\labelcolor,\diffthick,opacity=\arrowop] (f1label) -- (w1112label) node[midway, opacity=\diffop] {\tiny $1$};
\draw[->,\labelcolor,\diffthick,opacity=\arrowop] (f1label) -- (w0112label) node[midway, opacity=\diffop] {\tiny $-1$};

\draw[->,\labelcolor,\diffthick,opacity=\arrowop] (f2label) -- (w1011label) node[midway, opacity=\diffop] {\tiny $x_1$};
\draw[->,\labelcolor,\diffthick,opacity=\arrowop] (f2label) -- (w1021label) node[midway, opacity=\diffop] {\tiny $-x_3$};
\draw[->,\labelcolor,\diffthick,opacity=\arrowop] (f2label) -- (w2122label) node[midway, opacity=\diffop] {\tiny $1$};
\draw[->,\labelcolor,\diffthick,opacity=\arrowop] (f2label) -- (w1122label) node[midway, opacity=\diffop] {\tiny $-1$};

% 1-cell to 0-cell
\draw[->,\labelcolor,\diffthick,opacity=\arrowop] (w0001label) -- (v00label) node[midway, opacity=\diffop] {\tiny $x_2x_3$};
\draw[->,\labelcolor,\diffthick,opacity=\arrowop] (w0001label) -- (v01label) node[midway, opacity=\diffop] {\tiny $-1$};

\draw[->,\labelcolor,\diffthick,opacity=\arrowop] (w1011label) -- (v10label) node[midway, opacity=\diffop] {\tiny $x_2x_3$};
\draw[->,\labelcolor,\diffthick,opacity=\arrowop] (w1011label) -- (v11label) node[midway, opacity=\diffop] {\tiny $-1$};

\draw[->,\labelcolor,\diffthick,opacity=\arrowop] (w0102label) -- (v01label) node[midway, opacity=\diffop] {\tiny $-x_3$};
\draw[->,\labelcolor,\diffthick,opacity=\arrowop] (w0102label) -- (v02label) node[midway, opacity=\diffop] {\tiny $x_4$};

\draw[->,\labelcolor,\diffthick,opacity=\arrowop] (w1112label) -- (v11label) node[midway, opacity=\diffop] {\tiny $-x_3$};
\draw[->,\labelcolor,\diffthick,opacity=\arrowop] (w1112label) -- (v12label) node[midway, opacity=\diffop] {\tiny $x_4$};

\draw[->,\labelcolor,\diffthick,opacity=\arrowop] (w2122label) -- (v21label) node[midway, opacity=\diffop] {\tiny $-x_3$};
\draw[->,\labelcolor,\diffthick,opacity=\arrowop] (w2122label) -- (v22label) node[midway, opacity=\diffop] {\tiny $x_4$};

\draw[->,\labelcolor,\diffthick,opacity=\arrowop] (w0010label) -- (v00label) node[midway, opacity=\diffop] {\tiny $x_1$};
\draw[->,\labelcolor,\diffthick,opacity=\arrowop] (w0010label) -- (v10label) node[midway, opacity=\diffop] {\tiny $-x_3$};

\draw[->,\labelcolor,\diffthick,opacity=\arrowop] (w0212label) -- (v02label) node[midway, opacity=\diffop] {\tiny $x_1$};
\draw[->,\labelcolor,\diffthick,opacity=\arrowop] (w0212label) -- (v12label) node[midway, opacity=\diffop] {\tiny $-x_3$};

\draw[->,\labelcolor,\diffthick,opacity=\arrowop] (w1222label) -- (v12label) node[midway, opacity=\diffop] {\tiny $x_1$};
\draw[->,\labelcolor,\diffthick,opacity=\arrowop] (w1222label) -- (v22label) node[midway, opacity=\diffop] {\tiny $-x_3$};

\draw[->,\labelcolor,\diffthick,opacity=\arrowop] (w2232label) -- (v22label) node[midway, opacity=\diffop] {\tiny $x_1$};
\draw[->,\labelcolor,\diffthick,opacity=\arrowop] (w2232label) -- (v32label) node[midway, opacity=\diffop] {\tiny $-x_3$};

\draw[->,\labelcolor,\diffthick,opacity=\arrowop] (w0011label) -- (v00label) node[midway, opacity=\diffop] {\tiny $x_1x_2$};
\draw[->,\labelcolor,\diffthick,opacity=\arrowop] (w0011label) -- (v11label) node[midway, opacity=\diffop] {\tiny $-1$};

\draw[->,\labelcolor,\diffthick,opacity=\arrowop] (w1021label) -- (v10label) node[midway, opacity=\diffop] {\tiny $x_1x_2$};
\draw[->,\labelcolor,\diffthick,opacity=\arrowop] (w1021label) -- (v21label) node[midway, opacity=\diffop] {\tiny $-1$};

\draw[->,\labelcolor,\diffthick,opacity=\arrowop] (w0112label) -- (v01label) node[midway, opacity=\diffop] {\tiny $-x_1$};
\draw[->,\labelcolor,\diffthick,opacity=\arrowop] (w0112label) -- (v12label) node[midway, opacity=\diffop] {\tiny $x_4$};

\draw[->,\labelcolor,\diffthick,opacity=\arrowop] (w1122label) -- (v11label) node[midway, opacity=\diffop] {\tiny $-x_1$};
\draw[->,\labelcolor,\diffthick,opacity=\arrowop] (w1122label) -- (v22label) node[midway, opacity=\diffop] {\tiny $x_4$};

\draw[->,\labelcolor,\diffthick,opacity=\arrowop] (w2132label) -- (v21label) node[midway, opacity=\diffop] {\tiny $-x_1$};
\draw[->,\labelcolor,\diffthick,opacity=\arrowop] (w2132label) -- (v32label) node[midway, opacity=\diffop] {\tiny $x_4$};

%%%%%% ADDING THE FULL COMPLEX HERE
\def\xshift{5}
\def\yshift{-2}
\def\compscale{2}
\def\xsep{2}
\def\ysep{0.5}
\def\compcolor{black}
\def\compsize{\small}

\coordinate (t20) at (\xshift, \yshift-\compscale*\ysep);
\coordinate (t21) at (\xshift, \yshift+0*\compscale*\ysep);
\coordinate (t22) at (\xshift, \yshift+1*\compscale*\ysep);
\coordinate (t10) at (\xshift+1*\compscale*\xsep, \yshift-1*\compscale*\ysep);
\coordinate (t11) at (\xshift+1*\compscale*\xsep, \yshift+0*\compscale*\ysep);
\coordinate (t12) at (\xshift+1*\compscale*\xsep, \yshift+1*\compscale*\ysep);
\coordinate (t00) at (\xshift+2*\compscale*\xsep, \yshift-0.5*\compscale*\ysep);
\coordinate (t01) at (\xshift+2*\compscale*\xsep, \yshift+0.5*\compscale*\ysep);

\coordinate (plus20) at (\xshift, \yshift-0.5*\compscale*\ysep);
\coordinate (plus21) at (\xshift, \yshift+0.5*\compscale*\ysep);
\coordinate (plus10) at (\xshift+1*\compscale*\xsep, \yshift-0.5*\compscale*\ysep);
\coordinate (plus11) at (\xshift+1*\compscale*\xsep, \yshift+0.5*\compscale*\ysep);
\coordinate (plus00) at (\xshift+2*\compscale*\xsep, \yshift);

\node[draw=none,fill=none,\compcolor] (t20node) at (t20) {\compsize $S(1,-1)$};
\node[draw=none,fill=none,\compcolor,inner sep = 5pt] (t21node) at (t21) {\compsize $S(0,-1)^{\oplus 2}$};
\node[draw=none,fill=none,\compcolor] (t22node) at (t22) {\compsize $S(-1,-1)^{\oplus 3}$};

\node[draw=none,fill=none,\compcolor] (t10node) at (t10) {\compsize $S(1,-1)^{\oplus 4}$};
\node[draw=none,fill=none,\compcolor,inner sep = 5pt] (t11node) at (t11) {\compsize $S(0,-1)^{\oplus 6}$};
\node[draw=none,fill=none,\compcolor] (t12node) at (t12) {\compsize $S(-1,0)^{\oplus 4}$};

\node[draw=none,fill=none,\compcolor] (t00node) at (t00) {\compsize $S(1,-1)^{\oplus 3}$};
\node[draw=none,fill=none,\compcolor] (t01node) at (t01) {\compsize $S(0,0)^{\oplus 6}$};

\node[draw=none,fill=none,\compcolor] (plus20node) at (plus20) {\compsize $\oplus$};
\node[draw=none,fill=none,\compcolor] (plus21node) at (plus21) {\compsize $\oplus$};

\node[draw=none,fill=none,\compcolor] (plus10node) at (plus10) {\compsize $\oplus$};
\node[draw=none,fill=none,\compcolor] (plus11node) at (plus11) {\compsize $\oplus$};

\node[draw=none,fill=none,\compcolor, inner sep = 0.9cm] (plus00node) at (plus00) {\compsize $\oplus$};

\draw[->,\compcolor] (t21node) -- (t11node);
\draw[->,\compcolor] (t11node) -- (plus00node);

\node[draw=none,fill=none] at (\xshift-0.6*\xsep*\compscale,\yshift) {$F_\bullet^D(D)\ \ \ = $};

\end{tikzpicture}\]
    \caption{Depicted is the chain complex $F_\bullet^D(D)$ from \Cref{example:H2_cellular_resolution} overlaid on top of the cellular complex $P_D$ of \Cref{fig:labeled_cellular_complex}. We denote the term $S(i,j)$ by $(i,j)$ for pictorial clarity. After placing the twists associated to the $i$-cells in homological degree $i$, one obtains the complex $F_\bullet^D(D)$, depicted below the cellular complex.}
    \label{fig:full_cellular_resolution}
\end{figure}
\end{example}

The following proposition is a crucial technical ingredient in understanding $W^{\Cox}_{\Gamma}$ and one of the main results in \cite{CH25}. Since the work of Cranton Heller deals with projective toric varieties, we reprove the author's claim in our more general setting of toric stacks. Nevertheless, the argument for varieties largely carries through to our setting. 

\begin{proposition}\label{prop:FMonAmplesGivesF}
For any $d \in \Gamma \cap \wh{G}$, the Fourier-Mukai transform of $\CO_{X_{\Gamma}}(d) \in D(\CX_{\Gamma})$, $\Phi_{\Gamma}(\CO_{X_{\Gamma}}(d)) \in K_{\Theta} \subset D([\A^{n}/G])$, is given by $F^{D}_{\bullet}(D)$ (viewed as a complex of $\wh{G}$-graded $S$-modules by \Cref{rmk:FisinKtheta}).  
\end{proposition}

\begin{proof}
First, note that this proposition is almost identical to the main theorems in \cite{CH25}. The only difference is that the Fourier-Mukai functor in \cite{CH25} is defined for $D(X)$, where $X$ a smooth, projective toric variety $X$. The kernel of the transform is the same HHL complex, but viewed as a complex of line bundles on $X \times X$ instead of the product stack $\CX_{\Gamma} \times [\A^{n}/G]$ as we consider here.

Recall that the FM functor $\Phi_{\Gamma}$ was defined as: 
\begin{equation*}
\Phi_\Gamma \coloneqq Rq^{G}_*\lp Lp^*(-)\otimes_L \HHL_\Gamma\rp: \Db(\CX_\Gamma)\to \Db([\A^n/G]).
\end{equation*}
Since we pull-back a line bundle,  $Lp^{*} = p^{*}$. We compute that
\begin{equation}\label{eq:FMcomplexb4pushing}
   Lp^{*}(\CO_{\CX_\Gamma}(d))\otimes_{L}H_{\Gamma} \simeq \cdots \to \bigoplus_{\delta\in E_{j}} \CO_{\CX_\Gamma}(\pi(\floor{\delta})+d)\boxtimes S(-\pi(\ceil{\delta}))\to \cdots. 
\end{equation}
Using the projection formula for stacks, we have that 
\[
R^{t}q^{G}_{*}(\CO_{\CX_\Gamma}(\pi(\floor{\delta})+d)\boxtimes S(-\pi(\ceil{\delta})))\cong H^{t}(\CX_{\Gamma}, \CO_{\CX_{\Gamma}}(\pi(\floor{\delta})+d))\otimes_{k}S(-\pi(\ceil{\delta})).
\]
By \Cref{prop:nefcohomologyvanishes}, we have that if $t\neq0$
\[
H^{t}(\CX_{\Gamma}, \CO_{\CX_{\Gamma}}(\pi(\floor{\delta})+d)) =0.
\]
Hence, the terms in \Cref{eq:FMcomplexb4pushing} are $q_{*}$-acyclic. Since the complex is bounded, we conclude that $\Phi_{\Gamma}(\CO_{\CX_{\Gamma}}(d))$ is given by
\[
 \to \bigoplus_{\delta\in E_{j}}H^{0}(\CX_{\Gamma}, \CO_{\CX_{\Gamma}}(\pi(\floor{\delta})+d))\otimes_{k}S(-\pi(\ceil{\delta})) \to \bigoplus_{\delta\in E_{j-1}}H^{0}(\CX_{\Gamma}, \CO_{\CX_{\Gamma}}(\pi(\floor{\delta})+d))\otimes_{k}S(-\pi(\ceil{\delta}))\to
\]
By \Cref{prop:nefcohomologyvanishes}, we can rewrite the complex using the graded pieces of the Cox ring:

\[
 \to \bigoplus_{\delta\in E_{j}}S_{\pi(\floor{\delta})+d}\otimes_{k}S(-\pi(\ceil{\delta})) \to \bigoplus_{\delta\in E_{j-1}}S_{\pi(\floor{\delta})+d}\otimes_{k}S(-\pi(\ceil{\delta}))\to
\]
Note that this is precisely the complex considered in \cite{CH25} (see Equation 3.2). Then, the identification between this complex and the cell complex on $P_{D}$ is proven in \cite[Theorem 5.1]{CH25}.
\end{proof}

\begin{corollary}\label{cor:wcoxisgeneratedbyf}
$W^\Cox_{\Gamma}$ is generated by complexes $F^{D}_{\bullet}(D)$, where $d \in \Gamma\cap \wh{G}$.    
\end{corollary}
\begin{proof}
 This immediately follows from \Cref{prop:amplesgenerate} and \Cref{prop:FMonAmplesGivesF}.   
\end{proof}

\begin{corollary} \label{cor: window equality}
Suppose for every $d \in \Gamma\cap \wh{G}$, $F^{D}_{\bullet}(D) \in W_{\bs{\lambda}}^{\bs{w}}$. Then, $W^{\Cox}_{\Gamma} = W_{\bs{\lambda}}^{\bs{w}}$.
\end{corollary}
\begin{proof}
If $F^{D}_{\bullet}(D) \in W_{\bs{\lambda}}^{\bs{w}}$ for any  $d \in \Gamma$, then \Cref{cor:wcoxisgeneratedbyf} implies $W^{\Cox}_{\Gamma} \subset W_{\bs{\lambda}}^{\bs{w}}$. To show the equality, note that the following diagram commutes by the definition of a window category:
\[\begin{tikzcd}[column sep = 0cm]
	{W^{\Cox}_{\Gamma}} && {W_{\bs{\lambda}}^{\bs{w}}} \\
	& {\Db([U_\Gamma/G])}
	\arrow[hook, from=1-1, to=1-3]
	\arrow["{r^{*}_1}"', from=1-1, to=2-2]
	\arrow["{r^{*}_{2}}", from=1-3, to=2-2]
\end{tikzcd}\]
Since $r^{*}_{1}$ and $r^{*}_{2}$ are equivalences and $W_{\Gamma}^{\Cox}$ is an isomorphism-closed subcategory, it must be $W^{\Cox}_{\Gamma} = W_{\bs{\lambda}}^{\bs{w}}$.
\end{proof}

In light of \Cref{cor: window equality}, to show equality of the Cox window and the grade restriction window, we must understand the $\lambda$-weights of the generators $F_\bullet^D(D)$. The following lemma provides valuable insight into this computation.

\begin{lemma}\label{lemma:weights_in_Tminus_Tplus}
For any $\delta \in \cells(P_D)$, 
we have the following for any 1-PS $\lambda$ such that $\lambda\notin -|\Sigma_{GKZ}|^\vee$:
\begin{enumerate}[(i)]
    \item\label{lemma:weights_in_Tminus_Tplus_EXT} If $\lambda\in |\Sigma_{GKZ}|^\vee$, then $\wh{\lambda}(\pi_\R(\ceil{\delta}))\in [0,T_\lambda^{>0})$.    \item\label{lemma:weights_in_Tminus_Tplus_INT} If $\lambda\notin |\Sigma_{GKZ}|^\vee$, then $\wh{\lambda}(\pi_\R(\ceil{\delta}))\in (T_\lambda^{<0},T_\lambda^{>0})$. 
\end{enumerate}
\end{lemma}
\begin{proof}
 By \Cref{rmk:FisinKtheta}, it follows that $\pi_\R(\ceil{\delta})\in Z_G = \pi_\R([0,1)^n)$, where $\delta \subset P_{D}$.  If $\lambda\in |\Sigma_{GKZ}|^\vee$, then $ I_\lambda^{\ge 0} = [n]$ so
\[\wh{\lambda}(\pi_\R(\ceil{\delta}))\in(\widehat{\lambda}(\pi_\R([0,1)^n)) = \left[0, \sum_{j\in I_\lambda^{>0}}\wh{\lambda}(\beta_j) \right) = [0,T_\lambda^{>0}).\]
If instead $\lambda\notin |\Sigma_{GKZ}|^\vee$, because $\lambda\in \Gamma^\vee$ we get that $I_\lambda^{>0}\ne \emptyset$ and $I_\lambda^{<0}\ne\emptyset$. Thus
\[\wh{\lambda}(\pi_\R(\ceil{\delta}))\in (\widehat{\lambda}(\pi_\R([0,1)^n)) = \lp \sum_{j\in I_\lambda^{<0}}\wh{\lambda}(\beta_j), \sum_{j\in I_\lambda^{>0}}\wh{\lambda}(\beta_j)\rp = (T_\lambda^{<0}, T_\lambda^{>0}).\]
\end{proof}
\begin{remark}\label{rmk:automatic_grade_restriction_rules_satisfied}
Since the terms of $F_\bullet^D(D)$ are direct sums of $\widehat{G}$-graded $S$-modules of the form $S(-\pi(\ceil{\delta}))$, \Cref{lemma:weights_in_Tminus_Tplus} implies that for any external wall $\lambda$ of $\Sigma_{GKZ}$,
\[
\mu(F_\bullet^D(D), \lambda, V) \subset [0, T_{\lambda}^{>0}).
\]
Hence, if $\lambda\in |\Sigma_{GKZ}|^\vee$ appears in an HKKN stratification $\bs{\lambda}$, $F_\bullet^D(D)$ automatically satisfies the associated grade restriction rule in the window category $W^{\bs{0}}_{\bs{\lambda}}$. If $\lambda\notin |\Sigma_{GKZ}|^\vee$, then instead
\[
\mu(F_\bullet^D(D), \lambda, V) \subset (T_{\lambda}^{<0}, T_{\lambda}^{>0}).
\]
Therefore, in order for $F_\bullet^D(D)$ to be in $W^{\bs{0}}_{\bs{\lambda}}$, for every internal wall of $\bs{\lambda}$, one needs to check that a finite number of negative weights do not belong to $\mu(F_\bullet^D(D), \lambda, V)$. These observations were part of the initial motivation for studying the comparison between $W^{\bs{w}}_{\bs{\lambda}}$ and $W^{\Cox}_{\Gamma}$ via the Cox generators $F_\bullet^D(D)$.
\end{remark}

\begin{example}\label{example:H2_grade_restriction_rules_geometrically}

In this example, we perform the computation of weights $\mu(F_\bullet^D(D), \lambda, V)$ in our running example following the procedure described in \ref{subsection:grade_restriction_windows}. Recall that $\CX_{\Gamma} = \CH_2$ and $F_\bullet^D(D)$ is the cellular complex of free $S$-modules associated to the ample divisor $D = (0,0,1,1)$ from  \Cref{example:H2_cellular_resolution}. We pick an HKKN stratification of $\bs{\lambda}$ with $\lambda_1 =\xi_1, \lambda_2 = \xi_3, \lambda_3 = \xi_2$ from \Cref{example:H2_stratifications} and consider the associated grade restriction window $W_{\bs{\lambda}}$. By  \Cref{def:grade_restriction_windows}, every 1-PS in $\bs{\lambda}$ has an associated grade restriction rule. As noted in \Cref{rmk:automatic_grade_restriction_rules_satisfied}, complex $F_\bullet^D(D)$ automatically satisfies grade restriction rule whenever $\lambda\in |\Sigma_{GKZ}|^\vee$. Since $\lambda_1,\lambda_2 \in |\Sigma_{GKZ}|^\vee$ and $\lambda_3 \notin |\Sigma_{GKZ}|^\vee$, the only grade restriction rule that is not automatically satisfied is the one associated to $\lambda_3$. Hence, we are interested in determining weights in $\mu(F_\bullet^D(D), \lambda_{3}, (X_3)_{\lambda_3}^{=0})$. One computes that
\[V\coloneqq (X_3)^{=0}_{\lambda_3} = (\A^4)^{=0}_{\lambda_3}\setminus\bigcup_{j=1,2}(\A^4)_{\lambda_j}^{\le 0} = V(x_1,x_2,x_3) \setminus {\{0\}}.\]
Pulling back $F_\bullet^D(D)$ to $V$ amounts to replacing free $S$-modules with the corresponding equivariant sheaves $\CO_{V}$. Since $V = V(x_1,x_2,x_3) \setminus {\{0\}}$, all differentials containing variables $x_1,x_2,x_3$ become 0 morphisms. From the perspective of commutative algebra, we took a quotient of the polynomial ring by the ideal $(x_1,x_2,x_3)$ and inverted the variable $x_4$. In other words, the pull-back is a complex of free $k[x^{-1}_4,x_4]$-modules. 

\Cref{fig:restricted_cellular_resolution} depicts the pullback of $F_\bullet^D(D)$ to $V$, overlaid on top of the cellular complex $P_D$, following the depiction of $F_\bullet^D(D)$ in \Cref{fig:full_cellular_resolution}. Highlighted in distinct colors are the fixed-weight summands $(i^*F_\bullet^D(D))_{\lambda_3=k}$ appearing in the splitting of $i^*F_\bullet^D(D)$ guaranteed by \eqref{eqn:splitting}. Observe that every fixed-weight summand further splits as a direct sum. There are only three non-acyclic complexes, all consisting of one term, which we have encircled for clarity. Two terms were summands of $(i^*F_\bullet^D(D))_{\lambda_3=0}$, and one term was a summand of $(i^*F_\bullet^D(D))_{\lambda_3=1}$, thus we conclude
\[
\mu(F_\bullet^D(D), \lambda_{3}, (X_3)^{=0}_{\lambda_3}) = \{0,1\}.
\]
It follows that $F^{D}_{\bullet}(D) \in W_{\bs{\lambda}}$.
\begin{figure}
\[\hspace{-1.8cm}\begin{tikzpicture}
\def\polyscale{6.5}
\def\cellcolor{cyan}
\def\fillcolor{cyan}
\def\componecolor{red}
\def\comptwocolor{darkGreen}
\def\compthreecolor{orange}
\def\cellop{0.1}
\def\labelop{1}
\def\twocellop{0.03}
\def\arrowop{0.3}
\def\diffop{1}
\def\nodewidth{2pt}
\def\xshift{0}
\def\yshift{0}
\def\thickness{ultra thick}
\def\nodesize{\tiny}
\def\diffthick{thick}
\coordinate (nodeshift) at (-0.02*\polyscale,0.02*\polyscale);
\def\labelradius{0.6}

\tikzset{every node/.style={inner sep=0pt}}

\coordinate (v00) at (\xshift+0*\polyscale,\yshift+0*\polyscale);
\coordinate (v01) at (\xshift+0*\polyscale,\yshift+0.5*\polyscale);
\coordinate (v02) at (\xshift+0*\polyscale,\yshift+1*\polyscale);
\coordinate (v10) at (\xshift+1*\polyscale,\yshift+0*\polyscale);
\coordinate (v11) at (\xshift+1*\polyscale,\yshift+0.5*\polyscale);
\coordinate (v12) at (\xshift+1*\polyscale,\yshift+1*\polyscale);
\coordinate (v21) at (\xshift+2*\polyscale,\yshift+0.5*\polyscale);
\coordinate (v22) at (\xshift+2*\polyscale,\yshift+1*\polyscale);
\coordinate (v32) at (\xshift+3*\polyscale,\yshift+1*\polyscale);
\coordinate (v31) at (\xshift+3*\polyscale,\yshift+0.5*\polyscale);

% 1 cell middle coords
\coordinate (w0001) at ($(v00)!0.5!(v01)$); % wijkl is between vij and vkl
\coordinate (w0102) at ($(v01)!0.5!(v02)$);
\coordinate (w0011) at ($(v00)!0.5!(v11)$);
\coordinate (w0112) at ($(v01)!0.5!(v12)$);
\coordinate (w0212) at ($(v02)!0.5!(v12)$);
\coordinate (w0010) at ($(v00)!0.5!(v10)$);
\coordinate (w1011) at ($(v10)!0.5!(v11)$);
\coordinate (w1112) at ($(v11)!0.5!(v12)$);
\coordinate (w1021) at ($(v10)!0.5!(v21)$);
\coordinate (w1122) at ($(v11)!0.5!(v22)$);
\coordinate (w1222) at ($(v12)!0.5!(v22)$);
\coordinate (w2122) at ($(v21)!0.5!(v22)$);
\coordinate (w2132) at ($(v21)!0.5!(v32)$);
\coordinate (w2232) at ($(v22)!0.5!(v32)$);

% 2 cell middle coords
\def\twocellx{0.25} % xshift, yshift of the centers of the two cell labels
\def\twocelly{0.38}
\coordinate (f1) at ($(v01)!0.5!(v11)$);
\coordinate (f2) at ($(v11)!0.5!(v21)$);
\coordinate (f3) at ($(v21)!0.5!(v31)$);
\coordinate (e1) at ($(f1)+(-1*\polyscale*\twocellx,\polyscale*\twocelly)$);
\coordinate (e2) at ($(f2)+(-1*\polyscale*\twocellx,\polyscale*\twocelly)$);
\coordinate (e3) at ($(f3)+(-1*\polyscale*\twocellx,\polyscale*\twocelly)$);
\coordinate (g1) at ($(f1)+(1*\polyscale*\twocellx,-1*\polyscale*\twocelly)$);

\draw[-,\cellcolor,\thickness,opacity=\cellop] (v00) -- (v02);
\draw[-,\cellcolor,\thickness,opacity=\cellop] (v10) -- (v12);
\draw[-,\cellcolor,\thickness,opacity=\cellop] (v21) -- (v22);
\draw[-,\cellcolor,\thickness,opacity=\cellop] (v00) -- (v10);
\draw[-,\cellcolor,\thickness,opacity=\cellop] (v02) -- (v32);
\draw[-,\cellcolor,\thickness,opacity=\cellop] (v01) -- (v12);
\draw[-,\cellcolor,\thickness,opacity=\cellop] (v00) -- (v22);
\draw[-,\cellcolor,\thickness,opacity=\cellop] (v10) -- (v32);

\draw[draw=none,fill=\fillcolor,opacity=\twocellop] (v00) -- (v02) -- (v32) -- (v10) -- cycle;

\node[inner sep = \nodewidth, draw=\cellcolor, circle, thick, fill=\cellcolor,opacity=\cellop] at (v00) {};
\node[inner sep = \nodewidth, draw=\cellcolor, circle, thick, fill=\cellcolor,opacity=\cellop] at (v01) {};
\node[inner sep = \nodewidth, draw=\cellcolor, circle, thick, fill=\cellcolor,opacity=\cellop] at (v02) {};
\node[inner sep = \nodewidth, draw=\cellcolor, circle, thick, fill=\cellcolor,opacity=\cellop] at (v10) {};
\node[inner sep = \nodewidth, draw=\cellcolor, circle, thick, fill=\cellcolor,opacity=\cellop] at (v11) {};
\node[inner sep = \nodewidth, draw=\cellcolor, circle, thick, fill=\cellcolor,opacity=\cellop] at (v12) {};
\node[inner sep = \nodewidth, draw=\cellcolor, circle, thick, fill=\cellcolor,opacity=\cellop] at (v21) {};
\node[inner sep = \nodewidth, draw=\cellcolor, circle, thick, fill=\cellcolor,opacity=\cellop] at (v22) {};
\node[inner sep = \nodewidth, draw=\cellcolor, circle, thick, fill=\cellcolor,opacity=\cellop] at (v32) {};

% 0 cell labels
\node[draw=none,fill=none,color=\componecolor,opacity=\labelop] (v00label) at ($(v00)+(nodeshift)$) {\nodesize $(0,0)$};
\node[draw=none,fill=none,color=\comptwocolor,opacity=\labelop] (v01label) at ($(v01)+(nodeshift)$) {\nodesize $(1,-1)$};
\node[draw=none,fill=none,color=\componecolor,opacity=\labelop] (v02label) at ($(v02)+(nodeshift)$) {\nodesize $(0,0)$};
\node[draw=none,fill=none,color=\componecolor,opacity=\labelop] (v10label) at ($(v10)+(nodeshift)$) {\nodesize $(0,0)$};
\node[draw=none,fill=none,color=\comptwocolor,opacity=\labelop] (v11label) at ($(v11)+(nodeshift)$) {\nodesize $(1,-1)$};
\node[draw=none,fill=none,color=\componecolor,opacity=\labelop] (v12label) at ($(v12)+(nodeshift)$) {\nodesize $(0,0)$};
\node[draw=none,fill=none,color=\comptwocolor,opacity=\labelop] (v21label) at ($(v21)+(nodeshift)$) {\nodesize $(1,-1)$};
\node[draw=none,fill=none,color=\componecolor,opacity=\labelop] (v22label) at ($(v22)+(nodeshift)$) {\nodesize $(0,0)$};
\node[draw=none,fill=none,color=\componecolor,opacity=\labelop] (v32label) at ($(v32)+(nodeshift)$) {\nodesize $(0,0)$};

% 1 cell labels
\node[draw=none,fill=none,color=\comptwocolor,opacity=\labelop] (w0001label) at ($(w0001)+(nodeshift)$) {\nodesize $(1,-1)$};
\node[draw=none,fill=none,color=\componecolor,opacity=\labelop] (w0102label) at ($(w0102)+(nodeshift)$) {\nodesize $(0,-1)$};
\node[draw=none,fill=none,color=\compthreecolor,opacity=\labelop] (w0010label) at ($(w0010)+(nodeshift)$) {\nodesize $(-1,0)$};
\node[draw=none,fill=none,color=\comptwocolor,opacity=\labelop] (w0011label) at ($(w0011)+(nodeshift)$) {\nodesize $(1,-1)$};
\node[draw=none,fill=none,color=\componecolor,opacity=\labelop] (w0112label) at ($(w0112)+(nodeshift)$) {\nodesize $(0,-1)$};
\node[draw=none,fill=none,color=\compthreecolor,opacity=\labelop] (w0212label) at ($(w0212)+(nodeshift)$) {\nodesize $(-1,0)$};
\node[draw=none,fill=none,color=\comptwocolor,opacity=\labelop] (w1011label) at ($(w1011)+(nodeshift)$) {\nodesize $(1,-1)$};
\node[draw=none,fill=none,color=\componecolor,opacity=\labelop] (w1112label) at ($(w1112)+(nodeshift)$) {\nodesize $(0,-1)$};
\node[draw=none,fill=none,color=\comptwocolor,opacity=\labelop] (w1021label) at ($(w1021)+(nodeshift)$) {\nodesize $(1,-1)$};
\node[draw=none,fill=none,color=\componecolor,opacity=\labelop] (w1122label) at ($(w1122)+(nodeshift)$) {\nodesize $(0,-1)$};
\node[draw=none,fill=none,color=\compthreecolor,opacity=\labelop] (w1222label) at ($(w1222)+(nodeshift)$) {\nodesize $(-1,0)$};
\node[draw=none,fill=none,color=\componecolor,opacity=\labelop] (w2122label) at ($(w2122)+(nodeshift)$) {\nodesize $(0,-1)$};
\node[draw=none,fill=none,color=\componecolor,opacity=\labelop] (w2132label) at ($(w2132)+(nodeshift)$) {\nodesize $(0,-1)$};
\node[draw=none,fill=none,color=\compthreecolor,opacity=\labelop] (w2232label) at ($(w2232)+(nodeshift)$) {\nodesize $(-1,0)$};

% 2-cell labels
\node[draw=none,fill=none,color=\compthreecolor,opacity=\labelop] (e1label) at ($(e1)+(nodeshift)$) {\nodesize $(-1,-1)$};
\node[draw=none,fill=none,color=\compthreecolor,opacity=\labelop] (e2label) at ($(e2)+(nodeshift)$) {\nodesize $(-1,-1)$};
\node[draw=none,fill=none,color=\compthreecolor,opacity=\labelop] (e3label) at ($(e3)+(nodeshift)$) {\nodesize $(-1,-1)$};
\node[draw=none,fill=none,color=\componecolor,opacity=\labelop] (f1label) at ($(f1)+(nodeshift)$) {\nodesize $(0,-1)$};
\node[draw=none,fill=none,color=\componecolor,opacity=\labelop] (f2label) at ($(f2)+(nodeshift)$) {\nodesize $(0,-1)$};
\node[draw=none,fill=none,color=\comptwocolor,opacity=\labelop] (g1label) at ($(g1)+(nodeshift)$) {\nodesize $(1,-1)$};

% differentials
% 2-cell to 1-cell

\draw[->,\comptwocolor,\diffthick,opacity=\arrowop] (g1label) -- (w0011label) node[midway, opacity=\diffop] {\tiny $1$};
\draw[->,\comptwocolor,\diffthick,opacity=\arrowop] (g1label) -- (w1011label) node[midway, opacity=\diffop] {\tiny $-1$};

\draw[->,\compthreecolor,\diffthick,opacity=\arrowop] (e1label) -- (w0212label) node[midway, opacity=\diffop] {\tiny $-x_4$};

\draw[->,\compthreecolor,\diffthick,opacity=\arrowop] (e2label) -- (w1222label) node[midway, opacity=\diffop] {\tiny $-x_4$};

\draw[->,\compthreecolor,\diffthick,opacity=\arrowop] (e3label) -- (w2232label) node[midway, opacity=\diffop] {\tiny $-x_4$};

\draw[->,\componecolor,\diffthick,opacity=\arrowop] (f1label) -- (w1112label) node[midway, opacity=\diffop] {\tiny $1$};
\draw[->,\componecolor,\diffthick,opacity=\arrowop] (f1label) -- (w0112label) node[midway, opacity=\diffop] {\tiny $-1$};

\draw[->,\componecolor,\diffthick,opacity=\arrowop] (f2label) -- (w2122label) node[midway, opacity=\diffop] {\tiny $1$};
\draw[->,\componecolor,\diffthick,opacity=\arrowop] (f2label) -- (w1122label) node[midway, opacity=\diffop] {\tiny $-1$};

% 1-cell to 0-cell
\draw[->,\comptwocolor,\diffthick,opacity=\arrowop] (w0001label) -- (v01label) node[midway, opacity=\diffop] {\tiny $-1$};

\draw[->,\comptwocolor,\diffthick,opacity=\arrowop] (w1011label) -- (v11label) node[midway, opacity=\diffop] {\tiny $-1$};

\draw[->,\componecolor,\diffthick,opacity=\arrowop] (w0102label) -- (v02label) node[midway, opacity=\diffop] {\tiny $x_4$};

\draw[->,\componecolor,\diffthick,opacity=\arrowop] (w1112label) -- (v12label) node[midway, opacity=\diffop] {\tiny $x_4$};

\draw[->,\componecolor,\diffthick,opacity=\arrowop] (w2122label) -- (v22label) node[midway, opacity=\diffop] {\tiny $x_4$};

\draw[->,\comptwocolor,\diffthick,opacity=\arrowop] (w0011label) -- (v11label) node[midway, opacity=\diffop] {\tiny $-1$};

\draw[->,\comptwocolor,\diffthick,opacity=\arrowop] (w1021label) -- (v21label) node[midway, opacity=\diffop] {\tiny $-1$};

\draw[->,\componecolor,\diffthick,opacity=\arrowop] (w0112label) -- (v12label) node[midway, opacity=\diffop] {\tiny $x_4$};

\draw[->,\componecolor,\diffthick,opacity=\arrowop] (w1122label) -- (v22label) node[midway, opacity=\diffop] {\tiny $x_4$};
\draw[->,\componecolor,\diffthick,opacity=\arrowop] (w2132label) -- (v32label) node[midway, opacity=\diffop] {\tiny $x_4$};

\draw[dashed, thick, black] (v00label) circle (\labelradius);
\draw[dashed, thick, black] (v10label) circle (\labelradius);
\draw[dashed, thick, black] (w0010label) circle (\labelradius);

%%%%% DRAWING THE SPLIT COMPLEX
\def\xshift{9}
\def\yshift{-2.6}
\def\compscale{1.7}
\def\xsep{2}
\def\ysep{1}
\def\compsize{\small}
\def\compthick{thick}
\def\labelshift{2};
\tikzset{every node/.style={inner sep=5pt}}

\coordinate (t20) at (\xshift, \yshift-\compscale*\ysep);
\coordinate (t21) at (\xshift, \yshift+0*\compscale*\ysep);
\coordinate (t22) at (\xshift, \yshift+1*\compscale*\ysep);
\coordinate (t10) at (\xshift+1*\compscale*\xsep, \yshift-1*\compscale*\ysep);
\coordinate (t11) at (\xshift+1*\compscale*\xsep, \yshift+0*\compscale*\ysep);
\coordinate (t12) at (\xshift+1*\compscale*\xsep, \yshift+1*\compscale*\ysep);
\coordinate (t00) at (\xshift+2*\compscale*\xsep, \yshift-1*\compscale*\ysep);
\coordinate (t01) at (\xshift+2*\compscale*\xsep, \yshift+0*\compscale*\ysep);

 \coordinate (label0) at (\xshift-\labelshift*\compscale, \yshift-1*\compscale*\ysep);
\coordinate (label1) at (\xshift-\labelshift*\compscale, \yshift+0*\compscale*\ysep);
\coordinate (label2) at (\xshift-\labelshift*\compscale, \yshift+1*\compscale*\ysep);

\coordinate (plus10) at (\xshift-\labelshift*\compscale, \yshift-0.5*\compscale*\ysep);
\coordinate (plus11) at (\xshift-\labelshift*\compscale, \yshift+0.5*\compscale*\ysep);

\node[draw=none,fill=none,\comptwocolor] (t20node) at (t20) {\compsize $\CO_V(1,-1)$};
\node[draw=none,fill=none,\componecolor] (t21node) at (t21) {\compsize $\CO_V(0,-1)^{\oplus 2}$};
\node[draw=none,fill=none,\compthreecolor] (t22node) at (t22) {\compsize $\CO_V(-1,-1)^{\oplus 3}$};

\node[draw=none,fill=none,\comptwocolor] (t10node) at (t10) {\compsize $\CO_V(1,-1)^{\oplus 4}$};
\node[draw=none,fill=none,\componecolor] (t11node) at (t11) {\compsize $\CO_V(0,-1)^{\oplus 6}$};
\node[draw=none,fill=none,\compthreecolor] (t12node) at (t12) {\compsize $\CO_V(-1,0)^{\oplus 4}$};

\node[draw=none,fill=none,\comptwocolor] (t00node) at (t00) {\compsize $\CO_V(1,-1)^{\oplus 3}$};
\node[draw=none,fill=none,\componecolor] (t01node) at (t01) {\compsize $\CO_V(0,0)^{\oplus 6}$};

\node[draw=none,fill=none,black] (plus10node) at (plus10) { $\bigoplus$};
\node[draw=none,fill=none,black] (plus11node) at (plus11) { $\bigoplus$};

\draw[->,\comptwocolor,\compthick] (t20node) -- (t10node);
\draw[->,\componecolor,\compthick] (t21node) -- (t11node);
\draw[->,\compthreecolor,\compthick] (t22node) -- (t12node);

\draw[->,\comptwocolor,\compthick] (t10node) -- (t00node);
\draw[->,\componecolor,\compthick] (t11node) -- (t01node);

\node[draw=none,fill=none,\comptwocolor] at (label0) {$(i^* F_\bullet^D(D))_{\lambda_3=-1}$};
\node[draw=none,fill=none,\componecolor] at (label1) {$(i^* F_\bullet^D(D))_{\lambda_3=0}$};
\node[draw=none,fill=none,\compthreecolor] at (label2) {$(i^* F_\bullet^D(D))_{\lambda_3=1}$};

\node[draw=none,fill=none,\comptwocolor] at ($(label0)!0.5!(t20)$) {$=$};
\node[draw=none,fill=none,\componecolor] at ($(label1)!0.5!(t21)$) {$=$};
\node[draw=none,fill=none,\compthreecolor] at ($(label2)!0.5!(t22)$) {$=$};

\node[draw=none,fill=none] at (\xshift-1.8*\xsep*\compscale,\yshift) {$i^* F_\bullet^D(D) \ \ = $};

\end{tikzpicture}\]
    \caption{Here we present a pictorial representation of the restriction to $V$ of $F_\bullet^D(D)$ mirroring that of \Cref{fig:full_cellular_resolution}. We drop the arrows corresponding to differentials which are sent to 0 by restriction to $V$ in order to illustrate the splitting which occurs upon restriction, as described in \Cref{example:H2_grade_restriction_rules_geometrically}. We let $(i,j)$ denote the line bundle $\CO_V(i,j)$. All non-acyclic complexes are encircled.}
    \label{fig:restricted_cellular_resolution}
\end{figure}

\end{example}

\section{A criterion for window comparison}\label{section:window_comparison}
In light of \Cref{cor: window equality}, the goal of this section is to develop machinery for computing the $\lambda$-weights of $F^{D}_{\bullet}(D)$. Concretely, we start with a projective, primitive geometric GIT problem defined by the inclusion $G\subset \Gm^n$. Let $\Gamma$ be a chamber of the associated secondary fan $\Sigma_{GKZ}$ and let $\bs{\lambda} = \lambda_1,\hdots, \lambda_k$ be an HKKN stratification of $V(B_\Gamma)$. Given this data, we consider the associated grade restriction window $W^{\bs{w}}_{\bs{\lambda}}$ as in \Cref{def:window}. To compare $W^{\bs{w}}_{\bs{\lambda}}$ with $W^{\Cox}_{\Gamma}$, it suffices to check whether for any $d\in \wh{G}\cap \Gamma$, the generators of $W^{\Cox}_{\Gamma}$ (specifically the complexes $F^{D}_{\bullet}(D)$) are in $W^{\bs{w}}_{\bs{\lambda}}$. 

In \S\ref{subsection:eq_class_splitting} we describe how to use \Cref{prop:derivednakayama} (a derived version of Nakayama's Lemma) to do the $\lambda$-weight computation by restricting to a finite set of distinguished points.  This reduces the calculation of $\lambda$-weights of $F^D_\bullet(D)$ to vector space homology. In \S\ref{subsection:top_interpretation_of_splitting} we give a topological interpretation for this vector space homology. Concretely, we show that this homology is the relative homology of a CW-pair of subcomplexes of $P_D$. We use this to describe a topological criterion for when $W_{\bs{\lambda}}^{\bs{w}} = W_\Gamma^{\Cox}$. In \S\ref{subsection:reduction_to_nerve}, we introduce a simplicial complex whose homology is equal to the relative homology of this CW-pair. We then provide a number of sufficient conditions for the contractibility of this simplicial complex.

\subsection{Restriction of cellular resolution to distinguished points}\label{subsection:eq_class_splitting}
We now explain how  \Cref{prop:derivednakayama} can be used to check whether $F^{D}_{\bullet}(D) \in W_{\bs{\lambda}}^{\bs{w}}$. Recall from \Cref{def:grade_restriction_windows},
\[
W_{\bs{\lambda}}^{\bs{w}} = \big\{\CF_\bullet \in \Db([\A^n/G]) \ | \ \mu(\CF_\bullet, \lambda_{j}, (X_j)^{=0}_{\lambda_{j}}) \subset [w_j, w_j+T_{\lambda_j}^{>0})\text{ for }j=1,\hdots, \ell\big\}.
\]
In this subsection we describe how to compute the weights $\mu(F^D_\bullet(D), \lambda_j, (X_j)_{\lambda_j}^{=0})$ for a particular $\lambda_j$ in the HKKN stratification $\bs{\lambda}$. Let $i:(X_j)_{\lambda_j}^{=0}\inc \A^n$ be the inclusion. Recall from \eqref{eqn:splitting} there is a splitting:
\begin{equation}\label{eq:splittingF}
i^*(F^{D}_{\bullet}(D)) \cong \bigoplus_{w\in \Z}(i^*F^{D}_{\bullet}(D))_{\lambda_{j} =w}.
\end{equation}
By definition, $w \in \mu(F^{D}_{\bullet}(D), \lambda_{j}, (X_j)_{\lambda_j}^{=0})$ if and only if the summand of \eqref{eq:splittingF} of $\lambda_{j}$-weight $w$, $(i^*F^{D}_{\bullet}(D))_{\lambda_{j} =w}$, is not acyclic. We demonstrate now that \Cref{prop:derivednakayama}  allows us to determine the acyclicity of $(i^*F^{D}_{\bullet}(D))_{\lambda_{j} =w}$.

\begin{proposition}\label{prop:derivednakayamaonF}
The complex $(i^*F^{D}_{\bullet}(D))_{\lambda_{j} =w}$ is acyclic if and only if for all distinguished points $p_{L}$ in $(X_j)_{\lambda_j}^{=0}$, the complex of vector spaces $p^{*}_{L}((i^*F^{D}_{\bullet}(D))_{\lambda_{j} =w})$ is acyclic.
\end{proposition}
\begin{proof}
Note that $F^{D}_{\bullet}$ was defined as a free $\Z^{n}$-graded complex of $S$-modules or, equivalently, a $T$-equivariant complex of free sheaves. Since $(X_j)_{\lambda_j}^{=0}$ is $T$-invariant, the pullback $i^*F^{D}_{\bullet}(D)$ admits a $T$-equivariant structure. Hence, we can apply \Cref{prop:derivednakayama} (Derived Nakayama's Lemma). Since line bundles are adapted to derived pullback, we may take the usual pullback to $p_L$ instead.
\end{proof}
In light of \Cref{prop:derivednakayamaonF}, we want to understand the complex 
\[p^{*}_{L}((i^*F^{D}_{\bullet}(D))_{\lambda_{j} =w}),\] where $p_{L}$ is any distinguished point of the $T$-action on $(X_j)^{=0}_{\lambda_{j}}$. Since coherent sheaves on $p_{L}$ are merely vector spaces, the complex $p^{*}_{L}((i^*F^{D}_{\bullet}(D))_{\lambda_{j} =w}$ is just a complex of vector spaces. Since $F^{D}_{\bullet}$ was constructed by lifting a cellular homology complex on $P_{D}$, it is not surprising that $p^{*}_{L}((i^*F^{D}_{\bullet}(D))_{\lambda_{j} =w})$ has a cellular structure as we now explain.  

We define
\[\cells(P_D)_{\lambda_j=w} \coloneqq \{\delta\in \cells(P_D) \ : \ \wh{\lambda}_j(\pi(\ceil{\delta})) = w\}.\]
We also define an equivalence relation on the cells in $\cells(P_D)_{\lambda_j=w}$ induced by the set $L \subset[n]$ from the distinguished point $p_{L}$. We say $\delta \sim_{L} \tau$ if and only if $\ceil{\delta}-\ceil{\tau} = \sum_{j\in L}c_{j}e_{j}$, where $c_{j} \in \Z$. Equivalently, $\delta \sim_L \tau$ if and only if $\ceil{\delta}$ and $\ceil{\tau}$ are equal except possibly in their $L$ coordinates. We denote all cells in the same equivalence class as cell $\delta$ by $\cells(P_D)_{\lambda_j=w,[\delta]}$. 

\begin{proposition}\label{prop:algebraic_splitting_of_F}
Set $ K_{\bullet} \coloneqq p^{*}_{L}((i^*F^{D}_{\bullet}(D))_{\lambda_{j} =w})$. Then,
\[
K_{i} = \bigoplus_{\delta \in \cells^{i}(P_D)_{\lambda_j=w}} k\cdot e_{\delta} 
\]
and
\[
\del_{i}(e_{\delta}) = \bigoplus_{\tau\in \cells^{i-1}(P_D)_{\lambda_j=w, [\delta]}}\ve(\delta,\tau)e_\tau.
\]
\end{proposition}
\begin{proof}

Note that pulling back $F^D_{\bullet}(D)_{\lambda_j = w}$ to $p_{L}$ is the same as pulling back the original complex, $F^D_{\bullet}(D)$, to $p_{L}$ and considering the associated summand of $\lambda_{j}$-weight $w$. Pulling back to $p_{L}$ replaces every free $S$-module in the complex with $k$. Since $\cells(P_D)_{\lambda_j=w}$ is precisely the set of all cells of $P_{D}$ of $\lambda_{j}$-weight $w$,  
\[
K_{i} = \bigoplus_{\delta \in \cells^{i}(P_D)_{\lambda_j=w}} k\cdot e_{\delta}. 
\]

We now justify the formula for the differentials of the pullback. The pullback evaluates the monomial differentials of $(i^*F^{D}_{\bullet}(D))_{\lambda_{j} =w})$ at $p_{L}$. Recall that the differentials in $F_{\bullet}(D)$ arise from boundary maps and are given by monomials of the form $\ve(\delta,\tau)x^{\lceil \delta \rceil - \lceil \tau \rceil}$ from $S(-\ceil{\delta})$ to $S(-\ceil{\tau})$ where $\tau$ is a codimension 1 boundary cell of $\delta$. 
\begin{equation*}%\label{eqn:lauren_differentials}
\begin{tikzcd}[column sep = 2cm] \del_i\coloneqq S(-(\ceil{\delta}+D))\arrow[r,"\varepsilon(\delta{,}\tau)x^{\ceil{\delta}-\ceil{\tau}}"] & S(-(\ceil{\tau}+D)).
\end{tikzcd}
\end{equation*}
If $\delta \sim_{L} \tau$, then, by definition, $\ceil{\delta}-\ceil{\tau} = \sum_{j\in L}c_{j}e_{j}$. Since $(p_{L})_{t}=1$ if $t \in L$ and $(p_{L})_{t}=0$ otherwise, the evaluation of $x^{\lceil \delta \rceil - \lceil \tau \rceil}$ at $p_{L}$ is equal to 1. In summary, 
\begin{equation}\label{eq:restrictedLauren'sdiff}
\begin{tikzcd}[column sep = 1.3cm] p^{*}_{L}(\del_i): k\cdot e_{\delta}\arrow[r,"\varepsilon(\delta{,}\tau)"] & k\cdot e_{\tau}.
\end{tikzcd}
\end{equation}
If $\delta \not \sim_{L} \tau$, $\ceil{\delta}-\ceil{\tau}$ must have a non-zero coordinate in the set of coordinates $[n]\setminus L$. Hence, the evaluation of $x^{\lceil \delta \rceil - \lceil \tau \rceil}$ at $p_{L}$ is $0$. Together, we get
\[
\del_{i}(e_{\delta}) = \bigoplus_{\tau\in \cells^{i-1}(P_D)_{\lambda_j=w, [\delta]}}\ve(\delta,\tau)e_\tau.
\]
\end{proof}

\begin{corollary}\label{cor:eq_class_splitting_of_F}
Let $K_\bullet$ be as in \Cref{prop:algebraic_splitting_of_F}. There is a splitting:
\[
K_{\bullet} = \bigoplus_{[\delta]\in\cells(P_{D})_{\lambda_j=w}/\sim_{L}}K^{[\delta]}_{\bullet},
\]
where $K^{[\delta]}_{\bullet}$ is the subcomplex of $K_{\bullet}$ of summands $k\cdot e_\tau$ for $\tau\in \cells(P_D)_{\lambda_j=w, [\delta]}$.
\end{corollary}
\begin{proof}
This follows immediately from the description of $\del_i$ in \Cref{prop:algebraic_splitting_of_F}.
\end{proof}

\subsection{A topological interpretation of the homology complexes}\label{subsection:top_interpretation_of_splitting}
Building on the cellular description from the previous section, we now give a topological interpretation to the complexes $K_\bullet^{[\delta]}$ of \Cref{cor:eq_class_splitting_of_F}. We start with the following topological subspace of $P_D$. 
\begin{definition}\label{def:S_regions} Let $L\subset [n]$ and $\delta\in \cells(P_D)$. We define
\[\CS_{\delta, L} \coloneqq P_D\cap \lp \ceil{\delta}+ \prod_{i\in [n]\setminus L}(-1,0]\times \prod_{t\in L}(-\infty,\infty)\rp.\]
\end{definition}
The regions $\CS_{\delta,L}$ precisely capture the $\sim_L$ equivalence classes in the following sense.
\begin{lemma}\label{lemma:S_is_union_of_cells} Let $\delta\in \cells(P_D)_{\lambda_j=w}$. Then
\[\CS_{\delta,L} = \bigcup_{\tau\in \cells(P_D)_{\lambda_j=w,[\delta]}}\tau.\]
\end{lemma}
\begin{proof} 
First, we show that $\CS_{\delta,L}$ is contained in the union. Let $x \in \CS_{\delta,L}$. Since by definition $\CS_{\delta,L} \subseteq P_{D}$, $x \in \tau \subset P_{D}$ for some cell $\tau$ of $P_D$. Moreover, by definition of $\CS_{\delta,L}$, 
\[x = \ceil{\delta}-\sum_{i\in [n]\setminus L}{b_{i}e_{i}}+\sum_{k\in L}{c_{k}e_{k}}\]
for some $b_{i}\in [0,1)$ and $c_k \in \R$. Then, 
\[\ceil{\tau}=\ceil{x}= \ceil{ \ceil{\delta}-\sum_{i\in [n]\setminus L}{b_{i}e_{i}}+\sum_{k\in L}{c_{k}e_{k}}} = \ceil{\delta}+ \sum_{k\in L}\ceil{c_{k}}e_{k} .\]
Because $p_L\in (X_j)_{\lambda_j}^{=0}$, $\wh{\lambda}_j(\pi_\R(e_k)) = 0$ for $k\in L$. So, the above equality further implies that 
\[\wh{\lambda}_j(\pi_\R(\ceil{\tau})) = \wh{\lambda}_j(\pi_\R(\ceil{\delta}))=w,\]
so $\tau\in \cells(P_D)_{\lambda_j=w}$. As $\tau\sim_L\delta$ we get $\tau\in \cells(P_D)_{\lambda_j=w,[\delta]}$ as desired.

For the opposite inclusion, take any $\tau \in \cells(P_D)_{\lambda_j=w,[\delta]}$. Then, by definition, $\ceil{\tau} - \ceil{\delta} = \sum_{k \in L}{c_{k}e_{k}}$ for some $c_k\in \Z$. Rewriting, we get $\ceil{\tau} = \ceil{\delta} + \sum_{k \in L}{c_{k}e_{k}}$. Thus
\[\tau \subset P_D\cap \lp \prod_{i\in [n]\setminus L}(\ceil{\delta}_i-1,\ceil{\delta}_i]\times \prod_{k\in L} (c_k-1,c_k] \rp\subset \CS_{\delta,L}.\]
\end{proof}
While $\CS_{\delta,L}$ is a union of cells of $P_D$ as evidenced by \Cref{lemma:S_is_union_of_cells}, it is only locally-closed and thus is not a CW subcomplex of $P_D$. However, the closure and boundary of $\CS_{\delta,L}$ are subcomplexes of $P_D$, which we describe in the following lemma.
\begin{lemma}\label{lemma:coordinate_description_of_closure_boundary} Let $L\subset [n]$ and $\delta\in \cells(P_D)$. Let $\overline{\CS}_{\delta,L}$ be the closure of $\CS_{\delta,L}$ and $\del \CS_{\delta,L}\coloneqq \overline{\CS}_{\delta,L}\setminus \CS_{\delta,L}$ be the boundary. Then
\begin{align*}
    \overline{\CS}_{\delta,L} &= P_D\cap \lp \prod_{i\in [n]\setminus L}[\ceil{\delta}_i-1,\ceil{\delta}_i]\times \prod_{t\in L}(-\infty,\infty)\rp, \\
    \del \CS_{\delta,L} &= P_D\cap \bigcup_{i \in [n]\setminus L}\lp \{\ceil{\delta}_{i}-1\}\times \prod_{j\in [n]\setminus (L\cup \{i\})}[\ceil{\delta}_j-1,\ceil{\delta}_j]\times \prod_{t\in L}(-\infty,\infty)\rp.
\end{align*}
Further, $\overline{\CS}_{\delta,L}, \del \CS_{\delta,L}$ are both subcomplexes of $P_D$.
\end{lemma}
\begin{proof} We will only prove the coordinate description of $\overline{\CS}_{\delta,L}$; the coordinate description of $\del \CS_{\delta, L}$ is an immediate consequence of this, and that $\overline{\CS}_{\delta,L}$ and $\del \CS_{\delta,L}$ are subcomplexes of $P_D$ follows directly from their coordinate descriptions.

To see that $\overline{\CS}_{\delta,L}$ agrees with the above coordinate description, first observe that $\overline{\CS}_{\delta,L}$ is contained in the right-hand side by virtue of the right-hand side being a closed set which contains $\CS_{\delta,L}$. To see the opposite inclusion, take any $y$ in the right-hand side and fix any $x\in \delta\subset \CS_{\delta,L}$ (see \Cref{lemma:S_is_union_of_cells}). Let $\gamma$ be a straight line path defined as:
\[
\gamma(t)\coloneqq x(1-t)+yt.
\]
Writing the coordinate expansion of $\gamma = (\gamma_1,\hdots, \gamma_n)$, we have that:
\[
\gamma_{i}(t)= x_{i}(1-t)+y_{i}t.
\]
We first show that $\gamma([0,1))\subset \CS_{\delta,L}$. Since $P_{D}$ is a convex polytope and $\gamma(0)=x,\gamma(1)=y \in P_{D}$, it follows that for any $t\in[0,1)$, $\gamma(t) \in P_{D}$. Therefore, using the coordinate description of $\CS_{\delta,L}$ (see \Cref{def:S_regions}), it suffices to show that for any $t\in[0,1)$ and $i \in [n]\setminus L$, $\gamma_{i}(t) \in (\ceil{\delta}_i-1,\ceil{\delta}_i]$. By the coordinate description of $\CS_{\delta,L}$, we have that $x_i > \ceil{\delta}_i-1$ and $y_i \ge \ceil{\delta}_i-1$ by our assumption on $y$. Together, this implies
\[\gamma_i(t) = x_i(1-t)+y_it > (\ceil{\delta}_i-1)(1-t)+(\ceil{\delta}_i-1)t = \ceil{\delta}_i-1.\]
Similarly, $x_i \le \ceil{\delta}_i$ and $y_i \le \ceil{\delta}_i$ together imply $\gamma_i(t) \le \ceil{\delta}_i$ for $t\in [0,1)$. Together, this implies

$\gamma([0,1))\subset \CS_{\delta,L}$ hence $\gamma([0,1])\subset \overline{\CS}_{\delta,L}$.
\end{proof}

\begin{example}\label{example:S_delta_L_regions} 
Let us continue with our running example of $\CH_2$. Recall the ample lift $D = (0,0,1,1)$ and $P_D$ from \Cref{example:H2_cellular_resolution}. In \Cref{fig:S_delta_L_regions} we identify a 2-cell $\delta$ and a 1-cell $\tau$ in $P_D$. Choosing $L = \{4\}$, we depict the regions $\CS_{\delta,L},\overline{\CS}_{\delta,L},\del \CS_{\delta,L}$ and the corresponding regions for $\tau$. As guaranteed by \Cref{lemma:S_is_union_of_cells}, we see that $\CS_{\delta,L}$ is a union of $\delta$ with another 1-cell, while $\CS_{\tau,L} = \tau$. While $\CS_{\delta,L}, \CS_{\tau,L}$ are not subcomplexes of $P_D$, we see that the closure and boundaries are subcomplexes of $P_D$ as guaranteed by \Cref{lemma:coordinate_description_of_closure_boundary}. 
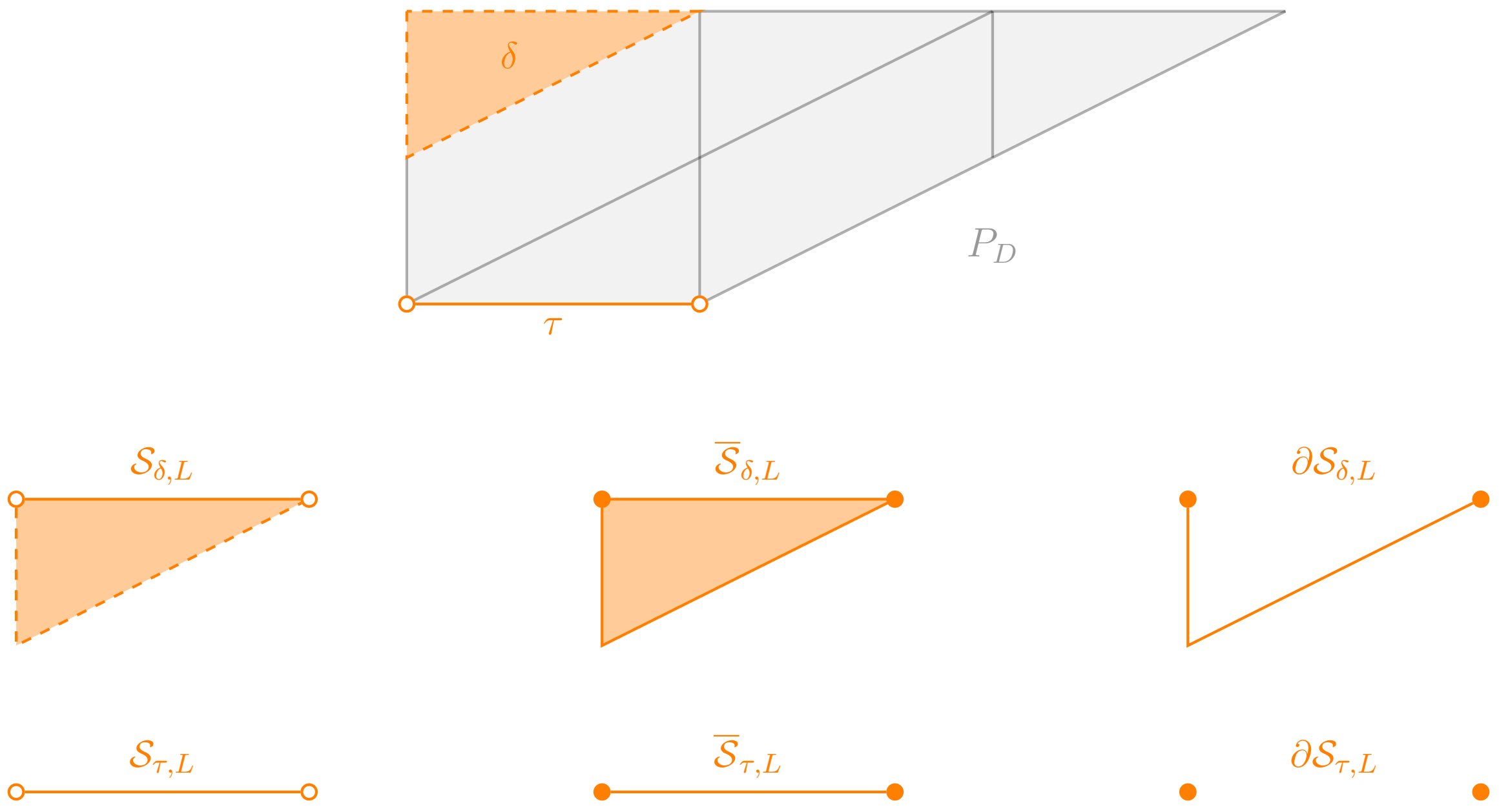
\begin{figure}[h!]
    \centering
\[\begin{tikzpicture}
\def\polyscale{3}
\def\xshift{0}
\def\yshift{0}
\def\polyop{0.3}
\def\cellop{0.4}
\def\cellOneColor{orange}
\def\cellTwoColor{orange}

\coordinate (v00) at (\xshift+0*\polyscale,\yshift+0*\polyscale);
\coordinate (v01) at (\xshift+0*\polyscale, \yshift+0.5*\polyscale);
\coordinate (v02) at (\xshift+0*\polyscale,\yshift+1*\polyscale);
\coordinate (v10) at (\xshift+1*\polyscale,\yshift+0*\polyscale);
\coordinate (v11) at (\xshift+1*\polyscale,\yshift+0.5*\polyscale);
\coordinate (v12) at (\xshift+1*\polyscale,\yshift+1*\polyscale);
\coordinate (v21) at (\xshift+2*\polyscale,\yshift+0.5*\polyscale);
\coordinate (v22) at (\xshift+2*\polyscale,\yshift+1*\polyscale);
\coordinate (v32) at (\xshift+3*\polyscale,\yshift+1*\polyscale);

% filling cells
\draw[draw=none,fill=\cellOneColor,opacity=\cellop] (v02) -- (v01) -- (v12) -- cycle;
\draw[draw=none,fill=gray,opacity=0.1] (v00) -- (v01) -- (v12) -- (v32) -- (v10) -- cycle;

\node[inner sep = 1.5pt, draw=\cellTwoColor, circle, thick, fill=white,opacity=1] (v00node) at (v00) {};
\node[inner sep = 1.5pt, draw=\cellTwoColor, circle, thick, fill=white,opacity=1] (v10node) at (v10) {};

\draw[-,black,thick,opacity=\polyop] (v00node) -- (v01);
\draw[-,black,thick,opacity=\polyop] (v10node) -- (v12);
\draw[-,black,thick,opacity=\polyop] (v21) -- (v22);
\draw[-,black,thick,opacity=\polyop] (v00node) -- (v10node);
%\draw[-,black,thick,opacity=\polyop] (v02node) -- (v12node);
\draw[-,black,thick,opacity=\polyop] (v12) -- (v32);
%\draw[-,black,thick,opacity=\polyop] (v01) -- (v12node);
\draw[-,black,thick,opacity=\polyop] (v00node) -- (v22);
\draw[-,black,thick,opacity=\polyop] (v10node) -- (v32);

% P_D
\node[draw=none,fill=none,black,opacity=\cellop] at (\xshift+2*\polyscale,\yshift+0.2*\polyscale) {$P_D$};

\draw[dashed,\cellOneColor,thick] (v02) -- (v01) -- (v12) -- (v02);
\node[draw=none,fill=none,\cellOneColor] at ($(v02)!0.5!(v12)!0.3!(v01)$) {\small $\delta$}; 

% cell 2
\draw[-,\cellTwoColor,thick] (v00node) -- (v10node);
\node[draw=none,fill=none,\cellTwoColor] at ($(v00)!0.5!(v10)$) [below]{\small $\tau$}; 

% S_delta,L drawings
\def\xshiftTwo{-4}
\def\yshiftTwo{-2}
\def\cellscale{3}
\def\cellhsep{2}
\def\cellvsep{1}

\coordinate (hsepvec) at (\cellhsep*\cellscale,0);
\coordinate (vsepvec) at (0,\cellvsep*\cellscale);

\coordinate (v02ul) at (\xshiftTwo,\yshiftTwo);
\coordinate (v01ul) at (\xshiftTwo,\yshiftTwo-0.5*\cellscale);
\coordinate (v12ul) at (\xshiftTwo+\cellscale,\yshiftTwo);
\coordinate (v02um) at ($(v02ul)+(hsepvec)$);
\coordinate (v01um) at ($(v01ul)+(hsepvec)$);
\coordinate (v12um) at ($(v12ul)+(hsepvec)$);
\coordinate (v02ur) at ($(v02ul)+2*(hsepvec)$);
\coordinate (v01ur) at ($(v01ul)+2*(hsepvec)$);
\coordinate (v12ur) at ($(v12ul)+2*(hsepvec)$);

\coordinate (v00ll) at (\xshiftTwo,\yshiftTwo-\cellvsep*\cellscale);
\coordinate (v10ll) at (\xshiftTwo+\cellscale,\yshiftTwo-\cellvsep*\cellscale);
\coordinate (v00lm) at ($(v00ll)+(hsepvec)$);
\coordinate (v10lm) at ($(v10ll)+(hsepvec)$);
\coordinate (v00lr) at ($(v00ll)+2*(hsepvec)$);
\coordinate (v10lr) at ($(v10ll)+2*(hsepvec)$);

\draw[draw=none,fill=\cellOneColor,opacity=\cellop] (v02ul) -- (v01ul) -- (v12ul) -- cycle;
\draw[draw=none,fill=\cellOneColor,opacity=\cellop] (v02um) -- (v01um) -- (v12um) -- cycle;

% cell 1 nodes
\node[inner sep = 1.5pt, draw=\cellOneColor, circle, thick, fill=white,opacity=1] (v02ulnode) at (v02ul) {};
\node[inner sep = 1.5pt, draw=\cellOneColor, circle, thick, fill=white,opacity=1] (v12ulnode) at (v12ul) {};
\node[inner sep = 1.5pt, draw=\cellOneColor, circle, thick, fill=\cellOneColor,opacity=1] (v02umnode) at (v02um) {};
\node[inner sep = 1.5pt, draw=\cellOneColor, circle, thick, fill=\cellOneColor,opacity=1] (v12umnode) at (v12um) {};
\node[inner sep = 1.5pt, draw=\cellOneColor, circle, thick, fill=\cellOneColor,opacity=1] (v02urnode) at (v02ur) {};
\node[inner sep = 1.5pt, draw=\cellOneColor, circle, thick, fill=\cellOneColor,opacity=1] (v12urnode) at (v12ur) {};

% cell 1 labels
\node[fill=none,draw=none, \cellOneColor] at ($(v02ul)!0.5!(v12ul)$) [above] {\small $\CS_{\delta, L}$};
\node[fill=none,draw=none, \cellOneColor] at ($(v02um)!0.5!(v12um)$) [above] {\small $\overline{\CS}_{\delta, L}$};
\node[fill=none,draw=none, \cellOneColor] at ($(v02ur)!0.5!(v12ur)$) [above] {\small $\del \CS_{\delta, L}$};

% cell 1 drawing
\draw[-,thick,\cellOneColor] (v02ulnode) -- (v12ulnode);
\draw[dashed,thick,\cellOneColor] (v02ulnode) -- (v01ul) -- (v12ulnode);
\draw[-,thick,\cellOneColor] (v02um) -- (v01um) -- (v12um) -- (v02um);
\draw[-,thick,\cellOneColor] (v02urnode) -- (v01ur) -- (v12urnode);

% cell 2 nodes
\node[inner sep = 1.5pt, draw=\cellTwoColor, circle, thick, fill=white,opacity=1] (v00llnode) at (v00ll) {};
\node[inner sep = 1.5pt, draw=\cellTwoColor, circle, thick, fill=white,opacity=1] (v10llnode) at (v10ll) {};
\node[inner sep = 1.5pt, draw=\cellTwoColor, circle, thick, fill=\cellTwoColor,opacity=1] (v00lmnode) at (v00lm) {};
\node[inner sep = 1.5pt, draw=\cellTwoColor, circle, thick, fill=\cellTwoColor,opacity=1] (v10lmnode) at (v10lm) {};
\node[inner sep = 1.5pt, draw=\cellTwoColor, circle, thick, fill=\cellTwoColor,opacity=1] (v00lrnode) at (v00lr) {};
\node[inner sep = 1.5pt, draw=\cellTwoColor, circle, thick, fill=\cellTwoColor,opacity=1] (v10lrnode) at (v10lr) {};

% cell 2 drawing
\draw[-,thick, \cellTwoColor] (v00llnode) -- (v10llnode);
\draw[-,thick, \cellTwoColor] (v00lmnode) -- (v10lmnode);

% cell 2 labels
\node[fill=none,draw=none, \cellTwoColor] at ($(v00ll)!0.5!(v10ll)$) [above] {\small $\CS_{\tau, L}$};
\node[fill=none,draw=none, \cellTwoColor] at ($(v00lm)!0.5!(v10lm)$) [above] {\small $\overline{\CS}_{\tau, L}$};
\node[fill=none,draw=none, \cellTwoColor] at ($(v00lr)!0.5!(v10lr)$) [above] {\small $\del \CS_{\tau, L}$};
\end{tikzpicture}\]
    \caption{\label{fig:S_delta_L_regions}We highlight a 2-cell $\delta$ and a 1-cell $\tau$ in the cellular complex $P_D$ of \Cref{example:S_delta_L_regions}. Setting $L = \{4\}$, we depict $\CS_{\delta,L}, \overline{\CS}_{\delta,L}$, and $\del\CS_{\delta,L}$. We draw the corresponding regions for $\tau$ below.}
\end{figure}
\end{example}

Let $C_\bullet(\overline{\CS}_{\delta,L}), C_\bullet(\del \CS_{\delta,L})$ be the chain subcomplexes of $C_\bullet(P_D)$ (see \eqref{eq:cellcomplexonPd}) associated to the CW subcomplexes $\overline{\CS}_{\delta,L}, \del \CS_{\delta,L}$ respectively. Let $\del^{\overline{\CS}}, \del^{\del \CS}$ be the respective differentials. By definition, $\del^{\overline{\CS}}, \del^{\del \CS}$ are exactly the differentials of $C_\bullet(P_D)$ restricted to the respective subcomplexes of $P_D$. With this in hand, we are in a position to describe $K_\bullet^{[\delta]}$ topologically.
\begin{proposition}\label{prop: relativehomology}
The complex $K^{[\delta]}_{\bullet}$ is the relative cellular chain complex of the pair
$(\overline{\CS}_{\delta,L},\del \CS_{\delta,L})$.
\end{proposition}

\begin{proof} 
First observe that as a consequence of \Cref{lemma:S_is_union_of_cells} and \Cref{lemma:coordinate_description_of_closure_boundary} we get that
\begin{equation}\label{eq:cell_disjoint_union}
    \cells^i(\overline{\CS}_{\delta,L}) = \cells^i(P_D)_{\lambda_j=w,[\delta]}\sqcup \cells^i(\del \CS_{\delta,L}).
\end{equation}
As the terms of $K_\bullet^{[\delta]}$ are indexed by $\cells(P_D)_{\lambda_j=w,[\delta]}$, we get well-defined morphisms $p:C_i(\overline{\CS}_{\delta,L})\to K_i^{[\delta]}$ via
\[p(e_\tau)\coloneqq \begin{cases} e_\tau &\text{if }\tau\in \cells^i(P_D)_{\lambda_j=w,[\delta]}\\
0 &\text{if }\tau\in \cells^i(\del \CS_{\delta,L})
\end{cases}\]
We now show that $p$ is a well-defined chain map. Let $\tau\in \cells^i(\overline{\CS}_{\delta,L})$ and, as in the first case, suppose $\tau \in \cells^i(P_D)_{\lambda_j=w,[\delta]}$. Then
\begin{align*}
    (p\circ \del_i^{\overline{\CS}})(e_\tau) &=  \sum_{\alpha\in \cells^{i-1}(\overline{\CS}_{\delta,L})}\ve(\tau,\alpha)p(e_\alpha) & \text{by definition of $\del_i^{\overline{\CS}}$}\\
    &= \sum_{\alpha\in \cells^{i-1}(P_D)_{\lambda_j=w,[\delta]}}\ve(\tau,\alpha) e_\alpha &\text{by \eqref{eq:cell_disjoint_union} and definition of $p$}\\
    &= \del_i(e_\tau) & \text{by \Cref{prop:algebraic_splitting_of_F}}\\
    &= (\del_i \circ p)(e_\tau). & \text{as $\tau\in \cells^i(P_D)_{\lambda_j=w,[\delta]}$}
\end{align*}
Now as in the second case, suppose $\tau\in \cells^i(\del \CS_{\delta,L})$ so that $(\del_i \circ p)(e_\tau) = 0$. By \Cref{lemma:coordinate_description_of_closure_boundary}, $\del \CS_{\delta,L}$ is a complex, so every boundary cell $\alpha$ of $\tau$ is also a cell of $\del \CS_{\delta,L}$. So we find
\[(p\circ \del_i^{\overline{\CS}})(e_\tau) = \sum_{\alpha\in \cells^{i-1}(\overline{\CS}_{\delta,L})}\ve(\tau,\alpha)p(e_\alpha) = 0 = (\del_i\circ p)(e_\tau).\]
Thus, $p$ is a well-defined chain map. 

Now let $\iota:C_\bullet(\del \CS_{\delta,L})\to C_\bullet(\overline{\CS}_{\delta,L})$ be the natural injective chain map induced by the subcomplex inclusion $\del \CS_{\delta,L}\subset \overline{\CS}_{\delta,L}$. It follows from \eqref{eq:cell_disjoint_union} that $p$ is surjective and $\iota(C_\bullet(\del \CS_{\delta,L})) = \ker p$ so that we get a short exact sequence
\[
0 \to   C_{\bullet}(\del \CS_{\delta,L})\xrightarrow{\iota} C_{\bullet}(\overline{\CS}_{\delta,L}) \xrightarrow{p} K_\bullet^{[\delta]} \to 0,
\]
completing the proof.
\end{proof}
From the short exact sequence guaranteed by \Cref{prop: relativehomology} we get the following topological description of the acyclicity of $K_\bullet^{[\delta]}$.
\begin{corollary}\label{cor:topological_acyclicity} The complex $K_\bullet^{[\delta]}$ is acyclic if and only if $\del \CS_{\delta,L}$ is acyclic. In particular, if $\del \CS_{\delta,L}$ is contractible, then $K_\bullet^{[\delta]}$ is acyclic.
\end{corollary}
\begin{proof} First observe that as $\overline{\CS}_{\delta,L}$ is a closed, convex set in $\R^n$, it is acyclic, so $\wt{H}_\ast(\overline{\CS}_{\delta,L}) \cong 0$. Then, from the long exact sequence guaranteed by \Cref{prop: relativehomology}, it follows that the homology of $K_\bullet^{[\delta]}$ is isomorphic to the reduced homology of $\del \CS_{\delta,L}$. 
\end{proof}

\begin{lemma}\label{lemma:lambda_weights_of_cox_generators} The $\lambda_j$-weights of $F_\bullet^D(D)$ along $(X_j)_{\lambda_j}^{=0}$ are given by
\[\mu(F_\bullet^D(D), \lambda_j, (X_j)_{\lambda_j}^{=0}) = \bigcup_{p_L\in (X_j)_{\lambda_j}^{=0}} \left\{ \wh{\lambda}_j(\pi(\ceil{\delta}))\ : \ \del\CS_{\delta,L}\text{ is not acyclic, for }\delta\in \cells(P_D)\right\}.\]
\end{lemma}
\begin{proof}
First recall the inclusion $i:(X_j)_{\lambda_j}^{=0}\inc \A^n$ and the associated splitting in \eqref{eq:splittingF}:
\[i^*F_\bullet^D(D) \cong \bigoplus_{w\in \Z} (i^*F_\bullet^D(D))_{\lambda_j=w}\]
With this we get the following chain of equalities:
\begin{align*}
    &\mu(F_\bullet^D(D), \lambda_j, (X_j)_{\lambda_j}^{=0})\\
    =& \{w\in \Z \ | \ (i^*F_\bullet^D(D))_{\lambda_j=w}\text{ is not acyclic}\} & \text{by definition}\\
    =& \bigcup_{p_L\in (X_j)_{\lambda_j}^{=0}} \{w\in \Z \ | \ p_L^*((i^*F_\bullet^D(D))_{\lambda_j=w})\text{ is not acyclic}\} & \text{by \Cref{prop:derivednakayamaonF}}\\
    =& \bigcup_{p_L\in (X_j)_{\lambda_j}^{=0}} \{w\in \Z \ | \ K_\bullet^{[\delta]}\text{ is not acyclic for some }[\delta]\in \cells(P_D)_{\lambda_j=w}/\sim_L\} & \text{by \Cref{cor:eq_class_splitting_of_F}}\\
    =& \bigcup_{p_L\in (X_j)_{\lambda_j}^{=0}} \{w\in \Z \ | \ \CS_{\delta,L}\text{ is not acyclic for some }[\delta]\in \cells(P_D)_{\lambda_j=w}/\sim_L\} & \text{by \Cref{cor:topological_acyclicity}}\\
    =& \bigcup_{p_L\in (X_j)_{\lambda_j}^{=0}}\left \{\wh{\lambda}_j(\pi(\ceil{\delta})) \ | \ \CS_{\delta,L}\text{ is not acyclic for some }\delta\in \cells(P_D)\right\}.
\end{align*}
The last equality follows from the fact that
\[\cells(P_D) = \bigcup_{w\in \Z}\cells(P_D)_{\lambda_j=w}\]
and that $\del \CS_{\delta,L} = \del \CS_{\tau,L}$ for $\delta \sim_L\tau$ by \Cref{lemma:S_is_union_of_cells}.
\end{proof}

\begin{figure}
\begin{tikzpicture}
\def\polyscale{5}
\def\componecolor{red}
\def\comptwocolor{darkGreen}
\def\compthreecolor{orange}
\def\thickness{thick}
\def\dashsty{dashed}
\def\twocellop{0.2}
\def\xshift{0}
\def\yshift{0}
%\tikzset{every node/.style={inner sep=0pt}}
\def\pert{0.2}
\def\slidepert{0.5}
\tikzset{opennode/.style={inner sep=1.5pt, fill= white,circle, thick},
closednode/.style={inner sep=1.5pt, circle, thick}}

\coordinate (e) at (\slidepert*\pert*\polyscale,0);
\coordinate (n) at (0,0.5*\pert*\polyscale);
\coordinate (ne) at (\slidepert*\pert*\polyscale,0.5*\slidepert*\pert*\polyscale);

\coordinate (nne) at ($(n)+(ne)$);
\coordinate (w) at (-\pert*\polyscale,0);

\coordinate (v00) at (\xshift+0*\polyscale,\yshift+0*\polyscale);
\coordinate (v01) at (\xshift+0*\polyscale,\yshift+0.5*\polyscale);
\coordinate (v02) at (\xshift+0*\polyscale,\yshift+1*\polyscale);
\coordinate (v10) at (\xshift+1*\polyscale,\yshift+0*\polyscale);
\coordinate (v11) at (\xshift+1*\polyscale,\yshift+0.5*\polyscale);
\coordinate (v12) at (\xshift+1*\polyscale,\yshift+1*\polyscale);
\coordinate (v21) at (\xshift+2*\polyscale,\yshift+0.5*\polyscale);
\coordinate (v22) at (\xshift+2*\polyscale,\yshift+1*\polyscale);
\coordinate (v32) at (\xshift+3*\polyscale,\yshift+1*\polyscale);

\coordinate (v00n) at ($(v00)+(n)$);
\coordinate (v00e) at ($(v00)+(e)$);
\coordinate (v00ne) at ($(v00)+(n)+(\pert*\polyscale,0)$);
\coordinate (v00nne) at ($(v00)+(nne)$);
\coordinate (v01n) at ($(v01)+(n)$);
\coordinate (v01ne) at ($(v01)+(ne)$);
\coordinate (v01nne) at ($(v01)+(nne)$);
\coordinate (v02e) at ($(v02)+(e)$);
\coordinate (v10w) at ($(v10)+(-\slidepert*\pert*\polyscale,0)$);
\coordinate (v10n) at ($(v10)+(n)$);
\coordinate (v10ne) at ($(v10)+(ne)$);
\coordinate (v10nne) at ($(v10)+(nne)$);
\coordinate (v11n) at ($(v11)+(n)$);
\coordinate (v11ne) at ($(v11)+(ne)$);
\coordinate (v11nne) at ($(v11)+(nne)$);
\coordinate (v12w) at ($(v12)+(w)$);
\coordinate (v12e) at ($(v12)+(e)$);
\coordinate (v21n) at ($(v21)+(n)$);
\coordinate (v21ne) at ($(v21)+(ne)$);
\coordinate (v21nne) at ($(v21)+(nne)$);
\coordinate (v22w) at ($(v22)+(w)$);
\coordinate (v22e) at ($(v22)+(e)$);
\coordinate (v32w) at ($(v32)+(w)$);

% filling cells first
\draw[draw=none,fill=\comptwocolor,opacity=\twocellop] (v00ne) -- (v10n) -- (v11) -- cycle;
\draw[draw=none,fill=\compthreecolor,opacity=\twocellop] (v02e) -- (v01nne) -- (v12w) -- cycle;
\draw[draw=none,fill=\compthreecolor,opacity=\twocellop] (v12e) -- (v11nne) -- (v22w) -- cycle;
\draw[draw=none,fill=\compthreecolor,opacity=\twocellop] (v22e) -- (v21nne) -- (v32w) -- cycle;
\draw[draw=none,fill=\componecolor,opacity=\twocellop] (v00nne) -- (v01ne) -- (v12) -- (v11n) -- cycle;
\draw[draw=none,fill=\componecolor,opacity=\twocellop] (v10nne) -- (v11ne) -- (v22) -- (v21n) -- cycle;

\node[opennode,draw =\comptwocolor] (v00nnode) at (v00n) {};
\node[opennode,draw =\comptwocolor] (v00nenode) at (v00ne) {};
\node[opennode,draw =\compthreecolor] (v00enode) at (v00e) {};
\node[opennode,draw =\componecolor] (v01nnode) at (v01n) {};
\node[opennode,draw =\componecolor] (v01nenode) at (v01ne) {};
\node[opennode,draw =\compthreecolor] (v02enode) at (v02e) {};
\node[opennode,draw =\compthreecolor] (v10wnode) at (v10w) {};
\node[opennode,draw =\comptwocolor] (v10nnode) at (v10n) {};
\node[opennode,draw =\comptwocolor] (v10nenode) at (v10ne) {};
\node[opennode,draw =\componecolor] (v11nnode) at (v11n) {};
\node[opennode,draw =\componecolor] (v11nenode) at (v11ne) {};
\node[opennode,draw =\compthreecolor] (v12wnode) at (v12w) {};
\node[opennode,draw =\compthreecolor] (v12enode) at (v12e) {};
\node[opennode,draw =\componecolor] (v21nnode) at (v21n) {};
\node[opennode,draw =\componecolor] (v21nenode) at (v21ne) {};
\node[opennode,draw =\compthreecolor] (v22wnode) at (v22w) {};
\node[opennode,draw =\compthreecolor] (v22enode) at (v22e) {};
\node[opennode,draw =\compthreecolor] (v32wnode) at (v32w) {};

\node[closednode,draw=\componecolor, fill = \componecolor] (v00node) at (v00) {};
\node[closednode,draw=\comptwocolor, fill = \comptwocolor] (v01node) at (v01) {};
\node[closednode,draw=\componecolor, fill = \componecolor] (v02node) at (v02) {};
\node[closednode,draw=\componecolor, fill = \componecolor] (v10node) at (v10) {};
\node[closednode,draw=\comptwocolor, fill = \comptwocolor] (v11node) at (v11) {};
\node[closednode,draw=\componecolor, fill = \componecolor] (v12node) at (v12) {};
\node[closednode,draw=\comptwocolor, fill = \comptwocolor] (v21node) at (v21) {};
\node[closednode,draw=\componecolor, fill = \componecolor] (v22node) at (v22) {};
\node[closednode,draw=\componecolor, fill = \componecolor] (v32node) at (v32) {};

\draw[\thickness,\comptwocolor] (v00nnode) -- (v01node);
\draw[\thickness,\comptwocolor] (v00nenode) -- (v11) -- (v10nnode);
\draw[\thickness,\comptwocolor] (v10nenode) -- (v21);
\draw[\thickness,\comptwocolor,\dashsty] (v00nenode) -- (v10nnode);

\draw[\thickness,\compthreecolor] (v02enode) -- (v12wnode);
\draw[\thickness,\compthreecolor] (v12enode) -- (v22wnode);
\draw[\thickness,\compthreecolor] (v22enode) -- (v32wnode);
\draw[\thickness,\compthreecolor] (v00enode) -- (v10wnode);
\draw[\thickness,\compthreecolor,\dashsty] (v02enode) -- (v01nne) -- (v12wnode);
\draw[\thickness,\compthreecolor,\dashsty] (v12enode) -- (v11nne) -- (v22wnode);
\draw[\thickness,\compthreecolor,\dashsty] (v22enode) -- (v21nne) -- (v32wnode);

\draw[\thickness,\componecolor] (v01nnode) -- (v02);
\draw[\thickness,\componecolor] (v01nenode) -- (v12) -- (v11nnode);
\draw[\thickness,\componecolor] (v11nenode) -- (v22) -- (v21nnode);
\draw[\thickness,\componecolor] (v21nenode) -- (v32);
\draw[\thickness,\componecolor,\dashsty] (v01nenode) -- (v00nne) -- (v11nnode);
\draw[\thickness,\componecolor,\dashsty] (v11nenode) -- (v10nne) -- (v21nnode);

\node[draw=none,fill=none, \componecolor] at (v00) [below left] {\small $\CS_{\delta_1,L}$};
\node[draw=none,fill=none, \componecolor] at (v10) [below right] {\small $\CS_{\delta_3,L}$};
\node[draw=none,fill=none, \compthreecolor] at ($(v00)!0.5!(v10)$) [below] {\small $\CS_{\delta_2,L}$};
\node[draw=none,fill=none, \comptwocolor] at ($(v00)!0.5!(v11)!0.33!(v10n)$){\small $\CS_{\tau_2,L}$};
\node[draw=none,fill=none, \comptwocolor] at ($(v00n)!0.5!(v01)$) [left]{\small $\CS_{\tau_1,L}$};
\node[draw=none,fill=none, \comptwocolor] at ($(v10ne)!0.5!(v21)$) [below right]{\small $\CS_{\tau_3,L}$};
\end{tikzpicture}
    \caption{Partition of $P_D$ into regions $\CS_{\delta,L}$. Cells $\delta$ in orange regions have $\wh{\lambda}_3(\pi_\R(\ceil{\delta}))=1$, cells $\delta$ in red regions have $\wh{\lambda}_3(\pi_\R(\ceil{\delta}))=0$, and cells $\delta$ in green regions have $\wh{\lambda}_3(\pi_\R(\ceil{\delta}))=-1$. See \Cref{example:H2_topological_weights} for a detailed explanation.}
    \label{fig:topological_weight_computation}
\end{figure}
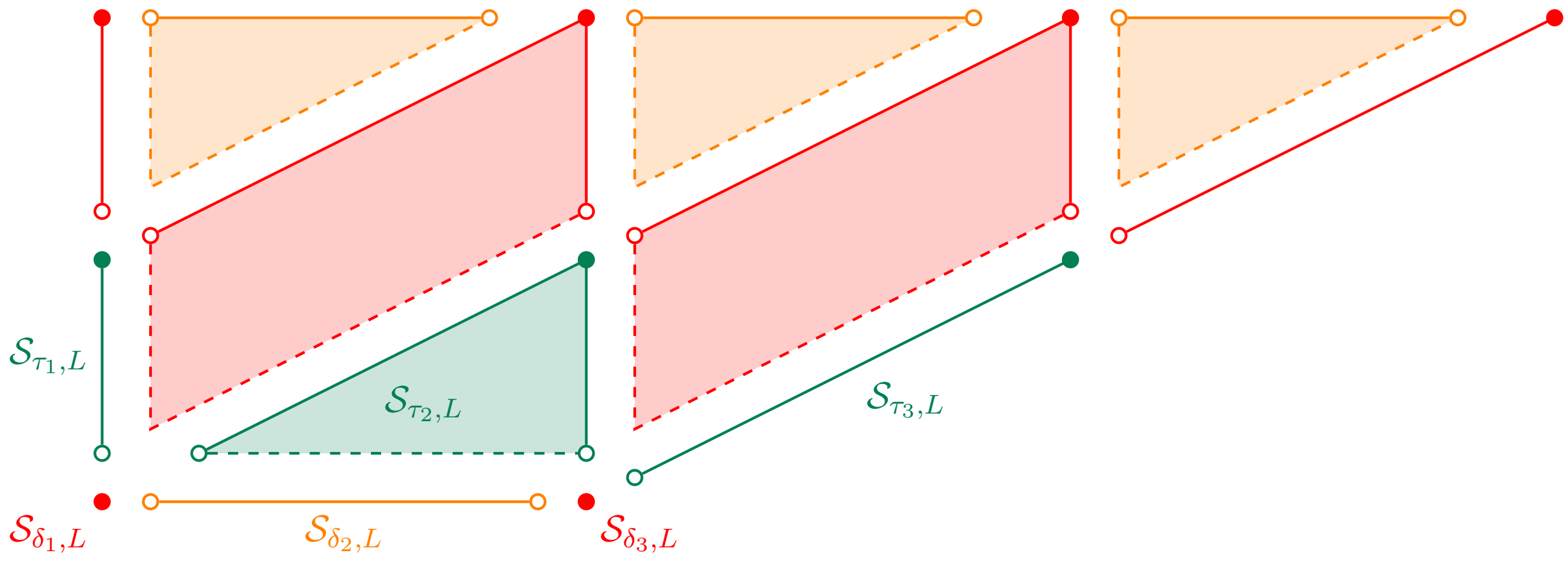

\begin{example}\label{example:H2_topological_weights}
We return to \Cref{example:H2_grade_restriction_rules_geometrically}. Fixing $F^{D}_{\bullet}(D)$ and $\bs{\lambda}$ as in the example, recall that we directly computed
\[
\mu(F_\bullet^D(D), \lambda_{3}, (X_3)_{\lambda_3}^{=0}) = \{0,1\}
\]
by pulling back $F^{D}_{\bullet}(D)$ to the fixed locus and inspecting the acyclicity of various fixed weight sub-complexes. Here, we demonstrate how to use \Cref{lemma:lambda_weights_of_cox_generators} to compute weights in $\mu(F_\bullet^D(D), \lambda_{3}, (X_3)_{\lambda_3}^{=0})$. First, recall that $(X_3)_{\lambda_3}^{=0} = V(x_1,x_2,x_3) \setminus {\{0\}}$, which can be written as a single orbit:
\[
(X_3)_{\lambda_3}^{=0} = T\cdot p_{L},
\]
where $p_{L}= (0,0,0,1)$.  \Cref{fig:topological_weight_computation} depicts the regions $\CS_{\delta,L}$ for $L=\{4\}$. There are only three regions $\CS_{\delta,L}$ with non-contractible boundary: $\CS_{\delta_1,L},\CS_{\delta_2,L},\CS_{\delta_3,L}$. Their associated $\lambda_{3}$-weights are $0, 1$, and $0$ respectively. Note that $\CS_{\delta,L}$'s with $\wh{\lambda}_3(\pi_\R(\ceil{\delta}))=-1$ (colored green) all have contractible boundaries. 
In summary, using \Cref{lemma:lambda_weights_of_cox_generators}, we recover that \[
\mu(F_\bullet^D(D), \lambda_{3}, (X_3)_{\lambda_3}^{=0}) = \{0,1\}.
\]
\end{example}

Throughout this section we have been translating the geometric condition of computing $\lambda_j$-weights of the Cox window generators $F_\bullet^D(D)$ into a topological condition. The following proposition is effectively a summary of this translation from geometry to topology. Recall $T_{\lambda_i}^{>0} = \sum_{j\in I_{\lambda_i}^{>0}}\wh{\lambda}_i(\beta_j)$.
\begin{proposition}\label{prop:cox_generator_grw_condition} Let $\Gamma$ be a chamber of $\Sigma_{GKZ}$ equipped with an HKKN stratification $\bs{\lambda} = \lambda_1,\hdots,\lambda_k$ and choose any shift vector $\bs{w} = (w_1,\hdots, w_k)\in \Z^k$. Let $\pi(D)\otimes 1\in  \Gamma$. Then $F^D_\bullet(D)\in W_{\bs{\lambda}}^{\bs{w}}$ if and only if for each $1 \le j \le k$, each $L\subset [n]$ such that $p_L\in (X_j)_{\lambda_j}^{=0}$, and each $\delta\in \cells(P_D)$ such that $\wh{\lambda}_j(\pi_\R(\ceil{\delta})) \notin [w_j,w_j+T_{\lambda_j}^{>0})$, $\del\CS_{\delta,L}$ is acyclic.
\end{proposition}
\begin{proof} This is an immediate consequence of the definition of $W_{\bs{\lambda}}^{\bs{w}}$ (see \Cref{def:grade_restriction_windows}) and \Cref{lemma:lambda_weights_of_cox_generators} which computes the $\lambda_j$-weights.
\end{proof}

As a consequence of this combinatorial characterization of the membership of the Cox generators $F_\bullet^D(D)$ in $W_{\bs{\lambda}}^{\bs{w}}$, we obtain the following.
\begin{theorem}\label{thm:GRW_is_Cox_iff} Let $\Sigma_{GKZ}$ be the secondary fan of a primitive geometric, projective GIT problem defined by $G\subset \Gm^n$. Let $\Gamma$ be a chamber of $\Sigma_{GKZ}$ equipped with an HKKN stratification $\bs{\lambda}$ of $V(B_\Gamma)$. Then, the associated grade restriction window $W_{\bs{\lambda}}$ is equal to the Cox window $W_\Gamma^{\Cox}$ if and only if for every $D$ with $\pi(D)\otimes 1\in \Gamma$, every $\lambda_j\notin |\Sigma_{GKZ}|^\vee$, every $L\subset [n]$ such that $p_L\in (X_j)_{\lambda_j}^{=0}$, and every $\delta\in \cells(P_D)$ such that $\wh{\lambda}_j(\pi_\R(\ceil{\delta})) < 0$, $\del \CS_{\delta,L}$ is acyclic.
\end{theorem}
\begin{proof} Recall that by \Cref{cor: window equality}, $W_\Gamma^{\Cox} = W_{\bs{\lambda}}$ if and only if for every $D$ such that $\pi(D)\otimes 1\in \Gamma$, $F_\bullet^D(D)\in W_{\bs{\lambda}}$. By \Cref{prop:cox_generator_grw_condition} we have that $F_\bullet^D(D)\in W_{\bs{\lambda}}$ if and only if for any $1 \le i \le k$, any $p_L\in (X_i)_{\lambda_i}^{=0}$, and any $\delta\in \cells(P_D)$ such that $\wh{\lambda}(\pi_\R(\ceil{\delta}))\notin [0, T_{\lambda_i}^{>0})$, that $\del \CS_{\delta,L}$ is acyclic. 

Fix some $i$ between $1$ and $k$. If $\lambda_i\in |\Sigma_{GKZ}|^\vee$, by \Cref{lemma:weights_in_Tminus_Tplus}\ref{lemma:weights_in_Tminus_Tplus_EXT}, $\wh{\lambda}_i(\pi_\R(\ceil{\delta}))\in [0,T_{\lambda_i}^{>0})$, so the condition to check is vacuous. If instead $\lambda_i\notin |\Sigma_{GKZ}|^\vee$, then $\wh{\lambda}_i(\pi_\R(\ceil{\delta}))\in (T_{\lambda_i}^{<0},T_{\lambda_i}^{>0})$ by \Cref{lemma:weights_in_Tminus_Tplus}\ref{lemma:weights_in_Tminus_Tplus_INT}. So, $\wh{\lambda}_i(\pi_\R(\ceil{\delta}))\notin [0,T_{\lambda_i}^{>0})$ if and only if $\wh{\lambda}_i(\pi_\R(\ceil{\delta})) < 0$.
\end{proof}
While showing that every Cox generator $F_\bullet^D(D)$ satisfies the conditions in \Cref{thm:GRW_is_Cox_iff} can be challenging, once we do so, we can obtain new equalities of grade restriction windows with the Cox window via the following.
\begin{corollary}\label{cor:two_stratifications_with_same_orbits_one_Cox} Let $\Gamma$ be as in \Cref{thm:GRW_is_Cox_iff} Let $\bs{\lambda}, \bs{\lambda}'$ be two HKKN stratifications of $V(B_\Gamma)$.  Suppose as unordered sets,
\[\bs{\lambda}'\setminus |\Sigma_{GKZ}|^\vee  \subset \bs{\lambda}\setminus |\Sigma_{GKZ}|^\vee.\]
and that for each $\lambda'_i\notin |\Sigma_{GKZ}|^\vee$ with corresponding $\lambda_j = \lambda'_i$ we have
\[ (X_i)_{\lambda'_i}^{=0}\subset (X_j)_{\lambda_j}^{=0}.\]
If  $W_{\bs{\lambda}} = W_\Gamma^\Cox$, then $W_{\bs{\lambda}'} = W_\Gamma^\Cox$ as well.  In particular, if $W_{\bs{\lambda}} = W_\Gamma^\Cox$ then $W_{\bs{\lambda}} = W_{\bs{\lambda'}}$  
\end{corollary}
\begin{proof} This is an immediate consequence of \Cref{prop:cox_generator_grw_condition} and \Cref{thm:GRW_is_Cox_iff}.
\end{proof}

\subsection{The nerve of the boundary subcomplex}\label{subsection:reduction_to_nerve}

In light of \Cref{thm:GRW_is_Cox_iff}, our goal is to understand when the topological spaces $\del \CS_{\delta, L}$ are acyclic. In this subsection we will introduce a cover of $\del \CS_{\delta,L}$ by closed convex sets whose nerve is homotopy equivalent to $\del \CS_{\delta,L}$. We will then provide explicit sufficient conditions for acyclicity of the nerve. For the remainder of the section we will fix some $L$ such that $p_L\in (X_j)_{\lambda_j}^{=0}$ and some $\delta\in\cells(P_D)$. 

In service of determining whether or not $\del\CS_{\delta,L}$ is contractible, we introduce the following definitions. 
\begin{definition}\label{def:nerve}
We define convex polytopes
\[
A_i \coloneqq P_D\cap \lp \{\ceil{\delta}_{i}-1\}\times \prod_{j\in [n]\setminus (L\cup \{i\})}[\ceil{\delta}_j-1,\ceil{\delta}_j]\times \prod_{t\in L}(-\infty,\infty) \rp
\]
so that
\[\del \CS_{\delta,L} = \bigcup_{i\in [n]\setminus L}A_i\]
by \Cref{lemma:coordinate_description_of_closure_boundary}. We define a family of closed, convex sets $\CC_{\delta,L} \coloneqq \{A_i\ | \ i\in [n]\setminus L\}$ and define its associated \emph{nerve} $N(\CC_{\delta,L})$ to be the abstract simplicial complex on $[n]\setminus L$ via
\[
N(\mathcal{C}_{\delta,L}) \coloneqq \left\{ J \subseteq [n]\setminus L \;:\; \bigcap_{j \in J} A_j \neq \varnothing\ \right\}.
\]
\end{definition}

We are in the exact situation which allows us to apply the nerve theorem to compare our topological space to its associated nerve.  Namely we have the following standard result.
\begin{proposition}\label{prop:nervecomputes}
$N(\mathcal{C}_{\delta,L})$ and $\del \CS_{\delta,L}$ are homotopy equivalent. 
\end{proposition}
\begin{proof}
This is a standard form of the Nerve Theorem, due to Leray \cite{Leray1945} and Borsuk \cite{Borsuk1948}; see also \cite{BjornerTopologicalMethods}. Our setting is precisely the regular-cell-complex setting in \cite[Theorem~10.6]{BjornerTopologicalMethods}: the $A_i$ are subcomplexes of a regular cell complex and every nonempty finite intersection of the $A_i$ is contractible.
\end{proof}
In light of the above proposition and \Cref{prop:cox_generator_grw_condition},
we must determine when the simplicial complex $N(\CS_{\delta,L})$ is acyclic. We do so by applying the following well-known observation.
\begin{lemma}\label{lemma:facet_common_element_implies_contractible}
Let $\Delta$ be a simplicial complex. If there exists a nonempty set $J'\in \Delta$ such that for any $J\in \Delta$, $J\cup J'\in \Delta$, then $\Delta$ is acyclic.
\end{lemma}
\begin{proof} This is a standard fact in algebraic topology (see 
\cite[Lemma 13.5]{Munkres1984}).
\end{proof} 
We make the following definition to capture the hypothesis of the previous lemma.
\begin{definition}\label{def:cone_set} $J' \in \Delta$ is a \emph{cone set} of a simplicial complex $\Delta$ if every vertex of $J'$ is a cone point. Equivalently, $J'$ is a cone set of $\Delta$ if for every $J\in \Delta$, $J\cup J'\in \Delta$.
\end{definition}
Recall our lift $D= \sum_{i\in [n]}a_ie_i\in \Z^n$ of $d\in \wh{G}\cap \Gamma$. In an effort to better understand $N(\CC_{\delta,L})$, we define \[\overline{K} \coloneqq \{i\in [n]\setminus L \ : \ \ceil{\delta}_i = -a_i\}.\]
As $\Kbar$ is defined in terms of the coordinates not in $L$, $\Kbar$ is well-defined for any choice of representative $\tau\in \cells(P_D)_{\lambda_j=w, [\delta]}$.

From the definitions of $\Kbar$ and $A_i$, it is clear that $i \in \overline{K} $ implies $ A_{i} = \emptyset$. In particular, $J\in N(\CC_{\delta,L})$ implies $J\subset [n]\setminus (L\cup \Kbar)$, so hereafter we will view $N(\CC_{\delta,L})$ as an abstract simplicial complex on $[n]\setminus (L\cup \Kbar)$.
Furthermore, when $i\notin \Kbar$, by unpacking $P_D = i_\R(M_\R)\cap \R^n_{\ge -D}$ in the definition of $A_i$ we get
\begin{equation}\label{eqn:Ai_coordinate_expression}
    A_i = B_i \cap i_\R(M_\R)
\end{equation}
where we define
\begin{equation}\label{eqn:Bi_coordinate_expression}
    B_i \coloneqq  \{\ceil{\delta}_{i}-1\}\times \prod_{j\in [n]\setminus ( \overline{K} \cup L \cup \{i\})}[\ceil{\delta}_j-1,\ceil{\delta}_j]\times \prod_{k \in \overline{K}}\{\ceil{\delta}_{k}\} \times \prod_{t\in L}[-a_{t},\infty).
\end{equation}
The nerve $N(\CC_{\delta,L})$ was defined in terms of whether or not a convex polyhedron in $\R^n$ was nonempty. Eventually we want to connect our analysis of these polyhedra to the secondary fan. The first step is to rephrase this nonemptyness as a convex geometric property in $\wh{G}_\R$ rather than $\R^n$. This is described in the following simple lemma.
\begin{lemma}\label{lemma:shortexactsequencenerve}
Given $J \subseteq [n] \setminus (L\cup \overline{K})$, $J \in N(\mathcal{C}_{\delta,L})$ if and only if
\[0\in \pi_{\R} \lp \bigcap_{j\in J} B_j\rp.\]
\end{lemma}
\begin{proof}
By definition, $J \in N(\CC_{\delta,L})$ if and only if $\bigcap_{j\in J}A_j\ne\emptyset$. However, by \eqref{eqn:Ai_coordinate_expression} it follows that
\[
\bigcap_{j \in J} A_j = \lp \bigcap_{j\in J} B_j\rp \cap i_\R(M_\R).
\]
By  the short exact sequence in \eqref{eqn:ses}, $i_\R(M_\R) = \ker(\pi_{\R})$, and the claim follows. 
\end{proof}
\begin{example}\label{example:H2_nerve_in_sceondary_fan}
Let us continue with our running example. Specifically, let $\CS_{\tau_2,L}$ be as in \Cref{fig:topological_weight_computation}. From \Cref{fig:topological_weight_computation}, it is evident that $\del \CS_{\tau_2,L}$ is a closed 1-cell, and thus contractible. By \Cref{prop:nervecomputes}, we could alternatively have determined the contractibility of $\del \CS_{\tau_2,L}$ via $N(\CC_{\tau_2,L})$. In \Cref{fig:topological_nerve_computation} we demonstrate the computation of $N(\CC_{\tau_3,L})$ using \Cref{lemma:shortexactsequencenerve}. In particular, for $J\subset [n]\setminus L = \{1,2,3\}$ we illustrate each region $B_J\coloneqq \bigcap_{j\in J}B_j$ and record when $0\in B_J$. In this way we determine
\[N(\CC_{\tau_3,L}) = \bigg\{\emptyset, \{1\}, \{2\}, \{3\}, \{1,2\}, \{2,3\}\bigg\}.\]
From this description it is evident that $J' = \{2\}$ is a cone set for $N(\CC_{\tau_3,L})$, thus $N(\CC_{\tau_3,L})\simeq \ast$ via \Cref{lemma:facet_common_element_implies_contractible}. 
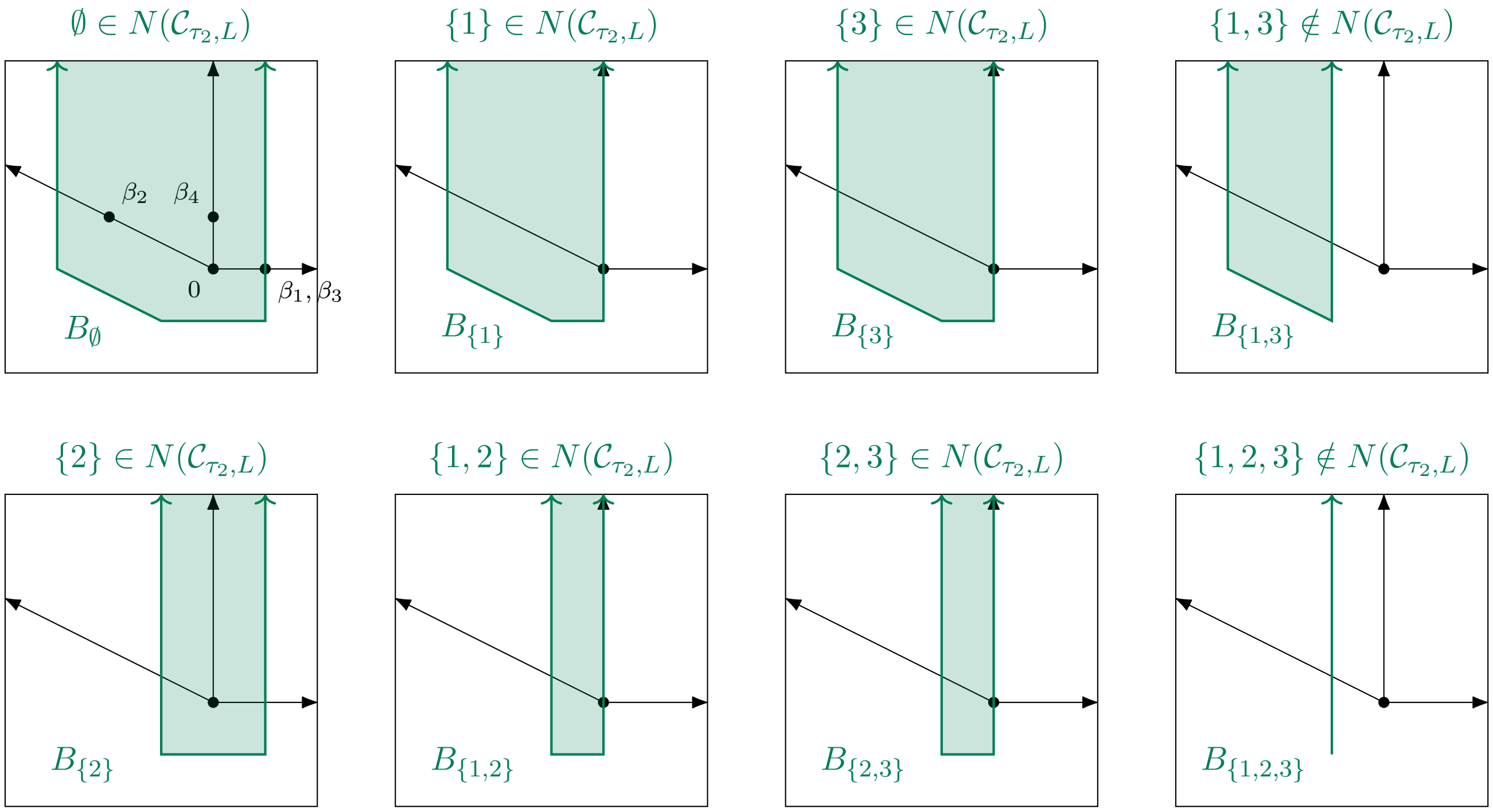
\begin{figure}
\begin{tikzpicture}
\def\fanscale{0.6}
\def\xmin{-4}
\def\xmax{2}
\def\ymin{-2}
\def\ymax{4}
\def\polycolor{darkGreen}
\def\polyop{0.2}
\def\thickness{thick};
\def\xsep{4.5}
\def\ysep{5}

\coordinate (shift00) at (0,0);
\coordinate (shift01) at (1*\xsep,0);
\coordinate (shift02) at (2*\xsep,0*\ysep);
\coordinate (shift03) at (3*\xsep,0*\ysep);
\coordinate (shift10) at (0*\xsep,-1*\ysep);
\coordinate (shift11) at (1*\xsep,-1*\ysep);
\coordinate (shift12) at (2*\xsep,-1*\ysep);
\coordinate (shift13) at (3*\xsep,-1*\ysep);

\node[draw=none, fill=none, \polycolor] at ($(shift00)+(-1*\fanscale, \ymax*\fanscale+0.4)$) {\small $\emptyset\in N(\CC_{\tau_2,L})$};
\node[draw=none,fill=none,\polycolor] at ($(shift00)+(-2.5*\fanscale,-1.2*\fanscale)$) {\small $B_{\emptyset}$};

\node[draw=none, fill=none, \polycolor] at ($(shift01)+(-1*\fanscale, \ymax*\fanscale+0.4)$) {\small $\{1\}\in N(\CC_{\tau_2,L})$};
\node[draw=none,fill=none,\polycolor] at ($(shift01)+(-2.5*\fanscale,-1.2*\fanscale)$) {\small $B_{\{1\}}$};

\node[draw=none, fill=none, \polycolor] at ($(shift02)+(-1*\fanscale, \ymax*\fanscale+0.4)$) {\small $\{3\}\in N(\CC_{\tau_2,L})$};
\node[draw=none,fill=none,\polycolor] at ($(shift02)+(-2.5*\fanscale,-1.2*\fanscale)$) {\small $B_{\{3\}}$};

\node[draw=none, fill=none, \polycolor] at ($(shift03)+(-1*\fanscale, \ymax*\fanscale+0.4)$) {\small $\{1,3\}\notin N(\CC_{\tau_2,L})$};
\node[draw=none,fill=none,\polycolor] at ($(shift03)+(-2.5*\fanscale,-1.2*\fanscale)$) {\small $B_{\{1,3\}}$};

\node[draw=none, fill=none, \polycolor] at ($(shift10)+(-1*\fanscale, \ymax*\fanscale+0.4)$) {\small $\{2\}\in N(\CC_{\tau_2,L})$};
\node[draw=none,fill=none,\polycolor] at ($(shift10)+(-2.5*\fanscale,-1.2*\fanscale)$) {\small $B_{\{2\}}$};

\node[draw=none, fill=none, \polycolor] at ($(shift11)+(-1*\fanscale, \ymax*\fanscale+0.4)$) {\small $\{1,2\}\in N(\CC_{\tau_2,L})$};
\node[draw=none,fill=none,\polycolor] at ($(shift11)+(-2.5*\fanscale,-1.2*\fanscale)$) {\small $B_{\{1,2\}}$};

\node[draw=none, fill=none, \polycolor] at ($(shift12)+(-1*\fanscale, \ymax*\fanscale+0.4)$) {\small $\{2,3\}\in N(\CC_{\tau_2,L})$};
\node[draw=none,fill=none,\polycolor] at ($(shift12)+(-2.5*\fanscale,-1.2*\fanscale)$) {\small $B_{\{2,3\}}$};

\node[draw=none, fill=none, \polycolor] at ($(shift13)+(-1*\fanscale, \ymax*\fanscale+0.4)$) {\small $\{1,2,3\}\notin N(\CC_{\tau_2,L})$};
\node[draw=none,fill=none,\polycolor] at ($(shift13)+(-2.5*\fanscale,-1.2*\fanscale)$) {\small $B_{\{1,2,3\}}$};

\coordinate (origin) at (0*\fanscale,0*\fanscale);
\coordinate (v0) at (-3*\fanscale,0*\fanscale);
\coordinate (v1) at (-1*\fanscale,-1*\fanscale);
\coordinate (v2) at (0*\fanscale,-1*\fanscale);
\coordinate (v3) at (1*\fanscale,-1*\fanscale);
\coordinate (t0) at (-3*\fanscale,\ymax*\fanscale);
\coordinate (t1) at (-1*\fanscale,\ymax*\fanscale);
\coordinate (t2) at (0*\fanscale,\ymax*\fanscale);
\coordinate (t3) at (1*\fanscale,\ymax*\fanscale);
\coordinate (beta13) at (1*\fanscale,0*\fanscale);
\coordinate (beta2) at (-2*\fanscale,1*\fanscale);
\coordinate (beta4) at (0*\fanscale,1*\fanscale);
\coordinate (ray13) at (\xmax*\fanscale,0*\fanscale);
\coordinate (ray2) at (\xmin*\fanscale,-0.5*\xmin*\fanscale);
\coordinate (ray4) at (0*\fanscale,\ymax*\fanscale);
\coordinate (boxll) at (\xmin*\fanscale,\ymin*\fanscale);
\coordinate (boxul) at (\xmin*\fanscale,\ymax*\fanscale);
\coordinate (boxur) at (\xmax*\fanscale,\ymax*\fanscale);
\coordinate (boxlr) at (\xmax*\fanscale,\ymin*\fanscale);

%%%%% FAN 00
\coordinate (shift) at (shift00);
\coordinate (sorigin) at ($(origin)+(shift)$);
\coordinate (sv0) at ($(v0)+(shift)$);
\coordinate (sv1) at ($(v1)+(shift)$);
\coordinate (sv2) at ($(v2)+(shift)$);
\coordinate (sv3) at ($(v3)+(shift)$);
\coordinate (st0) at ($(t0)+(shift)$);
\coordinate (st1) at ($(t1)+(shift)$);
\coordinate (st2) at ($(t2)+(shift)$);
\coordinate (st3) at ($(t3)+(shift)$);
\coordinate (sbeta13) at ($(beta13)+(shift)$);
\coordinate (sbeta2) at ($(beta2)+(shift)$);
\coordinate (sbeta4) at ($(beta4)+(shift)$);
\coordinate (sray13) at ($(ray13)+(shift)$);
\coordinate (sray2) at ($(ray2)+(shift)$);
\coordinate (sray4) at ($(ray4)+(shift)$);
\coordinate (sboxll) at ($(boxll)+(shift)$);
\coordinate (sboxul) at ($(boxul)+(shift)$);
\coordinate (sboxur) at ($(boxur)+(shift)$);
\coordinate (sboxlr) at ($(boxlr)+(shift)$);

\draw[draw=black,fill=none,thin] (sboxll) -- (sboxul) -- (sboxur) -- (sboxlr) -- cycle;

\draw[-arr] (sorigin) -- (ray13);
\draw[-arr] (sorigin) -- (ray2);
\draw[-arr] (sorigin) -- (ray4);

\node[circle, inner sep = 1pt, thick, fill = black, draw=black] at (sorigin) {};
\node[circle, inner sep = 1pt, thick, fill = black, draw=black] at (sbeta13) {};
\node[circle, inner sep = 1pt, thick, fill = black, draw=black] at (sbeta2) {};
\node[circle, inner sep = 1pt, thick, fill = black, draw=black] at (sbeta4) {};
\node[draw=none,fill=none] at (sorigin) [below left] {\tiny $0$};
\node[draw=none,fill=none] at (sbeta13) [below right] {\tiny $\beta_1,\beta_3$};
\node[draw=none,fill=none] at (sbeta2) [above right] {\tiny $\beta_2$};
\node[draw=none,fill=none] at (sbeta4) [above left] {\tiny $\beta_4$};

\draw[draw=none,fill = \polycolor, opacity=\polyop] (t0) -- (v0) -- (v1) -- (v3) -- (t3);
\draw[<->,draw = \polycolor, \thickness] (t0) -- (v0) -- (v1) -- (v3) -- (t3);

%%%%% FAN 01
\coordinate (shift) at (shift01);
\coordinate (sorigin) at ($(origin)+(shift)$);
\coordinate (sv0) at ($(v0)+(shift)$);
\coordinate (sv1) at ($(v1)+(shift)$);
\coordinate (sv2) at ($(v2)+(shift)$);
\coordinate (sv3) at ($(v3)+(shift)$);
\coordinate (st0) at ($(t0)+(shift)$);
\coordinate (st1) at ($(t1)+(shift)$);
\coordinate (st2) at ($(t2)+(shift)$);
\coordinate (st3) at ($(t3)+(shift)$);
\coordinate (sbeta13) at ($(beta13)+(shift)$);
\coordinate (sbeta2) at ($(beta2)+(shift)$);
\coordinate (sbeta4) at ($(beta4)+(shift)$);
\coordinate (sray13) at ($(ray13)+(shift)$);
\coordinate (sray2) at ($(ray2)+(shift)$);
\coordinate (sray4) at ($(ray4)+(shift)$);
\coordinate (sboxll) at ($(boxll)+(shift)$);
\coordinate (sboxul) at ($(boxul)+(shift)$);
\coordinate (sboxur) at ($(boxur)+(shift)$);
\coordinate (sboxlr) at ($(boxlr)+(shift)$);

\draw[draw=black,fill=none,thin] (sboxll) -- (sboxul) -- (sboxur) -- (sboxlr) -- cycle;

\draw[-arr] (sorigin) -- (sray13);
\draw[-arr] (sorigin) -- (sray2);
\draw[-arr] (sorigin) -- (sray4);

\node[circle, inner sep = 1pt, thick, fill = black, draw=black] at (sorigin) {};

\draw[draw=none,fill = \polycolor, opacity=\polyop] (st0) -- (sv0) -- (sv1) -- (sv2) -- (st2);
\draw[<->,draw = \polycolor, \thickness] (st0) -- (sv0) -- (sv1) -- (sv2) -- (st2);

%%%%% FAN 02
\coordinate (shift) at (shift02);
\coordinate (sorigin) at ($(origin)+(shift)$);
\coordinate (sv0) at ($(v0)+(shift)$);
\coordinate (sv1) at ($(v1)+(shift)$);
\coordinate (sv2) at ($(v2)+(shift)$);
\coordinate (sv3) at ($(v3)+(shift)$);
\coordinate (st0) at ($(t0)+(shift)$);
\coordinate (st1) at ($(t1)+(shift)$);
\coordinate (st2) at ($(t2)+(shift)$);
\coordinate (st3) at ($(t3)+(shift)$);
\coordinate (sbeta13) at ($(beta13)+(shift)$);
\coordinate (sbeta2) at ($(beta2)+(shift)$);
\coordinate (sbeta4) at ($(beta4)+(shift)$);
\coordinate (sray13) at ($(ray13)+(shift)$);
\coordinate (sray2) at ($(ray2)+(shift)$);
\coordinate (sray4) at ($(ray4)+(shift)$);
\coordinate (sboxll) at ($(boxll)+(shift)$);
\coordinate (sboxul) at ($(boxul)+(shift)$);
\coordinate (sboxur) at ($(boxur)+(shift)$);
\coordinate (sboxlr) at ($(boxlr)+(shift)$);

\draw[draw=black,fill=none,thin] (sboxll) -- (sboxul) -- (sboxur) -- (sboxlr) -- cycle;

\draw[-arr] (sorigin) -- (sray13);
\draw[-arr] (sorigin) -- (sray2);
\draw[-arr] (sorigin) -- (sray4);

\node[circle, inner sep = 1pt, thick, fill = black, draw=black] at (sorigin) {};

\draw[draw=none,fill = \polycolor, opacity=\polyop] (st0) -- (sv0) -- (sv1) -- (sv2) -- (st2);
\draw[<->,draw = \polycolor, \thickness] (st0) -- (sv0) -- (sv1) -- (sv2) -- (st2);

%%%%% FAN 03
\coordinate (shift) at (shift03);
\coordinate (sorigin) at ($(origin)+(shift)$);
\coordinate (sv0) at ($(v0)+(shift)$);
\coordinate (sv1) at ($(v1)+(shift)$);
\coordinate (sv2) at ($(v2)+(shift)$);
\coordinate (sv3) at ($(v3)+(shift)$);
\coordinate (st0) at ($(t0)+(shift)$);
\coordinate (st1) at ($(t1)+(shift)$);
\coordinate (st2) at ($(t2)+(shift)$);
\coordinate (st3) at ($(t3)+(shift)$);
\coordinate (sbeta13) at ($(beta13)+(shift)$);
\coordinate (sbeta2) at ($(beta2)+(shift)$);
\coordinate (sbeta4) at ($(beta4)+(shift)$);
\coordinate (sray13) at ($(ray13)+(shift)$);
\coordinate (sray2) at ($(ray2)+(shift)$);
\coordinate (sray4) at ($(ray4)+(shift)$);
\coordinate (sboxll) at ($(boxll)+(shift)$);
\coordinate (sboxul) at ($(boxul)+(shift)$);
\coordinate (sboxur) at ($(boxur)+(shift)$);
\coordinate (sboxlr) at ($(boxlr)+(shift)$);

\draw[draw=black,fill=none,thin] (sboxll) -- (sboxul) -- (sboxur) -- (sboxlr) -- cycle;

\draw[-arr] (sorigin) -- (sray13);
\draw[-arr] (sorigin) -- (sray2);
\draw[-arr] (sorigin) -- (sray4);

\node[circle, inner sep = 1pt, thick, fill = black, draw=black] at (sorigin) {};

\draw[draw=none,fill = \polycolor, opacity=\polyop] (st0) -- (sv0) -- (sv1) -- (st1);
\draw[<->,draw = \polycolor, \thickness] (st0) -- (sv0) -- (sv1) -- (st1);

%%%%% FAN 10
\coordinate (shift) at (shift10);
\coordinate (sorigin) at ($(origin)+(shift)$);
\coordinate (sv0) at ($(v0)+(shift)$);
\coordinate (sv1) at ($(v1)+(shift)$);
\coordinate (sv2) at ($(v2)+(shift)$);
\coordinate (sv3) at ($(v3)+(shift)$);
\coordinate (st0) at ($(t0)+(shift)$);
\coordinate (st1) at ($(t1)+(shift)$);
\coordinate (st2) at ($(t2)+(shift)$);
\coordinate (st3) at ($(t3)+(shift)$);
\coordinate (sbeta13) at ($(beta13)+(shift)$);
\coordinate (sbeta2) at ($(beta2)+(shift)$);
\coordinate (sbeta4) at ($(beta4)+(shift)$);
\coordinate (sray13) at ($(ray13)+(shift)$);
\coordinate (sray2) at ($(ray2)+(shift)$);
\coordinate (sray4) at ($(ray4)+(shift)$);
\coordinate (sboxll) at ($(boxll)+(shift)$);
\coordinate (sboxul) at ($(boxul)+(shift)$);
\coordinate (sboxur) at ($(boxur)+(shift)$);
\coordinate (sboxlr) at ($(boxlr)+(shift)$);

\draw[draw=black,fill=none,thin] (sboxll) -- (sboxul) -- (sboxur) -- (sboxlr) -- cycle;

\draw[-arr] (sorigin) -- (sray13);
\draw[-arr] (sorigin) -- (sray2);
\draw[-arr] (sorigin) -- (sray4);

\node[circle, inner sep = 1pt, thick, fill = black, draw=black] at (sorigin) {};

\draw[draw=none,fill = \polycolor, opacity=\polyop] (st1) -- (sv1) -- (sv3) -- (st3);
\draw[<->,draw = \polycolor, \thickness] (st1) -- (sv1) -- (sv3) -- (st3);

%%%%% FAN 11
\coordinate (shift) at (shift11);
\coordinate (sorigin) at ($(origin)+(shift)$);
\coordinate (sv0) at ($(v0)+(shift)$);
\coordinate (sv1) at ($(v1)+(shift)$);
\coordinate (sv2) at ($(v2)+(shift)$);
\coordinate (sv3) at ($(v3)+(shift)$);
\coordinate (st0) at ($(t0)+(shift)$);
\coordinate (st1) at ($(t1)+(shift)$);
\coordinate (st2) at ($(t2)+(shift)$);
\coordinate (st3) at ($(t3)+(shift)$);
\coordinate (sbeta13) at ($(beta13)+(shift)$);
\coordinate (sbeta2) at ($(beta2)+(shift)$);
\coordinate (sbeta4) at ($(beta4)+(shift)$);
\coordinate (sray13) at ($(ray13)+(shift)$);
\coordinate (sray2) at ($(ray2)+(shift)$);
\coordinate (sray4) at ($(ray4)+(shift)$);
\coordinate (sboxll) at ($(boxll)+(shift)$);
\coordinate (sboxul) at ($(boxul)+(shift)$);
\coordinate (sboxur) at ($(boxur)+(shift)$);
\coordinate (sboxlr) at ($(boxlr)+(shift)$);

\draw[draw=black,fill=none,thin] (sboxll) -- (sboxul) -- (sboxur) -- (sboxlr) -- cycle;

\draw[-arr] (sorigin) -- (sray13);
\draw[-arr] (sorigin) -- (sray2);
\draw[-arr] (sorigin) -- (sray4);

\node[circle, inner sep = 1pt, thick, fill = black, draw=black] at (sorigin) {};

\draw[draw=none,fill = \polycolor, opacity=\polyop] (st1) -- (sv1) -- (sv2) -- (st2);
\draw[<->,draw = \polycolor, \thickness] (st1) -- (sv1) -- (sv2) -- (st2);

%%%%% FAN 12
\coordinate (shift) at (shift12);
\coordinate (sorigin) at ($(origin)+(shift)$);
\coordinate (sv0) at ($(v0)+(shift)$);
\coordinate (sv1) at ($(v1)+(shift)$);
\coordinate (sv2) at ($(v2)+(shift)$);
\coordinate (sv3) at ($(v3)+(shift)$);
\coordinate (st0) at ($(t0)+(shift)$);
\coordinate (st1) at ($(t1)+(shift)$);
\coordinate (st2) at ($(t2)+(shift)$);
\coordinate (st3) at ($(t3)+(shift)$);
\coordinate (sbeta13) at ($(beta13)+(shift)$);
\coordinate (sbeta2) at ($(beta2)+(shift)$);
\coordinate (sbeta4) at ($(beta4)+(shift)$);
\coordinate (sray13) at ($(ray13)+(shift)$);
\coordinate (sray2) at ($(ray2)+(shift)$);
\coordinate (sray4) at ($(ray4)+(shift)$);
\coordinate (sboxll) at ($(boxll)+(shift)$);
\coordinate (sboxul) at ($(boxul)+(shift)$);
\coordinate (sboxur) at ($(boxur)+(shift)$);
\coordinate (sboxlr) at ($(boxlr)+(shift)$);

\draw[draw=black,fill=none,thin] (sboxll) -- (sboxul) -- (sboxur) -- (sboxlr) -- cycle;

\draw[-arr] (sorigin) -- (sray13);
\draw[-arr] (sorigin) -- (sray2);
\draw[-arr] (sorigin) -- (sray4);

\node[circle, inner sep = 1pt, thick, fill = black, draw=black] at (sorigin) {};

\draw[draw=none,fill = \polycolor, opacity=\polyop] (st1) -- (sv1) -- (sv2) -- (st2);
\draw[<->,draw = \polycolor, \thickness] (st1) -- (sv1) -- (sv2) -- (st2);

%%%%% FAN 13
\coordinate (shift) at (shift13);
\coordinate (sorigin) at ($(origin)+(shift)$);
\coordinate (sv0) at ($(v0)+(shift)$);
\coordinate (sv1) at ($(v1)+(shift)$);
\coordinate (sv2) at ($(v2)+(shift)$);
\coordinate (sv3) at ($(v3)+(shift)$);
\coordinate (st0) at ($(t0)+(shift)$);
\coordinate (st1) at ($(t1)+(shift)$);
\coordinate (st2) at ($(t2)+(shift)$);
\coordinate (st3) at ($(t3)+(shift)$);
\coordinate (sbeta13) at ($(beta13)+(shift)$);
\coordinate (sbeta2) at ($(beta2)+(shift)$);
\coordinate (sbeta4) at ($(beta4)+(shift)$);
\coordinate (sray13) at ($(ray13)+(shift)$);
\coordinate (sray2) at ($(ray2)+(shift)$);
\coordinate (sray4) at ($(ray4)+(shift)$);
\coordinate (sboxll) at ($(boxll)+(shift)$);
\coordinate (sboxul) at ($(boxul)+(shift)$);
\coordinate (sboxur) at ($(boxur)+(shift)$);
\coordinate (sboxlr) at ($(boxlr)+(shift)$);

\draw[draw=black,fill=none,thin] (sboxll) -- (sboxul) -- (sboxur) -- (sboxlr) -- cycle;

\draw[-arr] (sorigin) -- (sray13);
\draw[-arr] (sorigin) -- (sray2);
\draw[-arr] (sorigin) -- (sray4);

\node[circle, inner sep = 1pt, thick, fill = black, draw=black] at (sorigin) {};

\draw[->,draw = \polycolor, \thickness] (sv1) -- (st1);
\end{tikzpicture}
    \caption{Computing $N(\CC_{\tau_2,L})$ using \Cref{lemma:shortexactsequencenerve} where $L=\{4\}$ and $\CS_{\tau_2,L}$ is the same region as in \Cref{fig:topological_weight_computation}. See \Cref{example:H2_nerve_in_sceondary_fan} for a complete description.}
    \label{fig:topological_nerve_computation}
\end{figure}
\end{example}
We think about \Cref{lemma:shortexactsequencenerve} as stating that the membership of $J$ in $N(\CC_{\delta,L})$ can be determined by checking whether or not the polyhedron $\pi_\R\lp \bigcap_{j\in J}B_j\rp$ in $\wh{G}_\R$ contains $0$. When this polyhedron is in $\R$, it is an interval and thus determining whether or not it contains 0 is a matter of determining the bounds of the interval. To exploit this idea we will consider the image of this polyhedron under every 1-PS $\xi$. First we make the following definition.
\begin{definition}\label{def:xi_complex}
Take any $\wh{\xi}\in \wh{G}_\R^\ast$. We define $N(\CC_{\delta,L})_\xi$ to be the abstract simplicial complex consisting of $J\subset [n]\setminus (\Kbar\cup L)$ such that 
\[0\in (\wh{\xi}\circ \pi_{\R}) \lp \bigcap_{j\in J}B_j\rp\]
\end{definition}
\begin{lemma}\label{lemma:nerve_is_intersection_of_xi_posets}
\[N(\mathcal{C}_{\delta,L})= \bigcap_{\xi\in \wh{G}_\R^\ast} N(\CC_{\delta,L})_\xi\] 
\end{lemma}
\begin{proof}
Recall that every polyhedron is always equal to the intersection over all half spaces which contain it. Applied to $R\coloneqq \pi_\R\lp\bigcap_{j\in J}B_j\rp$, that is
\begin{equation}\label{eqn:intersection_over_strips}
    R = \bigcap_{\substack{(\xi, a)\in \wh{G}_\R^\ast \times \R, \\ R\subset H_{\xi}^{\ge -a}}} H_{\xi}^{\ge -a} = \bigcap_{\xi\in \wh{G}_\R^\ast} \wh{\xi}^{-1}\lp \wh{\xi}(R)\rp
\end{equation}
The second equality above comes from unraveling definitions. It follows that
\begin{align*}
 J \in N(\CC_{\delta,L}) &\iff 0 \in R & \text{ by \Cref{lemma:shortexactsequencenerve}}\\
& \iff \forall  \xi\in \wh{G}_\R^\ast,\ 0 \in \wh{\xi}(R) & \text{by \Cref{eqn:intersection_over_strips}} \\
& \iff J \in \bigcap_{\xi\in \wh{G}_\R^\ast} N(\CC_{\delta,L})_\xi. & \text{by \Cref{def:xi_complex}}
\end{align*}
\end{proof}

As a corollary to \Cref{lemma:nerve_is_intersection_of_xi_posets}, we obtain the following statement.
\begin{corollary}\label{cor:cone_set_for_xi_implies_nerve_contractible} If $J'\subset [n]\setminus (\Kbar\cup L)$ is a nonempty cone set for $N(\CC_{\delta,L})_\xi$ for every $\xi\in \wh{G}_\R^\ast$, then $J'$ is a cone set for $N(\CC_{\delta,L})$. In particular, $N(\CC_{\delta,L})$ is contractible by \Cref{lemma:facet_common_element_implies_contractible}.
\end{corollary}
In light of \Cref{lemma:nerve_is_intersection_of_xi_posets}, we must understand when exactly $J\in N(\CC_{\delta, L})_\xi$ for each $\xi\in \Lambda$. It turns out that there is an explicit combinatorial description for each $N(\CC_{\delta, L})_\xi$ which we will outline in the remainder of this subsection.
\begin{definition}\label{def:distinguished_character} We define the \textit{distinguished character} $\chi(D,\delta,L)\in \wh{G}_\R$ as 
\[\chi(D,\delta,L)\coloneqq \sum_{j\in [n]\setminus (\Kbar\cup L)} (1-\ceil{\delta}_j)\beta_j + \sum_{t\in \Kbar\cup L} a_t \beta_t.\]
When the context is clear, we will simply write $\chi \coloneqq \chi(D,\delta,L)$.
\end{definition}
\begin{lemma}\label{lemma:dist_character_and_interval} 
For any $J\subset [n]\setminus (L\cup \Kbar)$, $J\in N(\CC_{\delta,L})_\xi$ if and only if
\[
\wh{\xi}(\chi)\in \sum_{j\in [n]\setminus (\Kbar \cup L \cup J)}\Conv(0,\wh{\xi}(\beta_j))+\sum_{t\in L}\Cone(\wh{\xi}(\beta_t))\eqcolon \CI_{\xi,J}.
\]
Moreover, $\CI_{\xi,J}$ is a closed interval in $\R$. If $L\cap I_\xi^{>0}\ne \emptyset$, $\CI_{\xi,J}$ is unbounded above, and otherwise
the upper bound is $\sum_{i\in I_\xi^{>0}\setminus (J\cup \Kbar)}\wh{\xi}(\beta_i)$. Similarly, if $L\cap I_\xi^{<0}\ne \emptyset$, $\CI_{\xi,J}$ is unbounded below, and otherwise the lower bound is $\sum_{i\in I_\xi^{<0}\setminus (J\cup \Kbar)}\wh{\xi}(\beta_i)$.
\end{lemma}
\begin{proof} First observe that
\[\bigcap_{i\in J}B_i = \prod_{i\in J}\{\ceil{\delta}_i-1\}\times \prod_{j\in [n]\setminus (L\cup \Kbar)}[\ceil{\delta}_j-1,\ceil{\delta}_j]\times \prod_{k\in \Kbar} \{\ceil{\delta}_k\}\times \prod_{t\in L}[-a_t,\infty).\]
So by taking $\pi_\R$ and rearranging, one can see that
\[\pi_\R\lp \bigcap_{i\in J}B_i\rp  = -\chi+\sum_{j\in [n]\setminus (J\cup \Kbar \cup L)}\Conv(0,\beta_j)+\sum_{t\in L}\Cone(\beta_t).\]
Then, as $J\in N(\CC_{\delta,L})_\xi$ if and only if $0\in (\wh{\xi}\circ \pi_\R)\lp \cap_{i\in J}B_i\rp$, applying $\wh{\xi}$ to the above expression gives the desired result.
\end{proof}
The value in this explicit description of $N(\CC_{\delta,L})_\xi$ in terms of the containment of $\wh{\xi}(\chi)$ in $\CI_{\xi,J}$ is that the upper and lower bounds of $\CI_{\xi,J}$ (if they exist at all) are always non-negative and non-positive respectively. As we will see below, if $\wh{\xi}(\chi) \ge 0$, the lower bound is always trivially satisfied, so any $J'\subset I_\xi^{\le 0}$ is a cone set for $N(\CC_{\delta,L})_\xi$. We summarize a useful few of these elementary observations on cone sets of $N(\CC_{\delta,L})_\xi$ in the below corollary.

\begin{corollary}\label{cor:integer_poset_cor} 

Let $\chi$ be as above and $\wh{\xi}\in \wh{G}_\R^\ast$. Let $J'\subset [n]\setminus (\Kbar\cup L)$. If any of the following four conditions are met, then $J'$ is a cone set of $N(\CC_{\delta,L})_\xi$:
\begin{enumerate}[(a),itemsep=0.5cm]
    \item\label{integer_poset_cor:NEW_I_plus_minus}$I_\xi^{<0}\cap L \ne \emptyset$ and $I_\xi^{>0}\cap L\ne \emptyset$.
    \item\label{integer_poset_cor:NEW_I_plus_J_plus} $I_\xi^{>0}\cap L\ne \emptyset$ and $J'\subset I_\xi^{\ge 0}$.
    \item\label{integer_poset_cor:NEW_I_plus_chi_plus} $I_\xi^{>0} \cap L \ne \emptyset$ and $\wh{\xi}(\chi) \ge 0$.
    \item\label{integer_poset_cor:NEW_J_plus_chi_minus} $J' \subset I_\xi^{\ge 0}$ and $\wh{\xi}(\chi) \le 0$.
\end{enumerate}
\end{corollary}
\begin{proof} Suppose $J\in N(\CC_{\delta,L})_\xi$. We will show for each condition that $J\cup J'\in N(\CC_{\delta,L})_\xi$.

To show \ref{integer_poset_cor:NEW_I_plus_minus}, suppose $I_\xi^{>0}\cap L\ne \emptyset$ and $I_\xi^{<0}\cap L \ne \emptyset$ and let $J'$ be any subset of $[n]\setminus (\Kbar \cup L)$. Then by \Cref{lemma:dist_character_and_interval}, $\CI_{\xi,J\cup J'}$ is unbounded above and below, so $\wh{\xi}(\chi)\in \CI_{\xi, J\cup J'}$. Thus $J\cup J'\in N(\CC_{\delta,L})_\xi$ by \Cref{lemma:dist_character_and_interval}.

To show \ref{integer_poset_cor:NEW_I_plus_J_plus}, suppose $I_\xi^{>0}\cap L\ne\emptyset$ and $J'\subset I_\xi^{\ge 0}$. Because $I_\xi^{>0}\cap L\ne\emptyset$, $\CI_{\xi, J}$ and $\CI_{\xi, J\cup J'}$ have no upper bound. Because $J'\subset I_\xi^{\ge 0}$, the lower bound of $\CI_{\xi,J\cup J'}$ is equal to the lower bound of $\CI_{\xi, J}$. Thus, $\CI_{\xi, J} = \CI_{\xi, J\cup J'}$. It follows from \Cref{lemma:dist_character_and_interval} that $J\cup J'\in N(\CC_{\delta,L})_\xi$.

To show \ref{integer_poset_cor:NEW_I_plus_chi_plus}, suppose $I_\xi^{>0}\cap L\ne \emptyset$ and $\wh{\xi}(\chi) \ge 0$ and let $J'\subset [n]\setminus (\Kbar\cup L)$ be arbitrary. Then by \Cref{lemma:dist_character_and_interval}, $\CI_{\xi, J\cup J'}$ has no upper bound. As $\wh{\xi}(\chi) \ge 0$ and the lower bound of $\CI_{\xi ,J\cup J'}$ is non-positive, it follows that $\wh{\xi}(\chi)\in \CI_{\xi, J\cup J'}$. Thus $J\cup J'\in N(\CC_{\delta,L})_\xi$ by \Cref{lemma:dist_character_and_interval}. 

To show \ref{integer_poset_cor:NEW_J_plus_chi_minus}, suppose $J'\subset I_\xi^{\ge0}$ and $\wh{\xi}(\chi) \le 0$. Because $J'\subset I_\xi^{\ge 0}$, the lower bound of $\CI_{\xi,J\cup J'}$ is equal to the lower bound of $\CI_{\xi, J}$. As $\wh{\xi}(\chi) \le 0$, the upper bound of $\CI_{\xi, J\cup J'}$ is trivially satisfied. So $\wh{\xi}(\chi)\in \CI_{\xi, J}$ implies $\wh{\xi}(\chi)\in \CI_{\xi, J\cup J'}$. Thus $J\cup J'\in N(\CC_{\delta,L})_\xi$ by \Cref{lemma:dist_character_and_interval}.
\end{proof}

From \Cref{cor:integer_poset_cor}, it is clear that the sign of $\wh{\xi}(\chi)$ is important for determining when $J'$ is a cone set. The sign of $\wh{\xi}(\chi)$ is fully determined by which side of the hyperplane $H_{\xi}^{=0}$ contains $\chi$. As such, understanding the position of $\chi$ in $\wh{G}_\R$ will be useful in \S\ref{section:main_thm}.
\begin{lemma}\label{lemma:chi_eff} The distinguished character $\chi\coloneqq \chi(D,\delta,L)$ is contained in $|\Sigma_{GKZ}|$.
\end{lemma}
\begin{proof}
Because $\delta$ is a nonempty cell of $P_D$, there is a point $x\in \delta$. By \Cref{lemma:S_is_union_of_cells} this implies
\[x\in \CS_{\delta,L} =\lp \prod_{j\in [n]\setminus (\Kbar\cup L)}(\ceil{\delta}_j-1,\ceil{\delta}_j]\times\prod_{k\in \Kbar}\{\ceil{\delta}_k\}\times \prod_{t\in L}[-a_t,\infty)\rp\cap i_\R(M_\R).\]
Thus there are real numbers $y_j \ge 0$ for $j\in [n]\setminus \Kbar$ so that
\[\pi_\R(x) = -\chi+\sum_{j\in [n]\setminus\Kbar}y_j\beta_j.\]
However, as $x\in i_\R(M_\R)$, by \eqref{eqn:ses}, $\pi_\R(x) = 0$. Rearranging shows that $\chi\in \beta_{[n]\setminus \Kbar}\subset |\Sigma_{GKZ}|$.
\end{proof}
\section{Window comparison for arrangeable GIT problems}\label{section:main_thm}
In this section we will prove \Cref{thm:window_equality}, which states that for arrangeable GIT problems, an HKKN stratification can be constructed so that the grade restriction window exactly equals the Cox window. By \Cref{cor: window equality} this amounts to showing that the Cox window generators $F_\bullet^D(D)$ are in the grade restriction window. We showed in \Cref{prop:cox_generator_grw_condition} that $F_\bullet^D(D)$ being in the grade restriction window is equivalent to the acyclicity of a number of topological spaces $\del \CS_{\delta,L}$. By \Cref{prop:nervecomputes}, it suffices to show that the associated nerves $N(\CC_{\delta,L})$ are acyclic under the assumptions from \Cref{prop:cox_generator_grw_condition}.

Recall that by \Cref{rmk:automatic_grade_restriction_rules_satisfied}, $F_\bullet^D(D)$ satisfy many of the required grade restriction rules to be contained in $W_{\bs{\lambda}}$. In fact, the only rules that $F_\bullet^D(D)$ do not automatically satisfy are the acyclicity of negative $\lambda$-weights for $\lambda\notin |\Sigma_{GKZ}|^\vee$. In this section, we will first state a proposition that handles this additional required acyclicity. We will use this to prove our main theorem and we will then spend the remainder of the section proving the required acyclicity.

Throughout this section we will assume that $\Gamma$ is a chamber of a projective, primitive geometric GIT problem with secondary fan $\Sigma_{GKZ}$ and that $D = \sum_{\rho\in [n]}a_\rho e_\rho$ is a lift of a degree $d \in \Gamma\cap \wh{G}$. Crucially, we assume that $\Sigma_{GKZ}$ is arrangeable (see \Cref{def:arrangeable}). Let $\bs{\lambda} = \lambda_1,\hdots, \lambda_k$ be any HKKN stratification of $V(B_\Gamma)$ obtained by extending the $(\Gamma^\vee \cap \Lambda, \wrel)$ poset, guaranteed by \Cref{thm:existence_of_wall_stratifications}.

\begin{proposition}\label{prop:neg_weight_acyclicity} Let $\Sigma_{GKZ}, \Gamma,\bs{\lambda},$ and $D$ be as above. Suppose $\lambda_i\notin |\Sigma_{GKZ}|^\vee$ and let $p_L\in (X_i)_{\lambda_i}^{=0}$ be any distinguished point. Let $\delta$ be any cell of $P_D$. If $\wh{\lambda}_i(\pi_\R(\ceil{\delta})) < 0$, then $N(\CC_{\delta,L})$ is acyclic.
\end{proposition}
We postpone the proof of \Cref{prop:neg_weight_acyclicity} to after the proof of \Cref{thm:window_equality}. Together, \Cref{lemma:weights_in_Tminus_Tplus} and \Cref{prop:neg_weight_acyclicity} are the technical ingredients proving the desired nerve acyclicty required to show window equality. We critically use these statements in the below proof of our main theorem.
\begin{theorem}\label{thm:window_equality} Let $\Sigma_{GKZ}$ be the secondary fan of an arrangeable, primitive geometric GIT problem. Let $\Gamma$ be a chamber of $\Sigma_{GKZ}$ and let $\bs{\lambda} = \lambda_1,\hdots, \lambda_k$ be any HKKN stratification of $V(B_\Gamma)$ obtained by extending the poset $(\Gamma^\vee \cap \Lambda, \wrel)$. Then $W_\Gamma^{\Cox} = W_{\bs{\lambda}}$. 
\end{theorem}
\begin{proof} Recall that by \Cref{thm:GRW_is_Cox_iff}, to show $W_\Gamma^\Cox = W_{\bs{\lambda}}$, it suffices to show that  for any $d \in \Gamma \cap \wh{G}$, for any $i$ such that $\lambda_i\notin |\Sigma_{GKZ}|^\vee$, any $p_L\in (X_i)_{\lambda_i}^{=0}$, and any $\delta\in \cells(P_D)$ such that $\wh{\lambda}(\pi_\R(\ceil{\delta})) < 0$, that $\del \CS_{\delta,L}$ is acyclic. By \Cref{prop:neg_weight_acyclicity} for such a $p_L$ and $\delta$, $N(\CC_{\delta,L})$ is acyclic. Recall that $\del \CS_{\delta,L}$ is homotopy equivalent to $N(\CC_{\delta,L})$ by \Cref{prop:nervecomputes}, so we conclude that $\del \CS_{\delta,L}$ is acyclic.
\end{proof}

In the remainder of the section, we will prove \Cref{prop:neg_weight_acyclicity}. This will require a technical detour. To start, we will provide an explicit combinatorial description of the distinguished points of $(X_i)_{\lambda_i}^{=0}$ arising from arrangeable HKKN stratifications.
\begin{proposition}\label{prop:internally_simplicial_gives_ones} Let $\Gamma$ be a chamber of an arrangeable secondary fan $\Sigma_{GKZ}$ equipped with an HKKN stratification $\bs{\lambda}$ obtained from totalizing the poset $(\Gamma^\vee \cap \Lambda, \wrel)$.
Let $L\subset [n]$ be arbitrary. Then $p_L\in (X_i)_{\lambda_i}^{=0}$ if and only if $L\subset I_{\lambda_i}^{=0}$ and for every $\rho\in I_{\lambda_i}^{=0}$, $L\cap [\rho] \ne \emptyset$. In particular, for any $\xi\in \wh{G}_\R^*$, if $I_{\lambda_i}^{=0}\cap I_\xi^{>0}\ne\emptyset$, then $L\cap I_\xi^{>0}\ne\emptyset$.
\end{proposition}
\begin{proof} First recall
\[(X_i)_{\lambda_i}^{=0} = V( I_{\lambda_i}^{\ne 0})\setminus \bigcup_{j < i}V(I_{\lambda_j}^{>0}).\]
So $p_L\in (X_i)_{\lambda_i}^{=0}$ if and only if $L\subset I_{\lambda_i}^{=0}$ and for each $j < i$, $L\cap I_{\lambda_j}^{>0}\ne \emptyset$. In one direction suppose that indeed for each $j < i$, $L\cap I_{\lambda_j}^{>0} \ne \emptyset$. Now consider any $\rho \in I_{\lambda_i}^{=0}$. By \Cref{lemma:internally_simplicial_implies_close_previous_wall} there exists $\lambda_j\in \Gamma^\vee\cap \Lambda$ with $\lambda_j \wrel \lambda_i$ and
\[I_{\lambda_j}^{>0}\cap I_{\lambda_i}^{=0} = [\rho].\]
As the ordering on $\Gamma^\vee \cap \Lambda$ is an extension of $\wrel$, we have that $j < i$. So we get
\[ L \cap [\rho] = L \cap I_{\lambda_i}^{=0}\cap I_{\lambda_j}^{>0} = L\cap I_{\lambda_j}^{>0} \ne \emptyset.\]
In the other direction, suppose for every $\rho\in I_{\lambda_i}^{=0}$, that $L\cap [\rho]\ne \emptyset$. Now suppose for some $j$, $L\cap I_{\lambda_j}^{>0} = \emptyset$. By our initial observation, it suffices to show $i < j$. Because $L\cap [\rho]\ne \emptyset$ for every $\rho\in I_{\lambda_i}^{=0}$, $L\cap I_{\lambda_j}^{>0}= \emptyset$ implies $I_{\lambda_i}^{=0}\cap I_{\lambda_j}^{>0} = \emptyset$. Equivalently, $I_{\lambda_i}^{=0}\subset I_{\lambda_j}^{\le 0}$, so $\lambda_i\wrel \lambda_j$. As the total ordering extends $\wrel$, we have $i < j$ as desired.
\end{proof}
By \Cref{prop:internally_simplicial_gives_ones}, we know exactly where $L$ lives in relation to the secondary fan for any $\lambda\notin |\Sigma_{GKZ}|^\vee$, namely $L$ contains an element at each weight in $I_{\lambda}^{=0}$. This perspective on $L$ will allow us to exploit the convex geometric conditions from \Cref{cor:integer_poset_cor} to find common cone sets for $N(\CC_{\delta,L})_\xi$. This convex geometric intuition motivates the following proposition.

\begin{proposition}\label{prop:main_prop_general}
Let $\Sigma_{GKZ}, \Gamma, \bs{\lambda},$ and $D$ be as above. Let $\lambda\coloneqq \lambda_i\notin |\Sigma_{GKZ}|^\vee$ and let $p_L\in (X_i)_{\lambda}^{=0}$ be any distinguished point. Let $\pi(D)\in \wh{G} \cap \Gamma$ and $\delta$ be a cell of $P_D$. Let $\chi \coloneq \chi(D,\delta,L)$ (see \Cref{def:distinguished_character}). For any $\xi\in \wh{G}_\R^{\ast}$ we have the following:
\begin{itemize}
    \item[(i)] If $\wh{\lambda}(\chi) \ge 0$, then $I_\lambda^{<0}\setminus \Kbar$ is a cone set for $N(\CC_{\delta,L})_\xi$.
    \item[(ii)] If $\wh{\lambda}(\chi) \le 0$, then $I_{\lambda}^{>0}\setminus \Kbar$ is a cone set for $N(\CC_{\delta,L})_\xi$. 
\end{itemize}
\end{proposition}
\begin{proof}
Throughout we assume $\xi\ne 0$, otherwise $N(\CC_{\delta,L})_\xi = 2^{[n]\setminus (\Kbar \cup L)}$ and the result follows immediately. To prove the first statement, suppose $\wh{\lambda}(\xi) \ge 0$.

We will proceed by analyzing the various ways $I_{\lambda}^{=0}$ can intersect $I_{\xi}^{<0}, I_\xi^{=0},$ and $I_\xi^{>0}$. First, suppose $I_\lambda^{=0}\subset I_\xi^{=0}$. Then by \Cref{cor:schmitt}\ref{cor:schmitt_v}, $\xi$ is a nonzero scalar multiple of $\lambda$, so $N(\CC_{\delta,L})_\xi = N(\CC_{\delta,L})_\lambda$. Then $I_{\lambda}^{<0}\setminus \Kbar$ is a cone set of $N(\CC_{\delta,L})_{\lambda}=  N(\CC_{\delta,L})_\xi$ by \Cref{cor:integer_poset_cor}\ref{integer_poset_cor:NEW_J_plus_chi_minus} applied to $-\lambda$.

In the next case suppose $I_\lambda^{=0} \cap I_\xi^{>0}\ne\emptyset$ and $I_\lambda^{=0}\cap I_\xi^{<0}\ne\emptyset$. Then by \Cref{prop:internally_simplicial_gives_ones} this implies both $L\cap I_\xi^{>0}\ne\emptyset$ and $L\cap I_\xi^{<0}\ne\emptyset$. Thus $I_\lambda^{<0}\setminus \Kbar$ is a cone set of $N(\CC_{\delta,L})_\xi$ by \Cref{cor:integer_poset_cor}\ref{integer_poset_cor:NEW_I_plus_minus}.

If we are not in the first two cases, then $I_\lambda^{=0}\ne I_\xi^{=0}$ and either $I_\lambda^{=0}\subset I_\xi^{\le 0}$ or $I_\lambda^{=0}\subset I_\xi^{\ge 0}$. As $N(\CC_{\delta,L})_\xi = N(\CC_{\delta,L})_{-\xi}$, we may assume without loss of generality that \[I_\lambda^{=0}\subset I_\xi^{\le 0}.\]
As $I_\lambda^{=0}\ne I_\xi^{=0}$, this implies $I_\lambda^{=0}\cap I_\xi^{<0}\ne \emptyset$, which by \Cref{prop:internally_simplicial_gives_ones} implies $L\cap I_\xi^{<0}\ne\emptyset$. Because $I_\lambda^{=0}\subset I_{\xi}^{\le 0}$, we have $\beta_{I_{\lambda}^{=0}}\subset H_{\xi}^{\le 0}$. As $\Sigma_{GKZ}$ is arrangeable and $\lambda\in \Lambda^\circ$, $\beta_{I_{\lambda}^{=0}}= H_{\lambda}^{=0}\cap |\Sigma_{GKZ}|$. Together this implies $H_{\lambda}^{=0}\cap |\Sigma_{GKZ}|\subset H_{\xi}^{\le 0}$. So by \Cref{lemma:hyperplane_and_halfspace_containment}, either
\[H_{\lambda}^{\le 0}\cap |\Sigma_{GKZ}| \subset H_{\xi}^{\le 0}\quad \text{ or }\quad H_{\lambda}^{\ge 0}\cap |\Sigma_{GKZ}|\subset H_{\xi}^{\le 0}.\]
If $H_{\lambda}^{\le 0}\cap |\Sigma_{GKZ}|\subset H_\xi^{\le 0}$, then $I_\lambda^{<0}\setminus \Kbar\subset I_\xi^{\le 0}$. Together with the fact that $L\cap I_\xi^{<0}\ne\emptyset$, this implies $I_\lambda^{<0}\setminus \Kbar$ is a cone set of $N(\CC_{\delta,L})_\xi$ by \Cref{cor:integer_poset_cor}\ref{integer_poset_cor:NEW_I_plus_J_plus} applied to $-\xi$. 

If instead $H_{\lambda}^{\ge 0}\cap |\Sigma_{GKZ}|\subset H_\xi^{\le 0}$, then by \Cref{lemma:chi_eff}, $\chi\in |\Sigma_{GKZ}|$ so $\wh{\xi}(\chi) \le 0$. Together with the fact that $L\cap I_\xi^{<0}\ne \emptyset$, this implies $I_\lambda^{<0}\setminus \Kbar$ is a cone set of $N(\CC_{\delta,L})_\xi$ by \Cref{cor:integer_poset_cor}\ref{integer_poset_cor:NEW_I_plus_chi_plus} applied to $-\xi$.

To prove the statement in the case $\wh{\lambda}(\chi) \le 0$, one can apply nearly the same proof as above, swapping $\lambda$ for $-\lambda$. The only step in the above proof that does not pass directly to $-\lambda$ is the usage of \Cref{prop:internally_simplicial_gives_ones}, which requires that $\lambda$ come from a wall stratification of $V(B_\Gamma)$. However, as $I_\lambda^{=0} = I_{-\lambda}^{=0}$, the conclusion of \Cref{prop:internally_simplicial_gives_ones} is the same for both $\lambda$ and $-\lambda$, so we are fine to reuse the above proof.
\end{proof}

The conclusion of \Cref{prop:main_prop_general} implies there is a common cone set for every $N(\CC_{\delta,L})_\xi$, regardless of the value of $\wh{\lambda}(\pi_\R(\ceil{\delta}))$. However, to use \Cref{cor:cone_set_for_xi_implies_nerve_contractible} to show that $N(\CC_{\delta,L})$ is contractible, we have to show that the cone sets are nonempty. This turns out generally to only be true when $\wh{\lambda}(\pi_\R(\ceil{\delta})) < 0$. This is stated precisely in \Cref{prop:neg_weight_acyclicity}, which we are now in a position to prove.
\begin{proof}[Proof of \Cref{prop:neg_weight_acyclicity}]
By \Cref{prop:main_prop_general}, depending on the sign of $\wh{\lambda}(\chi)$, either $I_{\lambda}^{<0}\setminus \Kbar$ is a cone set for $N(\CC_{\delta,L})_\xi$ for every $\xi\in \wh{G}_\R^*$ or $I_{\lambda}^{>0}\setminus \Kbar$ is a cone set for $N(\CC_{\delta,L})_\xi$ for every $\xi\in \wh{G}_\R^*$. 
Then by \Cref{cor:cone_set_for_xi_implies_nerve_contractible}, it suffices to show that $\wh{\lambda}(\pi_\R(\ceil{\delta})) < 0$ implies that both $I_{\lambda}^{<0}\setminus \Kbar$ and $I_\lambda^{>0}\setminus \Kbar$ are nonempty.

We will first show that $I_\lambda^{>0}\setminus \Kbar$ is nonempty. We note in advance that this case does not actually require that $\wh{\lambda}(\pi_\R(\ceil{\delta})) < 0$. As $\delta$ is a nonempty cell of $P_D$ and $\delta\subset \CS_{\delta,L}$ by \Cref{lemma:S_is_union_of_cells}, we have
\[ \CS_{\delta,L} = \lp \prod_{j\in [n]\setminus (\Kbar\cup L)}(\ceil{\delta}_j-1,\ceil{\delta}_j]\times \prod_{k\in \Kbar}\{\ceil{\delta}_k\}\times \prod_{t\in L}[-a_t,\infty)\rp\cap i_\R(M_\R)\ne \emptyset.\]
In particular, as $\ceil{\delta}_k = -a_k$ for $k\in \Kbar$ and $\ceil{\delta}_j> -a_j$ for $j\in [n]\setminus (\Kbar \cup L)$, this implies
\begin{align*}
    \lp -D+\Cone(e_i \ | \ i\in [n]\setminus \Kbar)\rp\cap i_\R(M_\R)&\ne\emptyset.
\end{align*}
As $i_\R(M_\R) = \ker \pi_\R$ by \eqref{eqn:ses}, applying $\pi_\R$ to the above implies that $0\in -\pi_\R(D)+\beta_{[n]\setminus \Kbar}$. Rearranging, $\pi_\R(D)\in \beta_{[n]\setminus \Kbar}$. As $\pi_\R(D)\in \Int \Gamma$ by assumption of $D$ being a lift of an ample degree, $\beta_{[n]\setminus \Kbar}\cap \Int \Gamma \ne \emptyset$. Thus by \Cref{cor:gamma_hits_beta_J_implies_contained}, $\Gamma \subset \beta_{[n]\setminus \Kbar}$. If $I_\lambda^{>0}\setminus \Kbar = \emptyset$, then $\Gamma \subset \beta_{I_\lambda^{\le 0}\setminus \Kbar} \subset H_\lambda^{\le 0}$. This contradicts $\lambda\in \Gamma^\vee$ because $\Gamma$ is full-dimensional. We conclude $I_\lambda^{>0}\setminus \Kbar \ne \emptyset$. 

To show $I_\lambda^{<0}\setminus \Kbar$, we will use $\wh{\lambda}(\pi_\R(\ceil{\delta})) < 0$. Because $\delta\subset P_D\subset \R_{-D}^n$ and $\ceil{\delta}_k = -a_k$ for $k\in \Kbar$, we get that $\delta$ is contained in a smaller convex subset of $\R^n$:
\[\delta\subset \lp \ceil{\delta}+\prod_{j\in [n]\setminus \Kbar}(-1,0]\times \prod_{k\in \Kbar}\{0\}\rp\cap i_\R(M_\R).\]
In particular, the right hand side of the above inclusion is nonempty. As in the above proof, applying $\pi_\R$ to the above and noting $i_\R(M_\R) = \ker \pi_\R$ by \eqref{eqn:ses} implies that
$\pi_\R(\ceil{\delta})\in \beta_{[n]\setminus \Kbar}$. Applying $\wh{\lambda}$, we find
\[\wh{\lambda}(\pi_\R(\ceil{\delta}))\in \Cone\lp \wh{\lambda}(\beta_j) \ | \ j\in I_\lambda^{\ne 0}\setminus \Kbar)\rp\subset \R,\]
so $\wh{\lambda}(\pi_\R(\ceil{\delta})) < 0$ implies $I_\lambda^{<0}\setminus \Kbar\ne\emptyset$ as desired.
\end{proof}

\section{Examples and counterexamples}\label{section:examples}

In \Cref{thm:existence_of_wall_stratifications} we showed that for arrangeable problems, one can place an ordering on $\Gamma^\vee \cap \Lambda$ which induces an HKKN stratification. In \Cref{thm:window_equality} we showed that using this particular stratification, the grade restriction window exactly equals the Cox window. As discussed in \S\ref{section:stratifications}, there is generally great flexibility in choosing an HKKN stratification, so it is natural to wonder if for more general stratifications we obtain window equality. In this section we show this is not the case. In particular, we describe a chamber of a GIT problem and prove there are no HKKN stratifications yielding window equality. In addition, we show that even for arrangeable problems, for arbitrary HKKN stratifications we do not see window equality. As a consequence, we provide examples of what can occur when one alters an HKKN stratification by reordering or dropping some of the one-parameter subgroups.

\subsection{A GIT problem where the Cox window and GRW are never equal}\label{subsection:counter_example_theory} In this subsection, we prove the following theorem.
\begin{theorem}\label{thm:no_stratification_making_windows_equal_exist}
There exists a projective GIT problem with a chamber $\Gamma$ such that for any HKKN stratification $\bs{\lambda}$ of $\CX_\Gamma$ and any choice of $\bs{w}$,
\[W^{\bs{w}}_{\bs{\lambda}} \neq W^{\Cox}_{\Gamma}. \] 
\end{theorem}
To prove \Cref{thm:no_stratification_making_windows_equal_exist}, we will first provide a sufficient condition for when a chamber of a GIT problem has no grade restriction window matching the Cox window (\Cref{prop:non_convex_sufficient_condition}). Then we will construct such an example meeting this sufficient condition. Checking that the condition is satisfied is then a purely computational exercise and we include \lst{Macaulay2} code in the appendix which does this computation.

To describe this sufficient condition we will need to introduce some notation. As usual, let $G\subset \Gm^n$ induce a toric GIT problem with secondary fan $\Sigma_{GKZ}$. Fix a chamber $\Gamma$ of $\Sigma_{GKZ}$. Fix a degree $d\in \Gamma\cap \wh{G}$ and fix a lift $D$ under $\pi$, and a subset $L\subset [n]$. Then define
\begin{equation}\label{eqn:nonacyclic_cell_labels}
    N_{L} \coloneqq \{\ceil{\delta} \ | \ \delta\in \cells(P_D), \ H_*(\del \CS_{\delta, L})\ne 0\}.
\end{equation}
When $L$ is clear, we will drop it. In \Cref{prop:cox_generator_grw_condition}, we established a connection between $N_{L}$ and the membership of $F_\bullet^D(D)$ (see \S\ref{subsection:cox_generator_images}) in a given grade restriction window $W_{\bs{\lambda}}^{\bs{w}}$. Namely, $F_\bullet^D(D)\in W_{\bs{\lambda}}^{\bs{w}}$ if and only if  $\wh{\lambda}_i(\pi_\R(N_L)) \subset [w_i,w_i+T_{\lambda_i}^{>0})$ for certain choices of $i,L$ (see \Cref{prop:cox_generator_grw_condition} for the precise statement). By unpacking this statement, we can obtain a simple condition for when a 1-PS $\lambda$ cannot be in any HKKN stratification $\bs{\lambda}$ such that $F_\bullet^D(D)\in W_{\bs{\lambda}}^{\bs{w}}$.
\begin{lemma}\label{lemma:single_lambda_weights_too_wide_condition} 
Given an HKKN stratification $\bs{\lambda}$ of $V(B_\Gamma)$, if there exist $\bs{a},\bs{a}'\in N_{I_{\lambda_i}^{=0}}$ such that
\[\wh{\lambda}_i(\pi_\R(\bs{a}-\bs{a}')) \ge T_{\lambda_i}^{>0},\]
for some $i$, then there is no choice of weights $\bs{w}$ such that $W_\Gamma^{\Cox} = W_{\bs{\lambda}}^{\bs{w}}$.
\end{lemma}
\begin{proof} For ease of notation set $L\coloneqq I_{\lambda_i}^{=0}$. First recall that for any HKKN stratification $\bs{\lambda}$, for every $j < i$, $I_{\lambda_j}^{>0}\not\subset I_{\lambda_i}^{\ne 0}$ by \Cref{lemma:HKKN_strat_combo}\ref{lemma:HKKN_strat_combo_ii}. It follows that $p_L\in (X_i)_{\lambda_i}^{=0}$. Thus by \Cref{prop:cox_generator_grw_condition}, for $F_\bullet^D(D)$ to be in $W_{\bs{\lambda}}^{\bs{w}}$, we must have
\[\wh{\lambda}_i(\pi_\R(N_L)) \subset [w_i,w_i+T_{\lambda_i}^{>0}).\]
However, this containment directly contradicts the inequality assumption. We conclude that there are no weights $\bs{w}$ so that $F_\bullet^D(D)$. As the complexes $F_\bullet^D(D)$ generate $W_\Gamma^{\Cox}$ (\Cref{cor:wcoxisgeneratedbyf}), it follows that there are no weights $\bs{w}$ such that $W_{\bs{\lambda}}^{\bs{w}} = W_{\Gamma}^{\Cox}$.
\end{proof}
In light of the above, we call the half-space $H_{\pi(\bs{a}-\bs{a}')-
\sum_{i\in I}\beta_i}^{\ge 0}$  of 1-PSs $\lambda$ such that \[\wh{\lambda}\lp\pi_\R\lp\bs{a}-\bs{a}'\rp - \sum_{\rho \in I}\beta_\rho\rp \geq 0\] 
the \emph{violating half-space} of $(\bs{a},\bs{a}')\in N_J\times N_J$ for some fixed positive set $I\subset [n]$ and zero-set $J\subset [n]$. If $\lambda$ is a 1-PS in a violating half-space such that $I_\lambda^{>0} = I$ and $I_\lambda^{=0} = J$, then $\lambda$ cannot appear in an HKKN stratification such that the Cox window and grade restriction window agree by the above lemma. This motivates to consider the set of all such 1-PSs.

\begin{figure}
\begin{tikzpicture}
\def\fanscale{1.8}
\coordinate (lambda0) at (0,0);
\coordinate (lambda1) at (0,2*\fanscale);
\coordinate (lambda2) at (2*\fanscale, 2*\fanscale);
\coordinate (lambda3) at (3*\fanscale, 1*\fanscale);
\coordinate (lambda4) at (2*\fanscale, 0*\fanscale);

\draw[draw = none, fill = red, opacity = 0.2] (lambda0) -- (lambda1) -- (lambda2) -- (lambda3) -- (lambda4) -- cycle;

\draw[draw = red, thick] (lambda4) -- (lambda0) -- (lambda1) -- (lambda2);

\draw[black, dotted, thick, <->] ($(lambda2)!-0.2!(lambda3)$) -- ($(lambda2)!1.2!(lambda3)$);
\draw[black, dotted, thick, <->] ($(lambda3)!-0.2!(lambda4)$) -- ($(lambda3)!1.2!(lambda4)$);

\node[circle,fill,inner sep=1pt, red] at (lambda0) {};
\node[circle,fill,inner sep=1pt, red] at (lambda1) {};

\node[draw = none, fill = none, red] at ($(lambda0) + (-0.3,-0.3)$) {\small $\lambda_0$};
\node[draw = none, fill = none, red] at ($(lambda1) + (-0.3,0.3)$) {\small $\lambda_1$};
\node[draw = none, fill = none, red] at ($(lambda2) + (0.3,0.3)$) {\small $\lambda_2$};
\node[draw = none, fill = none, red] at ($(lambda3) + (0.4,0)$) {\small $\lambda_3$};
\node[draw = none, fill = none, red] at ($(lambda4) + (0.3,-0.3)$) {\small $\lambda_4$};

\draw[black, dotted, thick, <->] ($(lambda0)!-0.2!(lambda1)$) -- ($(lambda0)!1.2!(lambda1)$);
\draw[black, dotted, thick, <->] ($(lambda1)!-0.2!(lambda2)$) -- ($(lambda1)!1.2!(lambda2)$);
\draw[black, dotted, thick, <->] ($(lambda0)!-0.2!(lambda4)$) -- ($(lambda0)!1.2!(lambda4)$);

\node[draw = none, fill = none, black] at ($(lambda0)!0.5!(lambda1)+(-0.3,0)$) {\small $\beta_4$};
\node[draw = none, fill = none, black] at ($(lambda1)!0.5!(lambda2)+(0,0.3)$) {\small $\beta_2,\beta_3$};
\node[draw = none, fill = none, black] at ($(lambda2)!0.5!(lambda3)+(0.3,0.3)$) {\small $\beta_6$};
\node[draw = none, fill = none, black] at ($(lambda3)!0.5!(lambda4)+(0.3,-0.3)$) {\small $\beta_5$};
\node[draw = none, fill = none, black] at ($(lambda0)!0.5!(lambda4)+(0,-0.3)$) {\small $\beta_0,\beta_1$};

\node[draw = none, fill = none, red] at (lambda2) {};
\node[circle,inner sep=1pt, draw = red,fill=white, thick]  at (lambda3) {};
\node[draw = none, fill = none, red] at (lambda4) {};

\node[draw = none, fill = none, red] at (\fanscale*1.6,\fanscale*1) {\small $\overline{C}_{\{5,6\},\emptyset}$};
\node[draw = none, fill = none, red] at (\fanscale*1,\fanscale*1.85) {\small $\overline{C}_{\{5,6\},\{2,3\}}$};
\node[draw = none, fill = none, red] at (\fanscale*1,\fanscale*0.15) {\small $\overline{C}_{\{5,6\},\{0,1\}}$};
\node[draw = none, fill = none, red] at (\fanscale*0.45,\fanscale*1) {\small $\overline{C}_{\{5,6\},\{4\}}$};
\node[draw = none, fill = none, red] at (\fanscale*-0.7,\fanscale*2.4) {\small $\overline{C}_{\{5,6\},\{2,3,4\}}$};
\node[draw = none, fill = none, red] at (\fanscale*-0.7,\fanscale*-0.4) {\small $\overline{C}_{\{5,6\},\{0,1,4\}}$};

\node[draw = none, fill = none, black] at (\fanscale*3,\fanscale*2) {\small $\wh{G}_\R^\ast$};

%%%%%%% NEW FAN
\def\fanscale{1.8}
\coordinate (shift) at (8,0);

\coordinate (beta2) at ($(shift)+(0,2*\fanscale)$);
\coordinate (beta6) at ($(shift)+(3*\fanscale,0)$);
\coordinate (beta0) at (shift);
\coordinate (beta3) at (beta2);
\coordinate (beta4) at ($(shift)+(1*\fanscale,1*\fanscale)$);
\coordinate (beta1) at (beta0);
\coordinate (beta5) at ($(shift)+(2*\fanscale, 1*\fanscale)$);

\coordinate (alpha0) at ($(shift)+(1.5*\fanscale,1*\fanscale)$);
\coordinate (alpha1) at ($(shift)+(1.5*\fanscale,0.75*\fanscale)$);
\coordinate (alpha2) at ($(shift)+(1.714*\fanscale,0.857*\fanscale)$);

\draw[draw=none,fill=gray,opacity=.2] (beta1) -- (beta2) -- (beta5) -- (beta6)--cycle;

\node[circle,fill,inner sep=1pt] at (beta0) {};
\node[circle,fill,inner sep=1pt] at (beta1) {};
\node[circle,fill,inner sep=1pt] at (beta2) {};
\node[circle,fill,inner sep=1pt] at (beta3) {};
\node[circle,fill,inner sep=1pt] at (beta4) {};
\node[circle,fill,inner sep=1pt] at (beta5) {};
\node[circle,fill,inner sep=1pt] at (beta6) {};

\node[draw = none, fill = none, black] at ($(beta0)+(-0.6,0)$) {\small $\beta_0,\beta_1$};

\node[draw = none, fill = none, black] at ($(beta2)+(-0.6,0)$) {\small $\beta_2,\beta_3$};

\node[draw = none, fill = none, black] at ($(beta4)+(-0.5,0)$) {\small $\beta_4$};

\node[draw = none, fill = none, black] at ($(beta5)+(0.4,0)$) {\small $\beta_5$};

\node[draw = none, fill = none, black] at ($(beta6)+(0.4,0)$) {\small $\beta_6$};

\draw[red,dotted,thick, <->] ($(beta0)!-0.5!(beta4)$) -- ($(beta0)!1.5!(beta4)$);
\draw[red,dotted, thick, <->] ($(beta2)!-0.5!(beta4)$) -- ($(beta2)!1.5!(beta4)$);
\draw[red,dotted, thick, <->] ($(beta2)!-0.2!(beta6)$) -- ($(beta2)!1.2!(beta6)$);
\draw[red,dotted, thick, <->] ($(beta5)!-0.6!(beta6)$) -- ($(beta5)!1.6!(beta6)$);
\draw[red,dotted,thick, <->] ($(beta0)!-0.3!(beta5)$) -- ($(beta0)!1.3!(beta5)$);

\node[draw = none, fill = none, red] at ($(beta0)!0.5!(beta4)+(-0.2,0.2)$) {\small $\lambda_0$}; 
\node[draw = none, fill = none, red] at ($(beta2)!0.5!(beta4)+(-0.2,-0.2)$) {\small $\lambda_1$}; 
\node[draw = none, fill = none, red] at ($(beta2)!0.3!(beta6)+(0.2,0.2)$) {\small $\lambda_2$}; 
\node[draw = none, fill = none, red] at ($(beta5)!0.5!(beta6)+(0.2,0.2)$) {\small $\lambda_3$}; 
\node[draw = none, fill = none, red] at ($(beta0)!0.3!(beta5)+(0.2,-0.2)$) {\small $\lambda_4$}; 

\node[draw=none,fill=none,black] at ($(shift)+(2*\fanscale, 2*\fanscale)$) {$\wh{G}_\R$};
\end{tikzpicture}
    \caption{On the left is the cone $\overline{C}_{\{5,6\}, \emptyset}$ in $\wh{G}_\R^\ast$ as well as the cones $\overline{C}_{\{5,6\},\{0,1\}}$, $\overline{C}_{\{5,6\},\{4\}}$, $\overline{C}_{\{5,6\},\{2,3\}}$, $\overline{C}_{\{5,6\},\{0,1,4\}}$, and $\overline{C}_{\{5,6\},\{2,3,4\}}$ which are themselves faces of $\overline{C}_{\{5,6\},\emptyset}$. The facets of $\overline{C}_{\{5,6\},\emptyset}$ are labeled by the corresponding $\beta_\rho$. Rays of this cone are labeled by their generators $\lambda_i$. On the right is the corresponding image in $\wh{G}_\R$. The hyperplanes $H_{\lambda_i}^{=0}$ are labeled by their corresponding $\lambda_i$.}
    \label{fig:C56_cone_NEW}
\end{figure}
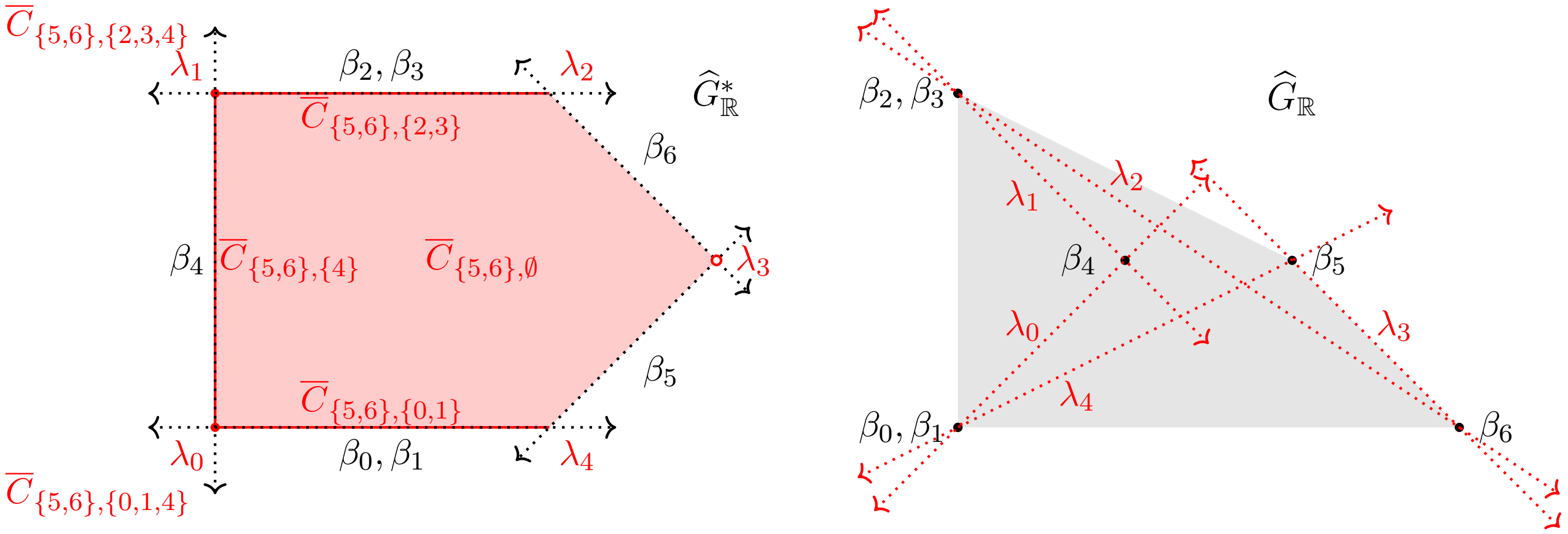

\begin{definition}\label{def:open_C_plus_zero_cone} Consider any two disjoint subsets $I,J \subset [n]$. We denote by 
\[C_{I,J} \coloneqq \{\lambda\in \wh{G}_\R^* \ | \ I_{\lambda}^{>0} = I \text{ and }I_{\lambda}^{=0} = J\}\]
the convex region which parametrizes all  1-PSs with positive set $I$ and zero-set $J$.
\end{definition}
We can now use \Cref{lemma:single_lambda_weights_too_wide_condition} to obtain a sufficient condition for window inequality.
\begin{proposition}\label{prop:non_convex_sufficient_condition} Suppose $\la I \ra$ is a component of the primary decomposition of $B_\Gamma$. If every $\lambda$ with positive set $I$ lies inside some violating half space,
\begin{equation*}
    C_{I,I^{=0}_{\lambda}} \subset \bigcup_{(\bs{a},\bs{a}')\in N_{I_\lambda^{=0}}\times N_{I_\lambda^{=0}}} H_{\pi(\bs{a}-\bs{a}')-\sum_{i\in I}\beta_i}^{\ge 0},
\end{equation*}
then $W_{\Gamma}^{\Cox} \neq W_{\bs{\lambda}}^{\bs{w}}$ for any choice of $\bs{\lambda}, \bs{w}$, i.e.\ there is no combination of HKKN stratification $\bs{\lambda}$ of $V(B_\Gamma)$ and weights $\bs{w}$ such that the Cox window and grade restriction window agree. 
\end{proposition}
\begin{proof} Let $\bs{\lambda}$ be any HKKN stratification of $V(B_\Gamma)$. Because $\la I\ra$ appears in the primary decomposition of $B_\Gamma$, by the ideal property of HKKN stratifications, \Cref{lemma:HKKN_strat_combo}\ref{lemma:HKKN_strat_combo_i}, there must be some $\lambda \coloneqq \lambda_i$ in the stratification such that $I_\lambda^{>0} = I$. As $\lambda\in C_{I,I_{\lambda}^{=0}}$, by the assumed inclusion there exist $\bs{a},\bs{a}'\in N_{I_{\lambda}^{=0}}$ such that
\begin{align*}
    \wh{\lambda}\lp \pi(\bs{a}-\bs{a}')-\sum_{i\in I}\beta_i\rp & \ge 0\\
    \wh{\lambda}(\pi(\bs{a}-\bs{a}')) & \ge \sum_{i\in I}\wh{\lambda}(\beta_i) = T_{\lambda}^{>0}.
\end{align*}
Thus by \Cref{lemma:single_lambda_weights_too_wide_condition}, there is no choice of weights $\bs{w}$ such that $W_\Gamma^{\Cox} = W_{\bs{\lambda}}^{\bs{w}}$.
\end{proof}
We are now in a position to describe the GIT problem and chamber appearing in the statement of \Cref{thm:no_stratification_making_windows_equal_exist}. Consider the GIT problem defined by $G\cong \Gm^3$ including into $\Gm^7$ via weights
\[\begin{array}{cccc} \beta_0 =\beta_1 = (0,0,1), & \beta_2 =\beta_3 = (0,2,1), & \beta_4 = (1,1,1),\\ \beta_5 = (2,1,1), & \beta_6 = (3,0,1).
\end{array}\]
\begin{figure}
\begin{tikzpicture}
%%% Fan 1

\def\fanscale{2}
\def\xshift{4}

\coordinate (beta2) at (0,2*\fanscale);
\coordinate (beta6) at (3*\fanscale,0);
\coordinate (beta0) at (0,0);
\coordinate (beta3) at ($(beta2)$);
\coordinate (beta4) at (1*\fanscale,1*\fanscale);
\coordinate (beta1) at ($(beta0)$);
\coordinate (beta5) at ($0.666*(beta6)+0.5*(beta2)$);

\coordinate (alpha0) at (1.5*\fanscale,1*\fanscale);
\coordinate (alpha1) at (1.5*\fanscale,0.75*\fanscale);
\coordinate (alpha2) at (1.714*\fanscale,0.857*\fanscale);

\draw[draw=none,fill=gray,opacity=.2] (beta1) -- (beta2) -- (beta5) -- (beta4)--(beta6)--cycle;

\draw[draw=none,fill=gray,opacity=.2] (beta5) -- (alpha0) -- (alpha2) -- (alpha1)--(beta6)--cycle;

\draw[draw=none,fill=red,opacity=.2] (beta4) -- (alpha0) -- (alpha2) -- (alpha1)--cycle;

\draw[-,black,thick] (beta1) -- (beta4);
\draw[-,black,thick] (beta2) -- (beta4);
\draw[-,black,thick] (beta0) -- (beta2);
\draw[-,black,thick] (beta1) -- (beta5);
\draw[-,black,thick] (beta4) -- (beta6);
\draw[-,black,thick] (beta0) -- (beta6);
\draw[-,black,thick] (beta0) -- (beta4);
\draw[-,black,thick] (beta5) -- (beta6);
\draw[-,black,thick] (beta2) -- (beta6);
\draw[-,black,thick] (beta4) -- (beta5);
\draw[-,black,thick] (beta2) -- (beta5);

\node[circle,fill,inner sep=1pt] at (beta0) {};
\node[circle,fill,inner sep=1pt] at (beta1) {};
\node[circle,fill,inner sep=1pt] at (beta2) {};
\node[circle,fill,inner sep=1pt] at (beta3) {};
\node[circle,fill,inner sep=1pt] at (beta4) {};
\node[circle,fill,inner sep=1pt] at (beta5) {};
\node[circle,fill,inner sep=1pt] at (beta6) {};

\node[draw = none, fill = none, black] at ($(beta0)+(-0.6,0)$) {\small $\beta_0,\beta_1$};

\node[draw = none, fill = none, black] at ($(beta2)+(-0.6,0)$) {\small $\beta_2,\beta_3$};

\node[draw = none, fill = none, black] at ($(beta4)+(-0.5,0)$) {\small $\beta_4$};

\node[draw = none, fill = none, black] at ($(beta5)+(0.4,0)$) {\small $\beta_5$};

\node[draw = none, fill = none, black] at ($(beta6)+(0.4,0)$) {\small $\beta_6$};
\node[draw = none, fill = none, red] at ($(beta4)+(0.9,-0.2)$) {\small $\Gamma$};
\end{tikzpicture}
    \caption{Secondary fan $\Sigma_{GKZ}$ associated to the weights described in \S\ref{subsection:counter_example_theory}.}
    \label{fig:GRW_ne_cox_GKZ}
\end{figure}
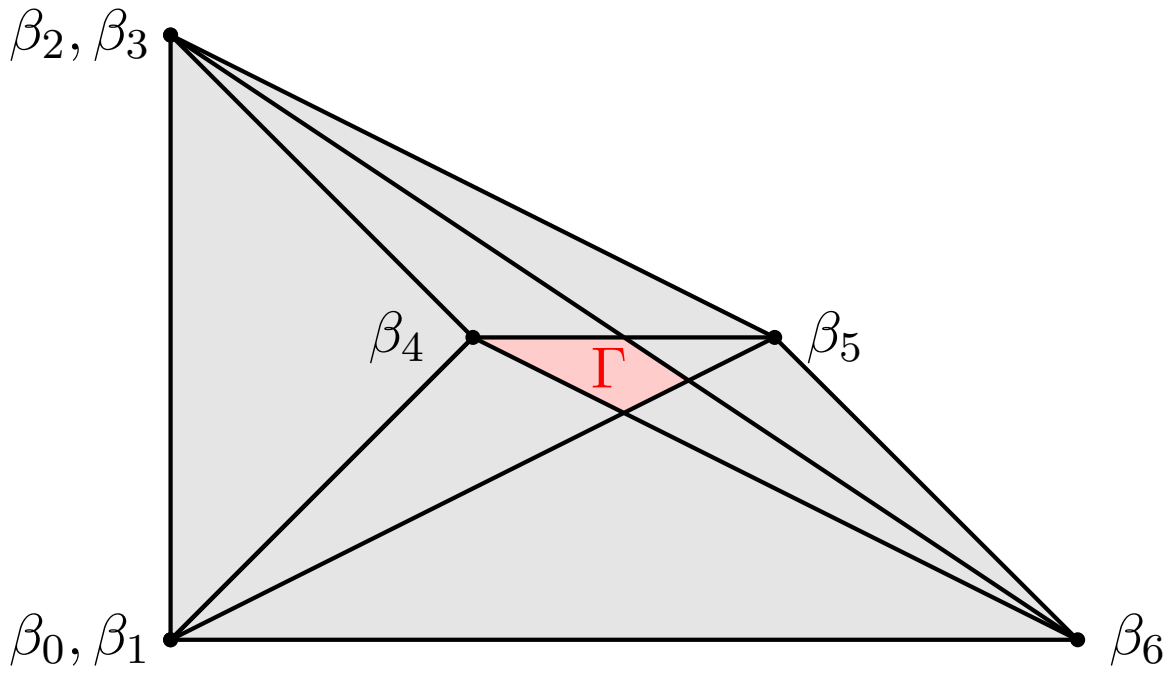
Each weight $\beta_\rho$ has third coordinate equal to 1, so we can draw the associated secondary fan $\Sigma_{GKZ}$ coming from a height 1 slice in the third coordinate (see \Cref{fig:GRW_ne_cox_GKZ}). Let $\Gamma$ be the red shaded chamber in the center of this secondary fan. We have previously encountered this chamber and secondary fan in \Cref{example:relation_not_a_poset} where we demonstrated that $(\Gamma^\vee \cap \Lambda, \wrel)$ was not a poset. $\CX_\Gamma$ is a smooth stack with simplicial, but non-smooth coarse moduli space $X_{\Sigma_\Gamma}$. The irrelevant ideal of this chamber can easily be computed (e.g. via \Cref{lemma:B_chi_via_betas}) as
\begin{align*}
B_\Gamma = (x_5,x_6)\cap(x_4,x_1,x_0)\cap(x_6,x_1,x_0)\cap (x_3,x_2,x_4)\cap (x_5,x_3,x_2).
\end{align*}
From this description it is evident that $(x_5,x_6)$ is a component of the primary decomposition of $B_\Gamma$. The cone $C_{\{5,6\},\emptyset}$ is depicted in \Cref{fig:C56_cone_NEW}. Its rays have generators
\[\begin{array}{ccc}
    \lambda_0 = (1, -1, 0), &\lambda_1 = (1, 1, -2), &\lambda_2 = (2, 3, -6),\\
    \lambda_3 = (1, 1 ,-3),&\lambda_4 = (1, -2, 0).
\end{array}\]
From the depiction in \Cref{fig:C56_cone_NEW}, we see that the only possible sets $I^{=0}$ for which there exist $\lambda$ such that $I_{\lambda}^{>0} = \{5,6\}$ and $I_{\lambda}^{=0} = I^{=0}$ are
\begin{equation*}
    I^{=0} = \emptyset, \{0,1\}, \{4\},\{2,3\}, \{0,1,4\},\text{ and }\{2,3,4\}.
\end{equation*}
In order to apply \Cref{prop:non_convex_sufficient_condition} to show that $W_\Gamma^{\Cox}$ is not a grade restriction window, we must additionally fix a lift of $d \in \Gamma\cap \wh{G}$. We set $D = (0,0,2,-5,7,5,-2)$. With this in hand we are in a position to prove \Cref{thm:no_stratification_making_windows_equal_exist}.
\begin{proof}[Proof of \Cref{thm:no_stratification_making_windows_equal_exist}] Let $\Gamma, \Sigma_{GKZ}, D$ be as above. As discussed above, $(x_5,x_6)$ is a component of the primary decomposition of $B_\Gamma$. As above, the only $I^{=0}$ for which $C_{\{5,6\}, I^{=0}}$ is nonempty are
\[I^{=0} = \emptyset, \{0,1\}, \{4\},\{2,3\}, \{0,1,4\},\text{ and }\{2,3,4\}.\]
Thus, by \Cref{prop:non_convex_sufficient_condition}, it suffices to show that for each of the six above choices for $I^{=0}$, that
\[C_{\{5,6\},I^{=0}}\subset \bigcup_{(\bs{a},\bs{a}')\in N_{I^{=0}}\times N_{I^{=0}}} H_{\pi(\bs{a}-\bs{a}')-\beta_5-\beta_6}^{\ge 0}.\]
To complete the proof we exhaustively verify each of these inclusions in \lst{Macaulay2} using a convex restatement of the above inclusion (\Cref{cor:convex_sufficient_condition}) using code introduced in \S\ref{subsection:appendix_methods}. See \Cref{example:no_stratifications_code} for this exact computation.
\end{proof}

\subsection{Grade restriction windows of closely related stratifications}\label{subsection:arrangeable_with_bad_stratification}  In this subsection, we prove Theorems \ref{intro thm: reordering} and \ref{intro thm: substratification} (see \Cref{cor:may_or_may_not})
by analyzing four different HKKN stratifications of a fixed chamber in the arrangeable GIT problem below.
Our first example is the most involved and provides a HKKN stratification of an arrangeable GIT problem whose Cox window and grade restriction window do not agree.  This illustrates the use of our special HKKN stratifications in the hypothesis of \Cref{thm:window_equality}. The stratification we use is very close to one in which the windows agree.  In fact, we go on to give three minor variations of this stratification all of which result in window equality.  Together, these examples illustrate that grade restriction windows may change under reordering and under restriction to an ordered subset.

\begin{proposition}\label{prop:arrangeable_with_bad_wall_strat}
There exists an arrangeable GIT problem with a chamber $\Gamma$ and an HKKN stratification $\bs{\lambda}$ of $V(B_\Gamma)$ with $\bs{\lambda} \subsetneq \Gamma^{\vee} \cap \Lambda$ such that for any choice of weights $\bs{w}$,
\[W^{\bs{w}}_{\bs{\lambda}}\neq W_{\Gamma}^{\Cox}.\]
\end{proposition}

We will now construct the GIT problem, chamber, and stratification which witnesses \Cref{prop:arrangeable_with_bad_wall_strat}. Let $\Sigma$ be the fan in $\R^3$ with rays $u_1,\hdots, u_6$ given by
\[\arraycolsep=1cm
\begin{array}{ccc} u_1 = (-1, -1, 0), &
    u_2 = (1, 0, 0), &
    u_3 = (-1, 0, 0), \\
    u_4 = (0, 1, 1),&
    u_5 = (0, -1 ,-1), &
    u_6 = (0,0,1),
\end{array}\]
and maximal cones
\[\arraycolsep=0.2cm
\Sigma_{\rm max} = \left\{\begin{array}{cccc} \Cone(u_2,u_5,u_6),&
    \Cone(u_3,u_5,u_6), &
    \Cone(u_2,u_4,u_6), &
    \Cone(u_1,u_3,u_4), \\
    \Cone(u_1,u_4,u_6),&
    \Cone(u_1,u_3,u_6),&
    \Cone(u_2,u_4,u_5),&
    \Cone(u_3,u_4,u_5),
\end{array}\right\}.\]
One can check that $X_\Sigma$ is the blowup of $\P^2\times \P^1$ at a torus fixed point. From our description of the rays of $\Sigma$, one computes that $\wh{G} = \Cl (X_\Sigma)\cong \Z^3$ with weights
\[\arraycolsep=1cm
\begin{array}{ccc} \beta_1 = (1,0,0),&
    \beta_2 = (1,1,0),&
    \beta_3 = (0,1,0),\\
    \beta_4 = (1,0,1),&
    \beta_5 = (0,0,1),&
    \beta_6 = (0,0,1).
\end{array}\]
\begin{figure}
\begin{tikzpicture}
\def\fanscale{6}
\def\trih{0.866} % height of equilateral triangle
\def\op{0.3}
\def\lineop{0.5}
\def\thickness{thick}
\def\dashsty{dashed}
\def\pslength{1.5cm}

\coordinate (beta1) at (-0.5*\fanscale, 0*\fanscale);
\coordinate (beta56) at (0.5*\fanscale, 0*\fanscale);
\coordinate (beta3) at (0*\fanscale, \trih*\fanscale);
\coordinate (beta2) at ($(beta1)!0.5!(beta3)$);
\coordinate (beta4) at ($(beta1)!0.5!(beta56)$);
\coordinate (center) at ($(beta4)!0.333!(beta3)$);

\draw[draw=none, fill=gray,opacity=\op] (beta1) -- (beta2) -- (beta4) -- cycle;
\draw[draw=none, fill=gray,opacity=\op] (beta2) -- (beta3) -- (beta56) -- (beta4) -- (center) -- cycle;
\draw[draw=none,fill = red, opacity=\op] (beta2) -- (center) -- (beta4) -- cycle;

\tikzset{beta/.style={circle,inner sep =1.5pt, fill = black, draw =black}}
\node[beta] at (beta1) {};
\node[beta] at (beta2) {};
\node[beta] at (beta3) {};
\node[beta] at (beta4) {};
\node[beta] at (beta56) {};

\draw[draw=black,\thickness,opacity=\lineop] (beta1) -- (beta3);
\draw[draw=black,\thickness,opacity=\lineop] (beta3) -- (beta56);
\draw[draw=black,\thickness,opacity=\lineop] (beta56) -- (beta1);
\draw[draw=black,\thickness,opacity=\lineop] (beta2) -- (beta56);
\draw[draw=black,\thickness,opacity=\lineop] (beta2) -- (beta4);
\draw[draw=black,\thickness,opacity=\lineop] (beta3) -- (beta4);

% 1-PS coordinates
\coordinate (lambda31) at ($(beta1)!-\pslength!(beta3)$);
\coordinate (lambda32) at ($(beta3)!-\pslength!(beta1)$);
\coordinate (lambda11) at ($(beta1)!-\pslength!(beta56)$);
\coordinate (lambda12) at ($(beta56)!-\pslength!(beta1)$);
\coordinate (lambda51) at ($(beta2)!-\pslength!(beta4)$);
\coordinate (lambda52) at ($(beta4)!-\pslength!(beta2)$);
\coordinate (lambda21) at ($(beta4)!-\pslength!(beta3)$);
\coordinate (lambda22) at ($(beta3)!-\pslength!(beta4)$);
\coordinate (lambda41) at ($(beta2)!-\pslength!(beta56)$);
\coordinate (lambda42) at ($(beta56)!-\pslength!(beta2)$);
% 1-PS
\draw[draw=red,\dashsty, <->,thick] (lambda11) -- (lambda12);
\draw[draw=red,\dashsty, <->,thick] (lambda21) -- (lambda22);
\draw[draw=red,\dashsty, <->,thick] (lambda31) -- (lambda32);
\draw[draw=red,\dashsty, <->,thick] (lambda41) -- (lambda42);
\draw[draw=red,\dashsty, <->,thick] (lambda51) -- (lambda52);

\node[draw=none,fill=none] at ($(beta1)+(-0.5,-0.3)$) {\small $\beta_1$};
\node[draw=none,fill=none] at ($(beta2)+(-0.5,-0.2)$) {\small $\beta_2$};
\node[draw=none,fill=none] at ($(beta3)+(-0.4,0)$) {\small $\beta_3$};
\node[draw=none,fill=none] at ($(beta4)+(-0.4,-0.3)$) {\small $\beta_4$};
\node[draw=none,fill=none] at ($(beta56)+(0.5,0.4)$) {\small $\beta_5,\beta_6$};

\node[draw=none,fill=none,red] at ($(beta2)!0.5!(beta4)!0.45!(center)$) {$\Gamma$};

\node[draw=none,fill=none,red] at ($(lambda11)+(-0.3,0)$) {\small $\lambda_1$};
\node[draw=none,fill=none,red] at ($(lambda21)+(0,-0.3)$) {\small $\lambda_2$};
\node[draw=none,fill=none,red] at ($(lambda31)+(-0.3,0)$) {\small $\lambda_3$};
\node[draw=none,fill=none,red] at ($(lambda41)+(-0.3,0.2)$) {\small $\lambda_4$};
\node[draw=none,fill=none,red] at ($(lambda51)+(-0.2,0.3)$) {\small $\lambda_5$};

\end{tikzpicture}
    \caption{The secondary fan $\Sigma_{GKZ}$ in $\wh{G}_\R\cong \R^3$ of the GIT problem described in \S\ref{subsection:arrangeable_with_bad_stratification}. The red dashed lines correspond to the hyperplanes of $\lambda_1,\hdots, \lambda_5$.}
    \label{fig:GKZ_for_arrangeable_bad_stratification}
\end{figure}
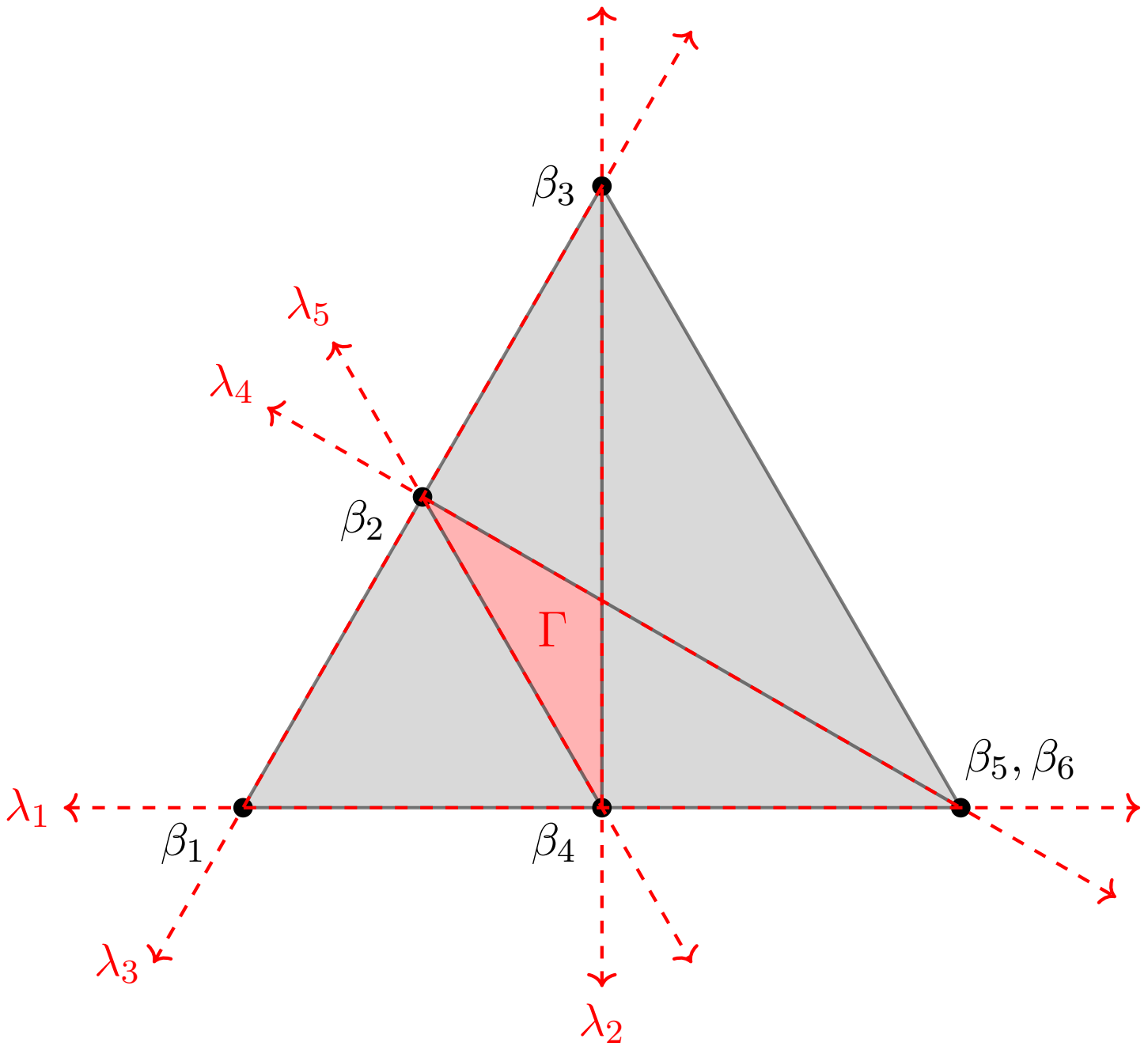
Let $\Sigma_{GKZ}$ be the secondary fan associated to these weights. The secondary fan is pictured in \Cref{fig:GKZ_for_arrangeable_bad_stratification} and the highlighted chamber $\Gamma$ is the chamber such that $\CX_{\Gamma} = X_\Sigma$. As $\Sigma_{GKZ}$ is a projective rank 3 GIT problem where no weights $\beta_i$ lie in $\Int |\Sigma_{GKZ}|$, we have that $\Sigma_{GKZ}$ is arrangeable (see \Cref{rmk:all_rank_2_problems_are_arrangeable}). 

We define the following five 1-PSs in $\wh{G}^\vee \cong \Z^3$:
\[\begin{array}{ccc} \lambda_1 = (0,1,0), & \lambda_2 = (1,0,-1), & \lambda_3 = (0,0,1), \\ \lambda_4 = (1,-1,0), & \lambda_5 = (-1,1,1).
\end{array}\]
The hyperplanes of these 1-PSs are drawn as dashed red lines in \Cref{fig:GKZ_for_arrangeable_bad_stratification}. By directly checking the conditions of \Cref{lemma:HKKN_strat_combo}, one can verify that $\bs{\lambda}\coloneqq \lambda_1,\hdots, \lambda_5$ is an HKKN stratification of $B_{\Gamma}$. Observe that each $\lambda_i\in \Gamma^\vee \cap \Lambda$, however $\lambda_1,\hdots,\lambda_5$ is not a totalization of $(\Gamma^\vee\cap \Lambda, \wrel)$. This is because the 1-PS associated to the wall $\beta_{\{3,5,6\}}$ is not in the ordered list $\bs{\lambda}$. Moreover, $\lambda_3\in |\Sigma_{GKZ}|^\vee$ so $\lambda_3 \wrel \lambda_2$.

The proof of \Cref{prop:arrangeable_with_bad_wall_strat} will proceed by computing the $\lambda_2$ weights of $F_\bullet^D(D)$ along $(X_2)^{\lambda_2}$ using \Cref{lemma:lambda_weights_of_cox_generators}, and then demonstrating that this weight range is wider that $T_{\lambda_1}^{>0}$. Let us, then, identify an ample lift $D = (0,0,2,3,0,-1)$. One checks that
\begin{align*}
    (X_2)^{\lambda_2} &= \lp \A^6\setminus V(I_{\lambda_1}^{>0})\rp^{\lambda_2}\\
    &= V(I_{\lambda_2}^{\ne 0})\setminus V(x_2,x_3)\\
    &= V(x_1,x_2,x_5,x_6)\setminus V(x_2,x_3)\\
    &= T\cdot p_{\{3\}}\sqcup T\cdot p_{\{3,4\}}.
\end{align*}
So, we can consider the case of $L = \{3\}$. 
We now define two 1-cells of $P_D$:
\begin{align*}
    \delta&\coloneqq \lp \{0\}\times \{2\}\times \{-2\}\times \{-2\} \times (0,1)\times (1,2)\rp \cap P_D,\\
    \tau&\coloneqq  \lp (0,1)\times (0,1)\times (-1,0)\times \{-1\}\times \{0\}\times \{1\}\rp\cap P_D.
\end{align*}
One can check that these are both indeed well-defined nonempty cells as $\delta$ contains the point $(0,2,-2,-2,1/2,3/2)$ and $\tau$ contains the point $(1/2,1/2,-1/2,-1,0,1)$.
\begin{proof}[Proof of \Cref{prop:arrangeable_with_bad_wall_strat}]
Let us continue from the above discussion. We will show that for $L = \{3\}$, both $\del \CS_{\delta,L}$ and $\del \CS_{\tau,L}$ are not acyclic. Let us assume for a moment that we have shown this. Observe that
\begin{align*}
    \wh{\lambda}_2(\pi_\R(\ceil{\delta})) &= -1,\\
    \wh{\lambda}_2(\pi_\R(\ceil{\tau})) &= 1.
\end{align*}
However, $T_{\lambda_2}^{>0} = \wh{\lambda}(\beta_1+\beta_2) = 2$. Thus there is no weight $w_2$ so that $[w_2,w_2+T_{\lambda_2}^{>0})$ contains each of the non-acyclic weights. So by \Cref{prop:cox_generator_grw_condition}, for any choice of weights $\bs{w}$, $F_\bullet^D(D)\notin W_{\bs{\lambda}}^{\bs{w}}$. As these complexes $F_\bullet^D(D)$ generate $W_\Gamma^{\Cox}$ by \Cref{prop:amplesgenerate}, we conclude that for any choice of weights $\bs{w}$, $W_{\bs{\lambda}}^{\bs{w}}\ne W_\Gamma^{\Cox}$. 

We will now show directly that both $\del \CS_{\delta,L}$ and $\del\CS_{\tau,L}$ are non-acyclic, completing the proof. By definition,
\begin{equation}\label{eqn:counter_example_2_S_delta_L}
\begin{aligned}
    \CS_{\delta, L} &= \lp \{0\}\times (1,2]\times (-2,\infty)\times (-3,-2] \times (0,1]\times (1,2]\rp \cap i_\R(M_\R)\\
    &= \lp \{0\}\times (1,2]\times (-2,\infty)\times (-3,-2] \times (0,1]\times (1,2]\rp \cap \ker(\pi_\R).
\end{aligned}
\end{equation}
We claim $\del \CS_{\delta,L}$ is homeomorphic to the union of two disjoint points. In particular, we will show that
\[\CS_{\delta,L} = \{(0,2,-2,-2,2-t,t) \ | \ t\in (1,2)\}.\]
One inclusion is immediate by computing $\pi_\R$ of the right side. To show the opposite inclusion, take $\bs{x}\in \CS_{\delta,L}$. That $\pi_\R(\bs{x}) = 0$ gives three equalities:
\begin{align*}
    x_1+x_2+x_4 &= 0,\\
    x_2+x_3 &= 0,\\
    x_4+x_5+x_6 &= 0.
\end{align*}
By eqn. \eqref{eqn:counter_example_2_S_delta_L} $x_1 = 0$, so the first equality implies $x_2 = -x_4$. Then, as $x_2\in (1,2]$ and $x_4\in (-3,-2]$ by eqn. \eqref{eqn:counter_example_2_S_delta_L}, these together imply $x_2 = 2$ and $x_4 = -2$. Then the second equality implies $x_3 = -2$. Now the third equality becomes $x_5 = 2-x_6$. As $x_5\in (0,1]$ eqn. \eqref{eqn:counter_example_2_S_delta_L}, this implies that $x_6\in [1,2)$. But eqn. \eqref{eqn:counter_example_2_S_delta_L} requires $x_6\in (1,2]$, so we get that $x_6\in (1,2)$ and conclude the desired equality. In particular we conclude that $\CS_{\delta,L}\cong (0,1)\subset \R$, so $\del \CS_{\delta,L}$ is homeomorphic to a union of two disjoint points and is thus not acyclic. By an exactly parallel method, one computes that
\begin{align*}
    \CS_{\tau,L} &= \lp (0,1]\times (0,1]\times (-\infty,\infty)\times (-2,-1]\times\{0\}\times \{1\}\rp\cap i_\R(M_\R)\\
    &= \{ (1-t, t,-t,-1,0,1) \ | \ t\in (0,1)\} \cong (0,1)\subset \R.
\end{align*}
Thus $\del\CS_{\tau,L}$ is homeomorphic to a disjoint union of two points and is thus non-acyclic. 
\end{proof}
The stratification used in \Cref{prop:arrangeable_with_bad_wall_strat} results in a grade restriction window distinct from the Cox window.  On the other hand, because the associated GIT problem is arrangeable, we know of stratifications where the windows agree by \Cref{thm:window_equality}.  We now compare these two stratifications $\bs{\lambda}^1, \bs{\lambda}^3$ along with two intermediate stratifications $\bs{\lambda}^2,\bs{\lambda}^4$ and find that slight variations in HKKN stratifications result in unpredictable behavior.
\begin{corollary}\label{cor:may_or_may_not} There exists a chamber $\Gamma$ of a projective toric GIT problem with four distinct HKKN stratifications $\bs{\lambda}^1,\hdots, \bs{\lambda}^4$ of $V(B_\Gamma)$ such that the following hold:
\begin{enumerate}[(a),itemsep=0.5cm]
    \item\label{cor:may_or_may_not_reorder_MAY} The stratifications $\bs{\lambda}^2,\bs{\lambda}^3$ are equal as unordered sets and $W_{\bs{\lambda}^2} = W_{\bs{\lambda}^3}$.
    \item\label{cor:may_or_may_not_reorder_NOT} The stratifications $\bs{\lambda}^1,\bs{\lambda}^4$ are equal as unordered sets and $W_{\bs{\lambda}^1} \ne W_{\bs{\lambda}^4}$.
    \item\label{cor:may_or_may_not_sub_MAY} The stratification $\bs{\lambda}^4$ is an ordered subset of $\bs{\lambda}^3$ and $W_{\bs{\lambda}^3} = W_{\bs{\lambda}^4}$.
    \item\label{cor:may_or_may_not_sub_NOT} The stratification $\bs{\lambda}^1$ is an ordered subset of $\bs{\lambda}^2$ and $W_{\bs{\lambda}^1} \ne W_{\bs{\lambda}^2}$.
\end{enumerate}
\end{corollary}
In the remainder of this section we will prove \Cref{cor:may_or_may_not}, continuing from the setup of the proof of \Cref{prop:arrangeable_with_bad_wall_strat}. We additionally define $\lambda_0$ to be the primitive 1-PS associated to $\beta_{\{3,5,6\}}$ oriented so that $\Gamma\subset H_{\lambda_0}^{\ge 0}$. We now define the four HKKN stratifications appearing in \Cref{cor:may_or_may_not}.
\begin{align*}
    \bs{\lambda}^1 &= \quad\: \:\lambda_1,\lambda_2,\lambda_3,\lambda_4,\lambda_5\\
    \bs{\lambda}^2 &= \lambda_0,\lambda_1,\lambda_2,\lambda_3,\lambda_4,\lambda_5\\
    \bs{\lambda}^3 &= \lambda_0,\lambda_1,\lambda_3,\lambda_2,\lambda_4,\lambda_5\\
    \bs{\lambda}^4 &= \quad\:\:\lambda_1,\lambda_3,\lambda_2,\lambda_4,\lambda_5
\end{align*}
By checking the conditions of \Cref{lemma:HKKN_strat_combo}, one can show that each $\bs{\lambda}^k$ indeed induces an HKKN stratification of $V(B_\Gamma)$. We note that $\bs{\lambda}^1$ was the stratification used in the proof of \Cref{prop:arrangeable_with_bad_wall_strat}, so we have $W_{\bs{\lambda}^1}\ne W_\Gamma^\Cox$. Additionally, observe that $\bs{\lambda}^3$ is a totalization of $(\Gamma^\vee \cap \Lambda, \wrel)$, so $W_{\bs{\lambda}^3} = W_\Gamma^\Cox$ by our main result, \Cref{thm:window_equality}. 

The internal walls in this GIT problem $\Lambda^{\circ}$ are $\pm\{\lambda_2,\lambda_4,\lambda_5\}$. So, for $i=2,4,5$ and $k=1,\hdots, 4$, denote by $(X^k)_{\lambda_i}^{=0}$ the fixed locus of $\lambda_i$ minus the strata from $\bs{\lambda}^k$ occurring before $\lambda_i$. That is,
\begin{equation}\label{eqn:fixed_loci_varied_stratifications}
(X^k)_{\lambda_i}^{=0} = V(I_{\lambda_i}^{\ne 0})\setminus \bigcup_{\lambda^k_j < \lambda_i}V(I_{\lambda^k_j}^{>0}).
\end{equation}
Recall then that each $(X^k)_{\lambda_i}^{=0}$ is a disjoint union of orbits $T\cdot p_L$. By directly computing each of the above set differences, we find each $L$ such that $T\cdot p_L\in (X^k)_{\lambda_i}^{=0}$ and summarize these in \Cref{table:orbits_for_may_or_may_not}. With this classification of orbits in hand, we can now prove \Cref{cor:may_or_may_not}.
\begingroup
\renewcommand{\arraystretch}{1.5}
\begin{table}[h]
    \centering
    \begin{tabular}{c|c|c|c}
         & $\lambda_2$ & $\lambda_4$ & $\lambda_5$ \\ \hline
        $\bs{\lambda}^1$ & $\{3,4\}, \{3\}$ & $\{2,5\}, \{2,5,6\}, \{2,6\}$ & $\{2,4\}$ \\ \hline
        $\bs{\lambda}^2$ & $\{3,4\}$ & $\{2,5\}, \{2,5,6\}, \{2,6\}$ & $\{2,4\}$ \\ \hline
        $\bs{\lambda}^3$ & $\{3,4\}$ & $\{2,5\}, \{2,5,6\}, \{2,6\}$ & $\{2,4\}$  \\ \hline
        $\bs{\lambda}^4$ & $\{3,4\}$ & $\{2,5\}, \{2,5,6\}, \{2,6\}$ & $\{2,4\}$  \\ \hline
    \end{tabular}
    \vspace{0.5cm}
    \caption{The sets $L$ for each orbit $p_L$ of the fixed locus of $\lambda_i$ in the stratification $\bs{\lambda}^k$. In particular, the list in the row labeled by $\bs{\lambda}^k$ and the column labeled by $\lambda_i$ correspond to the list of possible orbits $T\cdot p_L$ occurring in $(X^k)_{\lambda_i}^{=0}$ (see eqn. \eqref{eqn:fixed_loci_varied_stratifications}). The appearance of the additional orbit $L = \{3\}$ for $\bs{\lambda}^1$ is exactly what allows one to obtain non-acyclic regions $\CS_{\delta,L}$ in \Cref{prop:arrangeable_with_bad_wall_strat}.}
    \label{table:orbits_for_may_or_may_not}
\end{table}
\endgroup
\begin{proof}[Proof of \Cref{cor:may_or_may_not}] As mentioned in the above discussion, it was proven in \Cref{prop:arrangeable_with_bad_wall_strat} that $W_{\bs{\lambda}^1} \ne W_\Gamma^\Cox$ and by our main result, \Cref{thm:window_equality}, $W_{\bs{\lambda}^3} = W_\Gamma^\Cox$.  By construction,
\[\bs{\lambda}^2\setminus |\Sigma_{GKZ}|^\vee = \bs{\lambda}^3\setminus |\Sigma_{GKZ}|^\vee = \bs{\lambda}^4\setminus |\Sigma_{GKZ}|^\vee = \{\lambda_2,\lambda_4,\lambda_5\}.\]
By the description of the orbits of the fixed locus given in \Cref{table:orbits_for_may_or_may_not}, we conclude that for each $i = 2,4,5$ and any $k,k'\in \{2,3, 4\}$, that
\[(X^k)_{\lambda_i}^{=0} = (X^{k'})_{\lambda_i}^{=0},\]
(see \eqref{eqn:fixed_loci_varied_stratifications}). These observations, together with $W_{\bs{\lambda}^3}= W_\Gamma^\Cox$, imply
\[W_{\bs{\lambda}^2} = W_{\bs{\lambda}^3} = W_{\bs{\lambda}^4} = W_{\Gamma}^\Cox\]
by \Cref{cor:two_stratifications_with_same_orbits_one_Cox}.
\end{proof}

\appendix
\section{Macaulay2 Code}
In this appendix we provide the code that is used to verify \Cref{thm:no_stratification_making_windows_equal_exist}. In \S\ref{subsection:a_convex_condition} we provide a computationally tractable sufficient condition equivalent to the computationally intractable sufficient condition of \Cref{prop:non_convex_sufficient_condition}. In \S\ref{subsection:appendix_methods} we provide the \lst{Macaulay2} functions we wrote to compute nerves and to check this convex sufficient condition for window non-equality. Additionally, in \S\ref{subsection:appendix_methods} we provide the code used to implement these functions in our particular counterexample.

\subsection{A convex sufficient condition for window non-equality}\label{subsection:a_convex_condition}
Recall our usual setup. Namely, we start with a toric projective GIT problem defined by $G\subset \Gm^n$ with secondary fan $\Sigma_{GKZ}$, a chamber $\Gamma$, a lift $D$ of $d\in \Gamma \cap \wh{G}$ under $\pi$ from eqn. \eqref{eqn:ses_for_lattices}. In \Cref{prop:non_convex_sufficient_condition}, we showed that if $I^{>0}$ is such that $\la I^{>0}\ra$ is a component of the primary decomposition of $B_\Gamma$, then $W_\Gamma^{\Cox}$ is not a grade restriction window if for every $I^{=0}\subset [n]\setminus I^{>0}$,
\begin{equation}\label{eqn:non_convex_condition}\tag{$\ast$}
C_{I^{>0},I^{=0}}\subset \bigcup_{(\bs{a},\bs{a}')\in N_{I^{=0}}\times N_{I^{=0}}} H_{\pi(\bs{a}-\bs{a}')-\sum_{i\in I^{>0}}\beta_i}^{\ge 0}.
\end{equation}
While abstractly nice, this condition is computationally intractable in \lst{Macaulay2} for two reasons. The first is that the regions $C_{I^{>0},I^{=0}}$ are not themselves cones, but are instead partially open subsets of their closures which are indeed cones. The second is that a union of half-spaces is generally non-convex. Thus it is not straight-forward computationally to create the region on the right-hand side of \eqref{eqn:non_convex_condition} as an object in \lst{Macaulay2}. Instead, in the remainder of this subsection we develop a convex sufficient condition equivalent to the fundamentally non-convex condition \eqref{eqn:non_convex_condition} such that we can implement it in \lst{Macaulay2} in \S\ref{subsection:appendix_methods}. We first state the following elementary lemma without proof for the sole purpose of justifying our computation of $\overline{C}_{I^{>0},I^{=0}}$ as a \lst{Cone} in \S\ref{subsection:appendix_methods}.
\begin{lemma}\label{lemma:closure_of_C_plus_zero_cone} When $C_{I^{>0},I^{=0}}\ne\emptyset$, its closure $\overline{C}_{I^{>0},I^{=0}}$, is a cone in $\wh{G}_\R^\ast$ with the following half-space and dual cone descriptions:
\begin{equation}\label{eqn:lambda_cone_definition}
\begin{aligned}
    \overline{C}_{I^{>0},I^{=0}} &\coloneqq \bigcap_{i\in I^{>0}} H_{\beta_i}^{\ge 0} \cap \bigcap_{i\in I^{=0}} H_{\beta_i}^{=0} \cap \bigcap_{i\in [n]\setminus (I^{>0}\sqcup I^{=0})}H_{\beta_i}^{\le 0}\\
    &= \lp \beta_{I^{>0}\sqcup I^{=0}} - \beta_{[n]\setminus I^{>0}}\rp^\vee.
\end{aligned}
\end{equation}
\end{lemma}
We now decompose $\overline{C}_{I^{>0},I^{=0}}$ into a finite number of cones which control the width of $\wh{\lambda}(\pi_\R(N_{I_{\lambda}^{=0}}))$.
\begin{definition}\label{def:C_max_min_cones} Let $I^{>0},I^{=0}\subset [n]$ be disjoint sets such that $C_{I^{>0},I^{=0}}$ is nonempty. Then for any $\bs{a}\in N\coloneqq N_{I^{=0}}$ we define the following cones in $\wh{G}_\R^{\ast}$:
\begin{equation*}
\begin{aligned}
    C_{\bs{a}}^{\op{max}} &\coloneqq \overline{C}_{I^{>0},I^{=0}} \cap \bigcap_{\bs{b}\in N} H_{\pi(\bs{a}-\bs{b})}^{\ge 0},\\
    C_{\bs{a}}^{\op{min}} &\coloneqq \overline{C}_{I^{>0},I^{=0}} \cap \bigcap_{\bs{b}\in N} H_{\pi(\bs{a}-\bs{b})}^{\le 0}.
\end{aligned}
\end{equation*}
\end{definition}
The cones $C_{\bs{a}}^{\omax},C_{\bs{a}}^{\omin}$ are adorned with $\omax,\omin$ respectively in light of the following fact.
\begin{lemma}\label{lemma:max_min_cones} The cone $C_{\bs{a}}^{\op{max}}$ is precisely the cone of functionals $\lambda\in  \overline{C}_{I^{>0},I^{=0}}$ such that \[\max \wh{\lambda}(\pi(N)) = \wh{\lambda}(\pi(\bs{a})).\]
Similarly, the cone $C_{\bs{a}}^{\omin}$ is precisely the cone of functionals $\lambda\in \overline{C}_{I^{>0},I^{=0}}$ such that
\[\min \wh{\lambda}(\pi(N)) = \wh{\lambda}(\pi(\bs{a})).\]
As a consequence, for any $\lambda\in C_{\bs{a}}^{\omax}\cap C_{\bs{a}'}^{\omin}$, we have
\[\max_{(\bs{b},\bs{b}')\in N\times N}\wh{\lambda}(\pi_\R(\bs{b}-\bs{b}')) = \wh{\lambda}(\pi_\R(\bs{a}-\bs{a}')).\]
\end{lemma}
\begin{proof}
By construction,
\begin{align*}
    C_{\bs{a}}^{\omax} &= \{\lambda\in C_{I^{>0},I^{=0}} \ | \ \lambda\in H_{\pi(\bs{a}-\bs{b})}^{\ge 0} \text{ for all }\bs{b}\in N\}\\
    &= \{\lambda\in C_{I^{>0},I^{=0}} \ | \ \wh{\lambda}(\pi(\bs{a})) \ge \wh{\lambda}(\pi(\bs{b}))\text{ for all }\bs{b}\in N\}\\
    &= \{\lambda\in C_{I^{>0},I^{=0}} \ | \ \wh{\lambda}(\pi(\bs{a})) = \max \wh{\lambda}(\pi(N))\}.
\end{align*}
The claim for $C_{\bs{a}}^{\omin}$ follows similarly.
\end{proof}
With the above, we can now prove that the sufficient condition \eqref{eqn:non_convex_condition} is equivalent to a finite number of convex conditions.
\begin{lemma}\label{lemma:equivalent_convex_sufficient_condition} For any disjoint $I^{>0},I^{=0}$ such that $C_{I^{>0},I^{=0}}$ is nonempty, 
\begin{equation}\label{eqn:non_convex_sufficient_condition}\tag{$\ast$} C_{I^{>0},I^{=0}}\subset \bigcup_{(\bs{a},\bs{a}')\in N_{I^{=0}}\times N_{I^{=0}}} H_{\pi(\bs{a}-\bs{a}')-\sum_{i\in I^{>0}}\beta_i}^{\ge 0}
\end{equation}
if and only if for every $\bs{a},\bs{a}'\in N_{I^{=0}}$,
\begin{equation}\label{eqn:convex_sufficient_condition}\tag{$\ast\ast$} 
C_{\bs{a}}^{\omax}\cap C_{\bs{a}'}^{\omin}\subset H_{\pi(\bs{a}-\bs{a}')-\sum_{i\in I^{>0}}\beta_i}^{\ge 0}.
\end{equation}
\end{lemma}
\begin{proof} To show \eqref{eqn:convex_sufficient_condition} implies \eqref{eqn:non_convex_sufficient_condition}, take any $\lambda\in C_{I^{>0},I^{=0}}$. By \Cref{lemma:max_min_cones}, by taking a $\bs{a}$ among $\bs{b}\in N_{I^{=0}}$ which maximizes $\wh{\lambda}\circ \pi_\R$, we have that $\lambda\in C_{\bs{a}}^{\omax}$. Similarly there is a $\bs{a}'\in N_{I^{=0}}$ so that $\lambda\in C_{\bs{a}'}^{\omin}$, so $\lambda\in C_{\bs{a}}^{\omax}\cap C_{\bs{a}'}^{\omin}$. By our assumption of \eqref{eqn:convex_sufficient_condition}, it follows that $\lambda$ is contained in the union in \eqref{eqn:non_convex_sufficient_condition}. We conclude that \eqref{eqn:non_convex_sufficient_condition} holds.

To show that \eqref{eqn:non_convex_sufficient_condition} implies \eqref{eqn:convex_sufficient_condition}, first note that as the right-hand side of \eqref{eqn:non_convex_sufficient_condition} is closed, taking closures of both sides yields
\begin{equation}\label{eqn:non_convex_sufficient_condition_CLOSED} \tag{$\overline{\ast}$} \overline{C}_{I^{>0},I^{=0}}\subset \bigcup_{(\bs{a},\bs{a}')\in N_{I^{=0}}\times N_{I^{=0}}} H_{\pi(\bs{a}-\bs{a}')-\sum_{i\in I^{>0}}\beta_i}^{\ge 0}.
\end{equation}
Now take any $\bs{a},\bs{a}'\in N_{I^{=0}}$ and any $\lambda\in C_{\bs{a}}^{\omax}\cap C_{\bs{a}'}^{\omin}$. As $C_{\bs{a}}^{\omax}\cap C_{\bs{a}'}^{\omin}\subset \overline{C}_{I^{>0},I^{=0}}$, \eqref{eqn:non_convex_sufficient_condition_CLOSED} implies that there exists a pair $(\bs{c},\bs{c}')\in N_{I^{=0}}\times N_{I^{=0}}$ such that
\[\wh{\lambda}(\pi(\bs{c}-\bs{c}')) \ge \sum_{i\in I^{>0}}\wh{\lambda}(\beta_i).\]
However, as $\lambda\in C_{\bs{a}}^{\omax}\cap C_{\bs{a}'}^{\omin}$, \Cref{lemma:max_min_cones} implies
\[\max_{(\bs{b},\bs{b}')\in N_{I^{=0}}\times N_{I^{=0}}} \wh{\lambda}(\pi(\bs{b}-\bs{b}')) = \wh{\lambda}(\pi(\bs{a}-\bs{a}')).\]
Together, we conclude
\[\wh{\lambda}(\pi(\bs{a}-\bs{a}')) \ge \sum_{i\in I^{>0}}\wh{\lambda}(\beta_i).\]
as desired.
\end{proof}
The advantage of the sufficient condition \eqref{eqn:convex_sufficient_condition} is that it amounts to checking a finite number of inclusions of cones, which is readily implementable in \lst{Macaulay2}. We include this implementation in the following subsection. For the sake of completeness, we restate \Cref{prop:non_convex_sufficient_condition} using \Cref{lemma:equivalent_convex_sufficient_condition}.
\begin{corollary}\label{cor:convex_sufficient_condition} Suppose $\la I^{>0}\ra$ is a component of the primary decomposition of $B_\Gamma$. Suppose that for every $I^{=0}\subset [n]\setminus I^{>0}$ such that $C_{I^{>0},I^{=0}}\ne\emptyset$, and every $\bs{a},\bs{a}'\in N_{I^{=0}}$,
\[C_{\bs{a}}^{\omax}\cap C_{\bs{a}',\omin}\subset  H_{\pi(\bs{a}-\bs{a}')-\sum_{i\in I^{>0}}\beta_i}^{\ge 0}.\]
Then there is no combination of HKKN stratification $\bs{\lambda}$ of $V(B_\Gamma)$ and weights $\bs{w}$ such that $W_{\Gamma}^{\Cox} = W_{\bs{\lambda}}^{\bs{w}}$. 
\end{corollary}

\subsection{Methods}\label{subsection:appendix_methods} All code in this subsection can be found at \href{https://github.com/dashstevanovich/Windows}{github.com/dashstevanovich/Windows}. We start with a toric projective GIT problem defined by $G\subset \Gm^n$ which is encoded by the map $\pi: \Z^n \to \wh{G}$. In the below code, \lstinline{classMap} is a $\dim \wh{G}_\R$ by $n$ \lstinline{Matrix} of \lstinline{ZZ} equal to the map $\pi_\R$ from eqn. \eqref{eqn:ses}. We assume a chamber $\Gamma$ in the secondary fan has been fixed and that $D \in \Z^n$ is a lift of $d\in \Gamma \cap \wh{G}$. In the below code, \lstinline{divisorVector} is a $n$ by 1 \lstinline{Matrix} of \lstinline{ZZ} equal to $D$.
\begin{lstlisting}
getDivisorPolytope = method()
getDivisorPolytope(Matrix, Matrix) := (classMap, divisorVector) -> (
    zeroVector := matrix(mutableMatrix(ZZ,numrows classMap,1)); 
    imageOfM := affinePreimage(classMap,convexHull(zeroVector)); 
    positiveOctant := polyhedronFromHData(-1*map(ZZ^(numcols classMap)), divisorVector); 
    return intersection(imageOfM,positiveOctant)
    );
\end{lstlisting}
The method \lstinline{getDivisorPolytope} takes as inputs the \lstinline{Matrix classMap} and the \lstinline{Matrix divisorVector} and returns a \lstinline{Polyhedron} which is equal to $P_D\subset \R^n$ as described in \Cref{def:divisor_polytope}. This is done by realizing $M_\R$ as a \lstinline{Polyhedron} (line 4), realizing $\R^n_{\ge -D}$ as a \lstinline{Polyhedron} (line 5), and then taking their intersection (line 6).
\begin{lstlisting}
slicePolytopeByHyperplanes = method()
slicePolytopeByHyperplanes(Polyhedron,Matrix) := (P,f) -> (
    fImage := affineImage(f, P);
    intervalBounds := (entries vertices fImage)_0;
    minValue := min(apply(intervalBounds, vertex -> floor(vertex)));
    maxValue := max(apply(intervalBounds, vertex -> ceiling(vertex)));
    slicesInHyperplane := for i from minValue to maxValue list (
        hyperplane := polyhedronFromHData(f || -f, matrix({{i},{-i}}));
        sliceInHyperplane := intersection(P,hyperplane);
        if dim sliceInHyperplane == -1 then continue else sliceInHyperplane
        );
    slicesInBetween := for i from minValue to maxValue-1 list (
        region := polyhedronFromHData(f || -f, matrix({{i+1},{-i}}));
        intersection(P,region)
        );
    return join(slicesInHyperplane, slicesInBetween)
    );
\end{lstlisting}
The method \lst{slicePolytopeByHyperplane} takes as its input the \lst{Polyhedron P} representing a bounded polytope $P$ in $\R^n$ and a $1$ by $n$ \lst{Matrix f} of \lst{ZZ} representing $f:\R^n\to \R$ and returns a \lst{List} of \lst{Polyhedron} objects corresponding to all non-empty polytopes of the type $P\cap f^{-1}(\{i\})$ and $P\cap f^{-1}([i,i+1])$ for $i\in \Z$. We do this by computing the least integer upper bound (line 5) and greatest integer lower bound (line 6) of $f(P)$. Then for each integer $i$ between these values, we compute the associated slice $P\cap f^{-1}(\{i\})$ and list the non-empty slices (lines 8-10). We then similarly compute the non-empty regions $P\cap f^{-1}([i,i+1])$ for $i+1$ greater that the greatest integer lower bound and $i$ less than the least integer upper bound (lines 13, 14). 
\begin{lstlisting}
collectUnique = method()
collectUnique(List) := (inList) -> (
    outList := {};
    scan(inList, C -> if not any(outList, Q -> Q == C) then outList = append(outList, C));
    return outList
    );
\end{lstlisting}
The method \lst{collectUnique} takes as its input the \lst{List inList} of \lst{Polyhedron} objects and returns the \lst{List} of unique polyhedra occurring in \lst{inList}.
\begin{lstlisting}
getClosedCells = method()
getClosedCells(Polyhedron,List) := (P, fList) -> (
    allCells = {P};
    for f in fList do (
        cellsByf := for Q in allCells list (
            slicePolytopeByHyperplanes(Q,f)
            );
        allCells = collectUnique(flatten(cellsByf));
        );
    return allCells
    );
\end{lstlisting}
The method \lst{getClosedCells} takes as its input the \lst{Polyhedron P} representing a bounded polytope $P\subset \R^n$ and the \lst{List fList} representing a list of functionals $\{f_1,\hdots, f_k\}$, $f_i:\R^n\to \R$ and returns a \lst{List} of \lst{Polyhedron} objects representing all non-empty closed cells of the type
\[P\cap \bigcap_{i\in I}f_i^{-1}(\{a_i\}) \cap \bigcap_{j\in [n]\setminus I}f_j^{-1}([a_j,a_j+1]).\]
This is done by iterative application of \lst{slicePolytopeByHyperplanes}. In particular, when \lst{fList} is equal to the list of coordinate functionals $\{e_1^*,\hdots, e_n^*\}\in (\R^n)^*$ and \lst{P} is $P_D$, the output of \lst{getClosedCells} is exactly
\[\left\{\overline{\delta} \ | \ \delta\in \cells(P_D)\right\}.\]
See \S\ref{subsection:cox_generator_images} for a detailed description of $\cells(P_D)$. This method of computing cells by iteratively slicing a polytope by hyperplanes is essentially the same as which appears in Jay Yang's HHL package \cite{JayHHL}.
\begin{lstlisting}
getCellCeiling = method()
getCellCeiling(Polyhedron) := (P) -> (
    x := interiorPoint(P);
    xEntries := (entries transpose x)_0;
    xEntriesCeiling := apply(xEntries, xi -> ceiling(xi));
    return transpose matrix({xEntriesCeiling})
    );
\end{lstlisting}
The method \lst{getCellCeiling} takes as its input a \lst{Polyhedron P} representing a closed polytope $P\subset \R^n$ and returns a $n$ by 1 \lst{Matrix} representing the componentwise ceiling $\ceil{x}$ of some point $x\in \Relint P$. In particular, when the input polytope is $\overline{\delta}$ for some $\delta\in \cells(P_D)$, the output of \lst{getCellCeiling} is exactly $\ceil{\delta}$ (see \S\ref{subsection:cox_generator_images}).
\begin{lstlisting}
getNerve = method()
getNerve(Matrix, List, Polyhedron) := (cellCeiling, L, divisorPolytope) -> (
    n := ambDim divisorPolytope;
    Lcomp := select(toList(0..n-1), i -> not isMember(i,L));
    identityMatrix := map(ZZ^n);
    cellCeilingMinusOne := cellCeiling - matrix(table(n,1,(i,j) -> 1));
    rectangularRegionInequalities := submatrix(identityMatrix,Lcomp,) || -1*submatrix(identityMatrix,Lcomp,);
    rectangularRegionValues :=  submatrix(cellCeiling,Lcomp,) || -1*submatrix(cellCeilingMinusOne,Lcomp,);
    SdeltaLClosed := intersection(polyhedronFromHData(rectangularRegionInequalities,rectangularRegionValues),divisorPolytope);
    subsetsOfLcomp := subsets(Lcomp);
    lowerFacetRegions := new HashTable from apply(subsetsOfLcomp, J -> J => polyhedronFromHData(submatrix(identityMatrix,J,) || -1*submatrix(identityMatrix,J,), submatrix(cellCeilingMinusOne,J,) || -1*submatrix(cellCeilingMinusOne,J,)));
    return select(subsetsOfLcomp, J -> dim(intersection(SdeltaLClosed, lowerFacetRegions#J)) >= 0)
    );
\end{lstlisting}
The \lst{getNerve} method takes as its input a \lst{Matrix cellCeiling} corresponding to $\ceil{\delta}$ for some $\delta\in \cells(P_D)$, a \lst{List L} representing a subset $L\subset [n]$ which should be interpreted as associated to some distinguished point $p_L$ (see \S\ref{subsection:eq_class_splitting}), and a \lst{Polyhedron divisorPolytope} corresponding to $P_D$. The method returns a \lst{List} of \lst{List} objects representing subsets of $[n]\setminus L$ equal to the nerve $N(\CC_{\delta,L})$ (see \Cref{def:nerve}). This is done by first computing $\overline{\CS}_{\delta,L}$ (line 9) following \Cref{lemma:coordinate_description_of_closure_boundary}. Then, for each $J\subset [n]\setminus L$ we compute the \lst{Polyhedron}
\[C_J \coloneqq \prod_{j\in J}[\ceil{\delta}_i-1,\infty)\times\prod_{i\in [n]\setminus J} (-\infty,\infty)\]
and store $C_J$ in the \lst{HashTable lowerFacetRegions} with \lst{key} $J$ (line 11). Then we return each $J$ such that $\overline{\CS}_{\delta,L} \cap C_J \ne\emptyset$ (line 12). It is immediate that
\[\bigcap_{j\in J}A_j = \overline{\CS}_{\delta,L}\cap C_J.\]
where $A_j$ is as in \Cref{def:nerve}, so the output of \lst{getNerve} is indeed equal to $N(\CC_{\delta,L})$ as in \Cref{def:nerve}.
\begin{lstlisting}
loadPackage SimplicialComplexes
getNonacyclicLabels = method()
getNonacyclicLabels(Matrix,Polyhedron,List,List) := (classMap,divisorPolytope,allClosedCells,L) -> (
    n := numcols classMap;
    S := ZZ[x_0..x_(n-1)];
    SVars := vars S;
    allCellCeilings := unique(apply(allClosedCells, Q -> getCellCeiling(Q)));
    allNerves := apply(allCellCeilings, cellCeiling -> getNerve(cellCeiling,L,divisorPolytope));
    allNerveMonomials := apply(allNerves, nerve ->  apply(nerve, J -> promote(product(flatten entries(submatrix(SVars,J))),S)));
    allNerveComplexes := apply(allNerveMonomials, nerveMonomials -> simplicialComplex(nerveMonomials));
    nonacyclicLabelIndices := positions(allNerveComplexes, Delta -> HH Delta != 0);
    return apply(nonacyclicLabelIndices, i -> allCellCeilings_i)
    );
\end{lstlisting}
The method \lst{getNonacyclicLabels} takes as inputs an $n$ by $\dim \wh{G}_\R$ \lst{Matrix classMap} of \lst{ZZ} representing $\pi_\R$, a \lst{Polyhedron divisorPolytope} representing $P_D\subset \R^n$, a \lst{List allClosedCells} of \lst{Polyhedron} objects equal to $\{\overline{\delta} \ | \ \delta\in \cells(P_D)\}$, and a \lst{List L} corresponding to a distinguished point $p_L$ (see \S\ref{subsection:eq_class_splitting}). The output is a \lst{List} of $n$ by $1$ \lst{Matrix} of \lst{ZZ} objects representing
\begin{equation}\label{eqn:nonacyclic_labels_appendix}
    \{\ceil{\delta} \ | \ \delta\in \cells(P_D), \ \wt{H}_*(\del \CS_{\delta, L})\ne 0\}.
\end{equation}
This is done by computing $N(\CC_{\delta,L})$ for each $\delta\in \cells(P_D)$ using \lst{getNerve} (line 8). Then for each such nerve, using the \lst{SimplicialComplexes} package \cite{SimplicialComplexesArticle} we computed the \lst{SimplicialComplex} object associated to it by encoding faces $J\in N(\CC_{\delta,L})$ as monomials $\prod_{j\in J}x_j$ (line 10). Using the reduced homology function \lst{HH} in the \lst{SimplicialComplexes} package, we then record $\ceil{\delta}$ for each $\delta$ such that $\wt{H}_*(N(\CC_{\delta,L}))\ne 0$ (line 11). We note that $N(\CC_{\delta,L})\simeq \del \CS_{\delta,L}$ by the Nerve Lemma (\Cref{prop:nervecomputes}), so this computation indeed recovers the desired set in eqn. \eqref{eqn:nonacyclic_labels_appendix}.

Recall that the set in eqn. \eqref{eqn:nonacyclic_labels_appendix} directly determines whether or not $F_\bullet^D(D)$ satisfies a given grade restriction rule (\Cref{lemma:lambda_weights_of_cox_generators}, \Cref{prop:cox_generator_grw_condition}), so the above code can be used to determine exactly when exactly a Cox window generator $F_\bullet^D(D)$ lies in a given grade restriction window. This is done by checking the values of eqn. \eqref{eqn:nonacyclic_labels_appendix} under $\wh{\lambda}_i\circ \pi_\R$. 

We now move on to code used specifically in the proof of \Cref{thm:no_stratification_making_windows_equal_exist} to verify the sufficient condition established in \Cref{cor:convex_sufficient_condition}.
\begin{lstlisting}
getLambdaPlusMinusCone = method()
getLambdaPlusMinusCone(Matrix, List, List) := (classMap, Iplus, Izero) -> (
    n := numcols classMap;
    Iminus := select(toList(0..n-1), i -> not isMember(i, Iplus) and not isMember(i, Izero));
    greaterThanIneqs := transpose (submatrix(classMap, Iplus) | submatrix(classMap, Izero));
    lessThanIneqs := transpose (submatrix(classMap, Iminus) | submatrix(classMap, Izero));
    return coneFromHData(greaterThanIneqs || -1*lessThanIneqs)
    );
\end{lstlisting}
The method \lst{getLambdaPlusMinusCone} takes as inputs the \lst{Matrix classMap} representing $\pi_\R$, the \lst{List Iplus} representing $I^{>0}\subset [n]$, and the \lst{List Izero} representing $I^{=0}\subset [n]$ such that $C_{I^{>0},I^{=0}}\ne \emptyset$ (see \Cref{def:open_C_plus_zero_cone}). The output is a \lst{Cone} representing the cone $\overline{C}_{I^{>0},I^{=0}}$ from \Cref{lemma:closure_of_C_plus_zero_cone}. This is computed precisely via the half-space definition in \Cref{lemma:closure_of_C_plus_zero_cone}.
\begin{lstlisting}
weightsTooWide = method()
weightsTooWide(List,Cone,Matrix,List) := (nonacyclicImages,lambdaCone,classMap,Iplus) -> (
    CMaxCones := apply(nonacyclicImages, label1 -> intersection(lambdaCone, coneFromHData(matrix(apply(nonacyclicImages, label2 -> flatten entries transpose (label1-label2))))));
    CMinCones := apply(nonacyclicImages, label1 -> intersection(lambdaCone, coneFromHData(matrix(apply(nonacyclicImages, label2 -> flatten entries transpose (label2-label1))))));
    betaIplus := submatrix(classMap,Iplus)*matrix(table(length(Iplus),1,(i,j) -> 1));
    imageIterator := 0..(length(nonacyclicImages)-1);
    scan(imageIterator, i -> scan(imageIterator, j -> if not contains(coneFromHData(transpose(nonacyclicImages_i - nonacyclicImages_j-betaIplus)), intersection(CMaxCones_i,CMinCones_j)) then return false));
    return true
    );
\end{lstlisting}
The method \lst{weightsTooWide} takes as inputs the \lst{List nonacyclicImages} representing $\{\pi_\R(\bs{a}) \ | \ \bs{a}\in N_{I^{=0}} \}$, the \lst{Cone lambdaCone} representing $\overline{C}_{I^{=0},I^{>0}}$ (see \Cref{lemma:closure_of_C_plus_zero_cone}), the \lst{ZZ Matrix classMap} representing $\pi_\R$, and the \lst{List Iplus} representing $I^{>0}\subset [n]$. The output is a \lst{Boolean} describing whether or not the sufficient condition of \Cref{cor:convex_sufficient_condition} is met for the fixed $I^{=0}$. This is done by computing each cone $C_{\bs{a}}^{\omax}$ (line 3) and each cone $C_{\bs{a}}^{\omin}$ (line 3) as in \Cref{def:C_max_min_cones}. Then for each $\bs{a},\bs{a}'\in N$, we check that $C_{\bs{a}}^{\omax}\cap C_{\bs{a}'}^{\omin}\subset H_{\pi(\bs{a}-\bs{a}')-\sum_{i\in I^{>0}}\beta_i}^{\ge 0}$ (line 7).

We are now in a position to use the above methods to justify the computation performed in the proof of \Cref{thm:no_stratification_making_windows_equal_exist}. This is detailed in the below example.
\begin{example}\label{example:no_stratifications_code} In this example we provide the code to verify \Cref{thm:no_stratification_making_windows_equal_exist} below.
\begin{lstlisting}
classMap = transpose matrix({{0,0,1},{0,0,1},{0,2,1},{0,2,1},{1,1,1},{2,1,1},{3,0,1}})
divisorVector = transpose matrix({{0,0,2,-5,7,5,-2}})
Iplus = {5,6}
Ls = {{},{0,1},{4},{2,3},{0,1,4},{2,3,4}}
divisorPolytope = getDivisorPolytope(classMap, divisorVector)
fList = apply(entries map(ZZ^(numcols classMap)), rowVec -> matrix({rowVec}))
allClosedCells = getClosedCells(divisorPolytope, fList)
weightsTooWideWrapped = L -> (
    nonacyclicLabels = getNonacyclicLabels(classMap, divisorPolytope, allClosedCells, L);
    nonacyclicImages = unique(apply(nonacyclicLabels, c -> classMap*c));
    lambdaCone = getLambdaPlusMinusCone(classMap, Iplus, L);
    weightsTooWide(nonacyclicImages,lambdaCone,classMap, Iplus)
    );
apply(Ls, L -> weightstooWideWrapped(L))
\end{lstlisting}
See \S\ref{subsection:counter_example_theory} for the details of our particular GIT problem and chamber. The weights $\beta_\rho$ are encoded as $\pi_\R$ in line 1, the lift $D$ is encoded in line 2, the set $I^{>0} = \{5,6\}$ is encoded in line 3, and the possible $I^{=0}$ sets for the fixed $I^{>0}$ such that $C_{I^{>0},I^{=0}}\ne\emptyset$ are given in line 4. In line 5 we compute $P_D$ and compute $\cells(P_D)$ in line 7. Lines 8-13 wrap the method \lst{weightsTooWide} so that with input $I^{=0}$, the output is a \lst{Boolean} representing whether or not the sufficient condition of \Cref{prop:non_convex_sufficient_condition} is met. In line 14 we compute whether or not the sufficient condition is met for each possible $I^{=0}$ set. When one runs the above code, the output of line 14 is a \lst{List} with every entry \lst{true}. We conclude that the sufficient condition of \Cref{prop:non_convex_sufficient_condition} is met, thus verifying the conclusion of \Cref{thm:no_stratification_making_windows_equal_exist}.
\end{example}
\section*{Notational Index}\label{section:notational_index}
\renewcommand{\arraystretch}{1.5}
\begin{longtable}{p{0.2\textwidth} p{0.7\textwidth}}
\textbf{Notation} & \textbf{Meaning} \\
\hline

$i, i_{\R}, \pi, \pi_{\R}$, & the maps from \eqref{eqn:ses_for_lattices} and \eqref{eqn:ses}\\

$[n]$ & the set $\{1,2,\dots,n\}$. \\

$\beta_{\rho}$ & the weights of a GIT action defined as $\beta_{\rho}\coloneqq \pi_{\R}(e_{\rho})$. \\

$\beta_{J}$ & Cone in $\wh{G}_{\R}$ generated by $\beta_{\rho}$ with $\rho \in J$. \\

$[\rho]$ &  the subset of $\rho'\in [n]$ for which $\beta_{\{\rho'\}} = \beta_{\{\rho\}}$.\\

$I^{\star}_{\lambda}$ & subset of $i$ in $[n]$ for which $\wh{\lambda}(\beta_i)$ satisfies the inequality $\star$ for $\star\in \{=0,\ge 0, \le 0 , >0, <0\}$. \\

$H^{\star}_{f}$ & the subset of $V$ consisting of $v\in V$ so that for the functional $f:V\to \R$, $f(v)$ satisfies inequality $\star$ for any $a\in \R$ and $\star\in \{=a, \ge a, \le a, >a, <a\}$. \\

$B_{\Gamma}$ & the irrelevant ideal associated to a chamber $\Gamma$ of $\Sigma_{GKZ}$.  \\

$\Sigma_{\Gamma}$ & the fan in $N_{\R}$ associated to a chamber $\Gamma$.  \\

$\lambda$ & a one-parameter subgroup $\lambda:\G_{m} \to G$. \\

$\Lambda$ & the set of primitive 1-PSs $\lambda$ whose hyperplanes contain walls of $\Sigma_{GKZ}$ \\

$\Lambda^{\circ}$ & the subset of 1-PSs $\Lambda\setminus \lp |\Sigma_{GZK}|^\vee \cup -|\Sigma_{GKZ}|^\vee\rp$ \\

$\wh{(-)}$ & the functor $\wh{(-)}\coloneqq \Hom(-,\G_{m})$ either from linear algebraic groups to Abelian groups or vice-versa. \\

$S$ & the Cox ring of $G\subset \Gm^n$, i.e. the $\wh{G}$-graded polynomial ring $S\coloneqq k[x_1,\dots,x_n]$. \\

$x^J$ & the square-free monomial $\prod_{j\in J}x_j\in S$ for $J\subset [n]$. \\ 

$\CX_{\Gamma}$ & the quotient stack $[(\A^{n}\setminus V(B_{\Gamma}))/G]$. \\

$\CX_{\Sigma_\Gamma, i_\Gamma^\vee}$ & the quotient stack $\left[ \lp \A^{\Sigma(1)}\setminus V(B_\Sigma)\rp / G_{\Sigma,i_\Gamma^\vee}\right]$ associated with the stacky fan $(\Sigma_\Gamma, i_\Gamma^\vee)$.  \\

$D(\CX)$ & the bounded derived category of coherent sheaves on the smooth DM quotient stack $\CX$. \\

$\CO_{\CX_{\Gamma}}(\chi)$ & the $G$-equivariant line bundle on $\A^{n}\setminus V(B_{\Gamma})$ given a character $\chi \in \wh{G}$. \\

$D$ & a formal linear combination $D = \sum_{\rho\in I}a_\rho e_\rho$.  \\

$P_D$ & the polytope $\{x\in M_\R \ | \ \la x, u_\rho\ra \ge -a_\rho\text{ for each }\rho \in I\}$ or its embedding in $\R^n$ via $i_\R$ when clear.\\

$(X)_{\lambda}^{=0}$ & the fixed locus of $\lambda$ defined as $X_{\lambda}^{=0} \coloneqq \{x\in X \ | \ \lambda(t)\cdot x = x,\ \forall t\in \Gm\}$. \\

$(X)_{\lambda}^{\leq0}$ & the contracting locus of $\lambda$ defined as $X_{\lambda}^{\le 0} \coloneqq \{x\in X \ | \ \lim_{t\to \infty}\lambda(t)\cdot x\in X_{\lambda}^{=0}\}$.  \\

$X_{i}$ & the $G$-invariant open subvariety $\A^n\setminus \bigcup_{j<i}(\A^n)_{\lambda_j}^{\le 0}$ \\

$T^{>0}_{\lambda}, T^{<0}_{\lambda}$ & the positive and negative widths, i.e. $T_\lambda^{>0}\coloneqq \sum_{\rho\in I_{\lambda}^{>0}}\wh{\lambda}(\beta_\rho)$ and $ T_\lambda^{<0}\coloneqq \sum_{\rho\in I_{\lambda}^{<0}}\wh{\lambda}(\beta_\rho)$ \\

$W^{\bs{w}}_{\bs{\lambda}}$ & the grade restriction window induced by an HKKN stratification $\bs{\lambda}$ and a shift vector $\bs{w}$. \\

$Z_{G}$ & the image of $[0,1)^{n}$ under the map $\pi_{\R}$ in (\ref{eqn:ses}). \\

$\Theta$ & the Bondal-Thomsen collection, i.e. the set of characters $\Theta \subseteq \wh{G}$ that lie in $-Z_{G}$. \\

$K_{\Theta}$ & the subcategory of $D([\A^n/G])$ generated by line bundles $\CO_{[\A^n/G]}(\chi)$. \\

$W^\Cox_{\Gamma}$ & the Cox Window of $\CX_{\Gamma}$. \\

$F^{D}_{\bullet}(D)$ & the cellular complex of free, graded $S$-modules supported on $P_{D}$. \\

$\cells(P_{D}), \cells^{i}(P_{D})$ & the set of all cells of $P_D$ and the set of all $i$-dimensional cells of $P_D$. \\

$p_{L}$ & the distinguished point on $\A^{n}$ defined as $(p_{L})_t = 1$ if $t \in L$ and $(p_{L})_t = 0$ otherwise.\\

$\CS_{\delta,L}$ & the subspace of $P_D$ defined in \Cref{def:S_regions}.\\

$\chi(\delta,L,D)$ & the distinguished character defined in \Cref{def:distinguished_character}. \\ 

$\overline{K}$ & $\overline{K} \coloneqq \{i\in [n]\setminus L \ : \ \ceil{\delta}_i = -a_i\}$ \\ 
\end{longtable}

\bibliography{biblio}
\bibliographystyle{amsalpha}
\end{document}